\ifdefined\useorstylema

\documentclass[opre,nonblindrev]{./macros/informs3}
\usepackage{amsmath,amssymb,amsfonts}
\usepackage{amsthm}

\usepackage{./macros/packages-or}
\usepackage{./macros/editing-macros}%
\usepackage{./macros/statistics-macros-or}
\usepackage{./macros/queueing-macros-ms}

\usepackage[square, numbers]{natbib}
 \bibpunct[, ]{(}{)}{,}{a}{}{,}%
 	\def\bibfont{\small}%
 \def\bibsep{\smallskipamount}%

\usepackage{booktabs}
\usepackage{rotating}
\usepackage{fancyvrb}
\usepackage[colorlinks=true, linkcolor=blue, citecolor=blue, filecolor=blue, urlcolor=blue]{hyperref}
\usepackage{enumerate}
\usepackage{pgfplotstable}
\usepackage{graphicx}
\usepackage{subcaption}
\usepackage{float} 
\usepackage{tabularx}
\usepackage{overpic}
\usepackage{tikz}
\usepackage{rotating}
\usepackage{psfrag}
\usepackage{array}
\usepackage{appendix}
\usepackage{endnotes}
\usepackage{blindtext}
\usepackage{titlesec}
\usepackage{xcolor}

\DoubleSpacedXI 

\usepackage{endnotes}
\let\footnote=\endnote

\usepackage{natbib}
 \bibpunct[, ]{(}{)}{,}{a}{}{,}%
 \def\bibfont{\small}%
 \def\bibsep{\smallskipamount}%
\TheoremsNumberedThrough     
\ECRepeatTheorems

\EquationsNumberedThrough    

\usepackage{hyperref}

\usepackage{pgfplotstable}
\usepackage{graphicx}
\usepackage{subcaption}
\usepackage{float} 
\usepackage{tabularx}
\usepackage{overpic}
\usepackage{tikz}
\usepackage{rotating}
\usepackage{psfrag}
\usepackage{array}
  
\RUNAUTHOR{Lee, Namkoong, and Zeng}
\RUNTITLE{Design and Scheduling of an AI-based Queueing System}

\TITLE{Design and Scheduling of an AI-based Queueing System}

\ARTICLEAUTHORS{%

\AUTHOR{Jiung Lee}
\AFF{Coupang,
  \EMAIL{jiunglee28@gmail.com}} 

\AUTHOR{Hongseok Namkoong}
\AFF{Decision, Risk, and Operations Division, Columbia Business School, New York, NY 10027, \EMAIL{namkoong@gsb.columbia.edu}} 

\AUTHOR{Yibo Zeng}
\AFF{Department of Industrial Engineering and Operations Research,
  Columbia University, New York, NY 10027,
  \EMAIL{FILL IN}} 

} 

\ABSTRACT{\input{abstract}}


\maketitle

\else

\documentclass[11pt]{article}
\usepackage{amsmath,amssymb,amsfonts}
\usepackage{amsthm}
\usepackage[margin=1in]{geometry}
\usepackage[numbers]{natbib}
\usepackage{booktabs} 
\usepackage{./macros/packages-or}
\usepackage{./macros/editing-macros}%
\usepackage{./macros/statistics-macros-or}
\usepackage{./macros/queueing-macros-ms}
\usepackage[toc,page]{appendix} 
\usepackage{minitoc}
\newtheorem{assumption}{Assumption}
\newtheorem{definition}{Definition}
\newtheorem{lemma}{Lemma}
\newtheorem{proposition}{Proposition}
\newtheorem{theorem}{Theorem}
\newtheorem{claim}{Claim}

\begin{document}

\abovedisplayskip=8pt plus0pt minus3pt
\belowdisplayskip=8pt plus0pt minus3pt



\begin{center}
 {\LARGE Design and Scheduling of an AI-based Queueing System} \\ 
  \vspace{.5cm}
  {\Large Jiung Lee ~~~ Hongseok Namkoong$^1$ ~~~
    Yibo Zeng$^1$} \\ 
  \vspace{.2cm} 
  {\large Columbia University$^{1}$}  \\
  \vspace{.2cm}
  \texttt{jiunglee28@gmail.com}
  \hspace{1cm}
  \texttt{namkoong@gsb.columbia.edu}
    \hspace{1cm}
    \texttt{yibo.zeng@columbia.edu}
\end{center}


\begin{abstract}%
To leverage prediction models to make optimal scheduling decisions in service
systems, we must understand how predictive errors impact congestion due to
externalities on the delay of other jobs. Motivated by applications where
prediction models interact with human servers (e.g., content moderation), we
consider a large queueing system comprising of many single server queues where
the \emph{class} of a job is estimated using a prediction model.  By
characterizing the impact of mispredictions on congestion cost in heavy
traffic, we design an index-based policy that incorporates the predicted class
information in a near-optimal manner.  Our theoretical results guide the
design of predictive models by providing a simple model selection procedure
with downstream queueing performance as a central concern, and offer novel
insights on how to design queueing systems with AI-based triage. We illustrate
our framework on a content moderation task based on real online comments,
where we construct toxicity classifiers by finetuning large language models.
\end{abstract}

\fi

\section{Introduction}
\label{section:introduction}

Recent advances in predictive modeling present new opportunities for utilizing
rich features to inform allocation of scarce resources under stochastic
workloads.  Since prediction is not the end goal, effective decision-making
requires understanding complex endogenous interactions that noisy predictions
introduce. To crystallize how prioritizing a particular job based on a
misprediction incurs negative externalities that affect the congestion of
other jobs, we study a simple single server queue where the true class of a
job that governs its service time is unknown and predicted using a machine
learning model. Because solving for the optimal scheduling policy is
computationally intractable even when job classes are known, owing to large
state and policy spaces~\citep{PapadimitriouTs94}, we study highly congested
systems in the heavy traffic limit as is standard in the queueing
literature~\citep{Reiman84, VanMieghem95, HarrisonZe04, Whitt02,
  MandelbaumSt04}.

When the true class of a job is known, a simple index-based myopic
policy---the oracle \gcmu---that greedily prioritizes jobs with the highest
marginal cost of delay is known to be heavy traffic
optimal~\citep{VanMieghem95, MandelbaumSt04}.  When the true classes are
unknown and predicted using a classifier, a naive adaptation of this index
rule is
\begin{equation}
  \label{eq: naive gcmu}
  \argmax_{l\in[K]} \fmu_l(t) (C_l)'( \fa_l( t)) \qquad \qquad  \mbox{\naivemethod},
\end{equation}
where $\fmu_l$ is the service rate for a \emph{predicted} class $l$,
$C_k(\cdot)$ is a convex delay cost defined on the time between arrival and
service completion, and $\fa_l(\cdot)$ is the age of the oldest unfinished job
with predicted class $l$. This index policy does not consider
misclassifications and ignores the fact that delay cost depends on the true
class label instead of the predicted class. This mismatch leads to suboptimal
scheduling decisions as we show in Theorem~\ref{theorem: optimality of our
  policy} to come.

We propose and analyze a new index policy---the \ourmethod---that remains
simple to implement, and optimally accounts for the impact of prediction
errors on the overall cost of delay.  Instead of scheduling by a predicted
class's own cost as in the naive rule~\eqref{eq: naive gcmu}, the
\ourmethod~corrects for misclassification using the classifier's confusion
matrix.  We derive this precise form by analyzing the stochastic fluctuations
in the queue lengths of the \textit{unobservable} true class jobs, and
aggregating them to represent the fluctuation in the queue length of each
\textit{predicted class} (Section~\ref{section: convergence and lower bound}).
In particular, we characterize the optimal queueing cost
(Theorem~\ref{theorem: HT lower bound}) that captures these diffusion-scale
misclassification fluctuations and reveals the heavy-traffic target that an
optimal policy must approach.  The KKT conditions for the optimal workload
allocation across \emph{predicted} classes become a weighted average of
\emph{true} class costs, which motivates the particular form of \ourmethod~
(Section~\ref{section: heavy-traffic optimality of pcmu}).

Although our policy is a simple index rule, establishing its optimality is far
from routine. The central challenge of this paper is not the form of the
policy but its analysis: proving that this index is asymptotically optimal in
heavy traffic over all feasible policies. Misclassification makes the
heavy-traffic analysis substantially more subtle than the classical setting
suggests.

First, it is a priori unclear what target an optimal policy should even aim for.
Because the scheduler observes only predicted classes while delay costs are
incurred by the unobserved true classes, prioritizing a job no longer has a
transparent effect on cost: a predicted class queue mixes jobs from several true
classes, so serving it relieves congestion for some true classes more than
others. In the classical setting where classes are observed, the cumulative cost
of any policy is written directly in terms of the true delay costs $C_k$, which
makes the instantaneous optimization problem and the optimality condition
transparent~\citep{VanMieghem95, MandelbaumSt04}. Misclassification breaks this
correspondence: we must first identify how each predicted class queue is composed
of true class jobs at the diffusion scale, where the fluctuations of a predicted
class queue aggregate the \textit{unobservable} true class fluctuations split by
the confusion matrix. 

Second, showing that the \ourmethod~attains this target reduces to proving
that the largest gap among the predicted-class indices vanishes in heavy
traffic, so that the queues collapse onto the optimal workload allocation. The
indices are driven by the age processes, which evolve endogenously with the
scheduling decisions, so their gaps must be controlled along the policy's own
sample paths. This is precisely the convergence
that~\citet[Eq.~(54)]{VanMieghem95} states without proof, and that we found
nontrivial to justify. We establish it through an induction over
diffusion-scaled time intervals of size $O(n^{-1/2})$ that tracks the index
dynamics while controlling the approximation errors of the predicted-class
processes, and we find that strong convexity of the cost functions is needed
for the argument to close, as opposed to the strict convexity assumed
by~\citet{VanMieghem95}.

Finally, we resolve another analytic gap in the classical heavy traffic
theory. The \gcmu~\citep{VanMieghem95} and
$D\text{-}Gc\mu$~\citep{MandelbaumSt04} indices are defined on the \emph{ages}
of waiting jobs, but their optimality is argued through an idealized version
that schedules on sojourn times, which is unimplementable since it requires
information not yet realized at the time of the decision.  To overcome this,
our analytic machinery provides the first rigorous proof of heavy traffic
optimality for the \emph{age-based} \gcmu~rule.  For~\citet{MandelbaumSt04},
our proof recovers the key implications of their fluid-scale attraction
property in the single-server case, under weaker conditions on the cost
functions. Section~\ref{subsection: comparison to VM and MS} discusses how our
results complete and strengthen~\citet{VanMieghem95}
and~\citet{MandelbaumSt04}, with the proof-level details deferred to
Sections~\ref{subsection: comparison to the optimality result in
  VM},~\ref{subsection: comparison to the optimality result in MS},
and~\ref{section: adding to the exposition in classical heavy traffic queueing
  analysis}.

\paragraph{Empirical illustration}
We study an \emph{offline}, off-policy setting, in contrast to a growing body
of work on \emph{online} learning in queueing that designs exploration
algorithms and analyzes regret against a benchmark with known
parameters~\citep{ChoudhuryJoWaSh21, KrishnasamySeJoSh16, KrishnasamySeJoSh21,
  StahlbuhkShMo21, Walton14, FreundLyWe23}. Empirical $c\mu$ rules that plug
in online estimates of unknown service rates have already been studied in this
setting~\citep{KrishnasamyArJoSh18, ZhongBiWa22}. We instead take the
classifier as given and trained offline on previously collected data, which
reflects operational constraints of modern AI-based service systems where
online experimentation is risky or unwieldy. Since service times are
determined by the true classes, observed service times do carry information
about true labels that could in principle refine the classifier online; even
so, any prediction model must be thoroughly validated prior to deployment, and
the timescale for model development (weeks to months) is far longer than that
for scheduling decisions (hours to days). We thus view our offline heavy
traffic analysis as a useful analytic device for modeling AI-based queueing
systems that operate close to system capacity.

\begin{figure}[t]
  \vspace{-1em}
\centering
\includegraphics[width=1\columnwidth]{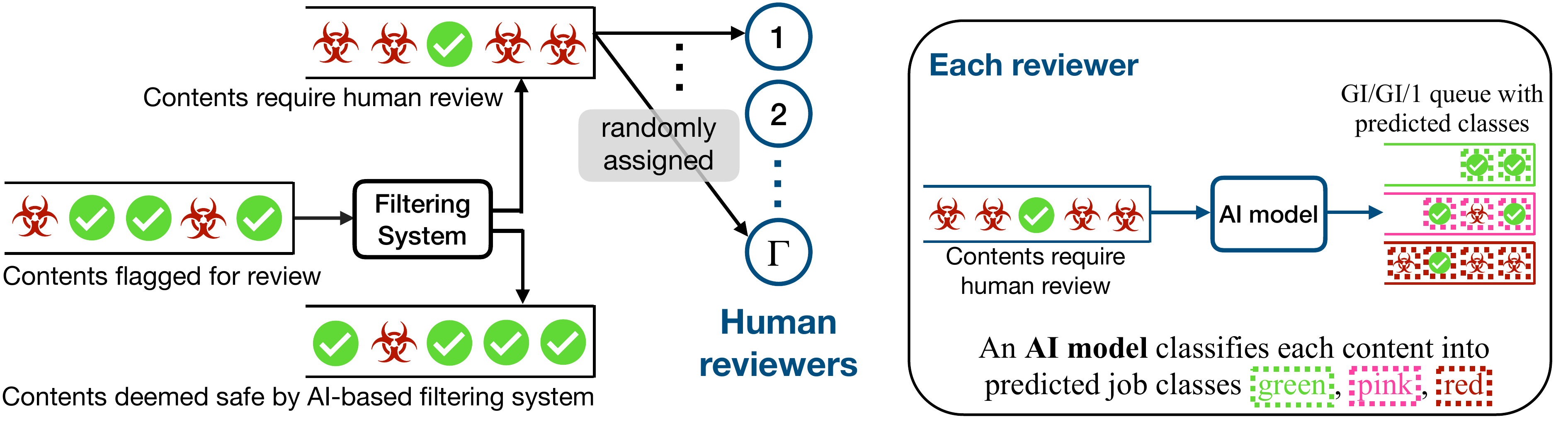}
\caption{\textbf{Schematic of a content moderation system as a triage
    system.}  Each content may be violating the user agreement (red toxicity
  symbol) or considered safe (green checkmark). This ground truth requires
  human review to uncover (``service'').  Contents are
  flagged for review by users or automated filters, which we view as
  ``entering'' the triage system. The online platform uses an initial AI
  model to filter out contents most deemed to be safe. Then,
  remaining jobs/contents are
  randomly assigned to the human reviewers, a common practice due to fairness
  considerations in terms of mental workload. An AI model classifies each
  content into different classes (e.g., hate speech on a protected group),
  placing them in the corresponding virtual queue for the predicted class. }
\label{fig:diagram}
\end{figure}

We illustrate our framework and its practical implications on content
moderation, a critical process for maintaining the health and sustainability of
online platforms. Delays in removing violating posts such as hate speech can
exacerbate their harm, and while clear-cut cases can be filtered out by an
initial AI-based filtering system, nuanced moderation requires human reviewers
who account for nonstationary social contexts and avoid unnecessary censorship
and violations of freedom of speech~\citep{AllouahKrZhAvBoDaGoPuShStTa23,
  MakhijaniShAvGoStMe21}. To ensure fairness and similar workload across
reviewers, jobs are typically assigned to human reviewers in an identically
random manner, which reduces the dynamic scheduling problem to a single-server
queue for each reviewer where job classes (such as toxicity, or whether content
targets a protected demographic group) are a priori unknown and predicted by an
AI model (Figure~\ref{fig:diagram}).

Our characterization of the optimal cost has two practical implications that we
develop in turn. First, it guides the design of the prediction model. Although
predictive performance is rarely the final goal, models are typically validated
on predictive measures such as precision or recall; yet overparameterized
models (e.g., neural networks) can attain the same predictive performance while
exhibiting very different downstream decision performance~\citep{DAmourEtAl20,
  BlackRaBa22}. We propose a model selection procedure based on the
\emph{cumulative queueing cost}, and demonstrate its advantages over
conventional predictive measures (Section~\ref{section: characterization of the
  optimal cost with q}). Second, it guides the design of the queueing system
itself. We design an AI-based triage system that determines staffing and
filtering levels by trading off filtering cost, predictive performance, hiring
cost, and congestion (Section~\ref{section: design of an AI triage system}).
Traditional prediction-based metrics accurately reflect the total cost only
when filtering or hiring costs dominate, and otherwise require computationally
expensive queueing simulations; our characterization instead determines the
optimal staffing and filtering levels by simulating a (reflected) Brownian
motion (Section~\ref{section: experiment triage system}).

\begin{figure}[t]
  \vspace{-2em}
\begin{minipage}[t]{0.485\textwidth}
\includegraphics[width=0.985\columnwidth]{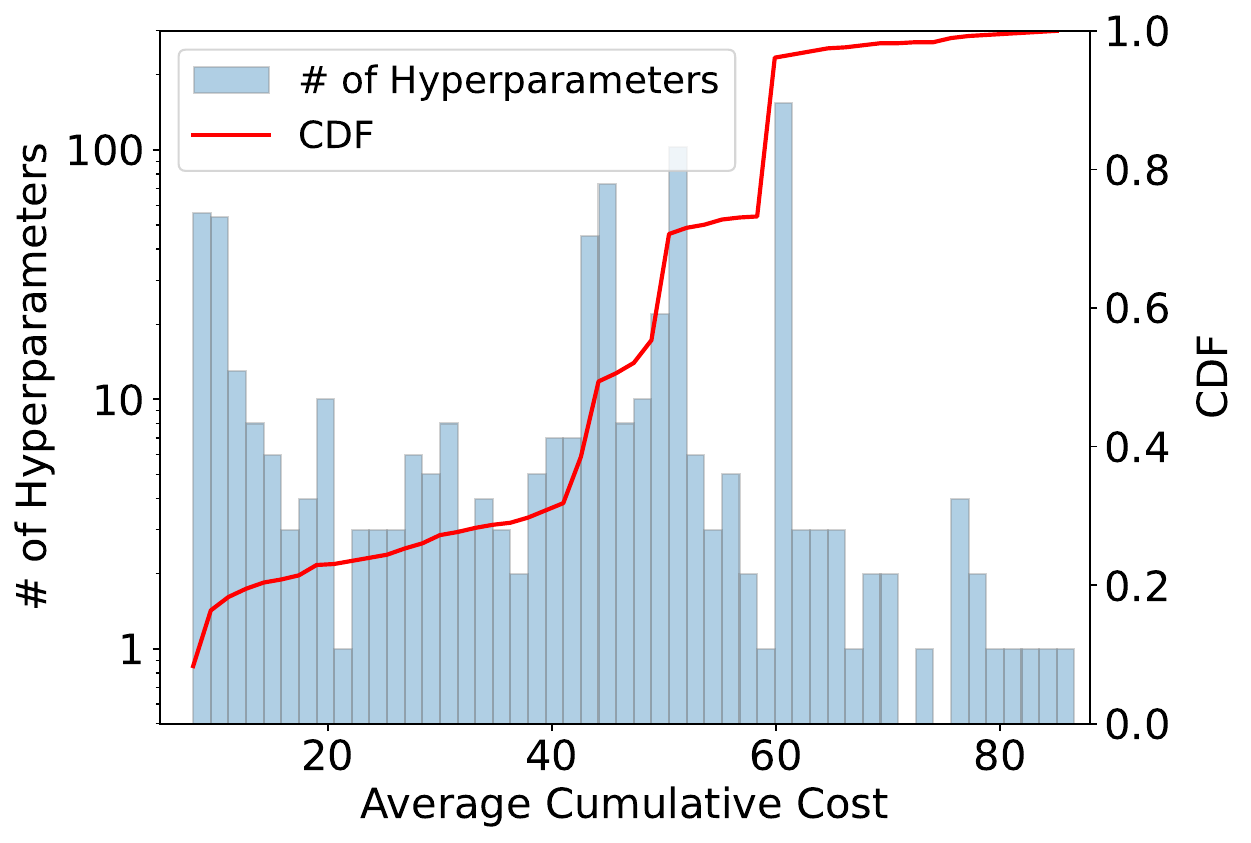}
\caption{Histogram of average cumulative queueing cost of deep Q-learning
  policies over 672  hyperparameter configurations. }
\label{fig: DRL_distribution}
\end{minipage}
\begin{minipage}[t]{0.485\textwidth}
\includegraphics[width=0.9\columnwidth]{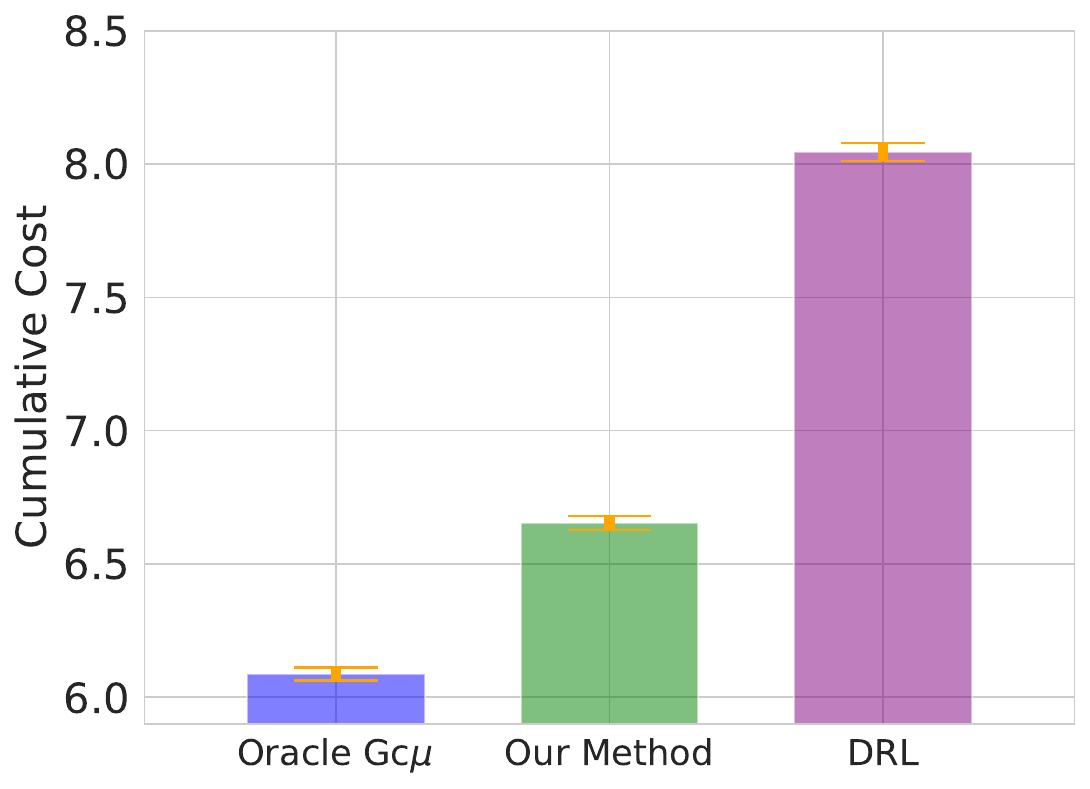}
\caption{Cumulative cost with 2$\times$ standard errors }
\label{fig: cumulative cost example}
\end{minipage} \vspace{-1em}
\end{figure}

We construct a simulated content moderation system from real online comments,
with toxicity classifiers obtained by finetuning large language models, and
compare the \ourmethod~against offline deep reinforcement learning (DRL). While
flexible, DRL methods require significant engineering effort to train reliably
and are highly sensitive to hyperparameters, implementation details, and even
random seeds~\citep{HendersonIsBa18, WaltonXu21, DaiGl21}. On a single-server
queue with 10 classes, deep Q-learning with experience replay exhibits
substantial variation in performance across hyperparameters even when using
identical instant reward functions, training and testing environments, and
random seeds across training runs (Figure~\ref{fig: DRL_distribution}), and the
simple index policy~\ourmethod~significantly outperforms the best-performing DRL
configuration (Figure~\ref{fig: cumulative cost example},
Section~\ref{section: numerical experiments}). See
Section~\ref{section:related-work} for a thorough literature review.

\section{Related work}
\label{section:related-work}


Our work builds on the large literature on queueing, as well as the more
nascent study of decision-making problems with prediction
models~\citep{ArgonZi09, MivsicPe20, SinghGuVa20, KotaryFiVaWi21, ChenDo21,
  SunArZi22}.  Unlike previous works that study relatively simple optimization
problems (e.g., linear programming~\citep{ElmachtoubGr21}) that take as input
predictions, our scheduling setting requires modeling the endogenous impact of
misclassifications.

\subsection{Heavy traffic analysis}
\label{subsection: comparison to VM and MS}

Heavy traffic analysis allows circumventing the complexity of state/policy
spaces via state-space collapse, thereby identifying asymptotically optimal
queueing decisions~\citep{HarrisonZe04, MandelbaumSt04, Reiman84,
  VanMieghem95, Whitt02}. The $c\mu$-rule was shown to be optimal among
priority rules in~\citep{CoxSm61}, and~\citet{VanMieghem95} laid out an
argument for heavy traffic optimality of the \gcmu~with convex delay costs in
single-server systems with general distributions of interarrival and service
times. Under a heavy traffic regime defined with a complete resource pooling
condition~\citep{HarrisonLo99}, \citet{MandelbaumSt04} extended the result to
multi-server and multi-class queues. In the many server Halfin-Whitt heavy
traffic regime (where the server pool is also
scaled~\citep{HalfinWh81}),~\citet{GurvichWh08} showed their state-dependent
policy that minimizes the holding cost reduces to a simple index-rule with
linear holding costs, and to the G$c\mu$-rule with convex costs. When
customers in queues can abandon systems, a similar index rule that accounts
for the customer abandonment rate was shown to minimize the long-run average
holding cost under the many-server fluid limit~\citep{AtarGiSh10}.

We focus on the single-server model and relax a common assumption that the
class of every job is known.  We study the impact of AI models in queueing
jobs, and use the heavy traffic limit to analyze the downstream impacts of
misclassifications. Our results provide a unified framework for evaluating and
selecting AI models for optimal queueing. Along the way, we also provide
rigorous proofs for the classical setting with known classes by establishing
steps that~\citet{VanMieghem95} states without proof and identifying the
conditions under which they hold.
The optimality of the generalized $c\mu$ rule turns on a single delicate step:
in heavy traffic, the largest gap between the class indices must vanish, so that
the queues collapse onto the cost-minimizing workload allocation.
\citet[Eq.~(54)]{VanMieghem95} \emph{states this convergence without proof} rather than
establishing it, and does so for an idealized index defined on workloads or
sojourn times rather than the implementable age-based index that the rule
actually uses. We found this step nontrivial to justify.

We close this gap directly. Under the \ourmethod, we prove that the index gaps
vanish through an induction over diffusion-scaled time intervals of size $O(n^{-1/2})$ that tracks
the index dynamics while controlling the approximation errors of the
predicted-class processes (Proposition~\ref{prop: convergence of max gcmu
  difference}), and we identify strong convexity of the cost functions as a
sufficient condition for the argument to close. Specialized to observed classes
($Q^n = I$), this already gives the first rigorous proof of heavy-traffic
optimality for the age-based \gcmu~itself. The optimality result
of~\citet{MandelbaumSt04} rests on a state-space-collapse attraction property,
established through a fluid-scale analysis of multi-server systems. Our induction
on the vanishing difference between the class indices shows that the key
implication of this property~\citep[Eq.~(66)]{MandelbaumSt04} can be established,
for the single-server case, under weaker conditions on the cost functions and
directly on the diffusion scale rather than the fluid scale. Details are provided
in Section~\ref{subsection: comparison to the optimality result in MS}.

\subsection{Broader literature on learning and queueing}

In contrast to our single-server setting, where stability is guaranteed under any work-conserving policy with nominal load below one and the main challenge is to quantify the endogenous impact of misclassification errors on congestion and delay costs in heavy traffic, a growing body of work studies information constraints in networked systems where stability itself is a first-order concern. In particular, \citet{MassoulieXu18} characterize the capacity region of information-processing networks and show how limited information can fundamentally restrict stabilizability, while \citet{XuZhong20} analyze how different forms of memory and information structures shape dynamic resource allocation and stability in stochastic networks. Extending our framework to such settings, where misclassification interacts with stability and capacity in the spirit of these works, is an important direction for future research.

The challenge of unknown traffic parameters was identified as early
as~\citet{Cox66}. Using an off-policy ML model in queueing systems was also
proposed for classifying jobs into different types or priority
classes~\citep{SinghGuVa20, SunArZi22} and predicting service
times~\citep{ChenDo21}.~\citet{ArgonZi09} focuses on minimizing the
mean \emph{stationary} waiting time with Poisson arrivals, while we allow
general distributions of arrival and service times in our heavy traffic
analysis. In the specific case of linear delay costs and steady-state waiting time as
the performance metric, the \ourmethod~bears surface-level resemblance
to~\citet[Section 8]{ArgonZi09}'s policy defined with conditional distributions
of the true classes given a job's signal. In contrast, we model the increasing
marginal cost of delay through strongly convex cost functions and prove heavy
traffic \emph{optimality} over all feasible policies, rather than dominance over
first-come-first-serve policies.~\citet{SunArZi22} consider a two-class setting (triage or not) where
classes can be inferred with additional time, and analyze when it is optimal
to triage all (or no) jobs.  Importantly, they assume service times follow
predicted classes.  ~\citet{ChenDo21} develop a two-class priority rule using
predicted service times and show the convergence of the queue length process
to the same limit as in the perfect information case when estimation error is
sufficiently small. In contrast, we characterize the optimal queueing cost
given a fixed classifier instead of aiming to match the performance of the
perfect classifier. Our approach allows us to provide guidance on model
selection for classifiers as we illustrate in Section~\ref{section:
  characterization of the optimal cost with q}.

Closest to our setting,~\citet{SinghGuVa20} also infer unobserved job types
from observable features through a classifier and likewise find that higher
classification accuracy need not lower queueing cost. They study a
static-priority M/G/1 system with linear delay costs, minimize the steady-state
average waiting cost, and compare a type-first design (predict the type, then
assign a queue) against a direct design that trains the classifier end to end
to map features to queues. Our framework differs along three axes. We analyze
the heavy-traffic diffusion limit with strongly convex delay costs, we establish
optimality over all feasible policies rather than over static-priority
assignments, and we take the classifier as fixed and select among given
classifiers (Section~\ref{section: characterization of the optimal cost with q})
rather than training one end to end.
  
Going beyond simple index policies, deep reinforcement learning (DRL)
algorithms can be used for queueing systems with unknown
parameters.~\citet{DaiGl21} develop a policy optimization approach for
multiclass Markovian queueing networks and proposes several variance reduction
techniques.~\citet{PavseChXiHa23} combine proximal and trust region-based
policy optimization algorithms~\citep{SchulmanLeAbJoMo15} with a
Lyapunov-inspired technique to ensure stability. Developing further approaches
to overcome the challenges of applying RL algorithms in queueing (e.g.,
infinite state spaces, unbounded costs) is a fertile direction of future
research.

There is a growing body of work on learning in queueing systems that focus on
\emph{online} learning and analyze regret, the performance gap between the
learning algorithm and the best policy in hindsight with the complete
knowledge of system parameters~\citep{DaiGl21, FreundLyWe22, FreundLyWe23,
  GaitondeTa23, KrishnasamySeJoSh16, KrishnasamyArJoSh18, KrishnasamySeJoSh21,
  SentenacBoPe21, Walton14, WaltonXu21, ZhongBiWa22}. Inspired by the
well-known static priority policies in queueing literature~\citep{AtarGiSh10,
  AtarGiSh11, MandelbaumSt04, PuhaWa19, PuhaWa21, VanMieghem95}, empirical
versions of such policies were proposed where plug-in estimates of unknown
parameters are used to compute static priorities. When service rates are
unknown,~\citet{KrishnasamyArJoSh18} propose an empirical $c\mu$ rule for
multi-server settings and show constant regret for linear cost functions,
and~\citet{ZhongBiWa22} develop an algorithm for learning service and
abandonment rates in time-varying multiclass queues with many servers and show
the empirical $c\mu/\theta$ rule achieves optimal regret. In the machine
scheduling literature, where a finite set of jobs are given (with no external
arrivals),~\citet{LeeVo21} studies settings where delay costs are unknown, and
show that a plug-in version of the $c\mu$-rule can achieve near-optimal regret
when coupled with an exploration strategy.


For more general queueing networks,~\citet{WaltonXu21} present a connection
between the MaxWeight policy~\citep{TassiulasEp92} and Blackwell
approachability~\citep{Blackwell56}, relating the waiting time regret to that
of a policy for learning service rates. Borrowing insights from the stochastic
multi-armed bandit literature~\citep{Slivkins19}, a body of
work~\citep{ChoudhuryJoWaSh21, KrishnasamySeJoSh16, KrishnasamySeJoSh21,
  StahlbuhkShMo21, Walton14} develops learning algorithms to minimize expected
queue length, addressing challenges in the queueing bandit model such as
ensuring stability until the parameters are sufficiently
learned~\citep{KrishnasamySeJoSh21}.~\citet{FreundLyWe23} propose a new
performance measure of time-averaged queue length, and show near-optimality of
the upper-confidence bound (UCB) algorithm in a single-queue multi-server
setting, as well as new UCB-type variants of MaxWeight and
BackPressure~\citep{TassiulasEp92} in multi-queue systems and queueing
networks, respectively. For queueing systems with multi-server multiclass
jobs, \citet{YangSrYi23} recently developed another UCB-type variant of the
MaxWeight algorithm. Learning service rates has also been studied in
decentralized queueing systems, where classes of jobs are considered as
strategic agents~\citep{FreundLyWe22, GaitondeTa23, SentenacBoPe21} and
stability is a primary concern. Motivated by content moderation, a concurrent
work~\citep{LykourisWe24} studies the joint decision of content classification,
and admission and scheduling for human review in an online learning framework.


\section{Model}
\label{section: model}

We begin by presenting our analytic framework in the heavy traffic regime.
There are two possible data generating processes we can study.  We could view
jobs as originating from a \emph{single} common arrival process, where
interarrival times are independent of job features, true classes labels, and
service times.  This single arrival stream allows us to disentangle the
arrival and service processes of predicted classes, and directly use the
diffusion limit to show optimality of the \ourmethod.

On the other hand, we may consider a more general generating process where the
arrival and service processes for different classes are exogenously given. In
this setting, we can still show similar mathematical guarantees as under the
single stream model using heavy traffic analysis techniques pioneered
by~\citet{MandelbaumSt04}.  However, this proof approach weakens our
optimality results: we can only establish optimality of the \ourmethod~over
p-FCFS policies (those serving jobs first-come-first-serve within each predicted
class), whereas the single stream model yields optimality over all feasible
policies. The gap is a consequence of misclassification. With class-specific
arrival streams, the arrival timing of a job is informative of its real class,
hence of its service law, so waiting jobs in a predicted class are no longer
identically distributed and the interchange argument behind the reduction to
p-FCFS breaks down. The single common arrival stream keeps interarrivals
independent of class, which preserves this reduction and is what yields
optimality over all feasible policies (see Section~\ref{subsection: convergence
of the classical queueing model}). We view the practical modeling capabilities
of the two data generating assumptions to be similar. The single arrival stream
is a good model of the content moderation system (as depicted in
Figure~\ref{fig:diagram}). Henceforth, we thus focus on the single common
arrival process for expositional clarity and crisp mathematical results.

In the remainder of this section we formalize the single-stream model and the associated
heavy-traffic framework. Our goal is to show that, under this modeling choice, the
heavy-traffic scaling and diffusion limits for the induced predicted-class processes
are entirely standard. To that end, we first lay out the indexing and notational
conventions, distinguishing true classes from predicted classes and prelimit processes
from their diffusion-scaled counterparts, which will be used consistently throughout
the paper. 

\paragraph{Indices and classes.}
We consider a sequence of single-server, multi-class queueing systems indexed by
$n \in \mathbb{N}$, where jobs are classified by a given classifier $\model$.
Within system $n$, jobs are indexed by $i = 1,2,\dots$, their true classes by
$k \in [K] := \{1,\dots,K\}$, and their predicted classes (the output of $\model$)
by $l \in [K]$. Throughout, $n$ always denotes the system index, $i$ a job index,
$k$ a true class, and $l$ a predicted class.

\paragraph{Data-generating process.}
System $n$ operates on a finite time horizon $[0,n]$ and starts empty. The system
is driven by a single exogenous arrival stream: interarrival times
$\{u_i^n\}_{i \ge 1}$ are i.i.d.\ with positive arrival rate $\tlambda^n$ (so
$\mathbb{E}[u_1^n] = 1/\tlambda^n$) and are \emph{independent of all job features,
labels, and service times}. For
$t\in [0, n]$, let $\tU_0^n (t):=\sum_{i=1}^{\lfloor t \rfloor} \tu_i^n$ be
the arrival time of the $\lfloor t \rfloor$th job in the system and
$A_0^n (t) = \max\{m: \tU_0^n (m) \leq t\}$ be the total number of jobs that
arrive up to time $t$.

For each arrival $i$ in system $n$, a feature vector $\feature_i^n \in \mathbb{R}^d$,
a one-hot true class label $Y_i^n = (Y_{i1}^n,\dots,Y_{iK}^n) \in \{0,1\}^K$, and a
service requirement $v_i^n > 0$ are realized. For each class $k$, let
$\tp_k^n:= \P^n[\tY_{1k}^n = 1]$ be the class prevalence and
$(\tmu_k^n)^{-1}:=\E^n[\tv_1^n \mid \tY_{1k}^n = 1]$ be the expected service
time in system $n$. Let $\tV_0^n (t) := \sum_{i=1}^{\lfloor t \rfloor} \tv_i^n,~t\in[0, n]$ be
the total service time required by the first $\lfloor t \rfloor$ jobs. 
The following shows underlying assumptions of our data generating process.

\begin{assumption}[Data-generating process]\label{assumption: data generating process}
For each system $n \in \mathbb{N}$:
\begin{enumerate}[(i)]
\item The sequence $\{(u_i^n, v_i^n, X_i^n, Y_i^n) : i \ge 1\}$ is i.i.d.
\item $\{u_i^n : i \ge 1\}$ is independent of $\{(v_i^n, X_i^n, Y_i^n) : i \ge 1\}$.
\item For any $i \ge 1$, $v_i^n$ and $X_i^n$ are conditionally independent given $Y_i^n$,
i.e., $v_i^n \perp X_i^n \mid Y_i^n$.
\end{enumerate}
\end{assumption}

We assume service time is conditionally independent of the covariates given the true class label $Y^n_i$ (Assumption~\ref{assumption: data generating process}(iii)), which simplies our anlysis by only
considering true class label's impact on service time. 
 In practice, if covariates directly influence service times (e.g., content length), we can mitigate such
dependency by creating more fine-grained true classes. 

\paragraph{True vs.\ predicted quantities.}
Each job $i$ is classified before entering the queue. Given features $X_i^n$, the
classifier $\model$ outputs a one-hot encoded predicted label $\fY_i^n = \model(\feature_i^n)
  = (\fY_{i1}^n,\dots,\fY_{iK}^n) \in \{0,1\}^K$.
We use a tilde ($\widetilde{\cdot}$) throughout to mark quantities indexed by predicted
classes, in contrast to their true-class counterparts; for example, $N_k^n(t)$ denotes
the number of true-class-$k$ jobs in the system at time $t$, whereas
$\widetilde N_l^n(t)$ denotes the number of jobs whose predicted class is $l$.

The classifier induces a confusion matrix
$Q^n = (\q_{kl}^n)_{k,l \in [K]}$, where $\q_{kl}^n := \P^n[\fY_{1l}^n = 1 \mid Y_{1k}^n=1]$ is
the probability of a class $k$ job being predicted as class $l,~k,l \in [K]$,
which describes how the single arrival stream is stochastically split into true and
predicted classes.

\paragraph{Heavy-traffic scaling and diffusion notation.}
We study a heavy-traffic regime in which system \(n\) operates on a time horizon
\([0,n]\). To obtain non-trivial stochastic limits, we apply the usual diffusion
scaling in time and space and use \textbf{boldface} to denote centered, diffusion-scaled
processes.

For a generic cumulative process \(X^n(\cdot)\) on \([0,n]\) with nominal drift
\(x \in \mathbb{R}\), we write its diffusion-scaled version on \([0,1]\) as
\[
  \mathbf{X}^n(t)
  \;:=\; n^{-1/2}\Big( X^n(nt) - x \, nt \Big),
  \qquad t \in [0,1].
\]
Whenever \(\mathbf{X}^n\) converges weakly as \(n\to\infty\), we denote
the limit by the same bold letter \(\mathbf{X}\). Here and throughout, weak convergence will be defined in the space of the right-continuous with left limits (RCLL) under the standard Skorohod $J_1$ topology; the precise definition and notation will be given shortly. 

\paragraph{Heavy traffic condition.}
Assumption~\ref{assumption: heavy traffic} below formalizes that all primitive parameters
converge at the rate of $n^{-1/2}$.
Intuitively, system $n$ becomes larger and more congested, with utilization
\(\rho^n := \tlambda^n \sum_k \tp_k^n / \tmu_k^n \uparrow 1\), while the centered fluctuations of key processes are of order $n^{1/2}$.

\begin{assumption}[Heavy Traffic Condition] \label{assumption: heavy traffic}
Given a classifier $\model$ and a sequence of queueing systems, there exist
$\tp_k, \q_{kl} \in [0,1]$ and $\tlambda, \tmu_k$ such that
$ \sum_{k=1}^K \tp_k \q_{kl} > 0$ for all $l\in [K]$,
$\tlambda\sum_{k=1}^K \frac{\tp_k}{\tmu_k} = 1$, and 
\begin{equation}\label{eq: convergence rate}
n^{1/2}\big(\tlambda^n - \tlambda \big) \rightarrow 0, \quad
n^{1/2}  \big(\tmu_k^n - \tmu_k\big) \rightarrow 0,\quad
n^{1/2} \big(\tp_k^n - \tp_k\big) \rightarrow 0, \quad 
n^{1/2} \big(\q_{kl}^n - \q_{kl} \big)
\rightarrow 0. 
\end{equation}
\end{assumption}
\noindent The $o(n^{-1/2})$-rate convergence of  $\tlambda^n$ and
$ \tmu_k^n$ aligns with the standard GI/GI/1
heavy-traffic scaling assumptions, e.g,~\cite[Eq.~(2)]{MandelbaumSt04}, and as usual we have that
\textit{traffic intensity} $\trho^n$ converges to $1$ at $o(n^{-1/2})$-rate
\begin{equation}
  \label{eq: heavy traffic condition}
  \begin{aligned}
    n^{1/2}[\trho^n - 1 ] 
    &=n^{1/2} \Big[\tlambda^n\sum\limits_{k=1}^K \frac{\tp_k^n}{\tmu_k^n} - 1\Big] \rightarrow 0.
  \end{aligned}
\end{equation}
The convergence rates in Assumption~\ref{assumption: heavy traffic} are
necessary for the results in Theorem~\ref{theorem: HT lower bound} and
Theorem~\ref{theorem: optimality of our policy} to come. We will show that the
predicted-class queues are GI/GI/1 and the first-order approximations of the 
arrival and service partial sum processes are linear (Section~\ref{subsection: appendix for arrival and service}), 
so the rates in~\eqref{eq: convergence rate} coincide with an $o(n^{-1/2})$ rate on the approximating processes. The classical G/G/1 framework
of~\citet{VanMieghem95} does not impose linearity, so the parameter rate no
longer determines the first-order approximations.
In Section~\ref{section: adding to the exposition in classical heavy traffic queueing analysis}, 
we provide the analogous condition on $\BtU_k^n, \BtV_k^n$ directly (per side, either
derivative convergence or $o(n^{-1/2})$-rate convergence) and use it to
complete~\citet[Propositions~3 and 7]{VanMieghem95}.

Similarly to~\citet{MandelbaumSt04}, we use an uniform integrability condition that allows us to apply the martingale functional limit theorem (FCLT). See Section~\ref{section: proof for fclt} for details.
\begin{assumption}[Uniform integrability]
  \label{assumption: second order moments}
  For any system $n\in \N$, we assume that
  \begin{enumerate}[(i)]
  \item $\E^n[(\tu_1^n)^2] < \infty$, $\E^n[(\tv_1^n)^2] < \infty$,
    $\E^n[(\feature_1^n)^2]< \infty$,  and for any $x \in \R$,
    \begin{equation*}
      \begin{aligned}
        & \E^n [(\tu_1^n)^2 \indic{\tu_1^n > x} ] \leq g_{\tu} (x),\quad
          \E^n [(\tv_1^n)^2 \indic{\tv_1^n > x} ] \leq g_{\tv} (x),
      \end{aligned}
    \end{equation*}
   where $g_u$ and $g_v$ are fixed functions such that $g_{\tu} (x) \rightarrow 0$, $g_{\tv}(x) \rightarrow 0$ as $x\rightarrow \infty$.
  \item There exist constants $\alpha_{\tu}\in (0, \infty)$ and $\alpha_{\tv, k} \in (0, \infty)$ for any $k\in [K]$ such that 
    \begin{equation*}
      \begin{aligned}
        & \alpha_{u}^n:=\E^n [(\tu_1^n)^2] 
          \rightarrow \alpha_{\tu},\quad
          \alpha_{v,k}^n:=\E^n [(\tv_1^n)^2| \tY_{1k}^n = 1] \rightarrow \alpha_{\tv, k}
      \end{aligned}
    \end{equation*}
    as $n\rightarrow \infty$.
  \end{enumerate}
\end{assumption}

Table~\ref{table: key notation} summarizes the notation of our model primitives.

\paragraph{Skorohod space and weak convergence.}
Let $\mathcal{C}$ be the space of continuous $[0, 1] \mapsto \mathbb{R}$
functions, $\mathcal{D}$ the set of the right-continuous with left limits (RCLL);
all stochastic processes will be RCLL.  Let $\mathcal{D}^k$ be its product
space and $\| \mathbf{x}(t)\| : = \max_{i\in [k]} | x_i(t)|$. Define
$d_{J_1}(\cdot, \cdot): \mathcal{D}\times \mathcal{D} \rightarrow
\mathbb{R}_+$ to be the $J_1$ (Skorohod) metric~\cite[Page 79]{Whitt02}. For
any vector-valued functions $\mathbf{x}(t), \mathbf{y}(t) \in \mathcal{D}^k$,
define
$d_p(\mathbf{x}, \mathbf{y}) = \sum_{i=1}^k d_{J_1}(x_i, y_i)$~\cite[Page
83]{Whitt02} and its topology $WJ_1$ (weak $J_1$ topology). Let $e$ be the identity function on $[0,1]$.

\begin{table}[t]
  \centering
  \small 
  \begin{tabular}{lp{11cm}} 
    \toprule
    Symbol & Description \\
    \midrule
    $n$ & System index in the heavy-traffic sequence \\
    $i$ & Job index within system $n$ \\
    $k \in [K]$ & True class index \\
    $l \in [K]$ & Predicted class index \\
    $\tu_i^n$ & Inter-arrival time of the $i\text{th}$ job arriving in system $n$ \\
    $\tv_i^n$ & Service time of the $i\text{th}$ job arriving in system $n$ \\
    $\feature_i^n$ & Feature vector of job $i$ in system $n$ \\
    $Y_i^n$ & One-hot vector of true class of job $i$ \\
    $\fY_i^n$ & One-hot vector of predicted class of job $i$ \\
    $\tp_k^n$ ($\tp_k$) & Prevalence of true class $k$ (limit) \\
    $\tmu_k^n$ ($\tmu_k$) & Service rate for true class $k$ (limit) \\
    $\tlambda^n$ ($\tlambda$) & Overall arrival rate (limit) \\
    $\q_{kl}^n$ ($\q_{kl}$) & Confusion matrix entry $\mathbb{P}(\fY_{1l}^n=1 \mid Y_{1k}^n=1)$ (limit) \\
    $\tA_0^n(t)$ & Total number of arrivals up to time $t\in [0, n]$ in system $n$ \\
    $\tU_0^n(t)$ & Arrival time of the $\lfloor t \rfloor\text{th}$ job in system $n$\\
    $\tV_0^n(t)$ & Total service time required by the first $\lfloor t \rfloor$ jobs in system $n$\\    
    \bottomrule
  \end{tabular}
  \caption{Key notation}
  \label{table: key notation}
\end{table}

\subsection{Discussion on Model Validity and Limitations}
\label{section: discussion model validity and limitations}

\paragraph{Model Validity}
We assume arrival and service rate are i.i.d., and that the scheduler has full
knowledge of the arrival and service rates, the confusion matrix, and the
convex cost function. While this assumption may not hold in certain scenarios,
such as cold-start queueing systems with limited data, it is a reasonable
assumption for applications like content moderation where the time scale of
the moderation is relatively short and abundant historical data is
typically available.

On large social media platforms such as Instagram and Facebook, historical
data on content---such as arrival times, content types, toxicity levels, and
the service time of human reviewers----often span months to years. On the other
hand, content moderation decisions are made on a much shorter time scale, such
as minutes to hours. This difference ensures the arrival and service processes
for content remain relatively stable during the short lifecycle of content
moderation. Consequently, we can reasonably model the arrival and service
rates as i.i.d. during the moderation period and the immediately preceding
timeframe. This stability allows us to use prior data to estimate key
parameters, such as the confusion matrix, arrival and service rates, and the
coefficients of the cost functions, by drawing on data from the recent past,
such as the previous few hours. For more on estimation techniques, refer to
this review paper:~\cite{AsanjaraniNaTa21}.

In imposing convex costs, we are implicitly focusing on a particular time scale
where significant delays are detrimental. This is a reasonable assumption in
content moderation problems where modeling the increasing marginal risk of
toxic contents going viral is first-order. For example, the terrorist attack
in two mosques in Christchurch was livestreamed on Facebook, and it quickly went
viral before the platform could intervene~\cite{TimbergHaShTrFu19}.  However, these costs should not
grow indefinitely in practice, although we believe the relatively short time
frame of content moderation problems makes convex costs a reasonable model for
real-world considerations. This assumption also aligns with a long line of
prior work in queueing~\cite{VanMieghem95, MandelbaumSt04} that consider
convex costs.

\paragraph{Model Limitations} While our model is a reasonable approximation to
interesting applications like content moderation, it is nevertheless not
without limitations. The assumption of i.i.d. arrival and service rates
typically holds over short time periods but can break under distribution
shifts. e.g., trending news can lead to spikes in arrival rates of related
contents. This increase may amplify further through social networks, rendering
hourly time-window estimates unreliable. Also, adversarial shifts,
such as coordinated efforts to flood the system with toxic content, can
increase the arrival rate of toxic contents, potentially overwhelming
moderation systems. These limitations are significant and we leave addressing
them as future work.

For classifiers, we assume the confusion matrix remains stable over short time
scales, but this assumption may also fail to hold under distribution
shifts. Additionally, estimating the confusion matrix from historical data may
introduce selection bias, especially if the historical data primarily reflects
those reviewed by human reviewers. This limitation can lead to estimation
errors but can be mitigated by incorporating diverse data sources, or applying
distributional robustness techniques to reduce bias in the classifier's
predictions.

\section{Lower bound on queueing cost}
\label{section: convergence and lower bound}

Our analysis relies on a diffusion limit for \emph{predicted classes} of the
model. Scheduling is based on predicted classes, but service times are determined by
the true classes. We characterize how misclassifications incur externalities
on other jobs, and derive the optimal queueing cost in Theorem~\ref{theorem:
  HT lower bound} to come. We begin by establishing the diffusion limits of the predicted-class arrival and service processes. In the presence of misclassifications, our setting requires a new method for handling diffusion scaling from the classical settings~\citep{VanMieghem95, MandelbaumSt04} that assume known job classes. 


\begin{table}[t]
  \centering
  \small 
  \begin{tabular}{ll} 
    \toprule
    Symbol & Description \\
    \midrule
    $\fA_l^n(t)$ & Number of arrivals whose predicted class is $l$ up to time $t\in [0,n]$ in system $n$\\
    $\fA_{kl}^n(t)$ & Number of arrivals of true-class-$k$ jobs that are predicted as class $l$ up to time $t$\\
    $\fS_l^n(t)$ & Service completions for predicted class $l$ up to cumulative service time $t$ \\
    $\fU_l^n(t)$ & Arrival time of the $\lfloor t \rfloor\text{th}$ job predicted as class $l$\\
    $\fV_l^n(t)$ & Cumulative service time of the first $\lfloor t \rfloor$ jobs predicted as class $l$\\
    $\fN_l^n(t)$ & Number of jobs with predicted class $l$ in the system at time $t$ \\
    $\fT_l^n(t)$ & Total service time devoted to predicted class $l$ up to time $t$ \\
    $\ftau_l^n (t)$ & Sojourn time of the last job predicted as class $l$ that remains in the system at time $t$\\
    $W_+^n(t)$ & Total remaining workload in the system at time $t$\\
    $\mathbf{X}^n$ & Diffusion-scaled, centered version of a process $X^n$ \\
    $C_k^n(\cdot)$, $C_k(\cdot)$ & Prelimit and limiting cost functions for true class $k$ \\
    $\fC_l^n (\cdot)$, $\fC_l(\cdot)$ & Prelimit and limiting predicted-class cost function defined in~\eqref{eq: pcmu cost} \\
    \bottomrule
  \end{tabular}
  \caption{Key notation on predicted-class processes and cost functions}
  \label{table: key processes}
\end{table}

\subsection{Convergence of predicted-class primitive processes}
\label{subsection: fundamental convergence results}

Define the counting processes for arrivals and service completions in the
predicted classes. 
For $l\in [K]$, let the $l$-th component of
$\VfA^n=(\fA_1^n, ..., \fA_K^n): [0,n] \rightarrow \mathbb{N}^K$ count the number of jobs that are
predicted as class $l$ until time $t\in [0,n]$, and for $k\in [K]$, let $\tA_{kl}^n (t)$ be the number of jobs that are of (true) class $k$ and predicted as class $l$.
Similarly, let
$\VfS^n=(\fS_1^n, ..., \fS_K^n): [0,n] \rightarrow \mathbb{N}^K$ count service completions as a
function of the total time that the server dedicates to each predicted
class.
We also define the associated partial sum processes. 
For $l\in [K]$, let the $l$-th component of $\fU^n=(\fU_1^n, ..., \fU_K^n): [0,n] \to \mathbb{R}_+^{K}$ give 
the arrival times of the jobs predicted as class $l$, and let $\fV^n=(\fV_1^n, ..., \fV_K^n): [0,n] \to \mathbb{R}_+^{K}$ denote the cumulative service times of predicted classes.
To streamline notation and ensure that all key processes are clearly
defined, Table~\ref{table: key processes} provides brief descriptions of the main stochastic processes
of the queueing system that are introduced in this section and used in our main results.
For ease of exposition, we defer a formal
discussion of diffusion limits to Section~\ref{section: proof for fclt} and defer
precise definitions to Section~\ref{subsection: approximation for W, T, N, and
  tau of the predicted classes}.

\paragraph{Feasible Policies} A scheduling policy $\policy_n$ is characterized
by an allocation process $\VfT^n=(\fT_1^n, ..., \fT_K^n): [0,n] \rightarrow \mathbb{R}^K$ whose $l$-th
coordinate denotes the total time dedicated to predicted class $l$ up to
$t\in [0, n]$. We use $\policy_n$ and $\VfT^n$ interchangably. Let
$\VfN^n=(\fN_1^n, ..., \fN_K^n): [0,n] \rightarrow \mathbb{N}^K$ be the queue length process; its
$l$-th coordinate denotes total jobs from predicted class $l$ remaining in
system at $t\in [0, n]$. Let $I^n(t): = t - \sum_l \fT^n_l(t)$ be the
cumulative idling time up to $t\in [0, n]$. The scheduler has full knowledge
of arrivals and the queue of \emph{predicted} classes.
\begin{definition}[Feasible Policies] \label{definition: feasible policies}
  The sequence of scheduling policies $\{\policy_n\}$ is feasible if the
  associated processes $\{\VfT^n(t), \VfN^n(t), I^n(t)\}$ satisfy for all
  $n\in \mathbb{N}$,
  \begin{enumerate}[(i)]
  \item $\VfT^n(0) = 0$, $\VfT^n$ is continuous and nondecreasing, $\VfN^n\geq 0$, and $I^n$ is nondecreasing;
  \item $\{\VfT^n(t), t\in [0, n]\}$ is adapted to the filtration $\sigma\{(\VfA^n (s), \VfN^n (s)): 0 \leq s < t\}$.
  \end{enumerate}
\end{definition}
\noindent
Condition (i) is natural, and condition (ii) ensures that $\{\pi_n\}$ only
relies on arrivals and queue lengths of predicted classes up to time $t$.  We
allow \textit{preemption} (preemptive-resume policy) so that the server can pause serving one job and switch to another in a \textit{different}
predicted class. Preemption is \textit{not} allowed between jobs from the same
predicted class, consistent with classical settings~\citep{MandelbaumSt04}.

\paragraph{Sojourn times and policy restrictions} Let $\ftau_{lj}^n$ be the
sojourn time of the $j$th job of predicted class $l$, and let
$\Vftau^n = \{\ftau^n_l\}_{l\in [K]}$ track the sojourn time of the most
recently arriving job in predicted class $l$, i.e.,
$\ftau^n_l(t) = \ftau^n_{l \fA_l^n(t)}$. Following the classical
settings~\citep{VanMieghem95, MandelbaumSt04}, we restrict attention to
\emph{p-FCFS} policies, which serve each predicted class in a
first-come-first-served manner, and to \emph{work-conserving} policies, under
which the server never idles while jobs are present. Any feasible policy is
stochastically dominated in cumulative queueing cost by a p-FCFS,
work-conserving counterpart (Lemmas~\ref{lemma: p fcfs} and~\ref{lemma: work
  conserving}), so we focus on such policies when establishing the lower bound.
We formalize the cost objective in Section~\ref{subsection: asymptotic lower
  bound of the cost functions}.

\paragraph{Diffusion limits of exogenous processes}
Our analysis begins by establishing diffusion limits of processes associated with predicted classes. To track the counts and cumulative service requirements of (mis)classified jobs, we define, for $k, l \in [K]$, $\fZ_{kl}^n (t)$ as the total number of true-class-$k$ jobs predicted as class $l$, and $\fR_l^n (t)$ as the total service requirement of jobs predicted as class $l$, among the first
$\lfloor t\rfloor$ jobs arriving in the system $n$:
\begin{equation*}
    \fZ_{kl}^n (t): = \sum_{i=1}^{\lfloor t \rfloor} \tY_{ik}^n \fY_{il}
^n, \quad \fR_{l}^n (t) := \sum_{i=1}^{\lfloor t \rfloor}\fY_{il}^n \tv_i^n, 
\quad t\in [0,n].
\end{equation*}
For $t\in [0, 1]$, let
$\TtU_0^n, \TtV_0^n, \VTfZ^n : = (\TfZ_{kl}^n)_{k, l \in [K]}, \VTfR^n:=
(\TfR^n_l)_{l\in [K]}$ be the diffusion-scaled processes, where
\begin{equation}\label{eq: def of tilde U and tilde V}
  \TtU_0^n (t) = n^{-1/2} [\tU_{0}^n (nt) -  
  (\tlambda^n)^{-1} \cdot nt], 
  \quad 
  \TtV_0^n (t) = n^{-1/2} [\tV_{0}^n (nt) -  \sum_{k=1}^n \frac{p^n_k}{\mu^n_k} \cdot nt], \quad t\in [0, 1];
\end{equation}
and formal definitions of $\VTfZ^n$ and $ \VTfR^n$ are deferred to
Definition~\ref{definition: U, Z, R, V}. Our analysis starts from the following joint weak convergence result, whose proof is provided in Section~\ref{section: proof for fclt}.
\begin{lemma}[Joint weak convergence] \label{lemma: joint weak convergence}
  Suppose that Assumptions~\ref{assumption: data generating
    process},~\ref{assumption: heavy traffic}, and~\ref{assumption: second
    order moments} hold.  Then, there exist Brownian motions
  $(\TtU_0, \VTfZ, \VTfR, \TtV_0)$ such that
\begin{equation*}
\begin{aligned}
  (\TtU_0^n, \VTfZ^n, \VTfR^n, \TtV_0^n) \Rightarrow (\TtU_0, \VTfZ, \VTfR, \TtV_0) \quad  \text{\quad in $(D^{K(K+1)+2}, WJ_1)$.} 
\end{aligned}
\end{equation*}
\end{lemma}

\paragraph{Sample path analysis}
Building off of our diffusion limit in Lemma~\ref{lemma: joint weak convergence}, we can strengthen the convergence to the 
\emph{uniform} topology using standard tools (e.g., see Lemma~\ref{lemma: Skorohod 
representation} and Lemma~\ref{lemma: equivalent to uniform convergence}), and 
conduct a \emph{sample path analysis} where we construct \emph{copies} of 
$(\TtU_0^n, \VTfZ^n, \VTfR^n, \TtV_0^n)$ and $(\TtU_0, \VTfZ, \VTfR, \TtV_0)$ that are identical in distribution 
to their original counterparts and converge almost surely under a common 
probability space.
With a slight abuse of notations, we use the same 
notation for the newly constructed processes.  

Sample path analysis lets us exploit uniform convergence, which substantially simplifies the arguments. All subsequent results and proofs in the Electronic Companions are established on the copied processes in a common probability space $(\Omega_\text{copy}, \mathcal{F}_\text{copy}, \mathbb{P}_\text{copy})$ with probability one, i.e., $\mathbb{P}_\text{copy} \text{-}a.s.$, and all convergence statements are taken in the \textit{uniform} norm $\|\cdot\|$. 
For example, Lemma~\ref{lemma: joint weak convergence} can be strengthened to $(\TtU_0^n, \VTfZ^n, \VTfR^n, \TtV_0^n) \rightarrow (\TtU_0,  \VTfZ, \VTfR, \TtV_0)$ in $(D^{K(K+1) + 2}, \| \cdot \|)$, $\mathbb{P}_{\text{copy}}$-a.s., as shown in Lemma~\ref{lemma: uniform convergence of UZRV}, and Lemma~\ref{lemma: convergence of the scaled primitive processes} shows convergence of the diffusion-scaled version of $\tA_0^n$ to $\TtA_0$ in $(\mathbb{D}, \|\cdot\|)$, $\mathbb{P}_{\text{copy}}$-a.s., 
where $\TtA_0$ is a function of $\TtU_0$. Since these copied processes are distributionally identical to the originals, almost-sure convergence results for the copied processes can be translated into weak convergence results for the originals; see Theorems~\ref{theorem: HT lower bound} 
and~\ref{theorem: optimality of our policy} for further discussion.

\paragraph{Convergence of the Endogenous processes}
Let $\VfW^n: [0,n] \to \mathbb{R}^K$ be the remaining workload process, with component $\VfW^n_l(t)$ denoting the service requirements of jobs---waiting or in service---predicted as class $l$ at time $t \in [0,n]$. Define the total remaining workload as $\sumtW^n(t) = \sum_{l} \fW_l^n(t)$, and let $\sumTtW^n$, $\VTfW^n$, $\VTfT$, $\VTftau^n$, and $\VTfN^n$ be the diffusion-scaled versions of $\sumtW^n$, $\VfW^n$, $\VfT$, $\Vftau^n$, and $\VfN^n$, respectively; see Section~\ref{subsection: approximation for W, T, N, and tau of the predicted classes}. We use $\sumtW^n$ without a tilde because, as shown in the following proposition, its limit $\sumTtW$ is invariant to both the classifier and the scheduling policy.

\begin{proposition}[Fundamental Convergence Results]
  \label{prop: convergence and approximation of predicted class N, tau, T, and W}
  Under Assumptions~\ref{assumption: data generating
    process},~\ref{assumption: heavy traffic}, and~\ref{assumption: second
    order moments}, and any work-conserving p-FCFS feasible policy
  \begin{enumerate}[(i)]
  \item (Invariant Convergence) 
    $\sumTtW^n \rightarrow \sumTtW := \phi\Big(\TtV_0 \circ \tlambda e +
    \sum_{k=1}^K \frac{\tp_k}{\tmu_k} \TtA_0 \Big)$, where $\phi$ is
    the reflection mapping as defined in \cite[Page 140, (2.5)]{Whitt02};
  \item (Equivalence of Convergence) For any predicted class $l\in [K]$,
    $\limsup_n \|\TfT_l^n\|$, $\limsup_n \|\TfN_l^n\|$,
    $\limsup_n \|\Tftau_l^n\|$, and $\limsup_n \|\TfW_l^n\|$ are all bounded. 
    If any of the processes
    $\TfT_l^n, \TfN_l^n, \Tftau_l^n,$ or $\TfW_l^n$ converges in 
    $(\mathcal{D}, \| \cdot \|)$, then all of $\TfT_l^n, \TfN_l^n, \Tftau_l^n,$ and $\TfW_l^n$ converge 
    in $(\mathcal{D}, \| \cdot \|)$. Moreover, if the limit of any one of these processes is continuous, then the limits of all four processes are continuous as well.
  \end{enumerate}           
\end{proposition}

Proposition~\ref{prop: convergence and approximation of predicted class N,
  tau, T, and W} extends the classical results of~\citet[Proposition
2]{VanMieghem95} by relaxing the assumption that true classes are known. When
true classes are known, convergence of the arrival and service processes of true
classes ($\VTtA^n$ and $\VTtS^n$) follows directly from the Functional Central
Limit Theorem (FCLT)~\cite[Assumption 1]{VanMieghem95}. Without this assumption,
establishing convergence of the diffusion-scaled arrival and service processes
of \textit{predicted classes} ($\VTfA^n$ and $\VTfS^n$) in
Proposition~\ref{prop: joint conv. of predicted class A, U, S, and V} demands a
different line of analysis. We exploit the joint convergence result in
Lemma~\ref{lemma: uniform convergence of UZRV} and characterize how
misclassification affects each subprocess, developing novel connections from the
primitives $\VTfZ^n$ and $\VTfR^n$ to $\VTfA^n$ and $\VTfS^n$ through random time
change and the continuous mapping theorem. We give the full proof in
Section~\ref{subsection: appendix for arrival and service}.

\subsection{Asymptotic lower bound of the scaled delay cost function}
\label{subsection: asymptotic lower bound of the cost functions}

The scheduler aims to minimize the cumulative queueing cost determined by the
true class labels. A true-class-$k$ job incurs queueing cost $\tC_k^n(\tau)$,
where $\tau$ is its sojourn time. Since the sojourn times $(\ftau_l^n)_{l\in
  [K]}$ are of order $n^{1/2}$ (Proposition~\ref{prop: convergence and
  approximation of predicted class N, tau, T, and W}), we impose a commensurate
scaling on the cost functions $\{C^n_k\}_{k\in[K]}$.
\begin{assumption}[Cost functions I] \label{assumption: on cost functions for
    showing the lower bound} For all $k\in [K]$, $\tC_k^n(\cdot)$ is
  differentiable, nondecreasing, and convex for all $n$.  There exists a
  continuously differentiable and strictly convex function $\tC_k$ with
  $\tC_k (0)=\tC_k'(0)=0$ such that $\tC_k^n(n^{1/2}\cdot) \rightarrow \tC_k(\cdot)$ and
  $n^{1/2}(C^n_k)'(n^{1/2}\cdot) \rightarrow C'_k(\cdot)$ uniformly on compact
  sets.
\end{assumption}
The scaled cumulative cost function incurred by $\policyn$ is
\begin{equation}
\label{eq: scaled cumulative cost function}
\begin{aligned}
    \TJ^n_\policyn (t; Q^n)
    =  n^{-1}\sum\limits_{l=1}^K \sum\limits_{k=1}^K
    \int_0^{nt}  C_k^n (\ftau_l^n (s))\mathrm{d}\fA_{kl}^n (s),
    ~\forall~t\in [0, 1],
\end{aligned}
\end{equation}
where $Q^n$ is the confusion matrix defined in Section~\ref{section: model}, and
$\mathrm{d}\fA_{kl}^n$ is the Lebesgue-Stieltjes measure induced by
$\fA_{kl}^n$. The cost $\TJ^n_\policyn (t; Q^n)$ depends on the scheduling policy
through the sojourn-time process $\{\ftau^n_l\}$.

We are now ready to present the main result of this section, the
asymptotic lower bound for the cumulative queueing cost in the heavy-traffic
limit. Our lower bound motivates the design of the \ourmethod~in
Section~\ref{section: heavy-traffic optimality of pcmu}. We let $\frho_l:=\sum_k \frac{\tlambda \tp_k \q_{kl}}{\mu_k} > 0$  (Assumption~\ref{assumption: heavy traffic}) for predicted class $l\in [K]$, and let $Q: = (\q_{kl})_{k,l\in [K]}$.

\begin{theorem}[Heavy-traffic lower bound] \label{theorem: HT lower bound}
Given a classifier $\model$ and a sequence of queueing systems, 
suppose that Assumptions~\ref{assumption: data
generating process},~\ref{assumption: heavy traffic},~\ref{assumption: second
order moments}, and~\ref{assumption: on
cost functions for showing the lower bound} hold.
Under any feasible scheduling policies
$\{\policyn\}$, the associated sequence of cumulative costs
$\{\TJ^n_\policyn (\cdot; Q^n): n \in \mathbb{N}\}$ satisfies
\begin{equation}\label{eq: definition of J*}
    \liminf _{n \rightarrow \infty}  \TJ^n_\policyn(t; Q^n)
    \geq  \TJ^*(t; Q)
    := \sum_{k=1}^K  \sum_{l=1}^K \int_0^t \tlambda \tp_k \q_{kl}
    C_k\Big(
    \frac{\big[h\big(\sumTtW(s)\big)\big]_l}{\frho_l}
    \Big)  \mathrm{d} s,~\forall t\in [0,1],
  \end{equation}
$\mathbb{P}_\text{copy}\text{-}a.s.,$ where $h(r)$ is an optimal solution to the following resource allocation
problem
\begin{equation}\label{eq: optimization problem}
\begin{aligned}
  \opt(r) := &\min_x 
  && \sum\limits_{l=1}^K \sum\limits_{k=1}^K \lambda p_k \q_{kl} C_k \Big(\frac{x_l}{\frho_l}\Big)\\
             &\text{s.t.} && \sum\limits_{l=1}^K x_l = r,\quad  x_l \geq 0,~\forall~l\in [K].
\end{aligned}
\end{equation}
Moreover, for the original processes under $\mathbb{P}^n$, under any feasible 
policies $\{\pi'_n\}$,
\begin{equation} \label{eq: HT lower bound in stochastic sense}
\liminf _{n \rightarrow \infty}  \mathbb{P}^n [\TJ^n_{\pi'_n} (t; Q^n) > x] 
\geq \mathbb{P}_\text{copy}[\TJ^* (t;Q) > x],~\forall~x\in \mathbb{R},~\forall~t\in [0,1].
\end{equation} 
\end{theorem}  
\noindent
According to Proposition~\ref{prop: convergence and approximation of predicted
  class N, tau, T, and W}, $\sumTtW$ is solely determined by the exogenous
processes $\TtA_0$ and $\TtV_0$. 
Consequently, the lower bound in Theorem~\ref{theorem: HT lower bound} depends on the prediction model $\model$ only through the optimally aligned workloads $h(\sumTtW)$.
In particular,~\eqref{eq: optimization problem} reveals that the predictive performance of the model, encoded by the confusion matrix $Q$, shapes the optimal cost through the normalized workload allocation across predicted classes
via the \emph{induced} intensities $(\frho_l)_{l\in [K]}$ and arrival rates $(\sum_k \lambda p_k \q_{kl})_{l\in [K]}$. In this sense,~\eqref{eq: optimization problem} makes explicit how misclassification errors ripple through the queueing system and determine the best achievable performance under any scheduling policy,
thereby providing a queueing-based benchmark for evaluating the model
$\model$ from an operations perspective.
Our proof is involved and deferred to
Section~\ref{section: proof for HT lower bound}, where we also contrast our
analytic approach to classical proof techniques.

\paragraph{Characterization of the optimal workload allocations}
For the convex optimization problem~\eqref{eq: optimization problem} with linear constraints,
its KKT conditions characterize the optimal workload allocation $h$. For
predicted class $l$, define the
\emph{P$c\mu$ cost function} $\fC_l$ as the weighted average of the true-class
costs $\tC_k$, with weights given by the posterior probability that a job
predicted as class $l$ is truly of class $k$,
\begin{equation}\label{eq: pcmu cost}
  \fC_l(t): = \sum_{k=1}^K\frac{ p_k \q_{kl}  }{\sum_{k'=1}^K p_{k'} \q_{k'l}} \cdot \tC_k(t),~~~ t \in [0, \infty),
\end{equation}
where $p_k$ is the prevalence of true class $k$ and $\q_{kl}$ is the probability
that a class-$k$ job is predicted as class $l$. The limiting service rate
$\fmu_l$ is given by the expected service time of a job predicted as class $l$,
$\fmu_l^{-1} = \sum_{k=1}^K \frac{ p_k \q_{kl} }{\sum_{k'=1}^K p_{k'} \q_{k'l}}
\cdot \frac{1}{\mu_k}$, a weighted average of the true-class expected service
times $1/\mu_k$ under the same posterior weights as $\fC_l$
(see Definition~\ref{def: concerned processes for arrival and service}).
\begin{lemma}[KKT conditions]
\label{lemma: KKT conditions}
  $\{x_l\}_{l\in[K]}$ is an optimal solution for $\opt(\sumTtW(t))$ if
  $x_l>0, ~\forall~l\in[K]$ and is a solution to
  \begin{equation}\label{eq: KKT artificial classes}
    \begin{aligned}
      \fmu_l \fC_l' \Big(\frac{x_l}{\frho_l}\Big)=\fmu_m \fC_m' \Big(\frac{x_m}{\frho_m}\Big),~\forall~l,m \in [K],\quad 
      \sum\limits_{l=1}^K x_l = \sumTtW (t).
    \end{aligned}
  \end{equation} 
\end{lemma}
\noindent
We also show that the KKT conditions~\eqref{eq: KKT artificial classes} have a
unique solution (Lemma~\ref{prop: properties of h and opt},
Section~\ref{section: proof for HT lower bound}) and thus $h(\sumTtW(t))$ is
well-defined.

\subsubsection{Illustrative two-class example}

To illustrate Theorem~\ref{theorem: HT lower bound} and the role of the
optimization problem~\eqref{eq: optimization problem}, we consider a simple
two-class setting ($K=2$) with limiting cost functions satisfying
Assumption~\ref{assumption: on cost functions for showing the lower bound}.
Let the overall arrival rate be $\lambda = 1$, with class prevalences
$p_1 = p$ and $p_2 = 1-p$, and service rates $\mu_1 > \mu_2$ so that class~2
jobs are harder to review. We assume per-class delay costs of the form
$  C_k(z) = \frac{1}{2} c_k z^2,~c_k > 0, \ k=1,2$, with $c_2 > c_1$ to reflect that delays for class~2 (e.g., harmful political misinformation) are more costly than for class~1 (e.g., benign content). In the two-class case, Lemma~\ref{lemma: KKT conditions} yields a closed-form
expression for the optimal predicted-class workloads:
\begin{equation}\label{eq: two-class-x-star}
  a_l(Q)
  := \frac{\sum_{k=1}^2 \lambda p_k q_{kl} c_k}
           {\big(\sum_{k=1}^2 \frac{\lambda p_k q_{kl}}{\mu_k}\big)^2},
  \quad
  x_1^*(r;Q) = r\,\frac{a_2(Q)}{a_1(Q)+a_2(Q)},\quad
  x_2^*(r;Q) = r\,\frac{a_1(Q)}{a_1(Q)+a_2(Q)}.
\end{equation}
The coefficient $a_l(Q)$ is the curvature of predicted class $l$'s cost in its
allocated workload: the per-class cost is $\tfrac12 a_l(Q)\,x_l^2$, so the
marginal cost of loading class $l$ is $a_l(Q)\,x_l$. The KKT conditions equalize
this marginal cost across classes, giving $x_l^*(r;Q)\propto 1/a_l(Q)$, so the
optimizer allocates workload inversely to curvature and makes a class a smaller
target precisely when its $a_l(Q)$ is large. Two competing forces shape
$a_l(Q)=\big(\sum_k \lambda p_k q_{kl} c_k\big)/\frho_l^2$: the cost-weighted
arrival rate in the numerator, and the predicted-class intensity $\frho_l$ in the
denominator, which enters squared because $x_l/\frho_l$ is the induced sojourn
time and the delay cost is quadratic. When $Q = I$, $a_l(Q)$ reduces to
$c_l \mu_l^2 / p_l$, recovering the classical dependence of optimal workloads on
costs and service rates.

\begin{figure}[t]
  \centering
  \begin{minipage}[t]{0.47\textwidth}
    \centering
    \includegraphics[
    width=\linewidth,
    height=0.7\linewidth,
  	keepaspectratio
    ]{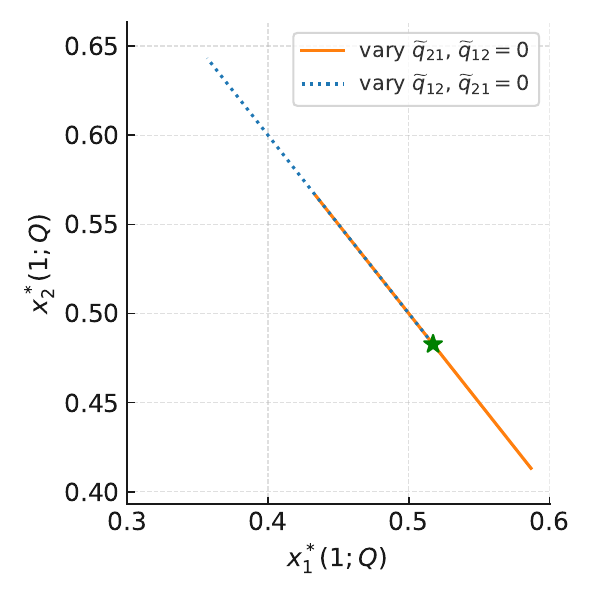}
    \vspace{0.3em}
    \textbf{(a)} Trajectories of $(x_1^*(1;Q),x_2^*(1;Q))$ as
    $\q_{21}$ (solid) or $\q_{12}$ (dotted) increase,
    with the green star marking $Q=I$.
  \end{minipage}
  \hfill
  \begin{minipage}[t]{0.47\textwidth}
    \centering
    \includegraphics[
    width=0.95\linewidth,
    height=0.7\linewidth,
  	keepaspectratio
    ]{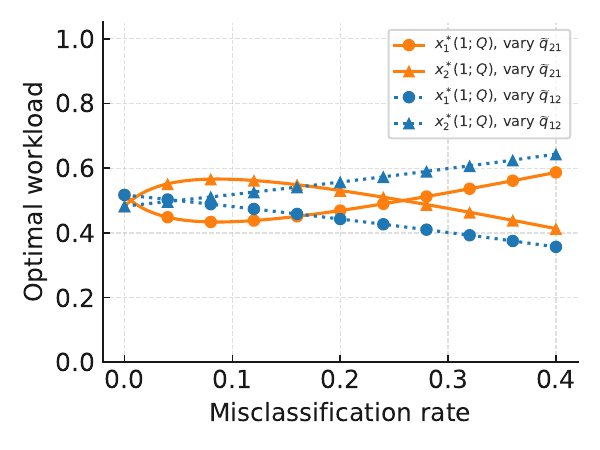}
    \vspace{0.3em}
    \textbf{(b)} Optimal workloads $x_1^*(1;Q)$ and $x_2^*(1;Q)$ as functions
    of the misclassification rate (solid: vary $\q_{21}$, $\q_{12}=0$;
    dotted: vary $\q_{12}$, $\q_{21}=0$).
  \end{minipage}
  \caption{Optimal predicted-class workloads in the two-class example under
  different misclassification patterns.}
  \label{fig:two-class-workload-geometry}
\end{figure}

To visualize the impact of misclassification, we consider a numerical example
with $\lambda = 1$, $p_1 = 0.3$, $p_2 = 0.7$, $\mu_1 = 2$, $\mu_2 = 1$, and
cost weights $c_1 = 1$, $c_2 = 10$, and examine the optimizer
$(x_1^*(r;Q),x_2^*(r;Q))$ at a representative
workload level $r = 1$.
Figure~\ref{fig:two-class-workload-geometry}(a) plots $(x_1^*(1;Q),x_2^*(1;Q))$ as
we vary $\q_{12}~(\q_{21})$ with $\q_{21}~(\q_{12})=0$, where the workload under $Q=I$ is marked by a green star, and Figure~\ref{fig:two-class-workload-geometry}(b) shows the individual trajectories of $x_1^*(1;Q)$ and $x_2^*(1;Q)$, with the dotted curves for the $\q_{12}$ sweep behaving symmetrically.

Figure~\ref{fig:two-class-workload-geometry}(a) also explains why the optimal
predicted-class workloads need not vary monotonically with the misclassification
rates. Since $x_l^*(1;Q)\propto 1/a_l(Q)$, the split is governed by how the two
curvatures move relative to each other as $\widetilde{q}_{21}$ sweeps from $0$ to
$1$, and they move in opposite phases. At first, a few costly, slow harmful jobs
mix into the small predicted-class-1 stream and raise its cost-weighted arrival
rate faster than its intensity $\frho_1$, so $a_1(Q)$ rises, overtakes $a_2(Q)$,
and $x_1^*(1;Q)$ decreases. As $\widetilde{q}_{21}$ grows further, the
predicted-class-2 stream thins out and its intensity $\frho_2$ falls toward zero. Since
$\frho_2$ enters $a_2(Q)$ squared, the curvature diverges and holding any
\emph{fixed} workload $x_2$ in this near-empty stream incurs an exploding quadratic delay
cost. The optimizer therefore drives $x_2^*(1;Q)$ toward zero and shifts
essentially all workload back onto predicted class~1. This crossover of the two
curvatures illustrates how misclassification errors can have a subtle,
endogenous effect on the optimal workload allocation even in a simple two-class
setting.
This perspective motivates the
empirical studies in Sections~\ref{section: numerical experiments}
and~\ref{section: characterization of the optimal cost with q} for tuning and selecting prediction models with the operational performance as a central criterion.

\paragraph{Comparison to the classical lower bound result}
When the job classifier is ``perfect'', i.e., $Q^n = I$, the setting reduces to the classical case in which the true classes are known. Even in this special case, our proof contributes to the classical heavy traffic theory by filling in the
missing arguments in \citet{VanMieghem95}'s original result for \textit{G/G/1}
systems. We identify sufficient conditions under which their original claims hold, and
verify that these conditions indeed hold for GI/GI/1 systems. As one example
(gap A1 in Table~\ref{table: classical gaps} in Section~\ref{section: heavy-traffic optimality of pcmu}), our
proof of Lemma~\ref{lemma: relation between N and W} demonstrates how i.i.d.\ service times
satisfy one of the conditions and establishes a key relationship between queue length and workload
processes~\cite[Proposition 3]{VanMieghem95}. See Sections~\ref{subsubsection:
  detailed comparison to the lower bound result in VM} and~\ref{section: sufficient conditions for the lower bound proof in Van Mieghem} for details.
Table~\ref{table: classical gaps} points to two other missing gaps in the
classical lower bound proof that we close analogously.


\section{Heavy-traffic optimality of the \ourmethod}
\label{section: heavy-traffic optimality of pcmu}

We are ready to formally derive the \ourmethod, which is motivated by the
convex optimization problem~\eqref{eq: optimization problem}.  We prove heavy
traffic optimality of the \ourmethod~by showing that it attains the lower
bound in Theorem~\ref{theorem: HT lower bound}.

\subsection{Optimality of the \ourmethod}
\label{subsection: optimality of the pcmu via kkt}

We first characterize the limiting cumulative cost of a convergent policy. Let
let $\BfA_{kl}$ be the limit of $n^{-1}\fA^n_{kl}(n \cdot)$ (see Definition~\ref{def:
  concerned processes for arrival and service} for a formal statement). In the
following, $\TJ_\policy(t; Q)$ is dependent on
$\VTftau = \{\Tftau_l\}_{l\in [K]}$ through the subscript $\policy$.
\begin{lemma}[Convergence of
$\TJ^n_\policyn(\cdot; Q^n)$]\label{lemma: convergence of tildeJ}
Given a classifier $\model$, suppose that Assumption~\ref{assumption: data
  generating process},~\ref{assumption: heavy traffic},~\ref{assumption: second
  order moments}, and~\ref{assumption: on
  cost functions for showing the lower bound} hold. For feasible policies $\{\policyn\}$ satisfying
$\Tftau^n_l \rightarrow \Tftau_l,~\forall~l\in [K]$,
\begin{equation}\label{eq: convergence of tildeJ}
    \sup\limits_{t\in [0,1]} |\tilde{J}^n_\policyn (t; Q^n) -\TJ_\policy(t; Q)|
    \rightarrow 0,
\end{equation}
where the limiting cumulative cost $\TJ_\policy(t; Q)$ is defined by
\begin{equation*}
  \TJ_\policy(t; Q): = \sum\limits_{l=1}^K \sum\limits_{k=1}^K 
  \int_0^{t}  C_k (\Tftau_l (s))\mathrm{d}\BfA_{kl} (s) 
  =  \sum\limits_{l=1}^K \sum\limits_{k=1}^K 
  \int_0^{t}  \lambda p_k \q_{kl}C_k (\Tftau_l (s))~\mathrm{d} s.
\end{equation*}
\end{lemma}
\noindent See Section~\ref{subsection:proof-convergence-tildeJ} for the proof.

Combining our characterization of the cumulative cost with the lower bound in
Theorem~\ref{theorem: HT lower bound}, we conclude that $\{\policyn\}$ is
asymptotically optimal if the following conditions are satisfied: i) the
scaled sojourn time processes converge, i.e.,
$\Tftau_l^n \rightarrow \Tftau_l,~\forall~l\in[K]$, and ii) the limiting
sojourn time processes satisfy
$\Tftau_l(t) = [h(\sumTtW(t))]_l / \frho_l,~\forall~t\in[0, 1],~l\in [K]$,
where $h(\cdot)$ is an optimal solution to the optimization problem~\eqref{eq:
  optimization problem}. Recalling $\opt(\sumTtW(t))$, the optimization
problem~\eqref{eq: optimization problem}, is convex with linear constraints,
its KKT conditions characterize the optimal workload allocation $h$. For
predicted class $l$, recall its limiting service rate $\fmu_l$ and the P$c\mu$
cost function~\eqref{eq: pcmu cost}.
By Lemma~\ref{lemma: KKT conditions}, the optimal allocation $h(\sumTtW(t))$
is the unique solution of the KKT conditions~\eqref{eq: KKT artificial classes}
(Proposition~\ref{prop: properties of h and opt}, Section~\ref{section: proof for HT lower bound}),
and is therefore well-defined.

The cost function $\fC_l(t)$~\eqref{eq: pcmu cost} arises from the KKT
conditions of $\opt(\sumTtW(t))$ as a weighted average with weights
proportional to $p_k \q_{kl}$, reflecting how predicted class $l$ is composed of jobs from different true classes. As $p_k$ and $\q_{kl}$ rely on the
arrival rates and misclassification errors, $\fC_l(t)$ can be viewed as the
exogenous average cost function associated with predicted class $l$. We aim to develop a scheduling policy that induces the workload allocation to
align with the exogenous cost $\fC_l(t)$, in the sense that the
conditions~\eqref{eq: KKT artificial classes} are satisfied for all
$t\in[0, 1]$.

According to Proposition~\ref{prop: convergence and approximation of predicted
  class N, tau, T, and W}, convergence of the sojourn time process
$\VTftau^n \rightarrow \VTftau$ is equivalent to convergence of workload
$\VTfW^n \rightarrow \VTfW$. Moreover, if $\VTfW^n$ converges, then
$\Tftau_l = \TfW_l / \frho_l,\forall~l \in [K]$ (see Lemma~\ref{prop: little's law} and Lemma~\ref{lemma: relation between N and W}). Consequently, our goal is to develop a policy that satisfies
$\VTftau^n \rightarrow \VTftau$ and
\begin{equation}\label{eq: KKT condition for tau}
  \fmu_l \fC'_l(\Tftau_l) = \fmu_m \fC'_m(\Tftau_m),~\forall~l, m \in [K],
\end{equation}
in the heavy traffic limit. When the balance~\eqref{eq: KKT condition for tau}
is achieved, the limiting workload allocation $\TfW_l=x_l := \Tftau_l \frho_l$
satisfies the KKT conditions~\eqref{eq: KKT artificial classes} and both
conditions (i) and (ii) are met, which leads to the policy's optimality.

Since the sojourn time---time between job arrival and service completion---is
\emph{not} observable, we substitute $\Vftau^n=\{\ftau_l^n\}_{l\in [K]}$ with
the observable \textit{age} processes.
\begin{definition}[Age Process]\label{definition: age process}
  Given a classifier $\model$ and feasible policies $\{\policyn\}$, a
  predicted class $l$ and time $t\in[0, n]$, let $\fa_l^n (t)$ be the waiting
  time of the oldest job in predicted class $l$ at time $t$, where a job being
  served is defined to be {\it waiting in the system}.  Let $\fa_l^n$ be the
  \emph{age} process of the predicted class $l\in [K]$ in system $n$, and let
  $\Tfa_l^n (t):=n^{-1/2} \fa_l^n (nt),~t\in [0,1]$ be the corresponding
  diffusion-scaled process.
\end{definition}
\noindent If either $\{\Tfa_l^n\}_{l\in [K]}$ or $\{\Tftau_l^n\}_{l\in [K]}$
converges, then both of the processes converge to the same limit, i.e.,
$\Tftau_l (t) = \Tfa_l (t),~\forall~l\in[K], t\in[0, 1]$ (see
Proposition~\ref{proposition: relation between a and tau}).  Thus, we can
equivalently reformulate the optimality condition for sojourn time~\eqref{eq:
  KKT condition for tau} into that with observable age processes
\begin{equation}\label{eq: KKT condition for age}
    \fmu_l \fC'_l(\Tfa_l) = \fmu_m \fC'_m(\Tfa_m),~\forall~l, m \in [K].
  \end{equation}

\paragraph{Heavy-Traffic Optimality}

In the heavy-traffic limit, the \ourmethod~has the same intuitive form as the
oracle \gcmu:
it greedily serves the predicted class with the highest marginal cost of delay,
but measures that cost with the confusion-weighted P$c\mu$ cost
$\fC_l$~\eqref{eq: pcmu cost} rather than a true- or predicted-class cost,
\begin{equation}
  \label{eq: Pcmu rule}
  \argmax_{l\in [K]} \fmu_l (\fC_l)'(\fa_l(t)) \qquad \qquad \mbox{\ourmethod}.
\end{equation}
We design the \ourmethod~in the \emph{prelimit} systems to achieve~\eqref{eq: KKT
  condition for age} in the heavy traffic limit; it prioritizes predicted classes
with the highest prelimit \ourpolicy~index, defined as follows.
\begin{definition}[\ourmethod]\label{definition: modified gcmu rule}
  Given a classifier $\model$, for any system $n$ at time $nt$ with $t\in [0, 1]$,
  the \ourmethod~serves the oldest job in the predicted class having the maximum
  \ourmethod~index, i.e., $l \in \arg \max_{m\in [K]} \fIpre^n_m(t)$, with
  preemption, where
  \begin{equation}\label{eq: gcmu index prelimit}
    \fIpre^n_l(t):=  \fmu^n_l\cdot n^{1/2} (\fC^n_l)'(\fa_l^n (nt)),~\forall~t\in [0,1],
  \end{equation}
  is the \ourmethod~index for predicted class $l$ at time $nt$ in system $n$,
  and
  $\fC_l^n(t): = \frac{\sum_k p^n_k \q^n_{kl} \tC^n_k(t) }{\sum_{k'} p^n_{k'}
    \q^n_{k'l}},~t\in [0, \infty),~l\in [K]$, is the weighted average of
  $C^n_k$ and the prelimit counterpart of $\fC_l(t)$.
\end{definition}
While this rule requires more information than the usual \gcmu, the arrival
rates and misclassification probabilities can be efficiently estimated on past
observations.

The \ourmethod~is a work-conserving p-FCFS policy by definition, and the $n^{1/2}$ scaling ensures
a well-defined heavy traffic limit. The \ourmethod~naturally allows for
\emph{preemption}: since we consider jobs being served as waiting in the
system, the age process $\ta_l^n$ corresponds to the same job waiting in the
queue until its service completion. 
We adopt preemption for analysis purposes. In particular, we can develop a 
non-preemptive counterpart of the \ourmethod~and show its optimality using the
same analytic framework.


We are now ready to present our optimality result, which shows that the
cumulative queueing cost associated with the \ourmethod,
$\TJ^n_\ourpolicy(\cdot; Q^n)$, converges to the asymptotic lower bound
$\TJ^* (\cdot;Q)$. Our proof relies on the fact that the \ourmethod~is a greedy
method minimizing the largest difference of the P$c\mu$ indices
, $\sup_{t\in [0,1]} \, \max_{l, m \in [K]} |\fIpre^n_l(t) - \fIpre^n_m(t)|$,
which guarantees (Proposition~\ref{prop: convergence of max gcmu difference},
Section~\ref{subsection: proof overview of the optimality})
\begin{equation}
  \label{eq: convergence of max difference of pcmu}
  \sup_{t\in [0,1]}  \, \max_{l, m \in [K]}  |\fIpre^n_l(t) - \fIpre^n_m(t)| \rightarrow 0.
\end{equation} 
We develop novel analysis techniques to show the convergence~\eqref{eq:
  convergence of max difference of pcmu}, which requires strong convexity of
the cost function.
\begin{assumption}[Cost functions II]\label{assumption: on cost functions for optimality}
  The limiting cost $C_k$ is strongly convex for all $k\in [K]$.
\end{assumption}
\begin{theorem}[Optimality of \ourmethod]
\label{theorem: optimality of our policy}
Given a classifier $\model$ and a sequence of queueing systems, suppose that
Assumptions~\ref{assumption: data generating process},~\ref{assumption: heavy traffic},~\ref{assumption: second order moments},~\ref{assumption: on cost
  functions for showing the lower bound}, and~\ref{assumption: on cost functions
  for optimality} hold. Then, 
$\TJ^n_\ourpolicy (\cdot;Q^n) \rightarrow 
\TJ^* (\cdot;Q)$ in $(\mathcal{D}, \| \cdot \|)$ $\mathbb{P}_\text{copy}$-a.s..
For the original processes under $\mathbb{P}^n$, 
$\TJ^n_\ourpolicy (\cdot;Q^n) \Rightarrow \TJ^* (\cdot;Q)$ 
in $(\mathcal{D}, J_1)$, and in particular,
$\mathbb{P}^n [\TJ^n (t;Q^n) > x] \rightarrow 
\mathbb{P}_{\text{copy}} [\TJ^* (t;Q) > x],~\forall~x\in \mathbb{R},~t\in [0,1].$
\end{theorem}
\noindent Our proof is highly involved so we provide a brief overview in
Section~\ref{subsection: proof overview of the optimality} and defer detailed
arguments to Section~\ref{section: proof for optimality of our policy} and
~\ref{section: proof of convergence of max gcmu
  difference}.

\subsection{Completing the classical optimality proofs}
\label{subsection: completing the classical optimality proofs}

While not the main contribution of this work, our analytic framework adds to
the standard heavy traffic analysis~\citep{VanMieghem95, MandelbaumSt04} when
specialized to the classical setting of known true classes.
We identify and resolve two issues in the classical analysis.
\begin{enumerate}
\item The optimality conditions in the previous section depended on the
  \emph{sojourn times} that are only revealed when a job exits the system. On
  the other hand, the \gcmu~\citep{VanMieghem95} and
  $D\text{-}Gc\mu$~\citep{MandelbaumSt04}~use indices based on the \emph{ages}
  of waiting jobs, but their mathematical arguments for optimality only apply
  to an unimplementable version in which the indices are defined by workloads
  or sojourn times.
\item Moreover, the argument in~\citet{VanMieghem95} \textit{states without
    proof} a critical property, that the differences between the
  (unimplementable) class indices vanish, and does not establish it for the
  implementable age-based indices. 
\end{enumerate}
The gap between the proposed indices and the mathematical optimality
guarantees that the formal arguments show, as well as the absence of a proof
for the vanishing difference between the indices~\citep{VanMieghem95} indicate
the incompleteness of the classical proofs for the optimality
in~\citet{VanMieghem95,MandelbaumSt04}, which we summarize in
Section~\ref{subsection: comparison to VM and MS} and detail in
Sections~\ref{subsection: comparison to the optimality result in
  VM},~\ref{subsection: comparison to the optimality result in MS},
and~\ref{section: adding to the exposition in classical heavy traffic queueing
  analysis}.

Concretely, the proofs in Sections~\ref{section: proof for HT lower
  bound},~\ref{section: proof for optimality of our policy}, ~\ref{section:
  proof of convergence of max gcmu difference}, and~\ref{section: adding to
  the exposition in classical heavy traffic queueing analysis} close five gaps
in the classical analysis of~\citet{VanMieghem95} listed below.  Three are
expository gaps in the lower bound result (A1--A3), supplying justifications
the original proofs omit, though we establish A3 only for the reflected
Brownian motion workload $\sumTtW$ of this paper and leave the general
reflected-process case open. The other two are substantive steps in the
optimality argument itself that we close directly: B1, the most critical,
establishes the index convergence on which the optimality proof rests, and B2
reconciles the index as defined with the quantity actually analyzed. For B2,
we establish the asymptotic equivalence of the age and sojourn time processes
in our model, and use it to give sufficient conditions for the same
equivalence in~\citet{VanMieghem95,MandelbaumSt04}, alongside the analogous
condition for A1 that completes the queue-to-workload translation
of~\citet[Proposition~3]{VanMieghem95} for the general G/G/1 systems studied
there. The matching subsection headings in the Electronic Companions carry
these tags for cross-reference.
\begin{table}[ht]
\centering
\small
\caption{Gaps in the classical analysis of~\cite{VanMieghem95,MandelbaumSt04} and where we address them}
\label{table: classical gaps}
\begin{tabular}{lll} 
\toprule
Gap & Where in \cite{VanMieghem95,MandelbaumSt04} & Addressed in \\
\midrule
A1 & \cite[Prop.~3 (Eq.~36)]{VanMieghem95}: queue-to-workload translation & Sections~\ref{section:proof-relation between N and W}, \ref{section: sufficient conditions for the lower bound proof in Van Mieghem} \\
\addlinespace 
A2 & \cite[Prop.~6 (Eq.~95)]{VanMieghem95}: continuity of the optimal allocation & Section~\ref{subsection: continuity of the optimal workload allocation} \\
\addlinespace
A3 & \cite[Prop.~6]{VanMieghem95}: vanishing mesh of the workload time-partition & Section~\ref{subsubsection: complementary proof for the partition size in VM} \\
\addlinespace
B1 & \cite[Prop.~7]{VanMieghem95}: max index-difference convergence stated without proof & Section~\ref{section: proof of convergence of max gcmu difference} \\
\addlinespace
B2 & $Gc\mu$~\cite{VanMieghem95}, $D\text{-}Gc\mu$~\cite{MandelbaumSt04}: defined on ages, analyzed on sojourn times & Sections~\ref{subsection: proof of the relationship betwen age and sojourn time}, \ref{section: sufficient conditions for the optimality proof in Van Mieghem} \\
\bottomrule
\end{tabular}
\end{table}

\noindent For B1, we provide a rigorous proof for the convergence of the age
processes in Proposition~\ref{prop: convergence of max gcmu difference}, using a
novel and involved induction argument that analyzes the dynamics of the indices.
As we overview in Section~\ref{subsection: proof overview of the optimality},
this proof simultaneously reveals a previously unstated sufficient condition
(strong convexity of the cost functions) while relaxing regularity assumptions
made in~\citep{MandelbaumSt04}.


\section{Empirical demonstration of the P$c\mu$-rule}
\label{section: numerical experiments} 
We demonstrate the effectiveness of \ourmethod~on a content moderation
problem using real-world user-generated text comments with the data generating
process in Section~\ref{section: model}. To operate at a massive scale, online
platforms use AI models to provide initial toxicity predictions. However,
these models are imperfect due to the inherent nonstationarity in the system.
For example, they cannot reliably detect context related to hate speech
following a recent terrorist attack. As a result, platforms must rely on human
reviewers as the final inspectors~\citep{MakhijaniShAvGoStMe21}, especially
since they bear the cost of mistakenly removing non-violating comments. Our
goal is to analyze the downstream impact of prediction errors on scheduling
decisions in the content moderation queueing system.

Different comments incur varying levels of negative impact on the
platform. If not removed in a timely manner, toxic comments attacking
historically marginalized or oppressed groups can have particularly harmful
effects. We model this using heterogeneous delay costs based on the level of toxicity and the demographic group targeted by the comment. These factors also affect processing time. For instance, reviewing comments about an ethnic minority group in a foreign nation is more challenging and time-consuming compared to domestic content.


We use real user-generated text comments on online articles from the
CivilComments dataset~\citep{BorkanDiSoThVa19}. Each comment has been labeled
by at least ten crowdsourcing workers with binary toxicity labels and whether
it mentions one of the 8 demographic identities: \emph{male, female, LGBTQ,
  Christian, Muslim, other religions, Black, White}.  For simplicity, we focus
on comments that mention one and only one of the common groups \emph{white,
  black, male, female, LGBTQ}.  By crossing them with binary toxicity labels,
we derive 10 job classes.  We assume the system has exact knowledge of target
group (using simple rule-based logic), but can only predict the toxicity
through an AI model.

The toxicity predictor, which can also be viewed as the job class predictor
$f_\theta$, utilizes the same neural network architecture and training
approach as described in~\citet{KohSaEtAl20}.  To showcase the versatility of
our scheduling algorithm regardless of the underlying prediction model, we
study three models fine-tuned based on a pre-trained language model
(DistilBERT-base-uncased~\citep{SanhDeChWo19}): empirical risk minimization
(ERM), reweighted ERM that upsamples toxic comments, and a simple
distributionally robust model trained to optimize worst-group performance over
target demographic groups (GroupDRO~\citep{SagawaKoHaLi19}). We observe
significant variation in predictive performance across the 10 job classes
defined by $\{$toxicity $\times$ target demographic$\}$. Across the three
models (ERM, Reweighted, GroupDRO), the worst-class accuracy (55\%, 68\%, 67\%) is
significantly lower than the mean accuracy (88\%, 84\%, 84\%), leading to diverse
patterns in the confusion matrix $Q$.
In particular, for each of the 5 identity groups, the confusion matrix $Q$
has the following entries:
\[
q_{\text{toxic, toxic}} = \text{TPR}, \quad
q_{\text{toxic, non-toxic}} = \text{FNR}, \quad
q_{\text{non-toxic, toxic}} = \text{FPR}, \quad
q_{\text{non-toxic, non-toxic}} = \text{TNR},
\]
resulting in a $10 \times 10$ confusion matrix across the 5 groups $\times$ 2 classes. We assume the scheduler can estimate $Q$ using a validation dataset. See Section~\ref{appendix: experiment details} for details.

\paragraph{Queueing system} We assume jobs are assigned to reviewers randomly
to ensure fairness, as mentioned in Section~\ref{section:introduction}, and view each reviewer as a single-server queueing system.  For
simplicity, we consider a queueing model operating in a finite time interval
$[0,1]$ with 10 job classes. New jobs/comments arrive with i.i.d.  exponential
interarrival times with rate 100 (uniformly drawn from the test set). Toxic
comments have a lower service rate and toxic comments mentioning minority
groups have an even lower service rate. The service times follow exponentially
distributions that solely depend on the true class label $Y$: for \emph{white,
  black, male, female, LGBTQ}, respectively,
$\mu_{\text{toxic}} = [100, 30, 110, 25, 15]$ and
$\mu_{\text{non-toxic}} = [150, 150, 150, 150, 150]$. (If the service rate
depends on the covariate $X$, e.g., length of the comment, we can create
further classes by splitting on relevant covariates.) Our queueing system
operates in heavy traffic with overall traffic intensity $\approx 1$, aligning
with Assumption~\ref{assumption: heavy traffic}. We set higher delay
costs for toxic comments and comments targeting historically marginalized or
oppressed groups. Specifically, for each demographic group $i$, we set the delay
cost as $C_{i, \cdot}(t) = c_{i, \cdot} t^2/2,$ with
$c_{i, \text{toxic}} = [10, 22, 12, 20, 25]$ and
$c_{i, \text{non-toxic}} = [1, 1, 1, 1, 1]$ for toxic and nontoxic comments
mentioning the aforementioned demographic groups, respectively. 
We adopt a quadratic cost function to maintain a unified experimental and analytical framework. While Theorem~\ref{theorem: optimality of our policy} applies to any cost function satisfying Assumption~\ref{assumption: on cost functions for optimality}, the quadratic form enables us to use a single queueing setting consistently across all sections, allowing closed-form derivations and explicit quantification of the model selection criteria in Section~\ref{section: characterization of the optimal cost with q}. This focus provides clarity and coherence while the main results (Theorem~\ref{theorem: optimality of our policy}) remain valid for the broader class of cost functions covered by Assumption~\ref{assumption: on cost functions for optimality}.


\paragraph{Queueing policies} We compare our proposed P$c\mu$-rule against
three scheduling approaches.
Unless otherwise stated, we employ the ERM-finetuned model with threshold $0.5$ as the predictor $\model$. All scheduling policies utilize the predicted class labels provided by the language model for jobs in the queue.
First, we consider the \naivemethod~\eqref{eq: naive gcmu} that treats the predicted classes as true, and employs the usual
Gc$\mu$-rule. For both \ourmethod~and \naivemethod, we assume the scheduler
has complete knowledge of the arrival/service rates of the predicted classes,
and use the confusion matrix $Q$ computed on the validation dataset. Second,
we study a black-box approach scheduling using deep reinforcement learning
methods (DRL), where we use a Q-learning method to estimate the value
function using a feedforward neural network (deep
Q-Networks~\citep{MnihKaSiGrAnWiRi13}),
where the value function uses predicted class label from
the language model as its state.
Finally, we consider the Oracle Gc$\mu$
policy, which knows the true class as well as associated arrival/service rates.
All policies are evaluated in the aforementioned setup, where the
scheduler predicts the class label using the AI model $\model$.

To train our DRL policy, we use~\citet{NamkoongMaGl22}'s discrete event
queueing simulator.
The DRL policy also uses the predicted class label from the language model, consistent with the other scheduling methods.
In particular,
we use $\{(\text{queue length, age of the oldest job})\}$ of all predicted
classes as our state space and let the predicted classes
$\{1, \ldots, K = 10\}$ be the action space. We learn a Q-function
parameterized by a three-layer fully connected network, and serve the oldest
job in the predicted class that maximizes the Q-function. As instantaneous
rewards, we use the sum of cost rates, $c_{i, \cdot}$ times the age, for all classes. We employ a similar training procedure as described in~\citep[Section D.1]{NamkoongMaGl22}, and
impose a large penalty to discourage the policy from serving empty queues.

\begin{figure}[t]
\centering  
\begin{minipage}[t]{0.325\textwidth}
\includegraphics[width=0.95\columnwidth]{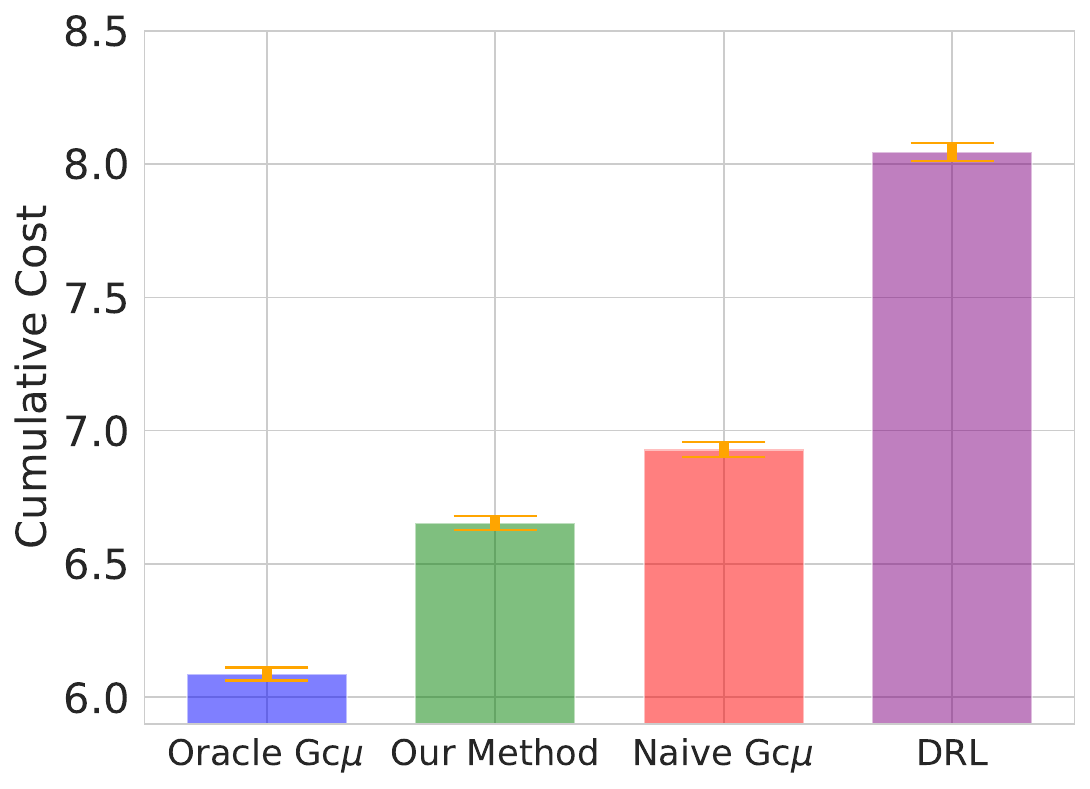} \\
\includegraphics[width=0.95\columnwidth]{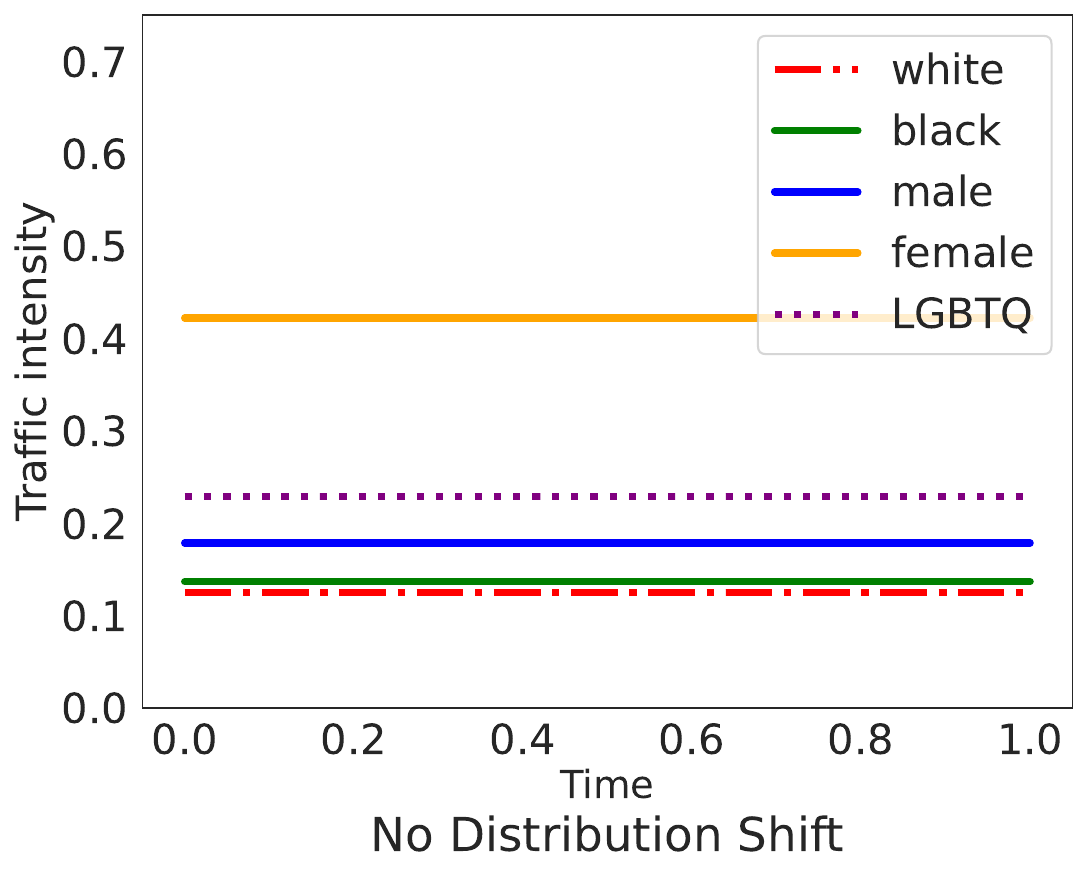} 
\end{minipage}
\begin{minipage}[t]{0.325\textwidth}
\includegraphics[width=0.95\columnwidth]{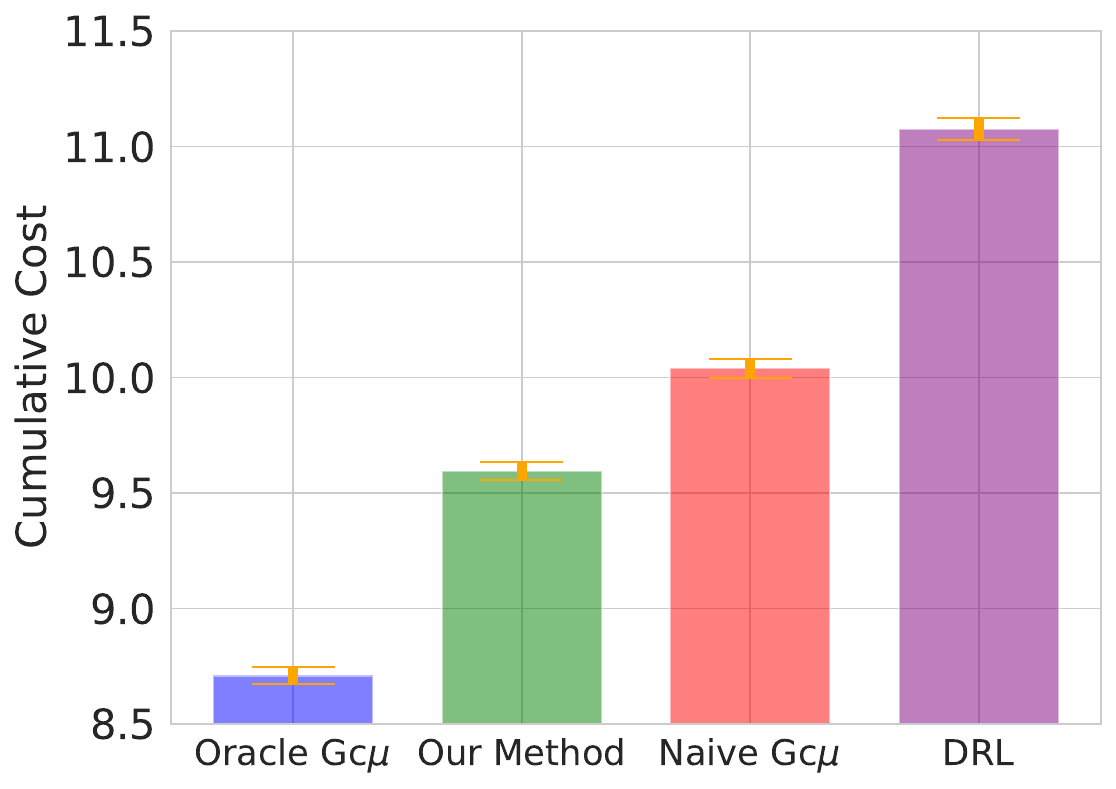} \\
\includegraphics[width=0.95\columnwidth]{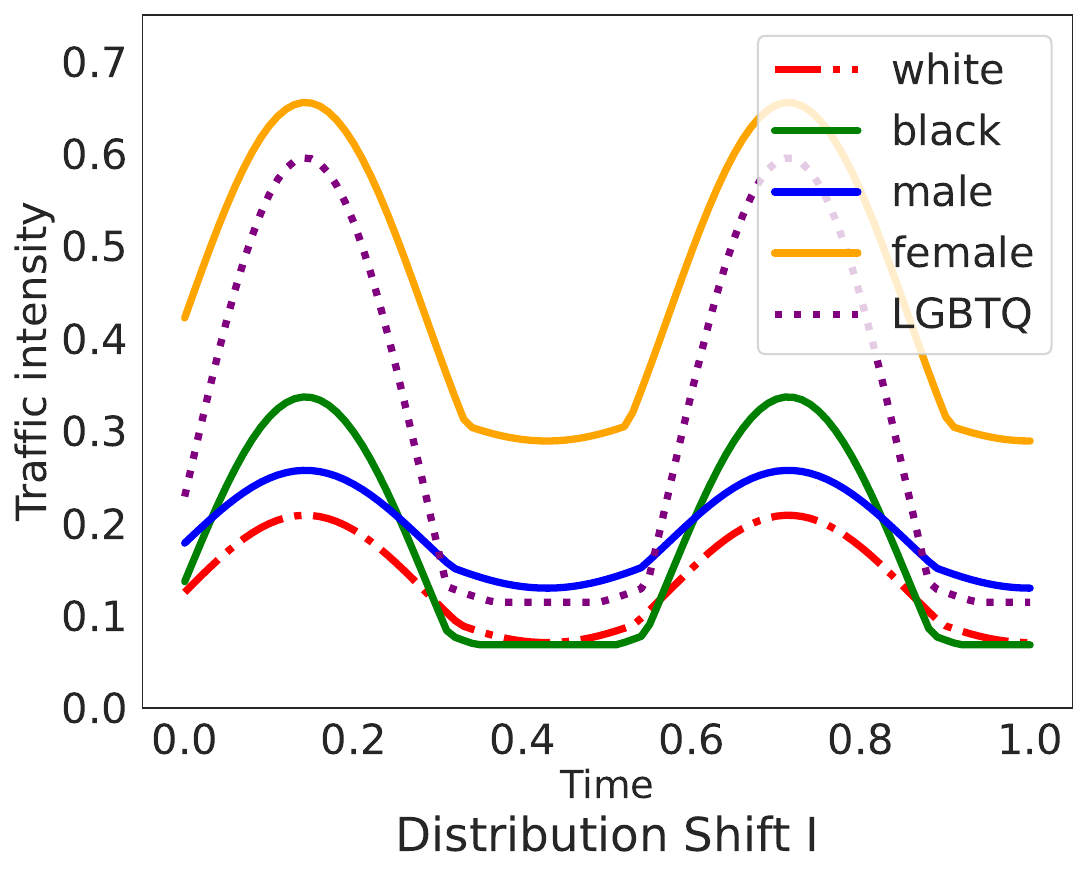} 
\end{minipage}
\begin{minipage}[t]{0.325\textwidth}
\includegraphics[width=0.95\columnwidth]{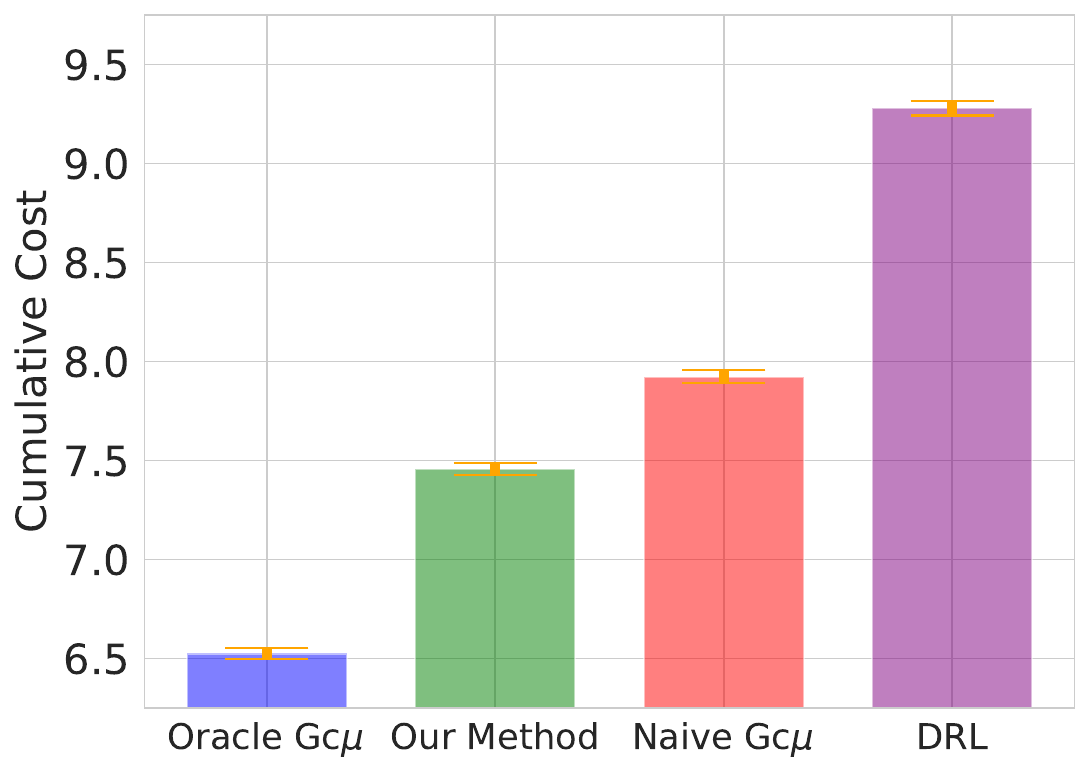} \\
\includegraphics[width=0.95\columnwidth]{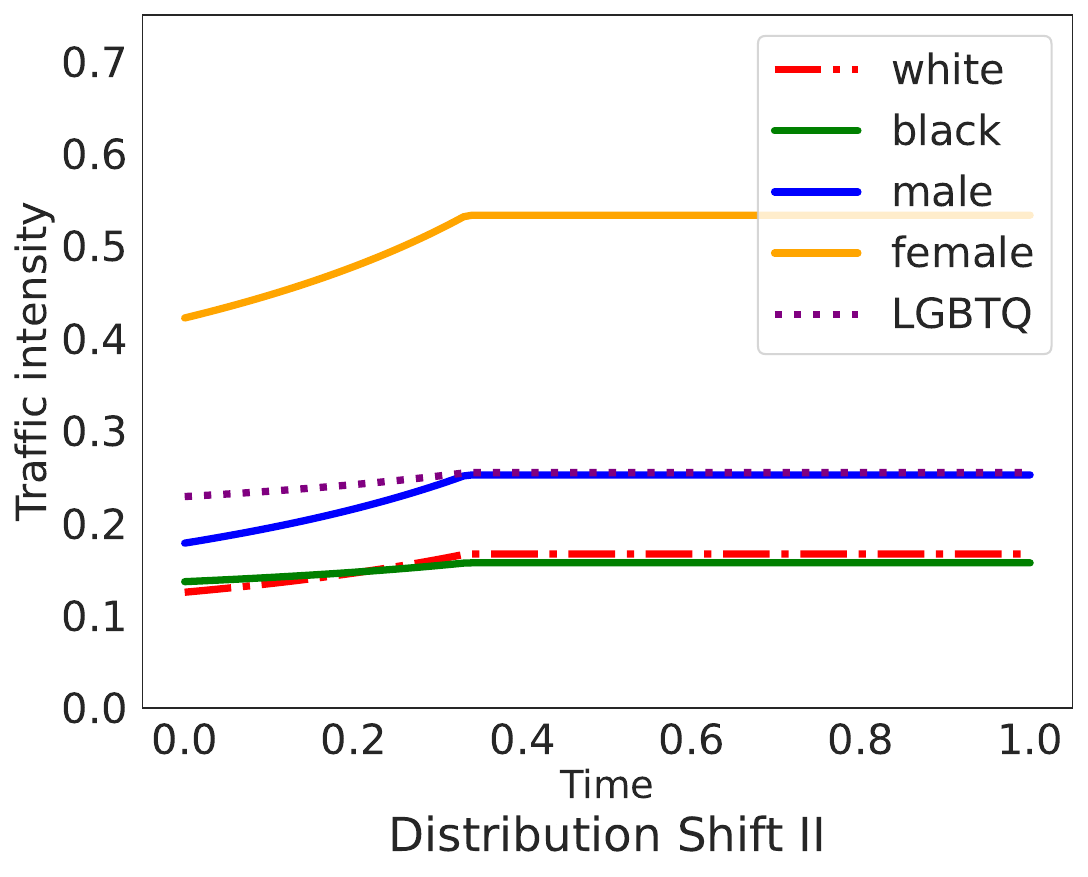} 
\end{minipage}
\caption{We present the cumulative cost for different policies under different
  testing environments (with 2$\times$ the standard error encapsulated in the
  orange bracket).} 
\label{fig: cumulative cost}
\vspace{-1em}
\end{figure}

\paragraph{Instability of reinforcement learning}
We run the deep Q-learning method with experience replay over 672 distinct
sets of hyperparameters and evaluate them based on the average cumulative
queuing cost over $5000$ independent sample paths simulated from the testing
enviroment. In Figure~\ref{fig: DRL_distribution}, we observe substantial
variation in queueing performance across hyperparameters even when using
identical instant reward functions and training/testing enviroments. (We also
use the same random seed across training runs.)  In particular, the minimum,
bottom \& top deciles of cumulative costs are 7.98, 9.79, and 60.46.  Our
empirical observation highlights the significant engineering effort required
to apply DRL approaches to scheduling and replicates previous findings in the
RL literature (e.g.,~\citep{HendersonIsBa18}). In the rest of the experiments,
we select the best hyperparameter based on average queueing costs reported in
Figure~\ref{fig: DRL_distribution}.

\paragraph{Main results}
In Figure~\ref{fig: cumulative cost}, we present cumulative cost averaged over
$50K$ sample paths. In the first column of Figure~\ref{fig: cumulative cost},
we test scheduling policies under the environment they were designed for:
constant arrival/service rates as we described above, with traffic intensity
$\approx$ 1. The P$c\mu$-rule outperforms \naivemethod~by $\sim30\%$ and 
DRL by $60-70\%$ in terms of the cost gap towards the Oracle G$c\mu$-rule.
While we expect the DRL policy can be further improved by additional
engineering (reward shaping, neural network architecture search etc), we view
the simplicity of our index-based policy as a significant practical advantage.
Next, we assess the robustness of the scheduling policies against
nonstationarity in the system. We consider two additional testing enviroments
with heavy traffic conditions that differ from that the policies were designed
for. See further details in Section~\ref{appendix: experiment details}. We observe the performance gains of the \ourmethod~hold over
nonstationarities in the system.

\begin{figure}[t]
\centering  
\includegraphics[width=0.5\columnwidth]{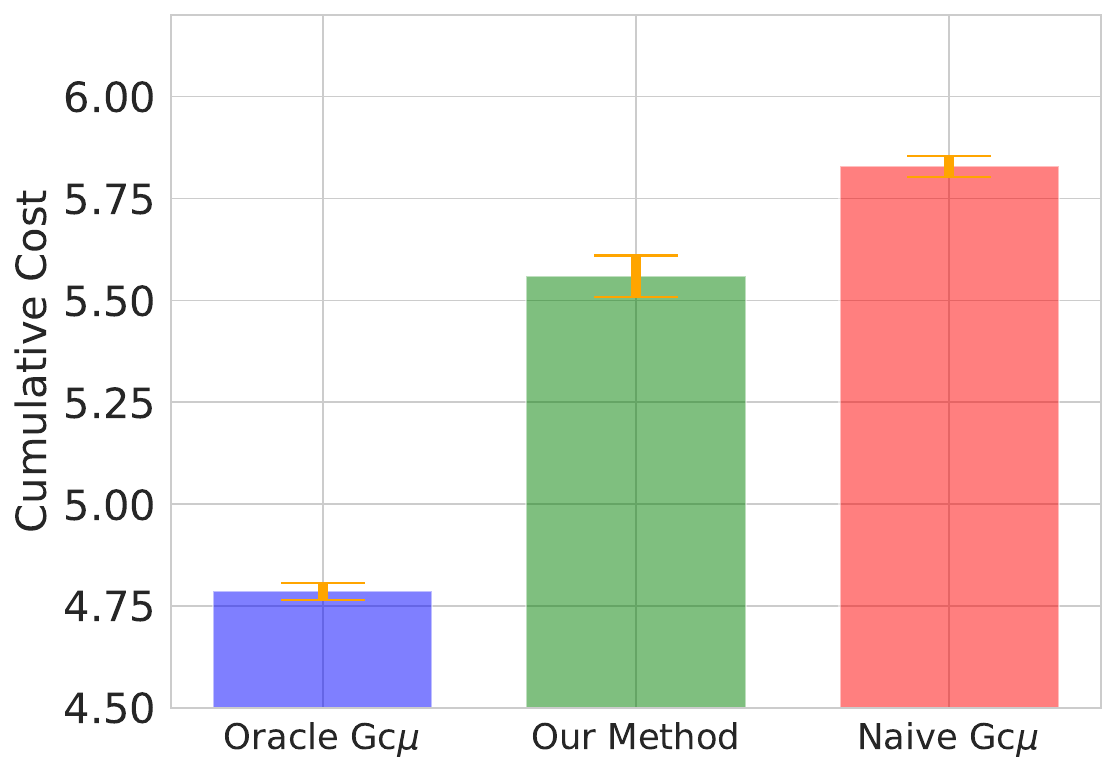} 
\caption{We present the cumulative cost for different policies where toxic classes have quadratic costs while nontoxic classes have cubic costs (with 2$\times$ the standard error encapsulated in the orange bracket).} 
\label{fig: different cost}
\vspace{-1em}
\end{figure}

\paragraph{Extension to other cost functions}
We also test heterogeneous cost forms: for toxic classes set
\(C_{i,\text{toxic}}(t)=c_{i,\text{toxic}}\,t^{2}\), and for non-toxic classes set
\(C_{i,\text{non-toxic}}(t)=c_{i,\text{non-toxic}}\,t^{3}\). This choice models
higher immediate cost for toxic items while allowing non-toxic items to incur
rapidly increasing penalties if delayed. The cubic cost is strictly but not
strongly convex, as its second derivative vanishes at the origin, so it lies
outside Assumption~\ref{assumption: on cost functions for optimality} required
by Theorem~\ref{theorem: optimality of our policy}. This experiment therefore
probes the empirical robustness of \ourmethod~beyond the regime covered by our
theory. All other experimental settings
(arrival/service rates, confusion matrix \(Q\), predictor, etc.) remain unchanged as in the first column of~Figure \ref{fig: cumulative cost}. We focus on Oracle Gc$\mu$, \ourmethod, and \naivemethod~for this experiment due to instability of reinforcement learning methods as shown above.
Figure~\ref{fig: different cost} shows that \ourmethod{} continues to outperform
\naivemethod{}, demonstrating its adaptability to different cost structures in content moderation tasks.

\section{Model selection based on queueing cost}
\label{section: characterization of the optimal cost with q}

Predictive models with similar accuracy levels can exhibit significant
differences in queueing performance.  By explicitly deriving the optimal
queueing cost under misclassification, our theoretical results allow designing
AI models with queueing cost as a central concern.  Since the \ourmethod~is
optimal in the heavy traffic limit, the corresponding $\TJ^*(t;Q)$ represents
the best possible cost when employing the given classifer, $f_\theta$, and the
relative regret $\TJ^*(t;Q)/\TJ^* (t;I)$ serves as an evaluation metric with
queueing performance as the central consideration. We empirically demonstrate
that this simple model selection criteria based on our theory can provide
substantial practical benefits in our content moderation simulator.

For quadratic cost functions, we can explicitly solve the optimization
problem~\eqref{eq: optimization problem} and derive analytic expressions for
$\TJ^* (\cdot; Q)$ and $\TJ^*(\cdot; I)$.
\begin{assumption}[Quadratic Cost Functions]\label{assumption: quadratic cost
    functions}
  For all $k\in [K]$, the cost functions are defined as
  $C^n_k(t) = \frac{1}{2n}c^n_kt^2,~t\in [0,n],~n \in \mathbb{N}$, and
  $C_k(t) = \frac{1}{2}c_kt^2,~t\in [0,1]$, where
  $\{c^n_k\}_{n\in \mathbb{N}}$ and $c_k$ are positive constants such that
  $c_k^n \rightarrow c_k$ as $n\rightarrow \infty.$
\end{assumption}

Under Assumption~\ref{assumption: quadratic cost functions}, \ourmethod~(Definition~\ref{definition: modified gcmu rule}) is defined using a linear cost rate $\fC_l' (t) = \fc_l t$ where $\fc_l := \sum_k \tp_k \q_{kl} c_k / (\sum_{k'} \tp_{k'} \q_{k'l})$. The following formulas are easy to approximate since the confusion matrix $Q$ can be effectively estimated on held-out data.
\begin{proposition}[Cumulative Cost Rate of the \ourmethod]
\label{proposition: information gap; quadratic cost}
Given a classifier $\model$ and a sequence of queueing systems, suppose that 
Assumptions~\ref{assumption: data generating process},
\ref{assumption: heavy traffic},~\ref{assumption: second order moments}, and~\ref{assumption: quadratic cost functions} hold. Then, we have that 
\begin{equation*}
	\TJ^*(t;Q) = \frac{1}{\sum_{m=1}^K (\beta_m (Q))^{-1}} \cdot \frac{1}{2} \int_0^t \sumTtW (s)^2 \mathrm{d}s,~\forall~t\in[0, 1],       
\end{equation*}
where $\beta_l (Q):= \fmu_l \fc_l / \frho_l$.
Under the \naivemethod,
the scaled cumulative queueing cost $\TJ^n_\text{Naive} (\cdot;Q^n)$ has the limit
\begin{equation*}
        \TJ_\text{Naive} (t;Q) = \sum\limits_{l=1}^K \frac{\beta_l (Q)}{\Big(\sum_m \frac{\beta_{l, \text{Naive}}(Q)}{\beta_{m, \text{Naive}}(Q)}\Big)^2} \cdot \frac{1}{2} \int_0^t \sumTtW (s)^2 \mathrm{d}s,
\end{equation*}
where $\beta_{l,\text{Naive}}(Q)=\fmu_l \tc_l / \frho_l$.
\end{proposition}
\noindent See Section~\ref{section: proof of cumulative cost rate} for the
proof of Proposition~\ref{proposition: information gap; quadratic cost}.
$\TJ^* (\cdot; Q)$ is dominated by small values of $\beta_m (Q)=\fmu_m \fc_m / \frho_m$, as is the
case for the limiting workload $\TfW_m$ under the \ourmethod~(see
Section~\ref{section: proof of cumulative cost rate}). 
Tuning the classifier to reduce the values of $\beta_m(Q)$, 
that is, ensuring that the resulting \ourmethod~indices remain small 
relative to the adjusted intensities $\{\frho_l\}_{l\in[K]}$, will drive 
the optimal cost downward.
In what follows, we heavily rely on the independence between $\sumTtW$ and misclassification
errors from Proposition~\ref{prop: convergence and approximation of predicted class N, tau, T, and W}.

\begin{figure}[t]
\vspace{-1em}
\centering  
\begin{minipage}[t]{0.485\textwidth}
\centering
\includegraphics[width=\columnwidth]{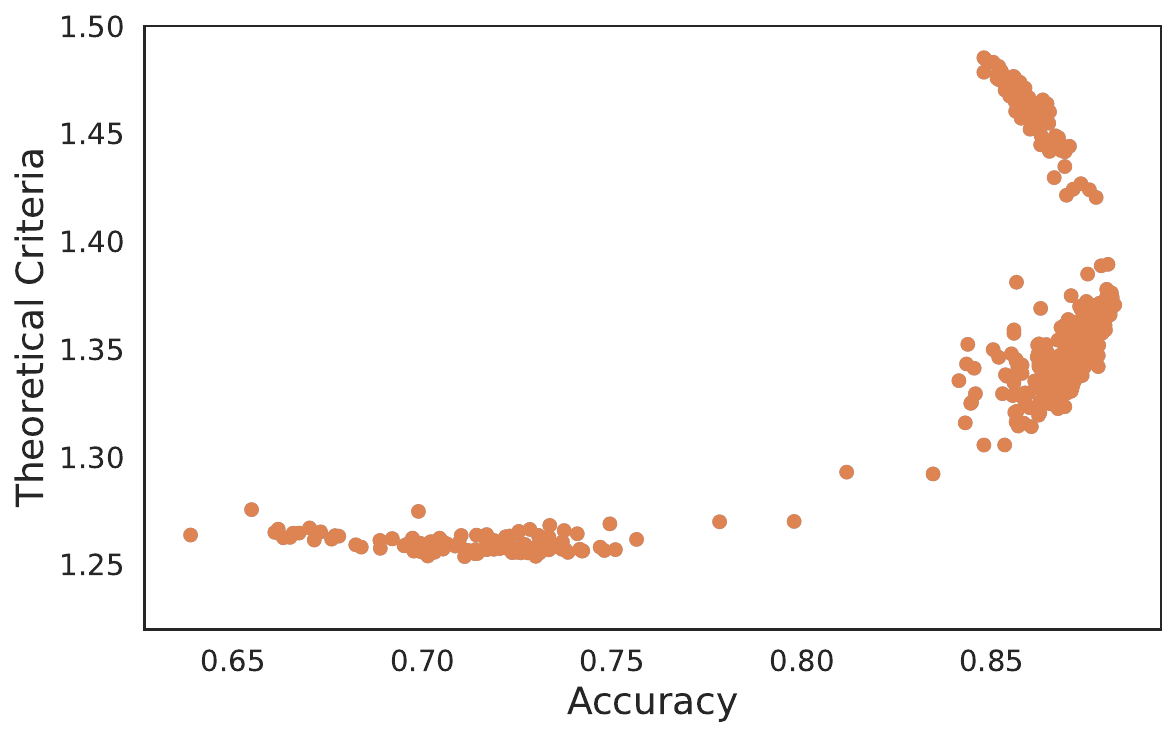} 
\end{minipage}
\begin{minipage}[t]{0.485\textwidth}
\centering
\includegraphics[width=\columnwidth]{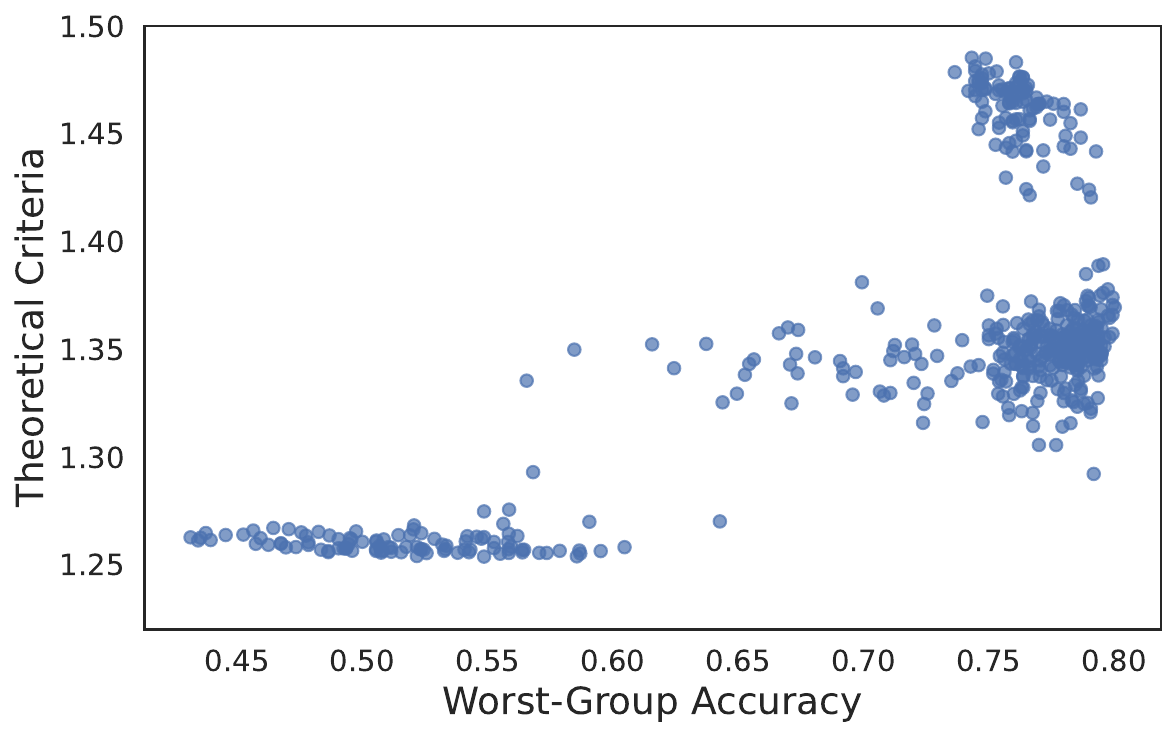} 
\end{minipage}
\caption{Distribution of relative regret $\TJ^*(t;Q)/\TJ^* (t;I)$ in the
  heavy traffic limit (``theoretical criteria'').
  Each dot represents a model obtained via
  last-layer finetuning with varying group loss weights.}
  \label{fig: model multiplicity}
\vspace{-1em}
\end{figure}
\paragraph{Model Multiplicity}
It is well known that models of equal prediction accuracy can perform
differently in downstream decision-making tasks~\citep{DAmourEtAl20,
  BlackRaBa22}. This phenomenon, known as model
multiplicity~\citep{BlackRaBa22}, is particularly important in our setting,
since prediction errors over different classes can have disproportionate
impacts on downstream queueing performance.  
We consider \ourmethod~in the testing environment from
Section~\ref{section: numerical experiments}
 to showcase that models with high accuracy levels can still exhibit
significant differences in queueing performances.

Given fixed costs, arrival rates, and service rates from
Section~\ref{section: numerical experiments}, we can explicitly quantify
the relative regret $\TJ^*(t;Q)/\TJ^* (t;I)$ for different classifiers through
the confusion matrix $Q$, considering the maximum and minimum possible
relative regret given a fixed accuracy or worst-group accuracy level.
We consider classifiers generated by further finetuning the ERM-based 
language models introduced in Section~\ref{section: numerical experiments}
on subsets of the {CivilComments} dataset's training data,
containing samples from one and only one of the common groups
\emph{white, black, male, female, LGBTQ}.
To make the procedure practical, we perform \emph{last-layer finetuning}, 
where the CLS embeddings of comments remain fixed while the classifier layer is updated. The language model parameters are not randomly reinitialized but continue from the ERM-trained weights in Section~\ref{section: numerical experiments}.
The finetuning process continues to use ERM with a learning rate of $10^{-3}$, batch size of $512$, and $10$ training epochs. To simulate different training objectives that balance performance across groups, we randomly sample $500$ sets of group loss weights from $[0, 100]$ for the 10 groups.

In Figure~\ref{fig: model multiplicity}, we plot the distribution of
the theoretical criteria for the 500 finetuned models.
The results show substantial variation in the criterion across different overall and worst-group accuracy levels, up to $15\%-20\%$ within the $80\%-85\%$ accuracy range and the $70\%-75\%$ worst-group accuracy range. 
Interestingly, models with lower accuracy or worst-group accuracy may exhibit better queueing performance according to the theoretical criterion. 
This experiment demonstrates that model multiplicity is a genuine concern for queueing systems as well: classifiers with comparable accuracy can lead to markedly different queueing outcomes, underscoring the value of our proposed evaluation metric in guiding model selection.

\begin{figure}[t]
\vspace{-1em}
\centering  
\begin{minipage}[t]{0.485\textwidth}
\centering
\includegraphics[width=0.985\columnwidth]{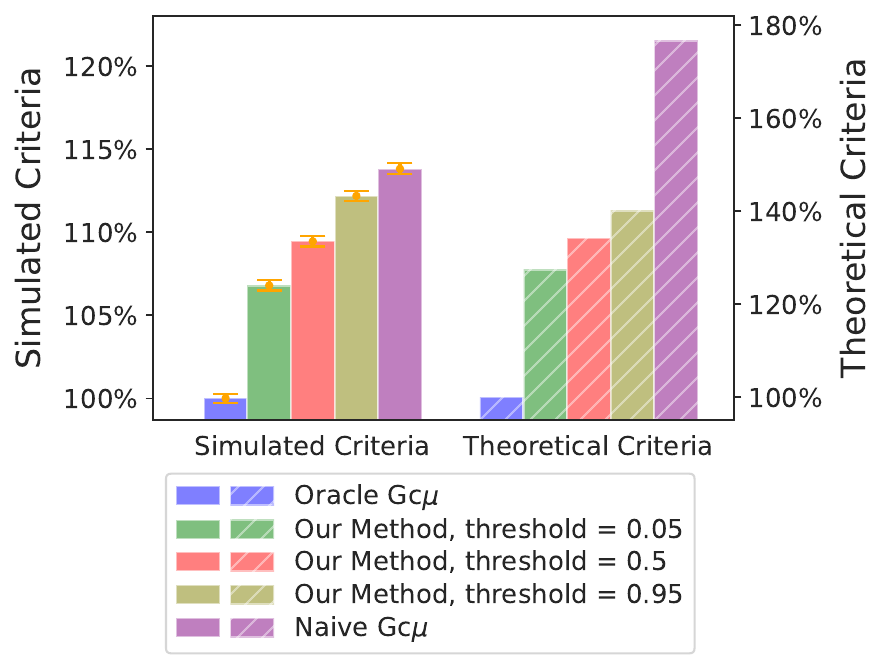} 
(a) Threshold selection
\end{minipage}
\begin{minipage}[t]{0.485\textwidth}
\centering
\includegraphics[width=0.985\columnwidth]{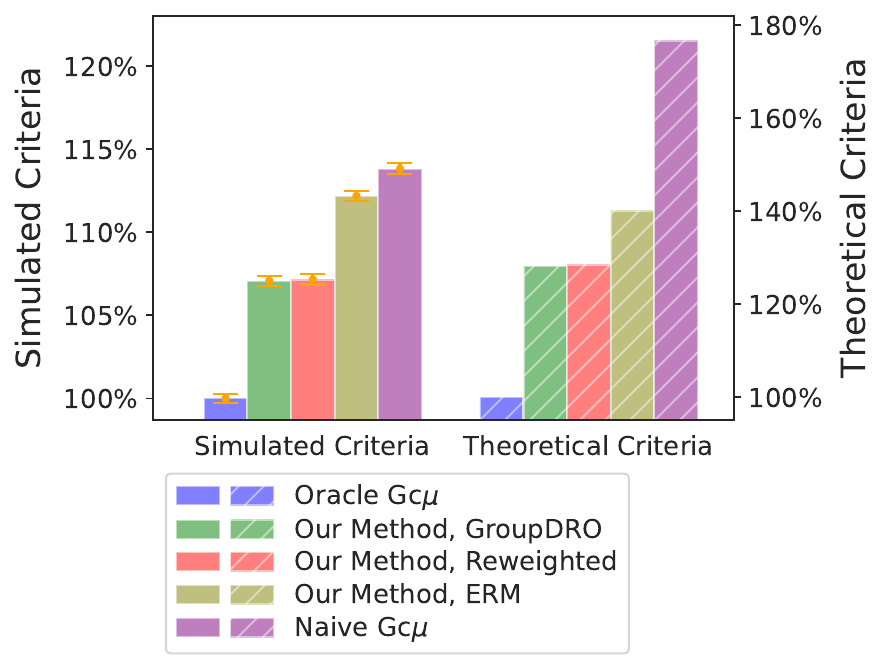} 
(b) Model selection
\end{minipage}
\caption{For different policies, we present our proposed model selection
criteria (theoretical criteria) based on the relative regret
$\TJ^*(t;Q)/\TJ^* (t;I)$ in the heavy traffic limit. To test its validity,
we plot the simulated/actual counterpart (simulated criteria) in the left.
The relative ranking of policies based on our theoretical criteria exactly
matches that given by the simulated quantities. 
}\label{fig: model_selection}
\vspace{-.5em}
\end{figure}

\paragraph{Model Selection between Fixed Candidates}
To further understand the validity of our proposed criteria,
we continue to use the same testing environment
from Section~\ref{section: numerical experiments} to conduct model selection
between fixed candidates.
We study the performance of the \ourmethod, Oracle
Gc$\mu$, and \naivemethod, using the cumulative queueing cost at $T=1$
across $5\times10^4$ independent sample paths. As the queueing cost of the
Oracle Gc$\mu$-rule converges to $\TJ^* (t;I)$, we normalize all simulated
cumulative cost over each sample path by the average cumulative cost of the
Oracle Gc$\mu$-rule. We refer to this quantity as ``simulated'' relative
regret.

We demonstrate the utility of selecting and evaluating classifiers based on
$\TJ^*(t, Q)/\TJ^*(t, I)$ using two tasks: (i) threshold selection for a fixed
classifier, and (ii) model selection for a given collection of classifiers.
In both cases, the ranking according to our proposed criteria aligns with
simulated counterparts, illustrating how an analytic characterization of
queueing cost can provide an effective comparison between ML models without
extensive queueing simulation. 
In Figure~\ref{fig: model_selection},
we present simulated relative regret using solid bars, with 2$\times$ the
standard error encapsulated in the orange brackets. The shaded bar depicts our
proposed model selection (theoretical criteria) given by the relative regret
in the heavy traffic limit. 
For threshold selection, we consider the \ourmethod~using
the aforementioned ERM predictor and compare its queueing performance with
thresholds being $[0.05, 0.5, 0.95]$, positioned from left to right in
Figure~\ref{fig: model_selection} (a). For model selection, we consider the
aforementioned three different classifiers: GroupDRO, Reweighted, and ERM with
thresholds 0.05, 0.05, and 0.95. (These thresholds are chosen to showcase
diverse queueing performances).  In Figure~\ref{fig: model_selection} (b), we
evaluate the \ourmethod~using these models. We also compare \ourmethod~to 
the \naivemethod, where the classifier is fixed to the aforementeioned ERM 
classifier with threshold 0.5 in Figure~\ref{fig: model_selection} (a) and (b).

We demonstrated that our proposed evaluation metric $\TJ^*(t, Q)/\TJ^*(t, I)$ effectively guides model selection by focusing on queueing performance. This approach ensures that the selected models optimize overall system performance, not just predictive accuracy, providing a robust basis for designing and selecting AI models in service systems.

\paragraph{Model Selection with Training}
Besides selecting between a few candidate thresholds or models, we can also directly train or
finetune models to achieve improved queueing performance. To illustrate this,
we consider the same experimental setup as before and consider last-layer finetuning of the ERM-finetuned model with threshold 0.5, introduced in
Section~\ref{section: numerical experiments}.

Traditionally, model training or finetuning typically focuses on
predictive performance, such as accuracy, precision, recall, or their
weighted combinations. We consider two straightforward methods for
last-layer finetuning:
(i) cost-weighted loss, where the each group's loss is weighted by the
queueing cost $c_{i, \text{toxic}}$ and $c_{i, \text{non-toxic}}$, and (ii)
cost/freq-weighted loss, where the each group's loss is weighted by its
queueing cost normalized by the group's frequency in the training data.
As we observe below, model rankings under the traditional methods
change across different queueing settings,
indicating their unreliability in queueing tasks.

We also consider last-layer finetuning using relative regret.
The relative regret $\TJ^*(t;Q)/\TJ^* (t;I)$ is not differentiable in the
classifier parameters, so it cannot serve directly as a training loss.
Because last-layer finetuning is fast, however, we can still optimize it
through random exploration: we derive 500 models by sampling group loss
weights from $[0, 100]$ for the 10 groups, as described in the
``Model Multiplicity" paragraph.
We select the model with the lowest relative regret as the final model and compare it with the aforementioned traditional methods as well as the base ERM model. For simplicity, all models use the same decision threshold of 0.5.

We consider two settings: (i) the setting from Section~\ref{section: numerical experiments}, and (ii) a modified setting that is nearly identical except for the class arrival rates. In setting (ii), arrivals remain i.i.d exponential with rate 100, but the probability of a toxic sample is set to $1/55$, while each non-toxic class has probability $10/55$.

\begin{figure}[t]
\vspace{-1em}
    \centering  
\begin{minipage}[t]{0.45\textwidth}
    \centering
    \includegraphics[width=0.95\columnwidth]{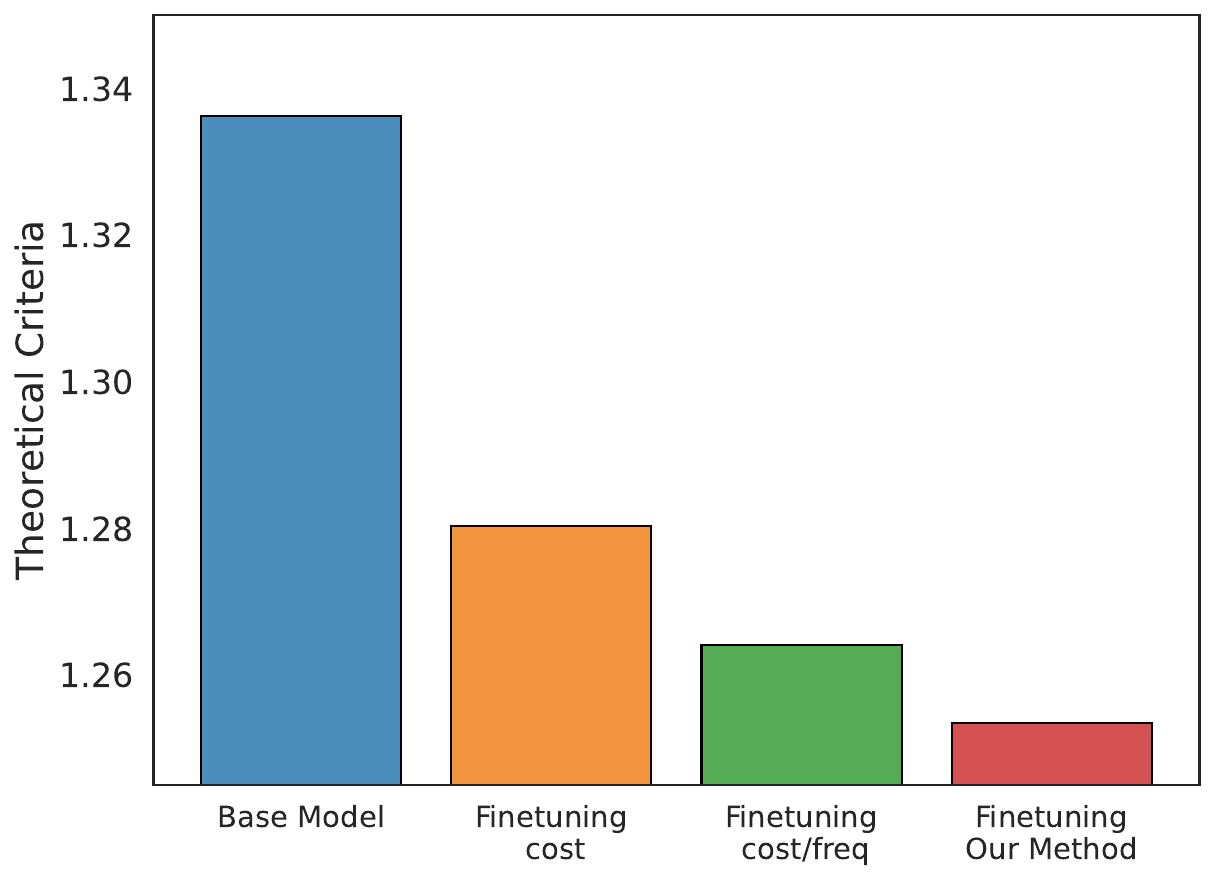} \\
    \text{(a) Queueing setting (i) in Section~\ref{section: numerical experiments}}
\end{minipage}
\begin{minipage}[t]{0.45\textwidth}
    \centering
    \includegraphics[width=0.95\columnwidth]{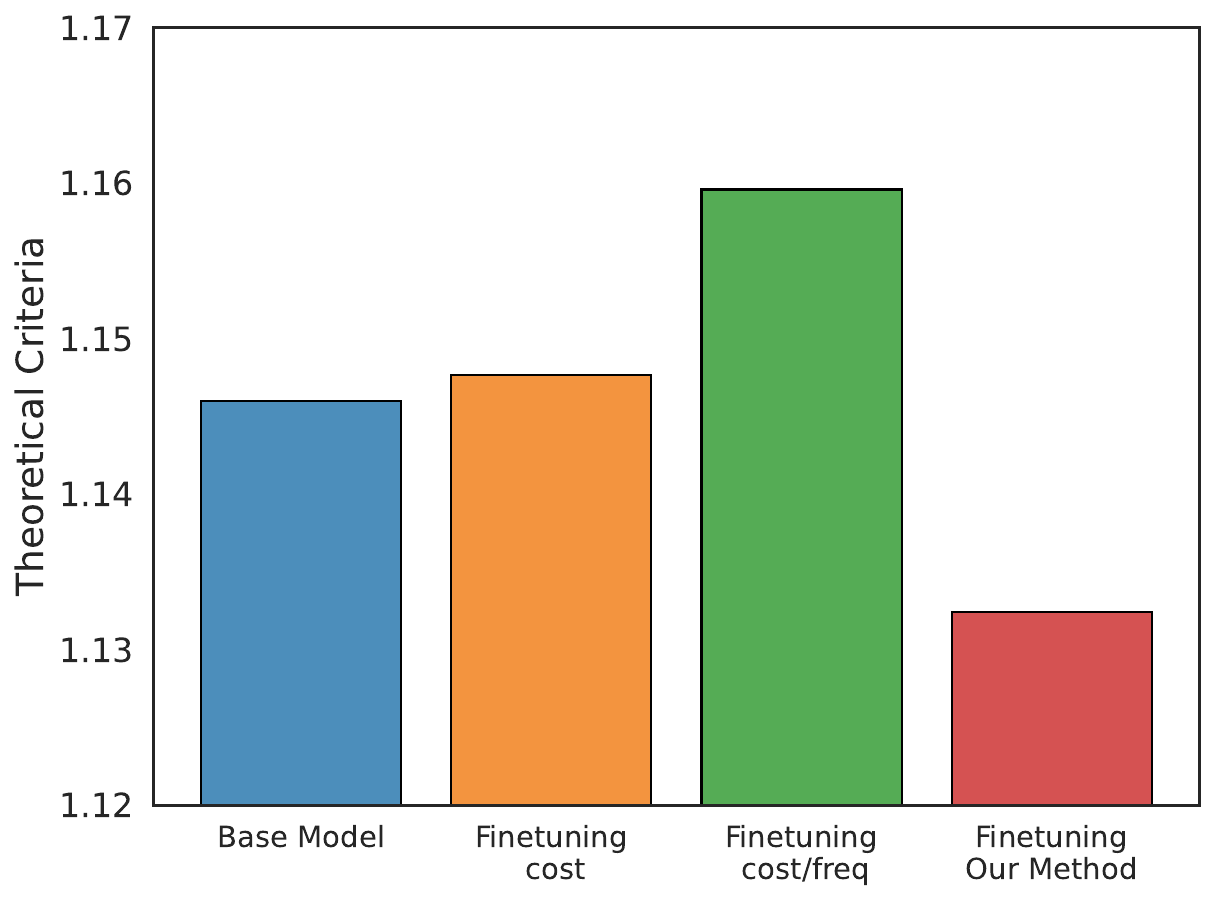} \\
    \text{(b) Queueing setting (ii)}
\end{minipage}
\vspace{0.75em}
\caption{We compare the base ERM model, models trained using traditional
  weighted-accuracy-based methods, and our proposed method based on the
  theoretical criteria. Our method reduces cumulative queueing costs, 
  while traditional methods exhibit varying rankings
  across different queueing settings.  }\label{fig: finetuning}
\vspace{-.75em}
\end{figure}

The corresponding theoretical criteria are shown in Figure~\ref{fig: finetuning}. We do not present the simulated criteria for two reasons. First, the finetuned models exhibit smaller variations in theoretical criteria compared to the more diverse models or thresholds considered previously, causing differences in the simulated criteria to be small and fall within the range of standard error, making them difficult to visualize within a reasonable computing budget. Second, as shown earlier, the ranking of the theoretical criteria exactly matches that of the simulated criteria, supporting the use of theoretical criteria for clarity.

As shown in Figure~\ref{fig: finetuning}, in this simplified toy example, our
method still outperforms traditional methods by $\sim 2-3\%$ in cumulative
queueing costs. This demonstrates the effectiveness of our evaluation metric when
queueing performance is the major concern.  Moreover, the rankings of
traditional methods and the ERM base model
vary across different queueing settings, indicating their
unreliability in queueing tasks.
While we conduct a brute force search over 500 finetuned models,
developing a practical and scalable training algorithm in more complex
scenarios is an important direction of future work.

\section{Design of an AI-based triage system}
\label{section: design of an AI triage system}

Our characterization of queueing cost can be further utilized to design
comprehensive job processing systems assisted by AI models. Motivated by
content moderation systems on social media
platforms~\citep{ChandakAb23,Vincent20,TiktokContentModeration24},
we study a triage system where an initial AI model filters out clear-cut
cases, after which the queueing system serves the remaining jobs
(Figure~\ref{fig:diagram}).  Standard triage systems in online platforms
determine the filtering level using simple metrics such as maximizing recall
subject to a fixed high precision level (e.g.,~\citep{ChandakAb23}). These
designs~\citep{MetaHarmfulContent21,Schroepfer19,WangFaKhMaMa21,FacebookHatespeechDetection20}
do not directly account for the downstream operational and queueing costs.

In this section, we provide a novel framework for designing AI-assisted triage
systems that jointly optimize the filtering and queueing systems, taking into
account all four types of costs: filtering costs, hiring costs,
misclassification costs, and queueing congestion costs.  Our objective can be
easily estimated using a small set of validation data and a simple simulation
of a (reflected) Brownian motion, allowing us to find the optimal filtering
level through methods like line search. In Section~\ref{section: experiment
  triage system}, we conduct numerical experiments to demonstrate
the effectiveness of this approach.  We find that prediction-based metrics, which are the norm
in practice, may align with the total cost when either filtering cost or
hiring cost dominates, but they fail to do so in more complicated settings with
trade-offs between different types of costs. Our method avoids computationally
expensive queueing simulations and consistently identifies the optimal
filtering and staffing levels in all of these scenarios by simply
simulating a (reflected) Brownian motion.

\subsection{Model of the AI-based triage system}
\label{section: model of the AI based system}

We consider a sequence of single-stream incoming jobs that arrive at the triage 
system. We assume the $n$th system operates on a finite time horizon $[0, n]$, starts
empty, and has i.i.d. interarrival times with an arrival rate of $\Lambda^n$. 
With a slight abuse of notation, we let $u^n_i$ be the interarrival time of the 
$i$th job in system $n$, $\tU^n_0(t)$ be the arrival time of the $\lfloor t\rfloor$th 
job in system $n$, and $\tA^n_0(t)$ be the total number of jobs arriving in the 
triage system $n$ up to time $t$.
For simplicity, we consider a two-class setting, with class 1 representing toxic content and class 2 representing non-toxic content. 
For each job, a tuple of (observed features, true class label, service time), 
denoted as $(X^n_i, Y^n_i, v^n_i)$, is generated \emph{identically} and
\emph{independently} of its arrival time $u^n_i$. Similar to our model in 
Section~\ref{section: model}, $v^n_i$ and $X^n_i$ are conditionally
independent given $Y^n_i$.

We use a binary classifier $\model$ for the filtering procedure across all
systems $n$. With a slight abuse of notation, let $\model(\cdot) \in [0, 1]$
now be the toxicity score instead of the predicted class label.  Specifically,
the classifier outputs $\model(X^n_i) \in [0,1]$ based on the observed
features $X^n_i$ for each job $i$.  The system designer is tasked with
choosing a threshold $\fl$ that affects the (triage) filtering level: an
arriving job in system $n$ can pass the filtering system and enter the
queueing system if and only if $\model(X^n_i) \geq \fl$.\footnote{For
simplicity, we only consider filtering out clearly safe content in this
section, though in practice, the system designer can choose another
threshold to filter out clearly toxic content from the human review
process and directly take further actions.}  Content that is filtered
out is not reviewed and can remain on the platform.  A higher filtering level
$\fl$ filters more jobs out, resulting in higher false negative rate (more
filtering and misclassification costs), fewer human reviewers required (lower
hiring cost), and a complex effect on the downstream queueing cost.

Each job that passes the filtering system is subsequently sent to human
reviewers (queueing system).  Given the filtering level $\fl$, we use the same
number of reviewers $\Gamma(\fl)$ across all systems, where $\Gamma(\fl)$ is a
predetermined decision variable, fixed in advance and not subject to
randomness.  We assume that for each system $n$, all reviewers have the same
service rate for class $k$ jobs, i.e.,
$\mu_{k, r}^n = \mu_{k}^n,~\forall~k\in \{1, 2\}$ for each reviewer
$r\in[\Gamma(\fl)]$.  To ensure workload equality and fairness among
reviewers, we assume jobs passing through the filtering system are assigned to
one and only one human reviewer with equal probability $1/\Gamma(\fl)$,
\emph{independently} of any other random objects.  Each human reviewer
operates their own single-server queueing system. The jobs allocated to the
$r$th reviewer correspond to their arrival process, denoted as
$\tApsr^n(t)$, which is split from the common arrival process after
filtering, denoted as $\tAps^n(t)$.  For the $j$th job passing through the
filtering system, let
$B^n_j: = (B_{j1}^n, \ldots, B_{j\Gamma(\fl)}^n)$ be the one-hot
encoded representation of the reviewer it is assigned to. Then,
$\tApsr^n (t) : = \sum_{j=1}^{\tAps^n(t)} B_{jr}^n,~\forall~t\in[0, n],
r\in[\Gamma(\fl)]$.

For simplicity, we assume all reviewers utilize $\model$ to predict the class
labels of incoming jobs.  The system designer must decide another threshold
$\tx \geq \fl$ that affects the toxicity classification. In particular, for
the $s$th job assigned to reviewer $r$, it is predicted to be toxic (class 1),
i.e., $\fY^n_{s1, r} = 1$, if and only if $\model(X^n_{s, r}) \geq \tx$.  We
assume all human reviewers adopt the same scheduling policy.  
Reviewers use the predicted class $\VfY^n_{s,r}=(\fY^n_{s1, r}, \fY^n_{s2, r})$ and feasible scheduling policies must satisfy a variant of
Definition~\ref{definition: feasible policies}.

Throughout this section, we use $i$ when counting jobs that arrive at the
triage system, $j$ for jobs that pass the filtering system and arrive at the
queueing system, $s$ for jobs assigned to a human reviewer, and $r$ for human
reviewers (servers). For a stochastic process $\tAps^n(t)$, the subscript
$\text{ps}$ indicates processes associated with jobs passing through the
filter and arriving at the queueing system. We use the subscript $0$ to
indicate the total arrival process and $r$ to indicate processes associated
with the reviewer $r$.  We denote our decision variables $(\fl, \tx)$ by
$z$. We summarize our assumptions for the AI-based triage system below.
\begin{assumption}[Data generating processes for the AI-based triage system]
\label{assumption: data generating process for triage system}
For any system $n\in \N$, (i)
$\{(\tu_i^n, \tv_i^n, \feature_i^n, \tY_i^n): i\in \N\}$ is a sequence of
i.i.d. random vectors; (ii) $\{\tu_i^n: i \in \mathbb{N}\}$ and
$\{(\tv_i^n, \feature_i^n, \tY_i^n): i \in \mathbb{N}\}$ are independent;
(iii) for any $i \in \mathbb{N}$, $v^n_i$ and $X^n_i$ are conditionally
independent given $Y^n_i$;
(iv) $\{B^n_j: j \in \mathbb{N}\}$ is a sequence of i.i.d. random vectors;
(v)  $\{B^n_j: j \in \mathbb{N}\}$ is independent of 
$\{\tu_i^n: i \in \mathbb{N}\}$ and 
$\{(\tv_i^n, \feature_i^n, \tY_i^n): i \in \mathbb{N}\}$.
\end{assumption}

The data generating process for the AI-based triage system is crucial to our
analysis, because it enables the reduction of the scheduling problem for all
reviewers to stochastically identical single-server scheduling problems across
the reviewers. Assumption~\ref{assumption: data generating process for triage
  system} (i), (iv) and (v) ensure that each reviewer $r$ has a single stream
of jobs with i.i.d. tuples
$\{(v^n_{s,r}, X^n_{s, r}, Y^n_{s, r}): s\in \mathbb{N}\}$. More importantly,
the tuples associated with reviewer $r$ become \emph{independent} of those of
any other reviewers. This leads to the joint convergence of the
diffusion-scaled processes defined by
$\{(v^n_{s,r}, X^n_{s, r}, Y^n_{s, r}): s\in \mathbb{N}\}$ across all
reviewers $r\in [\Gamma(\fl)]$.  In addition, similar to Section~\ref{section:
  model}, Assumption~\ref{assumption: data generating process for triage
  system} (i), (iv), and (v) allow us to disentangle the interarrival times
$\{u^n_{s,r}: s\in \mathbb{N}\}$ from the filtering process, service
processes, and the covariates, ensuring the joint convergence of the
diffusion-scaled processes defined by $\{u^n_{s,r}: s\in \mathbb{N}\}$ across
the reviewers $r\in [\Gamma(\fl)]$. Since Assumption~\ref{assumption: data
  generating process for triage system} (ii) and (v) imply independence
between
$\{(v^n_{s,r}, X^n_{s, r}, Y^n_{s, r}): s\in \mathbb{N}\}_{r\in
  [\Gamma(\fl)]}$ and $\{u^n_{s,r}: s\in \mathbb{N}\}_{r\in [\Gamma(\fl)]}$,
we can derive the desired joint convergence (Lemma~\ref{lemma: joint weak
  convergence of triage system}) and apply our \emph{sample path analysis at
  the reviewer level}.  For further discussion, see EC.~\ref{section:
  joint convergence of triage system}.

\subsection{Heavy traffic conditions for the AI-based triage system}
\label{section: heavy traffic triage system}

In the sequel, we assume the triage system operates under heavy traffic conditions 
and analyze the limiting system. Denote the
conditional probability of a class $k$ job passing the filtering level $\fl$ as
$g_k^n (\fl):=\mathbb{P}^n [\model (\feature_i^n) \geq \fl \mid \tY_{ik}^n=1],
~\forall~\fl\in [0, 1],~k\in \{1, 2\}$.
Similar to Assumption~\ref{assumption: heavy traffic}, we adopt the following 
heavy traffic conditions for the AI-based triage system.

\begin{assumption}[Heavy traffic conditions for AI-based triage system]
\label{assumption: AI triage system}
Given a classifier $\model$ and a sequence of triage systems,  
we assume that there exist $\Lambda$, $\mu_k$, and $g_k: [0,1] \rightarrow [0,1]$
such that 
(i) for any filtering level $\fl\in [0, 1]$
and class $k\in \{1, 2\}$, 
we have that 
\begin{eqnarray*}
    n^{1/2} (\Lambda^n - \Lambda) \rightarrow 0, 
    \quad n^{1/2} (\mu_{k}^n - \mu_k)  \rightarrow  0, \quad 
    n^{1/2} (g_k^n(\fl) - g_k(\fl)) \rightarrow 0;
\end{eqnarray*}
(ii) given the filtering level $\fl$, the number of hired reviewers satisfies 
$\Gamma (\fl)=\Lambda \sum_{k=1}^2 \frac{\tp_k g_k (\fl)}{\tmu_k}.$
\end{assumption}
\noindent 
We adopt Assumption~\ref{assumption: AI triage system} to ensure that each 
reviewer aligns with Assumption~\ref{assumption: heavy traffic}. Specifically, 
according to Assumption~\ref{assumption: AI triage system} (i) 
and~\citep[Theorem 9.5.1]{Whitt02}, we can show that for each reviewer $r$, 
their class prevalence 
$\tp^n_{k, r}(\fl): = \P^n[\tY_{sk, r}^n=1 \mid \model(X^n_{s, r}) \geq \fl]$
and confusion matrix 
$\q^n_{kl, r}(z): = 
\P^n [ \fY^n_{sl, r} = 1 \mid \model(X^n_{s, r}) \geq \fl, \tY_{sk, r}^n=1]$
all converge to their limits $\tp_k(\fl)$ and $\q_{kl}(z)$
at the rate of  $o(n^{-1/2})$ 
(Lemma~\ref{lemma: conv of lambda p q for each reviewer}).
We use $Q^n(z)$ and $Q(z)$ to denote the prelimit and limiting 
confusion matrix for each reviewer.
In addition, by Assumption~\ref{assumption: AI triage system} (ii), we have that 
\begin{equation} \label{eq: heavy traffic for an individual reviewer}
\begin{aligned}
n^{1/2} \Big[ \frac{\Lambda^n}{\Gamma(\fl)}
  \sum\limits_{k=1}^2 \frac{\tp_k^n g_k^n (\fl)}{\tmu_k^n} - 1\Big] \rightarrow 0,
\end{aligned}
\end{equation}
which indicates that each reviewer operates under heavy traffic conditions and 
matches~\eqref{eq: heavy traffic condition}. 
According to Assumption~\ref{assumption: AI triage system} (ii), 
when all reviewers
operate under heavy traffic conditions, the number of reviewers is solely
determined by limiting traffic intensity and the filtering level $\fl$.
Thus, our decision variables are the filtering level $\fl$ and the toxicity level $\tx$,
with the number of reviewers determined accordingly. Intuitively, 
as the filtering level $\fl$ increases, the traffic intensity decreases and 
the number of reviewers hired also decreases.

Starting from Assumptions~\ref{assumption: data generating process for triage system}
and~\ref{assumption: AI triage system}, we first establish the joint convergence result 
in Lemma~\ref{lemma: joint weak convergence of triage system}. As 
Assumptions~\ref{assumption: data generating process for triage system}
and~\ref{assumption: AI triage system} are compatible with 
Assumptions~\ref{assumption: data generating process} 
and~\ref{assumption: heavy traffic}, we derive a common probability space 
$\mathbb{P}_\text{copy}$ in Lemma~\ref{lemma: uniform convergence of triage system} 
and apply the previous results on 
single-server queueing systems to each reviewer. This allows us to establish 
the limiting total cost of the triage system in
Section~\ref{section: total cost of triage system}.

\subsection{Total Cost of the AI-based triage system}
\label{section: total cost of triage system}

Motivated by content moderation problems, we divide the total cost into four components: 
filtering cost, hiring cost, misclassification cost, and queueing cost.  
Since the \ourmethod~is optimal for each single-server queueing system under heavy traffic 
conditions, we can explicitly 
quantify the best possible queueing cost of the limiting system under the quadratic
cost assumption (Assumption~\ref{assumption: quadratic cost functions}).
This enables us to determine the limiting total cost and minimize it to find the 
optimal filtering and classification levels $(\fl, \tx)$ for a fixed classifier
$\model$. In the following, we first define each cost component and then establish 
the limiting total cost in Theorem~\ref{theorem: total cost of triage queue system}.

\begin{definition}[Total cost of the AI-based triage system]
\label{definition: total cost of triage system}
Given a classifier $\model$, filtering level $\fl$, toxicity level $\tx$, the number
of hired reviewers $\Gamma(\fl)$, and a sequence of AI-based triage systems,
for a sequence of feasible policies $\{\policyn\}$,
 define the cost incurred as the
following.
\begin{enumerate}[(i)]
\item (Filtering cost) For each job that is filtered out, the unit costs for toxic
  and non-toxic jobs are $\ctrp > 0$ and $\ctrn < 0$. The total filtering cost up
  to time $t\in[0, n]$ is
\begin{equation*}
  G^n(t; \fl) :=  \ctrp \sum_{i = 1}^{\tA_0^n(t)} 
  \mathbb{I}(\model(X^n_{i}) < \fl)
  \cdot Y^n_{i1} 
  + \ctrn \sum_{i=1}^{\tA_0^n(t)} 
  \mathbb{I}(\model(X^n_{i}) < \fl) \cdot Y^n_{i2},
\end{equation*}
and $\TG^n(t;\fl): = n^{-1}G^n(nt; \fl)$ is the scaled filtering cost.
\item (Hiring Cost) Each reviewer costs $c_{\text{rev}} > 0$ per unit of time.
\item (Misclassification Cost) The per-job cost of false positive, false
  negative, true positive, or true negative are $\cfp, \cfn, \ctp, \ctn$,
  respectively.  The total misclassification cost up to time $t$ is
  $M^n(t; z)$, and its scaled counterpart is
  $\TM^n(t; z): = n^{-1} M^n(nt; z)$. 
\item (Queueing Cost) For each system $n$ and reviewer $r$,
$J^n_{\policyn, r} (t; Q^n (z))$ is the cumulative queueing cost as
defined in Section~\ref{subsection: fundamental convergence results}, and
$\TJ^n_{\policyn, r} (t; Q^n (z)): = n^{-1} J^n_{\policyn, r} (nt;
Q^n (z)),~t\in [0,1]$, is its scaled counterpart.
\end{enumerate}
The total cost incurred up to time $t$ is defined by
\begin{equation*}
  F_{\policyn}^n (t; z) = 
  \underbrace{G^n(t; \fl)}_{\text{filtering}} 
  + \underbrace{c_{\text{rev}} \Gamma(\fl)t}_\text{hiring}
  + \underbrace{M^n(t; z)}_\text{misclassification} + 
  \underbrace{\sum_{r=1}^{\Gamma(\fl)} 
  J^n_{\policyn, r} (t; Q^n (z))}_\text{queueing},
  ~\forall~t\in [0, n],
\end{equation*}
and
$\TF_{\policyn}^n (t; z): = n^{-1} F_{\policyn}^n (nt; z)$ 
is its scaled counterpart.

\end{definition}

For any filtering level $\fl$ and toxicity level $\tx$, we can easily 
establish the optimal total cost of the AI-based triage system 
under heavy traffic limits by extending 
Proposition~\ref{prop: joint conv. of predicted class A, U, S, and V}, 
Theorem~\ref{theorem: optimality of our policy}, and
Proposition~\ref{proposition: information gap; quadratic cost}. 
Such optimal cost can be achieved by applying the \ourmethod~to all reviewers 
as shown in~\eqref{eq: total cost of triage queue system} in the following.

\begin{theorem}[Total cost of the AI-based triage system]
\label{theorem: total cost of triage queue system}
Given a classifier $\model$, filtering level $\fl$, toxicity level $\tx$, 
the number of hired reviewers $\Gamma(\fl)$, and 
a sequence of AI-based triage systems, suppose that
Assumptions~\ref{assumption: second order moments},~\ref{assumption: quadratic cost functions},~\ref{assumption: 
data generating process for triage system}, and~\ref{assumption: AI triage system} hold. 
There exists a common probability space $\mathbb{P}_\text{copy}$ such that
\begin{enumerate}[(i)]
\item (Lower bound) under any 
feasible policies $\{\policyn\}$, the associated total cost 
$\TF_{\policyn}^n (t; z)$ 
satisfies 
$\liminf\limits_{n} \TF_{\policyn}^n (t; z) \geq 
\TF^*(t; z)$, $\forall~t\in [0,1]~\mathbb{P}_\text{copy}-a.s.$. 
For the original processes under $\mathbb{P}^n$, under any feasible policies
$\{\pi'_n\}$, 
\begin{equation*}
\liminf\limits_{n} \mathbb{P}^n[\TF_{\pi'_n}^n (t; z) > x] 
\geq \mathbb{P}_\text{copy}[\TF^* (t; z) > x],~\forall x\in \mathbb{R}, t\in [0, 1];
\end{equation*}
\item (Optimality) under the \ourmethod, we have that 
$\TF^n_\ourpolicy (\cdot; z) \rightarrow
\TF^*(\cdot; z)$ $\text{in } (\mathcal{D}, \|\cdot \|),~\mathbb{P}_\text{copy}-a.s.$. For the original processes under $\mathbb{P}^n$, 
$\TF^n_\ourpolicy (\cdot; z) \Rightarrow
\TF^*(\cdot; z)$ $\text{in } (\mathcal{D}, J_1)$, and in particular, 
$\mathbb{P}^n[\TF^n_\ourpolicy (t; z) > x]
\rightarrow \mathbb{P}_\text{copy}[\TF^* (t; z) > x], 
~\forall~x\in \mathbb{R}, t\in [0, 1]$.
\end{enumerate}
Here, the optimal total cost $\TF^*(t; z)$ is defined as
\begin{equation}\label{eq: total cost of triage queue system}
\TF^*(t; z): = 
\TG^*(t; \fl) + c_{\text{rev}} \Gamma(\fl) t + \TM^*(t; z) 
+ \sum_{r=1}^{\Gamma(\fl)}  \TJ_r^*(t; Q(z)),
\end{equation}
where 
\begin{equation*}
\begin{aligned}
&\TG^*(t; \fl) 
= \Lambda t \cdot [\ctrp p_1(1 - g_1(\fl))
+ \ctrn p_2(1 - g_2(\fl))] , \\
&\TM^*(t; z)
=  \Lambda t \cdot
\big[p_1 g_1(\fl) [\ctp\q_{11}(z) + \cfn \q_{12}(z)]
+ p_2 g_2(\fl) [\cfp \q_{21}(z) + \ctn \q_{22}(z)]\big]\\
&\TJ_r^*(t; Q(z))
= \frac{\beta_1 (Q(z)) \beta_2 (Q(z))}{2\big[\beta_1 (Q(z)) +  \beta_2 (Q(z)) \big]} 
\int_0^t \sumTtW (s; z, r)^2 \mathrm{d}s, 
\end{aligned}
\end{equation*}
for all $t\in [0, 1]$, and $\sumTtW (t; z, r)$ is the limiting remaining 
total workload process of reviewer $r$ as defined in 
Lemma~\ref{lemma: conv. of L and W triage system}. 
\end{theorem}

According to Theorem~\ref{theorem: total cost of triage queue system}, we can
minimize~\eqref{eq: total cost of triage queue system} to find the optimal
filtering and toxicity levels $(\fl, \tx)$ for a given classifier $\model$.
In particular,~\eqref{eq: total cost of triage queue system} depends solely on
limiting exogenous quantities such as $\Lambda$, $p_k$, $g_k(\fl)$, $\q_{kl}(z)$
that can be easily estimated given a small set of validation data.
$\sumTtW(t; z, r)$ is a reflected Brownian motion with a known
drift and covariance (see further discussion in EC.~\ref{section:
  simulation of the total cost of the AI-based Triage System}), so we can
estimate the total cost using a simulated (reflected) Brownian motion.  The
optimal level $z^*$ can then be found through a simple line search
over $[0, 1]$. Our approach avoids traditional queueing simulations, which can
be costly and time-consuming, making it practical and scalable for real-world
applications.


\section{Discussion}
\label{section:discussion}

We discuss limitations of our framework and pose future directions of
research. First, implementing the \ourmethod~needs more information than the
previous index-based policies. It necessitates arrival information, $\tlambda$ and
$\{\tp_k\}_{k\in [K]}$, as well as the misclassification probabilities
$Q^n=(\q_{kl}^n)_{k,l\in [K]}$.  In practice, such parameters need to be
estimated on a limited amount of data and estimation errors are unavoidable.

\paragraph{Extension of the queueing model} 
We identify conceptual and analytical challenges in extending our framework to
the multiserver setting. Modeling the extension after~\citet{MandelbaumSt04}
who consider known true classes, we can posit the complete resource pooling
(CRP) condition. This condition requires the limit of the arrival rates to be located in
the outer face of the stability region, and to be uniquely represented as a
maximal allocation of the servers' service capacity. For our setting in
Section~\ref{section: model}, the main challenge is that the CRP condition
will not necessarily hold on the arrival and service rates of the
\textit{predicted} classes.  Because the service rate of each predicted class
in prelimit will be a mixture of the original service rates as $\fmu_l^n$ in Definition~\ref{def: concerned processes for arrival and service}, the CRP
condition on true classes may not be preserved for predicted classes.

Understanding of how prediction error interacts with queueing performance
under general dynamics is an important direction of future research.  For
example, when jobs exhibit abandonment behavior, a suitable adjustment to the
c$\mu$-rule minimizes the long-run holding cost under many-server fluid
scaling~\citep{AtarGiSh10}. Policies that simultaneously account for
predictive error and job impatience may yield fruit.

\paragraph{Design of queueing systems under class uncertainty} 

While we focus on optimal scheduling, an even more important operational lever
is the \emph{design} of the queueing system~\citep{FeldmanMaMaWh08,
  JenningsMaMaWh96}.  For example, designing priority classes that account for
predictive error is a promising research direction~\citep{ChenDo21}.  In our
model, if we keep the limiting distributions of the interarrival and service
times the same, class designs satisfying the heavy traffic
condition~\eqref{eq: heavy traffic condition} should have an identical
limiting total workload $\sumTtW$ by Proposition~\ref{prop: convergence and
  approximation of predicted class N, tau, T, and W}, implying similar forms
of the lower bound in~\eqref{eq: definition of J*}. Given this observation, we
may investigate how class design interacts with the AI model's predictive
performance.

\paragraph{Combining the \ourmethod~with AI-based approaches}
The performance of RL algorithms degrade under distribution shift, and simple
index-based policies may offer robustness benefits. The two approaches may
provide synergies. For example, we can pre-train a policy to initially imitate
an index-based policy, and then fine-tune it to maximize performance in
specific environments.

For quadratic costs, we showed the effectiveness of using the relative regret
$\TJ^*(t;Q)/\TJ^* (t;I)$ to select the classification threshold. Alternatively, we could directly fine-tune the classifier to minimize this metric, in the spirit of the end-to-end queuing-loss training of~\citet{SinghGuVa20}, which may
further enhance downstream queueing performance, albeit at the cost of increased engineering complexity.





\bibliographystyle{abbrvnat}

\ifdefined\useorstyle
\setlength{\bibsep}{.0em}
\else
\setlength{\bibsep}{.7em}
\fi

\bibliography{bib}

\ifdefined\useorstyle

\ECSwitch


\ECHead{Appendix}

\else
\appendix
\newpage
\part*{Appendices} 



\fi

\section{Experiment Details}
\label{appendix: experiment details}

\paragraph{Confusion Matrices}
In Sections~\ref{section: numerical experiments} 
and~\ref{section: characterization of the optimal cost with q}, 
we consider several language models finetuned on the 
CivilComments dataset. They exhibit diverse toxicity prediction accuracy 
across different identity groups, and thus lead to different confusion matrices.
We list the details of these models below, where the actual accuracy is based on 
test set performance and not known to the scheduler, while the estimated accuracy 
is derived from the validation set and used by the scheduler.

\begin{table}[h]
  \centering
  \begin{tabular}{l|ccccc|ccccc}
    \toprule
    & \multicolumn{5}{c|}{Non-Toxic Group, Actual Accuracy} & \multicolumn{5}{c}{Toxic Group, Actual Accuracy} \\
    \cmidrule(lr){2-6} \cmidrule(lr){7-11}
    Model & White & Black & Male & Female & LGBTQ & White & Black & Male & Female & LGBTQ \\
    \midrule
    ERM, 0.05 & 0.601& 0.594& 0.891& 0.883& 0.649 & 0.883& 0.909& 0.853& 0.821& 0.868 \\
    ERM, 0.5 & 0.898& 0.869& 0.959& 0.961& 0.916 & 0.609& 0.677& 0.678& 0.631& 0.546 \\
    ERM, 0.95 & 0.975& 0.967& 0.986& 0.989& 0.986 & 0.394& 0.445& 0.506& 0.460& 0.295 \\
    GroupDRO, 0.05 & 0.597& 0.596& 0.896& 0.886& 0.689 & 0.893& 0.890& 0.817& 0.826& 0.817 \\
    Reweighted, 0.05 &0.661& 0.649& 0.891& 0.883& 0.723 & 0.843& 0.886& 0.832& 0.810& 0.804 \\
    \midrule
    \multicolumn{11}{c}{\rule{0pt}{2.5ex}} \\[-2.5ex]
  \end{tabular}
  \label{tab:model_identity_metrics}
  \caption{Actual Accuracy on toxic and non-toxic comments
  for different identity groups.}
\end{table}

\begin{table}[h]
  \centering
  \begin{tabular}{l|ccccc|ccccc}
    \toprule
    & \multicolumn{5}{c|}{Non-Toxic Group, Estimated Accuracy} & \multicolumn{5}{c}{Toxic Group, Estimated Accuracy} \\
    \cmidrule(lr){2-6} \cmidrule(lr){7-11}
    Model & White & Black & Male & Female & LGBTQ & White & Black & Male & Female & LGBTQ \\
    \midrule
    ERM, 0.05 & 0.593& 0.589& 0.888& 0.882& 0.668 & 0.857& 0.887& 0.859& 0.823& 0.862 \\
    ERM, 0.5 & 0.882& 0.860& 0.959& 0.960& 0.903 & 0.598& 0.657& 0.688& 0.670& 0.547 \\
    ERM, 0.95 & 0.974& 0.962& 0.983& 0.988& 0.973 & 0.402& 0.392& 0.526& 0.466& 0.331 \\
    GroupDRO, 0.05 & 0.581& 0.562& 0.902& 0.891& 0.705 & 0.845& 0.892& 0.842& 0.823& 0.845 \\
    Reweighted, 0.05 & 0.640& 0.621& 0.900& 0.895& 0.728 & 0.806& 0.877& 0.838& 0.809& 0.812 \\
    \midrule
    \multicolumn{11}{c}{\rule{0pt}{2.5ex}} \\[-2.5ex]
  \end{tabular}
  \label{tab:estimated_model_identity_metrics}
  \caption{Estimated accuracy on toxic and non-toxic comments
  for different identity groups.}
\end{table}
\paragraph{Queueing Setup in Figure~\ref{fig: cumulative cost}}
Figure~\ref{fig: cumulative cost} illustrates queueing systems where all jobs have interarrival and service times that are independently and exponentially distributed. Arrival and service rates can vary over time across three queueing settings.

In the base setting used for DRL policy training, the constant arrival rates for toxic comments are $\lambda_\text{toxic} = [4.2, 2.9, 3.4, 5, 2.7]$, and for non-toxic comments $\lambda_\text{non-toxic} = [12.4, 6.1, 22.1, 33.4, 7.8]$, corresponding to the groups: \emph{white, black, male, female, and LGBTQ}. Service rates are also constant, where $\mu_\text{toxic} = [100, 30, 110, 25, 15]$ and $\mu_\text{non-toxic} = [150, 150, 150, 150, 150]$ for the same groups.

In Distribution Shift I, arrival rates for all 10 classes are modeled as time-dependent using a sinusoidal function: $\lambda(t) = \max(\lambda_0/2, \lambda_0 + r \sin(\omega t))$, where $\lambda_0$ is the baseline rate (i.e., the base arrival rate as in the base setting), $r=5$ controls the amplitude, and $\omega=11$ sets the frequency. This captures realistic fluctuations in comment volume due to daily cycles or trending topics. Service rates are constant and the same as the base setting to isolate the impact of arrival rate variability.

In Distribution Shift II, for all 10 classes, arrival rates are held constant as in the base setting, while service rates decrease over time to simulate system degradation. Service rates follow $\mu(t) = \max(d, \mu_0 - kt)$, where $\mu_0$ is the baseline service rate (the base service rate as in the base setting), $k$ determines the rate of decline, and $d$ is the lower bound. For toxic comments, $k=0$ (no degradation); for non-toxic comments, $k=150$ and $d = (2/3)\mu_0$ (moderate degradation). This scenario reflects performance drops due to model drift, resource constraints, or increased workload.

\section{Diffusion limits}
\label{section: proof for fclt}

We consider the following processes: partial sum process of interarrival time
$\tU_{0}^n(t)$ that depends solely on $\{\tu_i^n: i \in \mathbb{N}\}$, partial
sum process of service time $\tV_{0}^n$ and two other processes
$\VfZ^n : = (\fZ_{kl}^n)_{k, l \in [K]}, \VfR^n:= (\fR^n_l)_{l\in [K]}$ that
solely relies on $\{(\feature_i^n, \tY_i^n, \tv_i^n): i \in \mathbb{N}\}$, as well as the diffusion-scaled processes, $\TtU_0^n, \TtV_0^n, \VTfZ^n : = (\TfZ_{kl}^n)_{k, l \in [K]}, \VTfR^n:=
(\TfR^n_l)_{l\in [K]}$, as in Section~\ref{section: convergence and lower bound} and Definition~\ref{definition: U, Z, R, V} to come.

As explained in Section~\ref{section: convergence and lower bound}, our analysis begins by establishing joint diffusion limits of the primitive processes in Lemma~\ref{lemma: joint weak convergence}.
Deferring a detailed proof to
Section~\ref{section:proof-joint-weak-convergence}, we highlight the main
ingredients of the joint convergence result. Our main observation is that the
diffusion-scaled processes admit a martingale central limit result when
$\{(u^n_i, v^n_i, X^n_i, Y^n_i): i \in \mathbb{N}\}$ are i.i.d.
(Assumption~\ref{assumption: data generating process} (i)) and
$v^n_i\perp X^n_i \mid Y^n_i$ (Assumption~\ref{assumption: data generating
  process} (iii)). This allows us to show the weak convergence
$\TtU_0^n \Rightarrow \TtU_0$ in $(\mathcal{D}, J_1)$ and
$(\VTfZ^n, \VTfR^n, \TtV_0^n) \Rightarrow (\VTfZ, \VTfR, \TtV_0)$ in
$(\mathcal{D}^{K(K+1) + 1}, WJ_1)$. Since $\{\tu^n_i\}$ and
$\{(\tv_i^n, \feature_i^n, \tY_i^n)\}$ are independent
(Assumption~\ref{assumption: data generating process} (ii)), we can obtain the desired joint convergence (e.g.,
see~\citet[Theorem 11.4.4]{Whitt02} which we give as Lemma~\ref{lemma: joint convergence under indep}).

Building off of our diffusion limit, we can strengthen the convergence to the uniform topology using standard tools (e.g., see Lemma~\ref{lemma: Skorohod representation} and Lemma~\ref{lemma: equivalent to uniform convergence}), and conduct a sample path analysis where
we construct \emph{copies} of $(\TtU_0^n, \VTfZ^n, \VTfR^n, \TtV_0^n)$ and
$(\TtU_0, \VTfZ, \VTfR, \TtV_0)$ that are identical in distribution with their
original counterparts and converge almost surely under a common probability
space. Abusing notation, we use the same notation for the newly construced
processes.
\begin{lemma}[Uniform convergence] \label{lemma: uniform convergence of UZRV}
  Suppose that Assumptions~\ref{assumption: data generating
    process},~\ref{assumption: heavy traffic}, and~\ref{assumption: second
    order moments} hold.  Then, there exist stochastic processes
  $(\TtU_0^n, \VTfZ^n, \VTfR^n, \TtV_0^n),~\forall~n\geq 1$ and
  $(\TtU_0, \VTfZ, \VTfR, \TtV_0)$ defined on a common probability space
  $(\Omega_\text{copy}, \mathcal{F}_\text{copy}, \mathbb{P}_\text{copy})$
  such that $(\TtU_0^n, \VTfZ^n, \VTfR^n, \TtV_0^n),~\forall~n\geq 1$ and
  $(\TtU_0, \VTfZ, \VTfR, \TtV_0)$ are identical in distribution with their
  original counterparts and
\begin{equation}\label{eq: uniform convergence}
\begin{aligned}
  (\TtU_0^n, \VTfZ^n, \VTfR^n, \TtV_0^n) \rightarrow (\TtU_0, \VTfZ, \VTfR, \TtV_0) \quad  \text{\quad in $(\mathcal{D}^{K(K+1)+2}, \| \cdot \|)$,} \quad
  \mathbb{P}_\text{copy}\text{-a.s.}.
\end{aligned}
\end{equation}
\end{lemma}
\noindent We defer a detailed proof to Section~\ref{subsection: proof of
  uniform convergence} since it is a basic consequence of the Skorohod
representation theorem~\citep[Theorem 6.7]{Billingsley99}. Since the diffusion
limits $(\TtU_0, \VTfZ, \VTfR, \TtV_0)$ are multidimensional Brownian motions,
the copied processes of $(\TtU_0^n, \VTfZ^n, \VTfR^n, \TtV_0^n)$ jointly
converge to a continuous limit almost surely. We obtain the
result by noting that convergence in $J_1$ to a deterministic and continuous
limit is equivalent to uniform convergence on compact intervals (e.g.,
see~\citet[Proposition 4]{Glynn90} stated in Lemma~\ref{lemma: equivalent
  to uniform convergence}).

Sample path analysis allows us to leverage properties of uniform convergence
and significantly simplifies our analysis. All subsequent results and their
proofs in the Electronic Companions, will be established on the copied processes in the
common probability space
$(\Omega_\text{copy}, \mathcal{F}_\text{copy}, \mathbb{P}_\text{copy})$  with probability one, i.e., $\mathbb{P}_\text{copy} \text{-}a.s.$, and
all of the convergence results will be understood to hold in the
\textit{uniform} norm $\| \cdot \|$. 
Moreover, since these newly constructed processes are
identical in distribution with their original counterparts, all subsequent
results regarding almost sure convergence for the copied processes can be 
converted into corresponding weak convergence results for the original processes; 
see more discussion in Theorems~\ref{theorem: HT lower bound} 
and~\ref{theorem: optimality of our policy}.


We use the uniform integrability condition in Assumption~\ref{assumption: second order moments} to apply the martingale FCLT. For completeness, we review the martingale FCLT and Skorohod representation result before proving the main results of Section~\ref{section: model}.

\subsection{Review of basic results}
\label{subsection:review}

We review classical results on the martingale FCLT and the Skorohod construction.

\subsubsection{Martigale Functional Central Limit Theorem}

Our proof of Lemma~\ref{lemma: joint weak convergence} primarily relies on the
martingale FCLT~\cite[Theorem 8.1]{PangTaWh07}. We define the maximum jump and
the optional quadratic variation of processes and review the martingale FCLT
in Lemma~\ref{lemma: Martingale FCLT}.

Let $\mathcal{D}_{[0, \infty)}:= \mathcal{D}([0,\infty), \R)$ be the set of
right-continuous with left limits (RCLL) functions $[0,\infty) \to \R$, and
$\mathcal{D}^k_{[0, \infty)}:= \mathcal{D}([0,\infty), \R^k)$ be the product
space $\mathcal{D}_{[0, \infty)}\times\cdots\times \mathcal{D}_{[0, \infty)}$
for $k\in \N$.  With a slight abuse of notations, we also use $J_1$ to denote
be the standard Skorohod $J_1$ topology on $\mathcal{D}_{[0, \infty)}$ and
$WJ_1$ to denote the product $J_1$ topology on $\mathcal{D}_{[0, \infty)}^k$.
\begin{definition}[Maximum jump]\label{def: max jump}
For any function $x\in \mathcal{D}_{[0, \infty)}$, the  maximum jump of $x$ up to time $t$ is represented as  
\begin{equation}
\label{eq: maximum jump}
J(x, t) := \sup\{|x(s) - x(s^-)|: 0< s \leq t\}, \quad t > 0.
\end{equation}
\end{definition}
\begin{definition}[Optional quadratic variation]\label{def: quadratic variation}
  Let $M_1$ and $M_2$ be two martingales in $\mathcal{D}_{[0, \infty)}$ with
  respect to a filtration
  $\mathcal{F} \equiv\{\mathcal{F}_{t}: t \geq 0\}$ satisfying
  $M_1(0)=M_2(0)=0$. The optional quadratic variation between $M_1$ and $M_2$
  is defined as
\begin{equation}
\label{eq: optional quadratic variation}
[M_1, M_2](t) = \lim\limits_{m\rightarrow \infty} \sum\limits_{i=1}^{\infty} \big(M_1 (t_{m,i}) - M_1 (t_{m, i-1}) \big) \big(M_2 (t_{m,i}) - M_2(t_{m, i-1}) \big),\quad t > 0,
\end{equation}
where $t_{m,i}= \min (t, i2^{-m})$. 
\end{definition}
\noindent
\citet[Theorem 3.2]{PangTaWh07} shows that $[M_1, M_2](t)$ is well-defined
for any martingales pairs $(M_1, M_2)$ satisfying conditions outlined in
Definition~\ref{def: quadratic variation}.

\begin{lemma}[Multidimensional martingale FCLT]
  \label{lemma: Martingale FCLT}
  For $n \geq 1$, let $\mathbf{M}^n \equiv(M^n_1, \ldots, M^{n}_k)$ be a
  martingale in $(\mathcal{D}_{[0, \infty)}^k, WJ_1)$ with respect to a
  filtration $\mathcal{F}_n \equiv\left\{\mathcal{F}_{n, t}: t \geq 0\right\}$
  satisfying $\mathbf{M}^n(0)=$ $(0, \ldots,
  0)$. 
  If both of the following conditions hold
  \begin{enumerate}[(i)]
  \item the expected maximum jump is asymptotically negligible:
    $\lim _{n \rightarrow \infty}\mathbb{E} [J (M^n_i, T)]=0,
    ~\forall~i\in[k],~\forall~T\geq 0;$
  \item there exists a positive semidefinite symmetric matrix
    $ \mathbf{A} = \{a_{ij}\}_{i,j \in [k]} \in \mathbb{R}^{k\times k}$ such
    that for any $1 \leq i , j \leq k$ and $t>0$,
    $[M^n_{i}, M^n_{j}](t) \Rightarrow a_{ij} t$ in $\mathbb{R}$ as
    $n\rightarrow \infty,$
  \end{enumerate}
  then, we have that 
  $$ \mathbf{M}^n \Rightarrow \mathbf{M} \quad \text { in } (\mathcal{D}_{[0, \infty)}^k, WJ_1) \quad \text { as } n \rightarrow \infty,$$
  where $\mathbf{M}$ is a $k$-dimensional Brownian motion with mean vector and
  covariance matrix being
  $$\mathbb{E}[\mathbf{M}(t)]=(0, \ldots, 0) \quad \text { and } \quad 
  \mathbb{E} [\mathbf{M}(t) \mathbf{M}^\top(t)]=\mathbf{A} t, \quad t \geq 0.$$
\end{lemma}

\subsubsection{Skorohod representation}
\label{section: Skorohod representation}
Recall the definition of random elements on a metric space $(S, m)$~\cite[Page
78]{Whitt02}.
\begin{definition}[Random Element]\label{definition: random element}
  For a separable metric space $(S, m)$, we say that $\mathbf{X}$ is a random
  element of $(S, m)$ if $\mathbf{X}$ is a measurable mapping from some
  underlying probability space $(\Omega, \mathcal{F}, \mathbb{P})$ to
  $(S, \mathcal{B}(S))$, where $\mathcal{B}(s)$ is the Borel $\sigma$-field
  induced by $(S, m)$.
\end{definition}
The well-known Skorohod representation theorem~\citep[Theorem
6.7]{Billingsley99} gives the following.
\begin{lemma}[Skorohod representation]
\label{lemma: Skorohod representation}
Let $\{\mathbf{X}^n\}_{n\geq 1}$ and $\mathbf{X}$ be random elements of a
separable metric space $(S, m)$. If $\mathbf{X}^n \Rightarrow \mathbf{X}$ 
in $(S, m)$, then there exists other random elements
$\{\mathbf{X}^n_\text{copy}\}_{n\geq 1}$ and $\mathbf{X}_\text{copy}$ of
$(S, m)$, defined on a common probability space
$(\Omega, \mathcal{F}, \mathbb{P})$, such that (i)
$\mathbf{X}^n_\text{copy} \stackrel{d}{=} \mathbf{X}^n,~\forall~n\geq 1$ and
$\mathbf{X}_\text{copy}\stackrel{d}{=} \mathbf{X}$; (ii)
$\lim_{n\rightarrow +\infty} m(\mathbf{X}^n_\text{copy},
\mathbf{X}_\text{copy}) = 0$ $\mathbb{P}$-almost surely.
\end{lemma}

Let $d_{J_1}(\cdot, \cdot)$ be the $J_1$ metric (Skorohod metric) defined on
$\mathcal{D}:=\mathcal{D}([0, 1], \mathbb{R})$, the set of RCLL functions $[0,1] \rightarrow \mathbb{R}$~\cite[Page 79]{Whitt02}. Moreover, for the
product space $\mathcal{D}^k:= \mathcal{D}\times \cdots \times \mathcal{D}$, let $d_p (\cdot, \cdot)$ be the
product metric defined by
$d_p(\mathbf{x}, \mathbf{y}) : = \sum_{i=1}^K d_{J_1}(x_i, y_i),
~\forall~\mathbf{x}, \mathbf{y} \in D^k$~\cite[Page 83]{Whitt02}. It is known
that both $(\mathcal{D}, d_{J_1}(\cdot, \cdot))$ and $(\mathcal{D}^k, d_p(\cdot, \cdot))$ are
separable metric spaces with $J_1$ topology and $WJ_1$ (weak $J_1$) topology, respectively~\cite[Sections 3.3, 11.4, and 11.5]{Whitt02}.  Then, according to Lemma~\ref{lemma: Skorohod representation}, for weakly converging random elements, we can obtain copies
that converges almost surely. This enables us to conduct sample path
analysis. Specifically, if the limiting random element is continuous almost
surely, we can utilize the following theorem from~\cite[Proposition 4]{Glynn90}
to conduct analysis under uniform norm convergence, which can streamline our
analysis significantly.
\begin{lemma}
  \label{lemma: equivalent to uniform convergence}
  For a sequence of functions $X^n \in \mathcal{D}$,
  convergence to a continuous function, say
  $X \in \mathcal{C}$, in the $J_1$ metric
  $d_{J_1}(\cdot, \cdot)$ is equivalent to convergence in uniform norm
  $\| \cdot \|$, i.e.,
\begin{equation*} 
\lim_{n\rightarrow \infty} d_{J_1}(X^n, X) = 0 \quad 
\iff \quad 
\lim_{n\rightarrow \infty}\| X^n - X \| = 0.
\end{equation*}    
\end{lemma}

\subsection{Proof of Lemma~\ref{lemma: joint weak convergence}}
\label{section:proof-joint-weak-convergence}

First, we define arrival and service processes of predicted classes on which
we apply the martigale FCLT to establish their weak convergence.
\begin{definition}[Arrival and service processes of predicted classes I]
  \label{definition: U, Z, R, V}
  Given a classifier $\model$ and a sequence of queueing systems, we define
  the following for a given system $n$ and time $t\in [0, n]$:
  \begin{enumerate}[(i)]
\item (Counting process) for any real class $k\in [K]$ and predicted class
$l\in [K]$, let $\fZ_{kl}^n (t)$ be the total number of jobs from real class
$k$ predicted as class $l$, among the first $\lfloor t\rfloor$ jobs
arriving in the system, i.e.,
\begin{equation*}
\fZ_{kl}^n (t): = \sum_{i=1}^{\lfloor t \rfloor} \tY_{ik}^n \fY_{il}
^n,~\forall~t\in [0, n];
\end{equation*}
Moreover, let $\VTfZ^n = \{\TfZ_{kl}^n\}_{k,l\in[K]}$ be the corresponding
diffusion-scaled process, defined as
\begin{equation*}
\TfZ_{kl}^n (t) = n^{-\frac{1}{2}} \Big[ \sum_{i=1}^{\lfloor nt \rfloor} 
Y_{ik}^n \fY_{il}^n - p_k^n \q_{kl}^n \cdot nt\Big],~\forall~t\in [0, 1];
\end{equation*}
\item (Cumulative service time) for any predicted class $l\in [K]$, let
$\fR_l^n$ be the total service time requested by jobs predicted as class
$l$, among the first $\lfloor t\rfloor$ jobs arriving in the system, i.e.,
\begin{equation*}
\fR_{l}^n (t) := \sum_{i=1}^{\lfloor t \rfloor}\fY_{il}^n \tv_i^n, 
~\forall~t\in [0,n].
\end{equation*}
Moreover, let $\VTfR = \{\TfR_{l}\}_{l\in[K]}$ be the corresponding
diffusion-scaled process, defined as
\begin{equation*}
\TfR_l^n (t) 
=   n^{-\frac{1}{2}} \Big[\sum_{i=1}^{\lfloor nt \rfloor}  \fY_{il}^n v_i^n
- \sum\limits_{k=1}^K \frac{p_k^n}{\mu_k^n} \q_{kl}^n \cdot nt \Big],
\quad ~\forall~t\in [0,1].
\end{equation*}
\end{enumerate}
\end{definition}

We define $\fZ_{kl}^n$ and $\fR^n_l$ on $[0, n]$, and
$\TfZ_{kl}^n$ and $\TfR^n_l$ on $[0, 1]$ for analysis simplicity, and 
these processes can be naturally extended to $[0, +\infty)$ to apply the martingale FCLT
in Lemma~\ref{lemma: Martingale FCLT}. In addition,
we introduce the following rescaled and centered processes $\BrtU^n_0(t)$ and
$(\VBrfZ^n, \VBrfR^n, \BrtV_0^n)$ to facilitate the analysis.
\begin{definition}[Arrival and service processes of predicted classes
II] \label{definition: breve U, Z, R, V} Given a classifier $\model$ and a
sequence of queueing systems, we define the rescaled and centered processes 
for a given system $n$ and time $t\in [0, 1]$ as followings:
\begin{equation*}
\begin{aligned}
& \BrtU_0^n (t) =  n^{-\frac{1}{2}} \sum_{i=1}^{\lfloor nt \rfloor} (u_i^n - (\lambda^n)^{-1}), \quad &&
\BrtV_0^n (t) =n^{-\frac{1}{2}}\sum_{i=1}^{\lfloor nt \rfloor} 
\Big[v_i^n - \sum_{k=1}^K \frac{p_k^n}{\mu_k^n}\Big],\\
&\BrfZ_{kl}^n (t) = n^{-\frac{1}{2}} \sum_{i=1}^{\lfloor nt \rfloor} 
[Y_{ik}^n \fY_{il}^n - p_k^n \q_{kl}^n],~\forall~k, l\in[K],
&&\BrfR_l^n (t) =n^{-\frac{1}{2}}\sum_{i=1}^{\lfloor nt \rfloor} 
\Big[\fY_{il}^n v_i^n - \sum\limits_{k=1}^K \frac{p_k^n}{\mu_k^n} \q_{kl}^n\Big],
~\forall~l\in[K].
\end{aligned}
\end{equation*}
\end{definition}
\noindent
One can check that $\BrtU_0^n , \BrtV_0^n, \VBrfZ^n, \VBrfR^n $ are closely
related to $\TtU_0^n , \TtV_0^n, \VTfZ^n, \VTfR^n$ by noting that for any
$t\in[0, 1]$,
\begin{equation*} 
\begin{aligned}
  &\TtU_{0}^n (t)= \BrtU_0^n (t) 
    + n^{-\frac{1}{2}}(\lambda^n)^{-1} (\lfloor nt \rfloor - nt),
  &&\TtV_0^n (t) =\BrtV_0^n (t) 
     + n^{-\frac{1}{2}}\sum\limits_{k=1}^K \frac{\tp_k^n}{\tmu_k^n} (\lfloor nt \rfloor - nt),\\
  &\TfZ_{kl}^n (t) = \BrfZ_{kl}^n (t) + n^{-\frac{1}{2}} p_k^n \q_{kl}^n(\lfloor nt \rfloor - nt),
  &&\TfR_l^n (t) =\BrfR_l^n (t) + n^{-\frac{1}{2}} \sum\limits_{k=1}^K \frac{p_k^n}{\mu_k^n}\q_{kl}^n(\lfloor nt \rfloor - nt).
\end{aligned}
\end{equation*}

Under Assumptions~\ref{assumption: data generating process}
and~\ref{assumption: second order moments}, we establish the weak convergence
of $\BrtU^n_0(t)$ and $(\VBrfZ^n, \VBrfR^n, \BrtV_0^n)$ using the martingale FCLT in Lemma~\ref{lemma: Martingale FCLT}.
\begin{lemma}[Individual weak convergence]
  \label{lemma: individual weak convergence}
  Suppose that Assumptions~\ref{assumption: data generating
    process},~\ref{assumption: heavy traffic}, and~\ref{assumption: second order
    moments} hold. Then, there exist Brownian motions $\BrtU_0$ and
  $(\VBrfZ, \VBrfR, \BrtV_0)$ such that (i) $\BrtU_0^n \Rightarrow \BrtU_0$ in
  $(\mathcal{D}, J_1)$; (ii)
  $(\VBrfZ^n, \VBrfR^n, \BrtV_0^n) \Rightarrow (\VBrfZ, \VBrfR, \BrtV_0)$ in
  $(\mathcal{D}^{K(K+1) + 1}, WJ_1)$.
\end{lemma}
\noindent We defer a detailed proof of the lemma to Section~\ref{subsection:
  proof of individual weak convergence}.

The following processes are all well-defined deterministic functions of 
$t \in [0, 1]$:
\begin{equation*}
  n^{-\frac{1}{2}}(\lambda^n)^{-1} (\lfloor nt \rfloor - nt),
  ~~ n^{-\frac{1}{2}}\sum\limits_{k=1}^K \frac{\tp_k^n}{\tmu_k^n}
  (\lfloor nt \rfloor - nt), 
  ~~ n^{-\frac{1}{2}} p_k^n \q_{kl}^n(\lfloor nt \rfloor - nt),
  ~~ n^{-\frac{1}{2}} \sum\limits_{k=1}^K \frac{p_k^n}{\mu_k^n}\q_{kl}^n(\lfloor nt \rfloor - nt).
\end{equation*}
Assumption~\ref{assumption: heavy traffic} and
$n^{-1/2} \sup_{t\in [0,1]}(\lfloor nt \rfloor - nt) \rightarrow 0$ imply that
all of them converge to $0$ 
in $(\mathcal{D}, J_1)$. Using the jointly weak convergence
with a deterministic limit~\cite[Theorem 11.4.5]{Whitt02}, continuity of
addition~\cite[Theorem 4.1]{Whitt02} by almost-sure continuity of all limits,
and the continuous mapping theorem, it follows that there exist Brownian
motions $\TtU_0$ and $(\VTfZ, \VTfR, \TtV_0)$ such that (i)
$\TtU_0^n \Rightarrow \TtU_0$ in $(\mathcal{D}, J_1)$; (ii)
$(\VTfZ^n, \VTfR^n, \TtV_0^n) \Rightarrow (\VTfZ, \VTfR, \TtV_0)$ in
$(\mathcal{D}^{K(K+1) + 1}, WJ_1)$.

Since $\{\tu^n_i\}$ and $\{(\tv_i^n, \feature_i^n, \tY_i^n)\}$ are independent
(Assumption~\ref{assumption: data generating process} (ii)), we can use the following
result~\cite[Theorem 11.4.4]{Whitt02} to obtain our desired jointly weak convergence in Lemma~\ref{lemma: joint weak convergence}.
\begin{lemma}[Joint weak convergence for independent random elements]
  \label{lemma: joint convergence under indep}
  Let $\mathbf{X}^n$ and $\mathbf{Y}^n$ be independent random elements of
  separable metric spaces $(S^{\prime}, m^{\prime})$ and
  $(S^{\prime \prime}, m^{\prime \prime})$ for each $n \geq 1$. Then, there is
  joint convergence in distribution
  $$
  (\mathbf{X}^n, \mathbf{Y}^n) \Rightarrow(\mathbf{X}, \mathbf{Y}) \text { in }
  S^{\prime} \times S^{\prime \prime}
  $$
  if and only if $\mathbf{X}^n \Rightarrow \mathbf{X}$ in $S'$ and
  $\mathbf{Y}^n \Rightarrow \mathbf{Y}$ in $S''$.
\end{lemma}


\subsubsection{Proof of Lemma~\ref{lemma: individual weak convergence}}
\label{subsection: proof of individual weak convergence}

To utilize the martingale FCLT (Lemma~\ref{lemma: Martingale FCLT}), we extend
$\BrtU^n_0$ and $(\VBrfZ^n, \VBrfR^n, \BrtV^n_0)$ to
$\mathcal{D}_{[0, \infty)}$ and $\mathcal{D}_{[0, \infty)}^{K(K+1) + 1}$,
respectively, and establish individual weak convergence for these extended
stochastic processes. We can get the desired result by restricting the
extended stochastic processes to the time interval $[0, 1]$.

We establish weak convergence of the extended $\BrtU^n_0$ and
$(\VBrfZ^n, \VBrfR^n, \BrtV^n_0)$ separately.  To show the former, note that
$\{u_i^n: i \geq 1\}$ are i.i.d. random variables with mean $(\lambda^n)^{-1}$
by Assumptions~\ref{assumption: data generating process}.  Evidently,
$\{\BrtU_0^n: n\in \mathbb{N}\}$ is a martingale with respect to the natural
filtration and satisfies $\BrtU_0^n (0)=0$.  It thus suffices to validate the
conditions (i) and (ii) of Lemma~\ref{lemma: Martingale FCLT}.  To verify
condition (i), use the shorthand
$\Delta^n_i \defeq |u^n_i - (\lambda^n)^{-1}|$ to write
\begin{equation*}
  \begin{aligned}
    \E^n[J(\BrtU_0^n, t)^2]
    = n^{-1} \E^n \big[
    \max\limits_{1\leq i \leq \lfloor nt\rfloor} (\Delta^{n}_i)^2\big]
    \le
    \mathbb{E}^n \big[ 
     (\Delta^{n}_1)^2
    \indic{(\Delta^{n}_1)^2 \ge \sqrt{n}}\big] + \frac{1}{\sqrt{n}}.
\end{aligned}
\end{equation*}
From uniform integrability (Assumption~\ref{assumption: second order
  moments}), we have
$\mathbb{E}^n[J(\BrtU_0^n, t)] \leq \big(\mathbb{E}^n \big[|J(\BrtU_0^n,
t)|^2\big]\big)^{1/2} \to 0$. To verify condition (ii), first truncate the triangular array $\{|\tu_{i}^n - (\lambda^n)^{-1}|^2: i \in \mathbb{N}, n\in \mathbb{N}\}$ uniformly with a constant using the uniform integrability (Assumption~\ref{assumption: second order moments}), and apply the triangular weak law of large numbers (WLLN)~\cite[Theorem 2.2.6]{Durrett10} on the truncated array, with a choice of $b_n:=n$ in that theorem, to obtain 
\begin{equation*}
  [\BrtU_0^n, \BrtU_0^n] (t) = n^{-1}\sum_{i=1}^{\lfloor nt
    \rfloor} (u_i^n - \tlambda_{n}^{-1})^2 \cp c_u t
  ~~\mbox{where}~~c_u:=\lim_{n\rightarrow \infty} \var(u_1^n) = \alpha_{u} - (\lambda)^{-2}.
\end{equation*}

We now show the weak convergnece of
$\mathbf{G}^n:=(\VBrfZ^n, \VBrfR^n, \BrtV_0^n)$. We have
$\mathbf{G}^n(0) = \mathbf{0}$ and by Assumption~\ref{assumption: data
  generating process},
$\{(\tY_{i}^n, \feature_i^n, \tv_i^n): i\in \mathbb{N}\}$ are i.i.d. and
$X^n_i$ is independent of $v^n_i$ given $Y^n_i$. Therefore, by conditioning on
$\tY_i^n$, we have that for all $i \geq 1$,
\begin{equation*}
  \begin{aligned}
    \mathbb{E}^n [Y_{ik}^n \fY_{il}^n] = p_k^n q_{kl}^n, \quad 
    \mathbb{E}^n [\fY_{il}^n v_i^n] = \sum\limits_{k=1}^K \frac{p_k^n}{\mu_k^n} q_{kl}^n, \quad
    \mathbb{E}^n[v_i^n] = \sum\limits_{k=1}^K \frac{p_k^n}{\mu_k^n},
  \end{aligned}
\end{equation*}
which indicates that $\mathbf{G}^n$ is a martingale with respect to the
natural filtration. To apply Lemma~\ref{lemma: Martingale FCLT} towards
$\mathbf{G}^n$, we now validate its conditions (i) and (ii). Using a similar argument as above, uniform integrability yields condition (i) of Lemma~\ref{lemma: Martingale FCLT}
\begin{equation}
\label{eq: convergence of corrected R}
\begin{aligned}
  \mathbb{E}^n[|J(\BrtV_0^n, t)|] \rightarrow 0,
  \qquad \mathbb{E}^n[J(\BrfZ_{kl}^n, t)] \rightarrow 0,
   \qquad \mathbb{E}^n[J(\BrfR_l^n, t)] \rightarrow 0
\end{aligned}
\end{equation}
for all $k, l \in [K]$.
Similarly, the triangular WLLN gives condition (ii)
\begin{equation*}
  \begin{aligned}
    &[\BrtV_0^n, \BrtV_0^n](t) \Rightarrow c_v t, \quad
    &&[\BrfZ_{kl}^n, \BrfZ_{rs}^n](t) \Rightarrow c_{(k, l), (r, s)} t,  \quad
    &&[\BrfR_{l}^n, \BrfR_{s}^n](t) \Rightarrow c_{l, s} t,\\
    &[\BrtV_{0}^n, \BrfZ_{kl}^n](t) \Rightarrow c_{0, k, l} t,  \quad
    &&[\BrtV_{0}^n, \BrfR_{l}^n](t) \Rightarrow c_{0, l} t, \quad
    &&[\BrfZ_{kl}^n, \BrfR_{s}^n](t) \Rightarrow c_{k, l, s} t,
  \end{aligned}
\end{equation*}
where
\begin{equation*}
  \begin{aligned}
    &c_v  \defeq \sum_{k=1}^K p_k \alpha_{v,k}
    - (\sum_{k=1}^K p_k/\mu_k)^2,
    &&c_{(k,l), (r,s)} \defeq \Big\{
    \begin{array}{lr}
      p_k \q_{kl} (1 - p_k \q_{kl})
      &\text{if } (k,l) = (r,s),\\
      -p_k \q_{kl} p_r \q_{rs}
      & \text{if }(k,l) \neq (r,s),
    \end{array}\\
    &c_{l,s}  \defeq 
    \begin{cases}
      \sum\limits_{k=1}^K p_k \q_{kl} \alpha_{v, k}
      - \big(\sum\limits_{k=1}^K \frac{p_k \q_{kl}}{\mu_k}\big)^2
      &\text{if } l = s, \\
      -\big(\sum\limits_{k=1}^K \frac{p_k \q_{kl}}{\mu_k} \big)
      (\sum\limits_{k=1}^K \frac{p_k \q_{ks}}{\mu_k} \big)
      & \text{if }l\neq s,
    \end{cases}
    && c_{0,k,l} \defeq \sum_{k=1}^K \frac{p_k \q_{kl}}{\mu_k} 
    -  \Big(\sum_{k=1}^K \frac{p_k}{\mu_k}\Big) \Big(\sum_{k=1}^K  p_k \q_{kl}\Big),\\
    & c_{0,l} \defeq \sum_{k=1}^K p_k \q_{kl}\alpha_{\tv,k} 
    -  \Big(\sum_{k=1}^K \frac{p_k}{\mu_k}\Big) \Big(\sum_{k=1}^K \frac{p_k \q_{kl}}{\mu_k}\Big),
    &&c_{k,l,s} \defeq 
    \begin{cases}
      \sum\limits_{k=1}^K  \frac{p_k \q_{kl} }{\mu_k} -  
      \big(\sum\limits_{k=1}^K p_k \q_{kl}\big)
      (\sum\limits_{k=1}^K \frac{p_k \q_{kl}}{\mu_k} \big)
      &\text{if } l = s,\\
      -\big(\sum\limits_{k=1}^K p_k \q_{kl}\big)
      (\sum\limits_{k=1}^K \frac{p_k \q_{ks}}{\mu_k} \big)
      & \text{if }l\neq s.
    \end{cases}
\end{aligned}
\end{equation*}

\subsection{Proof of Lemma~\ref{lemma: uniform convergence of UZRV}}
\label{subsection: proof of uniform convergence}

From the Skorohod representation (Lemma~\ref{lemma: Skorohod representation}), there exist stochastic processes defined on some common probability space
$(\Omega_\text{copy}, \mathcal{F}_\text{copy}, \mathbb{P}_\text{copy})$, 
$(\TtU_0^n, \VTfZ^n, \VTfR^n, \TtV_0^n),~\forall~n\geq 1$ and
$(\TtU_0, \VTfZ, \VTfR, \TtV_0)$, such
that $(\TtU_0^n, \VTfZ^n, \VTfR^n, \TtV_0^n)$ and
$(\TtU_0, \VTfZ, \VTfR, \TtV_0)$ are identical in distribution to their
original counterparts and
\begin{equation*}
  \label{eq: a.s. joint convergence of U_0, Z, R, V_0}
  \begin{aligned}
    (\TtU_0^n, \VTfZ^n, \VTfR^n, \TtV_0^n)
    \rightarrow (\TtU_0, \VTfZ, \VTfR, \TtV_0)
    \quad  \text{\quad in $(\mathcal{D}^{K(K+1)+2}, WJ_1)$,} \quad
    \mathbb{P}_\text{copy}\text{-a.s.}.
  \end{aligned}
\end{equation*}
Or equivalently,  with probability one
\begin{equation*}
  d_p\big((\TtU_0^n, \VTfZ^n, \VTfR^n, \TtV_0^n),
  (\TtU_0, \VTfZ, \VTfR, \TtV_0)\big) \rightarrow 0,
\end{equation*}
where $d_p(\cdot, \cdot)$ is the product $J_1$ metric. By definition of
$d_p(\cdot, \cdot)$, each coordinate of
$(\TtU_0^n, \VTfZ^n, \VTfR^n, \TtV_0^n)$ converges to the limiting process in
$(\mathcal{D}, J_1)$ $\mathbb{P}_\text{copy}$-a.s.. Since
$(\TtU_0, \VTfZ, \VTfR, \TtV_0)$ is a multidimensional Brownian motion,
$(\TtU_0, \VTfZ, \VTfR, \TtV_0)$ is continuous
$\mathbb{P}_\text{copy}$-a.s.. By Lemma~\ref{lemma: equivalent to uniform
  convergence}, $\mathbb{P}_\text{copy}$-almost surely, every coordinate of
$(\TtU_0^n, \VTfZ^n, \VTfR^n, \TtV_0^n)$ converges to the limiting process in
$(\mathcal{D}, \| \cdot \|)$. This completes our proof.


\section{Proofs of results in Section~\ref{subsection: fundamental convergence results}}
\label{section: fundamental results}

We show convergence diffusion-scaled versions of the exogenous processes associated with predicted classes in Section~\ref{subsection: appendix for arrival and service}. Then, we provide a sequence of interim results required for us to prove Proposition~\ref{prop: convergence and approximation of predicted class N, tau, T, and W} in Section~\ref{subsection: proof for convergence and approximation of predicted class N, tau, T, and W}.

We begin by extending Lemma~\ref{lemma: uniform convergence of UZRV} to include the arrival process.  For any system $n$, let $A_0^n (t) := \max\{m: \tU_0^n (m) \leq t\},~\forall~t\in[0, n]$ be the total number of jobs that arrive in the system up to time $t$, and
\begin{equation} \label{eq: definitions of the scaled A0}
   \TtA_0^n (t) = n^{-1/2} \big[\tA_{0}^n (nt) -  \tlambda^n nt\big],
   ~ t\in [0, 1].
\end{equation}
By Lemma~\ref{lemma: uniform convergence of UZRV},
$\TtU^n_0 \rightarrow \TtU_0 \in \mathcal{C}$; since the limiting function is
continuous, convergence in weak $M_2$ topology is equivalent to convergence in
uniform metric~\cite[Corollary 12.11.1]{Whitt02}. Using the asymptotic
equivalence between counting and inverse processes with
centering~\cite[Corollary 13.8.1]{Whitt02}, convergence of $\TtA_0^n$ follows
from convergence of $\TtU_0^n$.
\begin{lemma}[Uniform convergence II] 
\label{lemma: convergence of the scaled primitive processes}
Suppose that Assumptions~\ref{assumption: data generating process},~\ref{assumption: 
heavy traffic}, and~\ref{assumption:
second order moments} hold. Then, there exists a multidimensional Brownian motion $(\TtA_0, \TtU_0, \VTfZ, \VTfR, \TtV_0)$  such that
\begin{equation}\label{eq: a.s. conv of the primitives}
\begin{aligned}
    (\TtA_0^n, \TtU_0^n, \VTfZ^n, \VTfR^n, \TtV_0^n) \rightarrow (\TtA_0, \TtU_0, \VTfZ, \VTfR, \TtV_0) \text{\quad in $(\mathcal{D}^{K(K+1)+3}, \| \cdot \|) $,} \quad
    \mathbb{P}_\text{copy}\text{-a.s.}
 \end{aligned}
\end{equation}
\end{lemma}

\subsection{Convergence of arrival and service processes of predicted classes}
\label{subsection: appendix for arrival and service}
We formally define the arrival and service processes associated predicted
classes, and provide corresponding diffusion limits in Proposition~\ref{prop:
  joint conv. of predicted class A, U, S, and V}. Given a classifier $\model$,
suppose that Assumption~\ref{assumption: heavy traffic} holds and consider
system $n$ operating in $t\in [0, n]$. For $j \in \mathbb{N},$ let $\fu_{l, j}^n$ and
$\fv_{l, j}^n$ be the interarrival and service times of the $j$th arriving job in predicted class $l \in [K]$.
\begin{definition}[Arrival and service processes of predicted classes II]
  \label{def: concerned processes for arrival and service}
  \hspace{23pt}
\begin{enumerate}[(i)]
\item (Arrival Process) Let
  $ \fA_{kl}^n (t): = \sum_{i=1}^{A^n_0(t)} \tY_{ik}^n \fY_{il} ^n$ be the
  number of jobs from real class $k$ predicted as class $l$ among jobs
  arriving up to time $t \in [0, n]$,
   and
  $\TfA_{kl}^n (t) = n^{-1/2} [\fA_{kl}^n (nt) - n\lambda^n \tp_k^n \q_{kl}^n t ]~t \in [0,
  1]$ be the diffusion-scaled process.  Let
  $\fA_l^n (t) := \sum_{k=1}^K \fA_{kl}^n (t)= \sum_{i=1}^{A^n_0(t)}
  \fY^n_{il}$ be the number of jobs predicted as class $l$ among jobs arriving
  up to time $t\in [0,n]$, 
  and
  $\TfA_{l}^n (t) := \sum_{k=1}^K \TfA_{kl}^n (t) = n^{-1/2} \Big[\fA_{l}^n
  (nt) - n \flambda^n_l t \Big]$ with $t \in [0, 1]$ be the diffusion-scaled
  process. Here, $\flambda^n_l:= \tlambda^n \fp_l^n$ and $\flambda_l:=\tlambda \fp_l$, where the occurrence of predicted class $l$ and its limit are denoted by
  $\fp_l^n:=\sum_{k=1}^K \tp_k^n \q_{kl}^n$ and
  $\fp_l:=\sum_{k=1}^K \tp_k \q_{kl}$.
\item (Sum of Interarrival Time) Let
  $\fU_l^n (t):= \sum_{j=1}^{\lfloor t \rfloor} \fu_{l, j}^n, ~t\in [0,n]$ be
  the sum of interarrival times among the first $\lfloor t\rfloor$ jobs
  predicted as class $l$,
    and
  $\TfU_l^n (t) = n^{-1/2} \big[\fU_{l}^n (nt) - n(\flambda_l^n)^{-1} t \big],~t\in
  [0,1]$ be the corresponding diffusion-scaled process.
\item (Sum of Service Time) Let
  $\fV_l^n (t):= \sum_{j=1}^{\lfloor t \rfloor} \fv_{l, j}^n, ~t\in [0,n]$ be
  the sum of service times among the first $\lfloor t\rfloor$ jobs predicted
  as class $l$,
   and
  $\TfV_l^n (t) = n^{-1/2} \Big[\fV_{l}^n (nt) - n (\fmu^n_l)^{-1} t \Big],~t\in
  [0,1],$ be the corresponding diffusion-scaled process. Here, $(\fmu^n_l)^{-1}:= \sum_{k=1}^K \frac{\tp_k^n \q_{kl}^n}{\fp_l^n} \frac{1}{\tmu_k^n}$ and $(\fmu_l)^{-1}:= \sum_{k=1}^K \frac{\tp_k \q_{kl}}{\fp_l} \frac{1}{\tmu_k}$ are the expected service times of an arbitrary job predicted as class $l$ and its limit.
\item (Service Process) Let $\fS_l^n(t) := \max \{ j \in \mathbb{N} : \fV^n_l(j) \leq t \}, \; t \in [0,n],$ denote the number of predicted class $l$ jobs that have been served and departed after $t$ units of service have been allocated to that class, and let
  $\TfS_l^n := n^{-1/2} [\fS_l^n (nt) -n \fmu_l^n t], ~t\in [0,1]$ be the
  corresponding diffusion-scaled process.
\end{enumerate}
\end{definition}
\noindent
For simplicity, we also use the vector processes
$\VTfA^n=(\TfA_l^n)_l$, $\VTfU^n=(\TfU_l^n)_l$, $\VTfS^n=(\TfS_l^n)_l$, and $\VTfV^n=(\TfV_l^n)_l$ to denote the second-order/diffusion-scaled processes.

Proposition~\ref{prop: joint conv. of predicted class A, U, S, and V} plays a
major role in our analysis of the endogenous processes in
Section~\ref{subsection: approximation for W, T, N, and tau of the predicted
  classes}.  We use the little-o notation $o_n(1)$ to denote uniform
convergence over $t\in[0, 1]$ as $n\rightarrow +\infty$.
\begin{proposition}[Convergence of exogenous processes of predicted classes] 
  \label{prop: joint conv. of predicted class A, U, S, and V}
  Given a classifier $\model$, suppose Assumptions~\ref{assumption: data
    generating process},~\ref{assumption: heavy traffic} and~\ref{assumption:
    second order moments} hold. There is a Brownian motion
  $(\VTfA, \VTfU, \VTfS, \VTfV)$ such that
  \begin{equation}
    \label{eq: conv. of predicted class A, U, S, V}
    \begin{aligned}
      (\VTfA^n, \VTfU^n, \VTfS^n, \VTfV^n) \rightarrow (\VTfA, \VTfU, \VTfS, \VTfV),
    \end{aligned}
  \end{equation}
  and for any predicted class $l\in [K]$
  \begin{equation} \label{eq: approximation of predicted A, U, S, V}
    \begin{aligned}
      &\fA_l^n (nt) =n \flambda^n_l t + n^{1/2}\TfA_l^n (t) + o_n(n^{1/2}),\\
      &\fU_l^n (nt) =n (\flambda_l)^{-1} t + n^{1/2}\TfU_l^n (t) + o_n(n^{1/2}),\\
      &\fS_l^n (nt) =n \fmu_l^n t + n^{1/2}\TfS_l^n (t) + o_n(n^{1/2}),\\
      &\fV_l^n (nt) =n (\fmu^n_l)^{-1} t + n^{1/2}\TfV_l^n (t) + o_n(n^{1/2}),
    \end{aligned}
  \end{equation}
  \begin{equation} \label{eq: vanishing jump sizes for U and V}
    \begin{aligned}
      n^{-1/2} \sup\limits_{1\leq j \leq \fA_l^n (n)} \fu_{l,j}^n \rightarrow 0,
      \qquad  n^{-1/2} \sup\limits_{1\leq j \leq \fA_l^n (n)} \fv_{l,j}^n \rightarrow 0.
    \end{aligned}
  \end{equation}
\end{proposition}
\begin{proof}
Recalling that $\TtA_0$, $\TfZ_{kl}$ are Brownian motions~\eqref{eq: a.s. conv
  of the primitives}, we begin by showing the limit
\begin{equation*}
\begin{aligned} 
  &\TfA_{kl}^n \rightarrow \TfA_{kl}: = \TfZ_{kl} \circ \tlambda e + \tp_k
  \q_{kl}\TtA_0, \quad
  \TfA_l^n \rightarrow \TfA_l: = \sum\limits_{k=1}^K \TfA_{kl}=\sum\limits_{k=1}^K \TfZ_{kl} \circ \tlambda e + \fp_l \TfA_0,\\
  &\TfU_l^n \rightarrow \TfU_l: = -\Big(\tlambda \sum\limits_{k=1}^K \tp_k
  \q_{kl}\Big)^{-1} \TfA_l \Big( \big(\tlambda \sum\limits_{k=1}^K \tp_k
  \q_{kl}\big)^{-1} e \Big).
\end{aligned}    
\end{equation*}
Recall from Definition~\ref{definition: U, Z, R, V} and Definition~\ref{def:
  concerned processes for arrival and service} that
$\fA_{kl}^n = \fZ_{kl}^n \circ \tA_0^n$ and
\begin{equation*}
\begin{aligned}
  \TfA_{kl}^n (t) = n^{-1/2} [ \fA_{kl}^n (nt) - \tlambda^n \tp_k^n \q_{kl}^n
  nt ] = \TfZ_{kl}^n (n^{-1}\tA_0^n (nt)) + \tp_k^n \q_{kl}^n \TtA_0^n (t).
\end{aligned}
\end{equation*}
Since $n^{-1} A^n_0(n\cdot) \to \lambda e$ by Lemma~\ref{lemma: convergence of the scaled primitive processes}, continuity of the composition
function~\cite[Theorem 13.2.1]{Whitt02} and the continuous mapping theorem
yields
\begin{equation*}
  \TfA_{kl}^n  \rightarrow 
  \TfZ_{kl} \circ \lambda e + p_k \q_{kl} \TtA_0
  \qquad \mbox{and} \qquad
  \TfA_l^n \rightarrow \TfA_l.
\end{equation*}
Since all limit have continuous sample paths, convergence in weak $M_2$
topology is equivalent to uniform convergence~\cite[Corollary
12.11.1]{Whitt02}. Asymptotic equivalence of counting and inverse processes
(with centering) gives convergence of $\TfU^n_l$~\cite[Corollary
13.8.1]{Whitt02}.

Using a nearly identical argument, we show convergence of $\TfS_l^n$ and
$\TfV_l^n$
\begin{equation*}
\begin{aligned}
  &\TfV^n_l \rightarrow \TfV_l :=
  \TfR_l \circ (\fp_l)^{-1}e + \fp_l (\fmu_l)^{-1} \TfM_l,\\
  &\TfS_l^n\rightarrow \TfS_l: = -\fmu_l \TfV_l \circ \fmu_l e = -\fmu_l
  \TfR_l \circ (\fp_l)^{-1} \fmu_l e - \fp_l \TfM_l \circ \fmu_l e,
\end{aligned}
\end{equation*}
where $\fM_l^n(t)$ is the total number of job arriving in the system until
arrival of $\lfloor t \rfloor$ jobs predicted as $l\in [K]$ and $\TfM_l^n(t)$
is the corresponding diffusion-scaled process,
\begin{equation}
  \label{eq:jobs}
\begin{aligned}
    &\fM_l^n(t) := \max \Big\{m\geq 0: \sum\limits_{i=1}^m \fY_{il}^n \leq t \Big\}
    = \max \Big\{m\geq 0:\sum_{k=1}^K \fZ_{kl}^n(m) \leq t \Big\},~\forall~t\in [0,n],\\
    & \TfM_l^n: = n^{-1/2} [\fM_l^n (nt) - (\fp_l^n)^{-1} nt ],~\forall~t\in [0,1].
\end{aligned}
\end{equation}
(Recall $\fp_l^n = \sum_{k=1}^K \tp_k^n \q_{kl}^n$.) $\fM_l^n$ is closely
related to $\sum_{k=1}^K\fZ_{kl}$ and can be understood as a counting process
with ``interarrival times" being $\{\fY_{il}^n: i\geq 1\}$.

Represent the service partial sum $\fV_l^n$ as a composition of $\fR_l^n$ with
the counting process $\fM_l^n$ 
\begin{equation*}
  \fV_l^n (t) = \fR_l^n (\fM_l^n (t)) = \sum_{i=1}^{\fM_l^n (t)} \fY_{il}^n \tv_i^n,~\forall~t\in [0,n]
  ~~~~~~\mbox{(Definitions~\ref{definition: U, Z, R, V}
and~\ref{def: concerned processes for arrival and service})}.
\end{equation*}
Since $\TfZ^n_{kl} \rightarrow \TfZ_{kl}\in \mathcal{C}$, a similar argument
as before gives
\begin{equation*}
  \TfM_l^n \rightarrow \TfM_l: = -(\fp_l)^{-1} \Big(\sum_{k=1}^K\TfZ_{kl} \Big)\circ \fp_l^{-1} e \in \mathcal{C}.
\end{equation*}
From the continuous mapping theorem, we have the convergence of $\TfV^n_l$ and
$\TfS^n_l$.
\end{proof}

\subsection{Dominance of p-FCFS and work-conserving policies}
\label{subsection: dominance of p-FCFS and work-conserving policies}

The results on the endogenous processes in Proposition~\ref{prop: convergence and approximation of predicted class N, tau, T, and W} and the lower bound in Theorem~\ref{theorem: HT lower bound} will be obtained assuming p-FCFS and work-conserving policies. We justify focusing on the set of p-FCFS and work-conserving policies.

\paragraph{p-FCFS}
Given a queueing system $n$ and a feasible policy, we can derive an associated feasible p-FCFS policy by swapping the service orders within each predicted class when the class has no previously preempted job. We show that the latter policy has stochastically smaller cumulative cost function $\TJ^n (t)$ for all $t\in [0,1]$. 
To do so, we analyze the distribution of the cost function under a modified data generating process governed by a new probability measure $\mathbb{Q}^n$ such that the distribution of $\TJ^n_{\policyn}$ remains the same as the original one under $\mathbb{P}^n$. The idea is to define the classes of jobs that govern the cost functions and service time distributions to be invariant under permutation within each predicted class, and use the convexity argument on the cost functions.

We assume that under $\mathbb{Q}^n$, $\{(u_i^n, X_i^n, \tY_i^n, \fY_i^n): i\in \mathbb{N}\}$ are generated in the same way as in Section~\ref{section: model}, but service times are generated differently. We introduce $\{(\widehat{Y}_{jl}^n, \fv_{jl}^n): j \in \mathbb{N}\}$ that are indexed according to the \textit{order of being served} rather than the order of arrivals within each predicted class $l\in [K]$. For \textit{any} job that is served as the $j$th distinct job within predicted class $l$ in system $n$, the service time is realized as $\fv_{jl}^n$ in a tuple $(\VaY_{jl}^n, \fv_{jl}^n)$, where $\VaY^n_{jl}: = (\aY^n_{jl, 1}, \ldots, \aY^n_{jl, K})$ denotes the one-hot encoding that determines the distribution of $\fv_{jl}^n$ as well as the cost function. In the sequel, we employ the subscripts $i$ and $j$ to signify indexing according to the arrival and service order within each predicted class, respectively.

We assume that for any queueing system $n$ and predicted class $l\in [K]$, 
\begin{enumerate}[(i)]
  \item $\{\VaY_{jl}^n, \fv_{jl}^n: j\in \mathbb{N}\}$ are i.i.d. random variables;
  \item $\{\VaY_{jl}^n, \fv_{jl}^n: j\in \mathbb{N}\}$ are independent of $\{(u_i^n, X_i^n, \tY_i^n, \fY_i^n): i\in \mathbb{N}\}$.
\end{enumerate}
Note that when swapping service orders between jobs within each predicted class, $(\VaY_{jl}^n, \fv_{jl}^n)$ remains unchanged in each sample path in $\mathbb{Q}^n$, a key property to be utilized in our proof.
To connect with the original data generating process under $\mathbb{P}^n$, we define the distribution of $(\VaY_{1l}^n, \fv_{1l}^n)$ as 
\begin{align} \label{eq: service time distribution of a predicted class}
  \mathbb{Q}^n [\aY_{1l, k}^n = 1, \fv_{1l}^n \leq x ]
  &:= \mathbb{P}^n [ \tY_{1k}^n = 1, v_{1}^n \leq x 
  \mid  \fY_{1l}^n = 1],
\end{align}
for any $k\in [K],x\in \mathbb{R}$, where $\mathbb{P}^n [ \tY_{1k}^n = 1, v_{1}^n \leq x \mid \fY_{1l}^n = 1] = \frac{p_k^n \q_{kl}^n}{\sum_{r=1}^K p_r^n \q_{rl}^n} \mathbb{P}^n [ v_{1}^n \leq x | \tY_{1k}^n=1]$ by conditional independence between $v_1^n$ and $\fY_{1l}^n$ given $\tY_{1k}^n$ in Assumption~\ref{assumption: data generating process}. Moreover, we use a modified (scaled) cumulative cost function $\aJ^n_{\pi^n}(t; Q^n)$, where for any job that is served as the $j$th distinct job within predicted class $l$, the cost is incurred according to its ``analytical class label'' $\aY_{jl}^n$ and defined by $C^n_{k:\aY_{jl,k}^n=1}(\cdot).$

\begin{lemma}[p-FCFS]\label{lemma: p fcfs}
Given a classifier $\model$ and a sequence of feasible policies $\{\pi_{n}\}$, suppose that Assumptions~\ref{assumption: data generating process}~\ref{assumption: heavy traffic},~\ref{assumption: second order moments}, and~\ref{assumption: on cost functions for showing the lower bound} hold. Then, for any queueing system $n$, there exists a feasible p-FCFS policy $\pi_{n, \text{p-FCFS}}$ such that $\TJ^n_{\pi_{n, \text{p-FCFS}}}(t; Q^n) \leq_{\text{st}} \TJ^n_{\pi_{n}}(t; Q^n),~\forall~t\in [0,1].$
\end{lemma}
\begin{proof}
Our proof uses a similar idea alluded to in the proof of~\cite[Theorem 2]{MandelbaumSt04} and provides a rigorous justification. Given a feasible policy $\policy_n$ in system $n\in \mathbb{N}$, we can define an associated p-FCFS policy, say $\pi_n'$, by applying the following basic operation: if there exists the $j$th arriving job in predicted class $l\in [K]$ that starts to be served by $\policy_n$ before the $i$th arriving job in the same predicted class with $\fU^n_l(i) < \fU^n_k(j)$ and $i$ being the smallest such index, then we swap service orders of the two jobs. It suffices to show that $\mathbb{P}^n [\TJ^n_{\pi_n'} (t; Q^n) > x] \leq \mathbb{P}^n [\TJ^n_{\pi_n} (t; Q^n) > x]$ for all $t\in [0,1],~x\in \mathbb{R}$.

We first claim that for all $t\in [0, 1]$,  $\TJ^n_{\policy_n}(t; Q^n)$ under $\mathbb{P}^n$ has the same marginal distribution as that of $\aJ^n_{\policy_n}(t; Q^n)$ under $\mathbb{Q}^n$. That is, 
\begin{equation}\label{eq: marginal distribution of cost}
  \mathbb{P}^n [\TJ^n_{\policy_n}(t; Q^n) > x] = 
  \mathbb{Q}^n [\aJ_{\policy_n}(t; Q^n) > x],
  ~\forall~x\in \mathbb{R}.
\end{equation}
The reason is that under $\mathbb{P}^n$ and $\mathbb{Q}^n$, the actual service time and the true/analytical class label of a job that determines the cost function to be applied are not known until the job starts to be served. Moreover, given arrival times $\{\fU_l^n (i): i\in \mathbb{N}\}$ in predicted class $l\in [K]$, service times and true/analytical class labels of waiting jobs are i.i.d. as~\eqref{eq: service time distribution of a predicted class} under the two probability measures. The latter implies that given the same realization of $\{\fU_l^n (i): i\in \mathbb{N}, l\in [K]\}$, the conditional distributions of $\TJ_{\policy_n}(\cdot;Q^n)$ and $\aJ_{\policy_n}(\cdot; Q^n)$ are identical.

By~\eqref{eq: marginal distribution of cost}, it suffices to show that $\pi_n'$ induced from $\pi_n$ by the basic operation satisfies
\begin{equation} \label{eq: dominance of p-FCFS}
  \mathbb{Q}^n \big[\aJ^n_{\pi_n'} (t; Q^n) \leq \aJ^n_{\pi_{n}} (t; Q^n),~\forall~t\in [0,1] \big] = 1. 
\end{equation}
To prove~\eqref{eq: dominance of p-FCFS}, fix a sample path under $\mathbb{Q}^n$ in system $n$. Suppose that at some time $nt'\in [0, n]$, there exist two jobs, $i_1$ and $i_2$, that \textit{arrived} as the $i_1$th and $i_2$th job in predicted class $l\in [K]$, respectively, with $\fU_l^n (i_1) < \fU_l^n (i_2)$, and have not been served at all. Suppose that $\policy_n$ chooses to serve $i_2$ at time $nt'$ as the $j_2$th distinct job served in predicted class $l$, and starts to serve $i_1$ later as the $j_1$th distinct job in that class with $j_1 > j_2$. Let $\Delta\fv_l^n (j_2, j_1) := \sum_{r=j_2 + 1}^{j_1 - 1} \fv_{rl}^n$ be the summation of service times for jobs in predicted class $l$ that are  served between $i_1$ and $i_2$. Also, suppose $\aY_{j_1,l,k_1}^n=\aY_{j_2,l, k_2}^n = 1$ for some $k_1,k_2\in [K]$. Note that $\Delta\fv_l^n (j_2, j_1), \fv_{j_2,l}^n, \fv_{j_1, l}^n, \aY_{j_1,l}^n$, and $\aY_{j_2,l}^n$ are identical \textit{regardless of which job is chosen for service} at time $nt'$. Similarly, under the conditions on preemption in Section~\ref{subsection: fundamental convergence results}, waiting times of jobs incurred by preemption during their service also remain the same independently of the job chosen for service at time $nt'$. Let $\Delta\fw_l^n (j_2, j_1)$ be the summation of waiting times incurred by preemption on jobs served between $i_1$ and $i_2$ in predicted class $l$, and let $\fw_{j_1,l}^n$ and $\fw_{j_2,l}^n$ be the waiting times by preemption on $i_1$ and $i_2$, respectively.

Now we are ready to show~\eqref{eq: dominance of p-FCFS}. We first show that changing service orders of $j_1$ and $j_2$ improves the cumulative cost at $t=1$. Specifically, the change $\aJ^n_{\policy_n'} (1;Q^n) - \aJ^n_{\policy_n} (1;Q^n)$ would be
\begin{equation*}
\begin{aligned}
&\Big[C_{k_2}^n \big( t' - \fU_l^n (i_1) + \fw_{j_2,l}^n + \fv_{j_2, l}^n\big) 
- C_{k_2}^n \big( t' - \fU_l^n (i_2) + \fw_{j_2,l}^n + \fv_{j_2, l}^n\big)\Big]\\
&-\Big[C_{k_1}^n \big(t'  - \fU_l^n (i_1) + \fw_{j_2,l}^n  + \fv_{j_2, l}^n + \Delta\fw_l^n (j_2, j_1) + \Delta\fv_l^n (j_2, j_1)  + \fw_{j_1,l}^n +\fv_{j_1, l}^n  \big)\\
& - C_{k_1}^n \big(t'  - \fU_l^n (i_2) + \fw_{j_2,l}^n + \fv_{j_2, l}^n + \Delta\fw_l^n (j_2,j_1) + \Delta\fv_l^n (j_2, j_1)  + \fw_{j_1,l}^n +\fv_{j_1, l}^n  \big)\Big] \leq 0,
\end{aligned}
\end{equation*}
where the inequality follows from convexity of $C_{k_1}^n, C_{k_2}^n$ and $\fU_l^n (i_1) < \fU_l^n (i_2)$, similarly to the proof~\cite[Proposition 1]{VanMieghem95}. In fact, one can observe that the cost reduction holds true for all $t\in [0,1]$ such that both jobs $i_1$ and $i_2$ are present in the system at time $nt$, and thus~\eqref{eq: dominance of p-FCFS} follows. This completes our proof.
\end{proof}

\paragraph{Work-Conserving}
For all system $n$ and any feasible policy
$\policyn$, we can always create a work-conserving counterpart policy
by having the server during an idle time to serve any waiting job, if
available. Since preemption is allowed without incurring additional costs, the server
can pause service and come back to the preempted job later, ensuring that the
cumulative cost does not increase as stated in the following lemma; see further discussion in~\cite[Section
2]{VanMieghem95}.  
\begin{lemma}[Work-Conserving]
\label{lemma: work conserving}
Given a classifier $\model$ and a sequence of feasible policy $\{\pi_{n}\}$, suppose that Assumptions~\ref{assumption: data generating process} and~\ref{assumption: on cost functions for showing the lower bound} hold. Then, for any queueing system $n$, there exists a feasible work-conserving policy $\pi_{n, \text{work-conserving}}$ such that 
$$\TJ^n_{\pi_{n, \text{work-conserving}}}(t; Q^n)
\leq \TJ^n_{\pi_{n}}(t; Q^n),
~\forall~t\in[0, 1],~\mathbb{P}^n\text{-a.s.}.$$
\end{lemma}

\subsection{Convergence of the endogenous processes of predicted classes}
\label{subsection: approximation for W, T, N, and tau of the predicted classes}

To prove Proposition~\ref{prop: convergence and approximation of predicted
  class N, tau, T, and W}, we formally define processes that are endogenous to
scheduling policies for predicted classes.
\begin{definition}[Endogenous processes]\label{definition: concerned process}
\begin{enumerate}[(i)]
\item (Total workload process) Let $\fL_l^n (t)$ be the total service time
  requested by all jobs predicted as class $l$ and arriving by time
  $t\in[0, n]$, and $\TfL_l^n (t)$ be the corresponding diffusion-scaled
  process,
\begin{equation*}
\begin{aligned}
  \fL_l^n (t) =\sum\limits_{i=1}^{\tA_0^n (t)} \fY_{il}^n \tv_i^n,~t\in [0,n],
  ~~~\TfL_l^n (t) = n^{-1/2}
    \Big[\fL_l^n (nt) - \tlambda^n \sum_{k=1}^K \frac{p^n_k}{\mu^n_k}\q_{kl}^n \cdot nt \Big],~t\in [0,1].
\end{aligned}
\end{equation*}
\item (Cumulative total input process) Let
  $\sumtL^n(t) = \sum_l \fL_l^n (t),~t\in[0, n]$ be the cumulative total input
  process and $\sumTtL^n(t) := \sum_{l=1}^K \TfL_l^n (t),~t\in[0, 1]$ be the
  corresponding diffusion-scaled process, i.e.,
\begin{equation*}
    \sumTtL^n (t) = n^{-1/2} \Big[\sumtL^n (nt) - \tlambda^n
    \sum_{k=1}^K \frac{\tp_k^n}{\tmu_k^n} \cdot nt\Big],~\forall~t\in [0,1].
\end{equation*}
\item (Policy process) Let $\fT_l^n (t)$ be total amount of time during
  $[0, t]$ that the server allocates to jobs from predicted class $l$, and
  $\TfT^n_l(t)$ be the corresponding diffusion-scaled process,
\begin{equation*}
  \TfT_l^n(t) = n^{-1/2} \Big[\fT^n_l(nt) - 
  \lambda^n\sum_{k=1}^K \frac{p^n_k}{\mu^n_k}\q^n_{kl} \cdot nt\Big],~t\in[0, 1].
\end{equation*}
\item (Remaining workload process) Let $\fW_l^n (t)$ be the remaining service
  time requested by jobs predicted as class $l$ and present (waiting for
  service or being served) in the system at time $t \in [0, n]$,
  \begin{equation}\label{eq: W in terms of L and T}
    \fW_l^n (t) = \fL_l^n (t) - \fT_l^n (t),\quad t\in [0,n].
  \end{equation}
  and $\TfW^n_l(t): = n^{-1/2} \fW^n_l(nt),~\forall~t\in [0, 1]$ be the
  corresponding diffusion scaled process.
\item (Total remaining workload process) Let
  $\sumtW^n(t) = \sum_l \fW_l^n (t)$ be the total remaining workload process
  and $\sumTtW^n(t): = n^{-1/2} \sum_{l=1}^K \fW^n_l(nt),~\forall~t\in [0, 1]$
  be the corresponding diffusion scaled process.
\item (Queue length process) Let $\fN_l^n (t)$ be the total number of jobs
  that are predicted as class $l$ and present (waiting for service or being
  served) in the system at time $t \in [0, n]$, and $\fN_{kl}^n (t)$ be the
  total number of true class $k$ jobs that are predicted as class $l$ and
  present in the system at time $t \in [0, n]$. Let
  $\TfN_l^n(t): = n^{-1/2} \fN^n_l(nt),~\TfN_{kl}^n(t): = n^{-1/2}
  \fN^n_{kl}(nt),~\forall~t\in [0, 1]$ be the corresponding scaled processes.
\item (Sojourn time process) Let $\ftau_{lj}^n$ be the sojourn time, the time
  span between arrival and service completion, of the $j$th job of predicted
  class $l$. Let $\ftau^n_l(t) = \ftau^n_{l,\fA_l^n(t)},~\forall~t\in [0, n]$
  be the sojourn time process where $\ftau^n_l (t)$ denotes the sojourn time
  of the latest job predicted as class $l$ and arriving by time $t$ and
  $\Tftau_l^n(t): = n^{-1/2} \ftau^n_l(nt),~\forall~t\in [0, 1]$ be the
  corresponding scaled process.
\end{enumerate}
\end{definition}

Note that $\ftau_l^n(\fU^n_l(i)) = \ftau^n_{l,i}$ and $\ftau^n_l$ only
exhibits jumps at arrival times $\{\fU^n_l(i)\}_{i = 1}^\infty$ of jobs
predicted as class $l$. By definition, $\ftau^n_l$ is also RCLL. Since
$\fL^n_l$ is an exogenous process, according to Eq.~\eqref{eq: W in terms of L
  and T}, we can also characterize the policy process $\fT^n_l$, or
equivalently, the scheduling policies, by the remaining workload process
$\fW^n_l$. The following results hold under p-FCFS feasible policies:
\begin{equation}\label{eq: predicted class N, tau, T, and W}
\begin{aligned}
\fN_l^n (t) &= \fA_l^n (t) - \fS_l^n (\fT_l^n (t)),
&\forall~t\in [0,n];\\
\ftau_l^n (t) &= \inf \{s \geq 0: \fW_l^n (t) \leq \fT_l^n (t + s) - \fT_l^n(t)\},
&\forall~t\in [0,n];\\ 
\fW_l^n (t) &= \fT_l^n \left(t + \ftau_l^n (t)\right) - \fT_l^n (t),
&\forall~t\in [0,n].
\end{aligned}
\end{equation}

We show the convergence of the scaled input process $\TfL^n_l$, which will be used to prove convergence of the workload $\sumtW^n$ in Section~\ref{subsection: proof for convergence and approximation of predicted class N, tau, T, and W}. 

\begin{lemma}[Convergence of scaled total input processes] \label{lemma: conv. of predicted L} Given a classifier $\model$, a sequence of queueing
systems, and a sequence of feasible policies $\{\policyn\}$, suppose that
Assumptions~\ref{assumption: data generating process},~\ref{assumption:
heavy traffic} and~\ref{assumption: second order moments} hold.  Then, for
any predicted class $l\in [K]$, we have that
\begin{equation*}
\TfL_l^n \rightarrow \TfL_l:=\TfR_l \circ \tlambda e + \sum_{k=1}^K \frac{\tp_k}{\tmu_k} \q_{kl} \TtA_0,
\quad  
\sumTtL^n \rightarrow \sumTtL:=\TtV_0 \circ \tlambda e + \sum\limits_{k=1}^K \frac{\tp_k}{\tmu_k} \TtA_0
\end{equation*}
as $n\rightarrow \infty$, where $e$ is the identity function on $[0,1]$. Also, 
for any system $n$ and time $t\in[0,1]$, 
\begin{equation*}
\begin{aligned}
\fL_l^n (nt) =\tlambda^n \sum_{k=1}^K \frac{\tp_k^n}{\tmu_k^n} \q^n_{kl}\cdot n  t + n^{1/2} \TfL_l^n (t) + o(n^{1/2});\quad 
\sumtL^n (nt) = \tlambda^n \sum_{k=1}^K \frac{\tp_k^n}{\tmu_k^n} \cdot n  t + n^{1/2} \sumTtL^n (t) + o(n^{1/2}).
\end{aligned}
\end{equation*}
\end{lemma}
\begin{proof}
Note that  $\fL_l^n = \fR_l^n \circ \tA_0^n$ by Definition~\ref{definition: concerned process}. Therefore, we have
\begin{equation*}
\begin{aligned}
\TfL^n_l (t) 
=&~     n^{-1/2} \Big[\fR^n_l(A^n_0(nt)) - \sum_{k=1}^K \frac{p^n_k}{\mu^n_k}\q_{kl}^n \cdot  A^n_0(nt)\Big]
+ \sum_{k=1}^K \frac{p^n_k}{\mu^n_k}\q_{kl}^n \cdot n^{-1/2} \Big[A^n_0(nt) - \tlambda^n\cdot  nt\Big]\\
=&~ \TfR_l^n(n^{-1} A^n_0(nt)) + \sum_{k=1}^K \frac{p^n_k}{\mu^n_k}\q_{kl}^n \cdot  \TtA^n_0(t).
\end{aligned}
\end{equation*}
Recall $\TtA^n_0(t) \rightarrow \TtA_0(t)$ and $n^{-1} A^n_0(n\cdot) \rightarrow \lambda e$ by Lemma~\ref{lemma: convergence of
  the scaled primitive processes}. Since $\lambda e$ is continuous,
continuity of the composition mapping~\citep[Theorem 13.2.1]{Whitt02} and
the continuous mapping theorem yields
$\TfL^n_l (t) \rightarrow \TfR_l \circ \tlambda e + \fp_l (\fmu_l)^{-1}
\TtA_0.$
Convergence of $\sumTtL^n$ is a direct consequence of continuous mapping theorem and $\sumTtL^n = \sum_{l=1}^K \TfL^n_l$ by Definition~\ref{definition: concerned process}.
\end{proof}

\subsubsection{Proof of Proposition~\ref{prop: convergence and approximation of predicted class N, tau, T, and W}}
\label{subsection: proof for convergence and approximation of predicted class N, tau, T, and W}

We establish Proposition~\ref{prop: convergence and approximation of predicted
  class N, tau, T, and W} based on Proposition~\ref{prop: joint conv. of
  predicted class A, U, S, and V} and Lemma~\ref{lemma: conv. of predicted
  L}. Our approach is similar to the proof of~\cite[Proposition
2]{VanMieghem95}, and we complement the latter with additional details in the
proof.  Since $\{\policyn\}$ is work-conserving, the remaining workload
process $\sumtW^n$ can be written as
$\sumtW^n = \phi (\sumtL^n - e)$~\citep{Whitt02}, where $\phi$ is the
one-sided reflection mapping. By Lemma~\ref{lemma: conv. of predicted L} and
heavy-traffic conditions (Assumption~\ref{assumption: heavy traffic}), for
any $t\in [0, 1]$
\begin{equation*}
  \begin{aligned}
    (\sumtL^n - e) (nt)
    & = \Big(\tlambda^n \sum_{k=1}^K
    \frac{\tp_k^n}{\tmu_k^n}  - 1 \Big)\cdot nt
    + n^{1/2} \sumTtL^n (t) + o(n^{1/2})\\
    & = n^{1/2} \Big[\sumTtL^n (t)
    + n^{1/2}\big(\tlambda^n
    \sum_{k=1}^K \frac{\tp_k^n}{\tmu_k^n} - 1 \big)t \Big] + o(n^{1/2}) \\
    & = n^{1/2} \sumTtL^n (t)  + o(n^{1/2}),
  \end{aligned}
\end{equation*}
where we used
$n^{1/2}\big(\tlambda^n \sum_{k=1}^K \frac{\tp_k^n}{\tmu_k^n} - 1 \big) =
o_n(1)$ in the final line. Combining with the relation
$\sumtW^n = \phi (\sumtL^n - e)$, we get
\begin{equation*}
\begin{aligned}
n^{-1/2}\sumtW^n (nt) =&~ n^{-1/2}\phi (\sumtL^n - e) (nt) 
= \phi (n^{-1/2}(\sumtL^n - e) (nt))\\
=&~ \phi( \sumTtL^n(t) + o_n(1)) 
= \phi( \sumTtL^n(t) ) + o_n(1),
\end{aligned}
\end{equation*}
where the first line follows from definition of $\phi$, and the last line
results from Lipschitz property of $\phi$ with the uniform metric~\citep[Lemma
13.5.1]{Whitt02}. Since $\sumTtW^n (t): = n^{-1/2}\sumtW^n (nt)$ in
Definition~\ref{definition: concerned process}, the convergence
$\sumTtW^n \rightarrow \phi(\sumTtL)$ follows from the above and
Lemma~\ref{lemma: conv. of predicted L}.

Next, we consider $\TfW^n_l, \TfN^n_l$, and $\Tftau^n_l$. 
Notice that $\TfW^n_l \geq 0$,  $\sum_{l=1}^K \TfW^n_l \rightarrow \phi(\sumTtL)$, and $\phi(\sumTtL)$ is a continous function on $[0, 1]$. Therefore, it is clear that $\limsup_n \| \TfW^n_l \| < +\infty,~\forall~l\in [K]$.
For $\TfT^n_l$, by  Definition~\ref{definition: concerned process} and Lemma~\ref{lemma: conv. of predicted L}, we have that
\begin{equation}\label{eq: approximation of tilde T}
\begin{aligned}
    \TfT_l^n (t) 
=&~ n^{-1/2} \Big[ \fT_l^n (nt)  - 
    \lambda^n \sum_{k=1}^K \frac{p_k^n}{\mu_k^n} \cdot nt \Big]\\
=&~ n^{-1/2} \Big[ \fL_l^n (nt)  - 
\lambda^n \sum_{k=1}^K \frac{p_k^n}{\mu_k^n} \cdot nt \Big]
- n^{-1/2}\fW_l^n (nt) \\
=&\TfL^n_l(t) - \TfW^n_l(t),
\end{aligned}
\end{equation}
where the second line follows from $\fT_l^n (nt) = \fL_l^n (nt) - \fW_l^n (nt)$ by  Definition~\ref{definition: concerned process}.
Since $\TfL^n_l \rightarrow \TfL_l$ by Lemma~\ref{lemma: conv. of predicted L}, one can check that for any $l\in [K]$, $\TfT_l^n$ converges if and only if $\TfW_l^n$ converges. Also, $\limsup_n \| \TfT_l^n\| < +\infty,~\forall~l\in [K]$.

Recalling the relation~\eqref{eq: predicted class N, tau, T, and W}, we have
$\fN_l^n (nt) = \fA_l^n (nt) - \fS_l^n (\fT_l^n (nt))$.  Using
Proposition~\ref{prop: joint conv. of predicted class A, U, S, and V}, we can
rewrite $\TfN^n_l(t)$ as
\begin{equation*}
\begin{aligned}
    \TfN^n_l(t) 
=&~ n^{-1/2}[\fA_l^n (nt) - \fS_l^n (\fT_l^n (nt)) ]\\
=&~ n^{1/2} \lambda^n \fp_l^n t + \TfA^n_l(t) -
    n^{1/2} \fmu_l^n t \cdot n^{-1} \fT_l^n (nt)
    -  \TfS^n_l(n^{-1} \fT_l^n (nt))
    + o(1).\\
\end{aligned}
\end{equation*}
Note that
\begin{equation}\label{eq: approximation of T}
    n^{-1}\fT_l^n (nt) = \lambda^n \sum_{k=1}^K \frac{p^n_k \q_{kl}^n}{\mu^n_k} \cdot t + n^{-1/2} \TfT^n_l(t) + o(n^{-1/2}),
\end{equation}
where $ n^{-1/2} \TfT^n_l(t) = o_n (1)$ as $\limsup_n \|  \TfT^n_l \| < +\infty$.
Therefore, $\TfN^n_l(t) $ can be rewritten as 
\begin{align}
    \TfN_l^n(t) 
=&~ n^{1/2} \Big[\lambda^n \fp^n_l\cdot t 
    - \fmu^n_l  \lambda^n \sum_{k=1}^K \frac{p^n_k \q_{kl}^n}{\mu^n_k}  \cdot t
    \Big]
    - \fmu^n_l \TfT^n_l(t) + \TfA^n_l(t)
    - \TfS^n_l\Big( \lambda^n \sum_{k=1}^K \frac{p^n_k \q_{kl}^n}{\mu^n_k} \cdot t + o_n (1)\Big) + o_n (1)\nonumber \\
=&~ - \fmu^n_l \TfT^n_l(t) + \TfA^n_l(t)
- \TfS^n_l\Big( \lambda^n \sum_{k=1}^K \frac{p^n_k \q_{kl}^n}{\mu^n_k} \cdot t + o_n (1)\Big) + o_n (1),
\label{eq: approximation of predicted N}
\end{align}
by definition of $\fp^n_l$ and $\fmu^n_l$ in Definition~\ref{def: concerned processes for arrival and service}. Since 
$\lambda^n \sum_{k=1}^K \frac{p^n_k \q_{kl}^n}{\mu^n_k} \cdot e \rightarrow
\lambda \sum_{k=1}^K \frac{p_k \q_{kl}}{\mu_k} \cdot e  \in \mathcal{C}$, 
by continuity of composition~\cite[Theorem 13.2.1]{Whitt02} and continuous mapping theorem, for any $l\in [K]$, $\TfT_l^n$ converges if and only if $\TfN_l^n$ converges, and $\limsup_n \| \TfN_l^n\| < +\infty,~\forall~l\in [K]$.

Finally, for $\Tftau^n_l$, once again by~\eqref{eq: predicted class N, tau, 
T, and W}, we have that 
\begin{equation*}
    \TfW_l^n (t) = n^{-1/2}[\fT_l^n (nt + \ftau_l^n (nt)) - \fT_l^n (nt)],~\forall~t\in [0,n].
\end{equation*}
According to the previous result~\eqref{eq: approximation of T}, we have 
\begin{equation}\label{eq: relation between W, tau, T}
\begin{aligned}
\TfW_l^n (t) &= \tlambda^n 
\sum_{k=1}^K \frac{p^n_k \q_{kl}^n}{\mu_k^n}
\cdot  \Tftau_l^n (t) + \TfT_l^n (t + n^{-1} \ftau_l^n (nt)) 
- \TfT_l^n (t) + o_n (1).
\end{aligned}
\end{equation}
For any predicted class $l\in [K]$, $\limsup_n \| \TfW_l^n\| < +\infty$ and
$\limsup_n \| \TfT_l^n\| < +\infty$, so that
$\limsup_n \| \Tftau_l^n\| < +\infty$ and $n^{-1}\ftau_l^n (n \cdot) = o_n (1)$.
Moreover, for any predicted class $l\in [K]$, 
since $\TfT_l^n$ converges if and only if $\TfW_l^n$ converges, 
it is easy to show that if $\TfT_l^n$ and $\TfW_l^n$ converge, 
then $\Tftau_l^n$ converges. Also, 
for the case that $\Tftau_l^n$ converges in $(\mathcal{D}, \| \cdot \|)$, 
according to $\TfT_l^n (t) = \TfL^n_l(t) - \TfW^n_l(t)$ by~\eqref{eq: approximation of tilde T}, we have that
\begin{equation} \label{eq: relation between W, tau, L}
\TfW_l^n(t + n^{-1} \ftau_l^n (nt)) = 
\tlambda^n 
\sum_{k=1}^K \frac{p^n_k \q_{kl}^n}{\mu_k^n}
\cdot  \Tftau_l^n (t) + \TfL_l^n (t + n^{-1} \ftau_l^n (nt)) 
- \TfL_l^n (t)  + o_n(1).
\end{equation}
Since $\TfL^n_l \rightarrow \TfL_l \in \mathcal{C}$ by Lemma~\ref{lemma: conv. of predicted L} and $n^{-1}\ftau_l^n (n \cdot)$ vanishes, convergence of $\Tftau_l^n$ in $(\mathcal{D}, \| \cdot \|)$ and the right-continuity imply convergence of $\TfW_l^n$ in $(\mathcal{D}, \| \cdot \|)$. Finally, by Eqs.~\eqref{eq: approximation of tilde T},~\eqref{eq: approximation of predicted N},~\eqref{eq: relation between W, tau, L}, if any of $\TfT^n_l, \TfW^n_l, \TfN^n_l$ and $\Tftau^n_l$ has a continuous limit, so does the others.

\subsection{Diffusion limits of the classical queueing model}
\label{subsection: convergence of the classical queueing model}
We extend the classical queueing model in~\citet{VanMieghem95} and~\citet{MandelbaumSt04} in the presence of misclassification errors. Key convergence results analogous to Lemma~\ref{lemma: uniform convergence of UZRV}, Proposition~\ref{prop: joint conv. of predicted class A, U, S, and V}, and Proposition~\ref{prop: convergence and approximation of predicted class N, tau, T, and W} can be shown similarly to our proofs. However, in this framework we can only \emph{establish} optimality of \ourmethod~among \textit{p-FCFS} policies, a weaker guarantee than Theorem~\ref{theorem: optimality of our policy}, where optimality holds over all feasible policies. Whether \ourmethod~remains optimal over all feasible policies in the classical model is left open.

\subsubsection{Diffusion limit in the classical framework}
\label{subsubsection: diffusion limit in the classical framework} 
We explain a new data generating process given external arrivals from $K$ real classes as in~\citep{VanMieghem95, MandelbaumSt04}. For $k\in [K]$ and $n\in \mathbb{N}$, i.i.d.\ random vectors $\{(\tu_{ki}^n, \feature_{ki}^n, \tv_{ki}^n): i\in \mathbb{N}\}$ are generated, where $\tu_{ki}^n$ is the interarrival time of the $i$th arriving job of real class $k$ in system $n$ with a constant arrival rate $\tlambda_{k}^n$, $(\tlambda_k^n)^{-1}:=\E^n[\tu_{k1}^n] > 0$. The tuple $(\feature_{ki}^n, \tv_{ki}^n)$ is generated \textit{independently} of $\tu_{ki}^n$, where
$\feature_{ki}^n \in \mathbb{R}^d$ represents the feature vector of the job, and $\tv_{ki}^n$ indicates the time required to serve the job. Let
$(\tmu_k^n)^{-1}:=\E^n[\tv_{k1}^n]$ be the expected service
time of a class-$k$ job. Let $\tU_k^n (m):=\sum_{i=1}^m \tu_{ki}^n$ be the arrival epoch of the $m$th real class $k$ job, and $A_k^n (t) = \max\{m: \tU_k^n (m) \leq t\},~t\in [0,n],$ the corresponding arrival counting process of real class $k$.

For each $k\in [K]$, the predicted class of the $i$th real class $k$ job is defined by the one-hot vector $\fY_{ki}^n := \model (\feature_{ki}^n)=(\fY_{ki}^n (1) ,..., \fY_{ki}^n (K))$. The classification probabilities are defined as $\q_{kl}^n := \P^n[\fY_{k1}^n (l) = 1]$ for $k,l\in [K]$. We assume that $v^n_{ki}$ is independent of $\feature^n_{ki}$, implying $v_{ki}^n \perp \fY^n_{ki}$. The data generating processes and the corresponding heavy traffic conditions are summarized as the following.
\begin{assumption}[Alternative data generating processes]
\label{assumption: alternative data generating process}
For any system $n\in \N$, 
\begin{enumerate}[(i)]
\item the sequences of random vectors $\{(\tu_{ki}^n, \tv_{ki}^n, \feature_{ki}^n): i\in \N\}$ are independent over $k\in [K]$;
\item $\{(\tu_{ki}^n, \tv_{ki}^n, \feature_{ki}^n): i\in \N\}$ is a sequence of i.i.d random vectors for each class $k\in [K]$;
\item $\{\tu_{ki}^n: i \in \mathbb{N}\}$, $\{\tv_{ki}^n: i \in \mathbb{N}\}$, and $\{\feature_{ki}^n: i \in \mathbb{N}\}$ are independent for each class $k\in [K]$.
\end{enumerate}
\end{assumption}

\begin{assumption}[Heavy traffic condition] 
\label{assumption: heavy traffic in classical framework}
Given a classifier $\model$ and a sequence of queueing systems, there exist
$\tlambda_k, \tmu_k \in (0, \infty) $ and $\q_{kl} \in [0,1]$ for $k,l\in [K]$ such that
$ \sum_{k=1}^K \q_{kl} > 0,~\forall~l\in [K]$,  
$\sum_{k=1}^K \frac{\tlambda_k}{\tmu_k} = 1$, and as
$n\rightarrow \infty$, for all $k,l \in [K]$
\begin{equation}\label{eq: convergence rate of the classical model}
n^{1/2}\big(\tlambda_k^n - \tlambda_k \big) \rightarrow 0, \quad
n^{1/2}  \big(\tmu_k^n - \tmu_k\big) \rightarrow 0,\quad
n^{1/2} \big(\q_{kl}^n - \q_{kl} \big)
\rightarrow 0. 
\end{equation}
\end{assumption}

Under the identification $\tlambda_k = \tlambda \tp_k$, the critical-load
condition $\sum_{k=1}^K \tlambda_k/\tmu_k = 1$ coincides with
$\tlambda \sum_{k=1}^K \tp_k/\tmu_k = 1$ in Assumption~\ref{assumption: heavy
traffic}, so the two models share the same heavy-traffic scaling.

\paragraph{Diffusion limit}
To derive the diffusion limit in the classical model, the key processes in Definition~\ref{definition: U, Z, R, V} are modified to
\begin{equation*}
    \fZ_{kl}^n (t): = \sum_{i=1}^{\lfloor t \rfloor} \fY_{ki}^n(l)
    , \quad \fR_{kl}^n (t) := \sum_{i=1}^{\lfloor t \rfloor}\fY_{ki}^n (l) \tv_{ki}^n, 
\quad t\in [0,n],~\forall k,l\in [K].
\end{equation*}
Note that $\fR_{kl}^n$ is now defined for each pair of $k,l\in [K]$. Then, using Assumptions~\ref{assumption: second order moments},~\ref{assumption: alternative data generating process},~\ref{assumption: heavy traffic in classical framework}, the convergence results analogous to Lemma~\ref{lemma: joint weak convergence} and Lemma~\ref{lemma: uniform convergence of UZRV} can be obtained using the martingale FCLT (Lemma~\ref{lemma: Martingale FCLT}) as in Section~\ref{section:proof-joint-weak-convergence} and Section~\ref{subsection: proof of uniform convergence}. Building off of the initial diffusion limit, we can show convergence of the processes of predicted classes as in Proposition~\ref{prop: joint conv. of predicted class A, U, S, and V} and Proposition~\ref{prop: convergence and approximation of predicted class N, tau, T, and W} using similar techniques. Specifically, let arrival processes associated with predicted classes be $\fA_{kl}^n (t):= \sum_{i=1}^{ \tA_k^n (t) } \fY_{ki}^n (l),~\fA_l^n (t) :=  \sum_{k=1}^K \fA_{kl}^n (t),~t\in [0,n]$ for $k,l\in [K]$, and adapt the definitions of the other processes (Definition~\ref{def: concerned processes for arrival and service}, Definition~\ref{definition: concerned process}) and their characterizations analogously. For example, similarly to the proof of Proposition~\ref{prop: joint conv. of predicted class A, U, S, and V}, $\fV_l^n$ will have to be represented as a composition to apply the random time change technique:
\begin{equation*}
  \fV_l^n (t) = \sum\limits_{k=1}^K \fR_{kl}^n \big( (\tA_k^n \circ \fU_l^n) (t)\big) =
  \sum\limits_{k=1}^K \sum\limits_{i=1}^{\tA_k^n (\fU_l^n (t))} \fY_{ki}^n (l) \tv_{ki}^n,~t\in [0,n].
\end{equation*}
Apart from this per-pair decomposition of $\fR_{kl}^n$ and the composition above, every step coincides with the single-stream proofs, so the analogous convergence results hold \emph{mutatis mutandis}.

\subsubsection{Failure of the p-FCFS interchange argument under the classical queueing model}
\label{subsubsection: weaker lower bound result under the classical model}
We show that the interchange argument behind the p-FCFS reduction (Lemma~\ref{lemma: p fcfs}) breaks down in the classical queueing model, so that our proof no longer reduces feasible policies to p-FCFS policies without loss. The obstruction is specific to misclassification. Under perfect classification ($Q^n = I$), each predicted class consists of a single real class with i.i.d.\ service times, Lemma~\ref{lemma: p fcfs} applies verbatim, and optimality over all feasible policies follows as in Theorem~\ref{theorem: optimality of our policy}. Under misclassification, a predicted class mixes several real classes with different service laws, and with class-specific arrival streams the arrival timing of a job is informative of its real class. As a result, service times of waiting jobs in a predicted class are not generally i.i.d.\ with respect to the filtration that policies are adapted to (Definition~\ref{definition: feasible policies}), except for the special case of independent Poisson arrivals, where memoryless interarrivals render the timing uninformative. The proof of Lemma~\ref{lemma: p fcfs} is then inapplicable, and we can only establish the distributional lower bound~\eqref{eq: HT lower bound in stochastic sense} and the optimality in Theorem~\ref{theorem: optimality of our policy} over p-FCFS policies rather than all feasible policies. Whether \ourmethod~remains optimal over all feasible policies in this model is left open.

To be concrete, consider a two-class system $n$ where $\{\tu_{1i}^n\}$ and $\{\tu_{2i}^n\}$ take values of either $100$ or $150$ and $1$ or $3$, respectively. Let service times $\{\tv_{1i}^n\}$ and $\{\tv_{2i}^n\}$ be either $2$ or $6$ and $\frac{1}{2}$ or $\frac{3}{2}$, respectively, and $\q_{kl}^n=\frac{1}{2}$ for all $k,l=1,2$. The service supports $\{2,6\}$ and $\{\frac{1}{2},\frac{3}{2}\}$ are disjoint, so a revealed service time identifies the real class. Moreover, real class $k$ arrival epochs lie in the additive semigroup generated by the support of its interarrival law, namely $\langle 100,150\rangle$ for class $1$ and $\langle 1,3\rangle$ for class $2$. Suppose \emph{predicted class} $1$ has two waiting jobs with arrival times $\fU_{1,j}^n,~j=1,2$.

\emph{Not identically distributed.} Consider $\fU_{1,1}^n=100,~\fU_{1,2}^n=103$. The first job may be of real class $1$ or $2$, whereas $103\notin\langle 100,150\rangle$ forces the second job to be of real class $2$. The two waiting jobs thus have different service-time distributions.

\emph{Not independent.} Consider $\fU_{1,1}^n=100,~\fU_{1,2}^n=150$. If the first job's service time is $2$ or $6$, then it is of real class $1$, and the second job cannot then be of real class $1$ since $150-100=50\notin\{100,150\}$, so it is of real class $2$. If instead the first service time is $\frac{1}{2}$ or $\frac{3}{2}$, the first job is of real class $2$ and the second can be of real class $1$ with positive probability. The conditional law of the second job's service time therefore depends on the first.

Either failure invalidates the i.i.d.\ premise behind the proof of Lemma~\ref{lemma: p fcfs}.

\subsection{Proof of Lemma~\ref{lemma: convergence of tildeJ}}
\label{subsection:proof-convergence-tildeJ}

Using the shorthand
$f^n_{kl}(s) := C_k^n(\ftau_l^n(n s)) =C_k^n(n^{1/2}\Tftau_l^n(s))$, 
$f_{kl}(s) := C_k(\Tftau_l(s))$, 
$\mathrm{d}\xi^n_{kl}(\cdot): = 
\mathrm{d} \left(n^{-1}\fA^n_{kl}(n \cdot)\right) (\cdot)$,
and $\mathrm{d}\xi_{kl}(\cdot): = \mathrm{d} \left( \tlambda \tp_{k} \q_{kl} e (\cdot)\right)$,
triangle inequalities gives
\begin{equation}\label{eq: limit of J eq1}
\begin{aligned}
&~    \sup\limits_{t\in [0,1]} |\TJ^n_\policyn(t; Q^n) - \TJ_\policy(t; Q)| 
=     \sup\limits_{t\in [0,1]} \Big| \sum\limits_{k=1}^K \sum\limits_{l=1}^K
      \int_0^t  f^n_{kl} (s ) \mathrm{d}\xi^n_{kl}(s)
      - \sum\limits_{k=1}^K \sum\limits_{l=1}^K
      \int_0^t  f_{kl}(s) \mathrm{d}\xi_{kl}(s) \Big| \\
\leq&~ \sum\limits_{k=1}^K	 \sum\limits_{l=1}^K \left\{ \sup\limits_{t\in [0,1]}
      \int_0^t  \big|f^n_{kl}(s) - f_{kl}(s) \big| \mathrm{d}\xi^n_{kl}(s)
+    \sup\limits_{t\in [0,1]}\Big| 
      \int_0^t f_{kl}(s) \mathrm{d}\xi^n_{kl}(s)
      - \int_0^t f_{kl}(s) \mathrm{d}\xi_{kl}(s)\Big| \right\}.
\end{aligned}
\end{equation}

The first term of~\eqref{eq: limit of J eq1} $\rightarrow 0$ since $f^n_{kl}(s) 
\rightarrow f_{kl}(s) $ by Assumption~\ref{assumption: on cost functions for showing
the lower bound} and $ \lim_{n}\xi^n_{kl}([0, 1]) = \xi_{kl}([0,1]) <  +\infty$
by Proposition~\ref{prop: joint conv. of predicted class A, U, S, and V}.
For the second term of~\eqref{eq: limit of J eq1}, by Proposition~\ref{prop: joint conv. 
of predicted class A, U, S, and V} and generalized Lebesgue convergence 
theorem~\cite[Page 270]{Royden88}, it is clear that 
$\int_0^{t'} f_{kl}(s) \mathrm{d}\xi^n_{kl}(s) - \int_0^{t'} f_{kl}(s) 
\mathrm{d}\xi_{kl}(s) \rightarrow 0$ as $n\rightarrow +\infty$ for any fixed 
$t'\in[0, 1]$. To achieve uniform convergence, we partition $[0, 1]$ into $M$ intervals 
$0 = a_0 < a_1 < \cdots < a_M = 1$ with $a_i - a_{i-1} = 1/M$. Then, for any fixed $M$,
we have $\max_{1\leq i \leq M}  | \int_0^{a_i} f_{kl}(s) \mathrm{d}\xi^n_{kl}(s)
- \int_0^{a_i} f_{kl}(s) \mathrm{d}\xi_{kl}(s)  | \rightarrow 0$ 
as $n\rightarrow +\infty$. Using $\| f_{kl}(s) \| < +\infty$  and
\begin{equation*}
\sup_{ |t_1 - t_2| \leq \frac{1}{M}} \Big | 
      \int_{t_1}^{t_2} f_{kl}(s) \mathrm{d}\xi^n_{kl}(s)
  - \int_{t_1}^{t_2} f_{kl}(s) \mathrm{d}\xi_{kl}(s) \Big | \leq 
\|f_{kl}\| \sup_{ |t_1 - t_2| \leq  \frac{1}{M} } \left\{
\Big|\int_{t_1}^{t_2} \mathrm{d}\xi^n_{kl}(s)\Big| + 
\Big|\int_{t_1}^{t_2} \mathrm{d}\xi_{kl}(s)\Big| \right\} \rightarrow 0
\end{equation*}
as $M, n \rightarrow +\infty$ by Proposition~\ref{prop: joint conv. of predicted class 
A, U, S, and V}, we can show the second term of~\eqref{eq: limit of J eq1} 
also $\rightarrow 0$ as $n\rightarrow +\infty$. This completes our proof. 


\section{Proof of heavy traffic lower bound (Theorem~\ref{theorem: HT lower
    bound})}
\label{section: proof for HT lower bound}

The proof closes three exposition gaps in~\citet[Proposition 6]{VanMieghem95}:
(A1) the queue-to-workload translation underlying Eq.~(36)
(Section~\ref{section:proof-relation between N and W});
(A2) continuity of the optimal allocation, used implicitly in Eq.~(95)
to pass a Riemann sum to an integral
(Section~\ref{subsection: continuity of the optimal workload allocation});
and (A3) the vanishing partition mesh, justified here for $\sumTtW$
a reflected Brownian motion
(Section~\ref{subsubsection: complementary proof for the partition size in VM}).

\subsection{Overview} 
\label{subsection: proof overview for HT lower bound}

Since the queue based on the predicted classes contains a mixture of true
classes due to misclassification, we must characterize its asymptotic
compositions in order to analyze the queueing cost. For $k,l\in [K]$, let
$\fN_{kl}^n (t),~t\in [0,n]$ be the number of true class $k$ jobs that are
predicted as class $l$ and remain in system $n$ at time $t$, and let
$\TfN_{kl}^n(t): = n^{-1/2} \fN_{kl}^n (t)$ denote its the diffusion-scaled
version. (See Section~\ref{subsection: approximation for W, T, N, and tau of
  the predicted classes} for the formal definition.)
\begin{proposition}[Proportion of true class labels]\label{proposition: prop of true class}
Given a classifier $\model$ and a sequence of queueing systems, 
suppose that Assumptions~\ref{assumption: data
generating process},~\ref{assumption: heavy traffic} and~\ref{assumption:
second order moments} hold. Under any work-conserving p-FCFS policy, we have that
for any $k, l\in [K]$ and $t\in [0,1]$,
  \begin{equation}
    \label{eq: proportion of true classes in a predicted class queue}
    \begin{aligned}
      \TfN_{kl}^n (t)
      = \frac{\tp_k^n \q_{kl}^n}
      {\sum_{r=1}^K \tp_{r}^n \q_{rl}^n} \TfN^n_l (t) + o_n (1).
    \end{aligned}
  \end{equation}
\end{proposition}

For any predicted class $l\in[K]$, Proposition~\ref{proposition: prop of true
  class} states the \textit{unobservable} (scaled) queue length of true class
$k$ jobs, $\TfN_{kl}^n$, is proportional to the overall queue length
$\TfN_{l}^n$. Moreover, the proportion is asymptotically ``stable" in the
sense that $ \frac{\tp_k^n \q_{kl}^n}{\sum_{r=1}^K \tp_{r}^n \q_{rl}^n}$
converges to a constant under Assumption~\ref{assumption: heavy traffic}.
Since the actual cost incurred by a job is governed by the job's true class
label, the decomposition~\eqref{eq: proportion of true classes in a predicted
  class queue} enables to approximate the aggregated cost incurred by jobs in predicted
class $l\in [K]$ according to their true class labels (see Eq.~\eqref{eq: cost
  in terms of Nkl} to come for details).

Next, we use Proposition~\ref{prop: convergence and approximation of predicted
  class N, tau, T, and W} to reveal asymptotic relationships between
endogenous processes such as $\VTfW_l^n$ and $\VTfN_l^n$ (e.g., see
Lemma~\ref{lemma: relation between N and W} in Section~\ref{section:
  proof for HT lower bound}). Combining this with the decomposition~\eqref{eq:
  proportion of true classes in a predicted class queue}, we establish a link
between the actual cost incurred in the presence of misclassification errors
and the exogenous component $\sumTfW^n$ (see Eq.~\eqref{eq: Nkl and Wl for
  lower bound} to come). Our analysis allows us to idenfify a lower bound as a
workload allocation over the predicted classes as we characterize in
Proposition~\ref{theorem: HT lower bound}.

\paragraph{Discussion of proof}
The proof of Proposition~\ref{proposition: prop of true class} is nontrivial,
but once we arrive at the decomposition~\eqref{eq: proportion of true classes
  in a predicted class queue}, it sheds light on the construction of the
P$c\mu$ rule (Definition~\ref{definition: modified gcmu rule}).  The main challenge in deriving the cost
functions~\eqref{eq: pcmu cost} used in the P$c\mu$ rule (Definition~\ref{definition: modified gcmu rule})
is the proof of Proposition~\ref{proposition: prop of true class}. We
decompose the stochastic fluctuation $\TfN_{kl}^n$ into fluctuations of other
processes, including the service process $\TfS_l^n$ and the classification
partial sum process, $\TfZ_{kl}^n$. Since service times and the true/predicted
class labels are correlated in our model, it is not a priori clear how the
corresponding fluctuations in $\TfS_l^n$ and $\TfZ_{kl}^n$ jointly influence
that of $\TfN_{kl}^n$. The derivation of~\eqref{eq: proportion of true classes
  in a predicted class queue} requires articulating the stochastic fluctuation of
$\TfN_{kl}^n$. Toward this goal, we provide a novel characterization of the
service completion in the predicted classes from the perspective of the common
stream of arrivals in Eq.~\eqref{eq: characterization of service completions
  in predicted classes}. The proof of the proposition is provided in
Section~\ref{subsection: proportion of true classes}.


Assumption~\ref{assumption: heavy traffic} is essential to prove Theorem~\ref{theorem: HT lower bound}. Proposition~\ref{proposition: prop of true class} relies on Proposition~\ref{prop: convergence and approximation of predicted class N, tau, T, and W}, which builds on the $o(n^{-1/2})$ convergence rate in the heavy traffic condition~\eqref{eq: heavy traffic condition} following from Assumption~\ref{assumption: heavy traffic}. Importantly, the assumption implies $|\fmu_l^n - \fmu_l | = o_n (1)~\forall~l\in [K]$, which allows to establish a pivotal relationship between $\TfW_l^n$ and $\TfN_l^n$ in Lemma~\ref{lemma: relation between N and W}. In Section~\ref{subsection: proof overview of the optimality}, we explain how Assumption~\ref{assumption: heavy traffic} also leads to a crucial equivalence between
the age and sojourn time processes, laying the foundation of the optimality of the
\ourmethod~in Theorem~\ref{theorem: optimality of our policy} to come.

\subsubsection{Comparison to the analysis of~\citet{VanMieghem95}}
\label{subsubsection: detailed comparison to the lower bound result in VM}


Plugging $Q^n=I$ into Theorem~\ref{theorem: HT lower bound}, we recover the
classical result under perfect classification in~\citet[Proposition
6]{VanMieghem95}. In addition to the broader generality of Theorem~\ref{theorem: HT lower bound}, our proof addresses key conditions that are missing in~\citet[Proposition 6]{VanMieghem95} and incomplete justifications in their proof, even within the classical setting when all true classes are known.

\paragraph{Convergence of service rate functions}
As mentioned in Section~\ref{subsection: proof overview for HT lower bound}, the convergence of $\tmu_k^n, \tp_k^n, \q_{kl}^n~\forall~k,l\in [K]$, as stated in Assumption~\ref{assumption: heavy traffic}, plays a critical role in establishing Theorem~\ref{theorem: HT lower bound} by guaranteeing the convergence of the \textit{service rate} $\fmu_l^n $, which is essential for translating the queue length $\TfN_l^n$ into the workload $\TfW_l^n$ (Lemma~\ref{lemma: relation between N and W}). Similarly, we identify that a comparable convergence condition, e.g., $(\bar{V}_k^n)' - (\bar{V}_k^*)' = o_n(1)$ in their notation, is one of the two alternative sufficient conditions that can complement the argument for the lower bound in~\citet{VanMieghem95}. 

Without assuming i.i.d. service times, as in our framework, we found that~\citet[Assumption 1]{VanMieghem95} alone does not suffice to derive~\citet[Proposition 3]{VanMieghem95}, a key intermediate result analogous to our Lemma~\ref{lemma: relation between N and W}, which underpins their lower bound proof in~\citet[Proposition 6]{VanMieghem95}. In Section~\ref{section: adding to the exposition in classical heavy traffic queueing analysis}, we provide detailed analysis and identify conditions that can be incorporated into~\citet{VanMieghem95} to rigorously establish the results claimed in that work.

\paragraph{Continuity of the optimal workload allocation}
In both the classical and our setting, 
the optimal workload allocation $h$, which solves~\eqref{eq: optimization problem}, must be continuous with respect
to the total workload $\sumTtW (t),~t\in [0,1]$. 
This crucial argument was not proven in~\citet{VanMieghem95}, making 
their results incomplete. We address this gap by proving continuity of $h$ in Proposition~\ref{prop: properties of h and opt}.

\paragraph{Parition of the time interval}
In the proof of Theorem~\ref{theorem: HT lower bound}, we partition the time
interval $[0,1]$ to bound the accrued cost over each small subinterval, and
the approximation errors due to the finite partitioning is handled accordingly
(see Eq.~\eqref{eq: using mean value theorem}). In the proof
of~\citet[Proposition 6]{VanMieghem95}, however, the partition is chosen by a
different method than ours, and the author claims that the partition size,
hence the approximation error, can be arbitrarily small without justification.
When the workload is a \textit{general reflected process} as
in~\citet{VanMieghem95}'s setting, we found this claim to be challenging to
prove. As a result, we provide a rigorous justification for their claim in the
particular case when $\sumTtW$ is a reflected Brownian motion in
Section~\ref{subsubsection: complementary proof for the partition size in VM}.

\subsection{Detailed proof of heavy traffic lower bound 
(Theorem~\ref{theorem: HT lower bound})}

We begin by proving~\eqref{eq: definition of J*}. 
We analyze 
$\TJ^n_\policyn (t; Q^n)$ for a fixed $t\in [0,1]$. By definition,
\begin{equation*}
  \TJ^n_\policyn (t; Q^n) 
  = n^{-1}\sum\limits_{l=1}^K \sum\limits_{k=1}^K 
  \int_0^{nt}  C_k^n \left(\ftau_l^n (s)\right)\mathrm{d}\fA_{kl}^n (s).
\end{equation*}
For a fixed $\varepsilon>0$, partition $[0, 1]$ into
$ 0 = t_0 < t_1 < \ldots < t_{M} = 1$ such that
$\sup_i(t_{i+1} - t_i) = \varepsilon$, where $M$ is a constant dependent on
$\varepsilon$. Let
$\mathrm{d} \xi_{kl, i}^n: = \frac{\mathrm{d} \fA_{kl}^n}{\fA_{kl}^n(nt_{i+1})
  - \fA_{kl}^n(nt_i)}$ be a probability measure over $[nt_i, nt_{i+1}]$,
convexity of $C^n_k$ and Jensen's inequality yields
\begin{equation}\label{eq: Jensen's ineq for lower bound}
\begin{aligned}
\TJ^n_\policyn (t; Q^n) 
=&~ n^{-1}\sum\limits_{l=1}^K \sum\limits_{k=1}^K \sum_i 
    \int_{nt_i}^{nt_{i+1}}
    C_k^n (\ftau_l^n (s))\mathrm{d}\fA_{kl}^n (s)\\
=&~ n^{-1}\sum_k \sum_l \sum_i [\fA_{kl}^n(nt_{i+1}) - \fA_{kl}^n(nt_i)]
\mathbb{E}_{\xi_{kl, i}^n} [C_k^n(\ftau_l^n)]\\
\geq&~ n^{-1} \sum_k \sum_l \sum_i [\fA_{kl}^n(nt_{i+1}) - \fA_{kl}^n(nt_i)]
C_k^n( \mathbb{E}_{\xi_{kl, i}^n}[\ftau_l^n]).
\end{aligned}
\end{equation}
By connecting $\mathbb{E}_{\xi_{kl, i}^n}[\ftau_l^n]$ with the workload
process, we can show the following claim. Recall $o_n(1) \rightarrow 0$
uniformly over $t\in [0, 1]$.
\begin{claim}
  \label{claim:cost-lb}
  \begin{align}
    & \TJ^n_\policyn (t; Q^n)
      \geq  n^{-1} \sum_k \sum_l \sum_i 
      [\fA_{kl}^n(nt_{i+1}) - \fA_{kl}^n(nt_i)]
      C_k^n( \mathbb{E}_{\xi_{kl, i}^n}[\ftau_l^n])
      \label{eq: limit of cumulative cost function eq3} \\
    & =  \sum_k \sum_l \sum_i 
      [\lambda p_k \q_{kl} (t_{i+1} - t_i) + o_n(1)]
      \cdot C_k^n\Big(n^{1/2} \Big[
      [\frho_l(t_{i+1} - t_i)]^{-1}
      \int_{t_i}^{t_{i+1}} \TfW^n_l(s) \mathrm{d}s 
      + o_n(1) \Big]\Big). \nonumber 
  \end{align}
\end{claim}
Since $C_k^n(n^{1/2}\cdot)\rightarrow C_k(\cdot)$ and $C_k'$ is bounded on the
compact set $[0, 2\lim\sup_n \| \sumTtW^n \|/\frho_l ]$, the right hand side of
inequality~\eqref{eq: limit of cumulative cost function eq3} can be rewritten as
\begin{equation}
  \label{eq: limit of cumulative cost function eq5}
  \begin{aligned}
    & \sum_i (t_{i+1} - t_i) \sum_k \sum_l \lambda p_k \q_{kl}
      C_k\Big(
      \frac{1}{t_{i+1} - t_i}\int_{t_i}^{t_{i+1}} \TfW^n_l(s) / \frho_l ~\mathrm{d}s \Big) 
      + o_n(1) \\
    & \geq   \sum_i (t_{i+1} - t_i) \sum_k \sum_l
      \lambda p_k \q_{kl} C_k \Big(
      [h(y_i^n)]_l  / \frho_l 
      \Big) + o_n(1) 
  \end{aligned}
\end{equation}
where $h(\cdot)$ is the solution to $\opt(r)$~\eqref{eq: optimization problem}
and
\begin{equation*}
  y_i^n: =  \sum_l\frac{1}{t_{i+1} - t_i} \int_{t_i}^{t_{i+1}} \TfW^n_l(s)~\mathrm{d}s 
  =  \frac{1}{t_{i+1} - t_i} 
  \int_{t_i}^{t_{i+1}} \sumTtW^n(s) \mathrm{d}s.
\end{equation*}
By $\sumTtW^n \rightarrow \sumTtW$ and the continuity of $\sumTtW$ in
Proposition~\ref{prop: convergence and approximation of predicted class N,
  tau, T, and W}, applying the mean value theorem for integrals yields the
existence of $\zeta_i\in[t_i, t_{i+1}]$ such that
\begin{equation}\label{eq: using mean value theorem}
  y^n_i
  = \frac{1}{t_{i+1} - t_i} \int_{t_i}^{t_{i+1}} \TfW^n_+(s) \mathrm{d}s
  = \frac{1}{t_{i+1} - t_i} \int_{t_i}^{t_{i+1}} \sumTtW (s) \mathrm{d}s + o_n(1)
  = \sumTtW (\zeta_i) + o_n(1).
\end{equation}

We use continuity of $h(\cdot)$ to complete the proof of~\eqref{eq: definition 
of J*}.  For any $r\geq 0$,
although $\opt(r)$ can potentially have multiple optimal solutions, it
suffices to study properties of one specific optimal solution.
\begin{lemma}[Properties of the optimal allocation]
  \label{prop: properties of h and opt}
  Given a classifier $f_\theta$, suppose Assumptions~\ref{assumption: data
    generating process},~\ref{assumption: heavy traffic},~\ref{assumption: second
    order moments}, and~\ref{assumption: on cost functions for showing the lower bound} hold.  Let $\underline{h}(0) = \mathbf{0}$ and for any
  $r > 0$, let $\underline{h}(r)$ be the solution to the following equations
  \begin{equation}\label{eq: linear systems for KKT}
    \fmu_l \fC_l' \Big(\frac{x_l}{\frho_l}\Big)=\fmu_m \fC_m' \Big(\frac{x_m}{\frho_m}\Big),~\forall~l,m \in [K];\quad 
    \sum\limits_{l=1}^K x_l =r; \quad 
    x_l \geq 0, ~\forall~l\in [K].
  \end{equation}
  Then, i) for any $r > 0$, there exists a unique solution, ii)
  $\underline{h}: [0, \infty) \rightarrow \R^{K}$ is continuous, iii) for any
  $r\geq 0$, $\underline{h}(r)$ is an optimal solution to $\opt(r)$~\eqref{eq:
    optimization problem}.
\end{lemma}
\noindent See Section~\ref{subsection: continuity of the optimal workload
  allocation} for the proof.

By Lemma~\ref{prop: properties of h and opt} and uniform continuity of $\underline{h}$ and $C_k$
on compact sets (Assumption~\ref{assumption: on cost functions for showing the
  lower bound}),
\begin{equation*}
\begin{aligned}
        \lim\inf_n \TJ^n_\policyn (t;Q^n) 
\geq&~  \liminf_n
        \sum_i  (t_{i+1} - t_i) 
        \sum_k\sum_l   \lambda p_k \q_{kl}
        C_k\Big([h(y_i^n)]_l / \frho_l\Big)\\
=&~  \sum_i  (t_{i+1} - t_i) 
        \sum_k\sum_l   \lambda p_k \q_{kl}
        C_k\Big(\big[h\big( \sumTtW (\zeta_i)\big)\big]_l / \frho_l\Big).
\end{aligned}
\end{equation*}
Note that the function
$\lambda p_k \q_{kl} C_k([h(\sumTtW(\cdot))]_l/ \rho_l)$ is continuous and
thus Riemann integrable. Letting $\varepsilon \rightarrow 0$ results in~\eqref{eq: definition of J*}:
\begin{equation*}
    \lim\inf_n \TJ^n_\policyn (t;Q^n)  \geq 
    \sum_{k=1}^K  \sum_{l=1}^K \int_0^t \tlambda \tp_k \q_{kl} C_k\Big(
    \frac{\big[h \big( \sumTtW(s)\big)]_l}{\frho_l} \Big)  \mathrm{d} s.
\end{equation*}
To show~\eqref{eq: HT lower bound in stochastic sense}, consider feasible p-FCFS policies $\{\pi_n'\}$. For all $n\in \mathbb{N}$, the original processes under $\mathbb{P}^n$ satisfy $\mathbb{P}^n [\TJ^n_{\policy_n'}(t; Q^n) > x] 
= \mathbb{P}_\text{copy} [\TJ^n_{\policy_n'}(t; Q^n) > x],~\forall~x\in\mathbb{R}, t~\in [0, 1]$, according to 
the Skorohod representation. By Fatou's lemma, for any 
$x\in \mathbb{R}, t\in [0, 1]$, we have that 
\begin{equation*}
\begin{aligned}
\lim\inf_n \mathbb{P}^n [\TJ^n_{\pi_n'}(t; Q^n) > x]  
= \lim\inf_n \mathbb{P}_\text{copy} [\TJ^n_{\policy_n'}(t; Q^n) > x]
\geq \mathbb{E}_{\mathbb{P}_\text{copy}}
[ \lim\inf_n \mathbb{I}\{\TJ^n_{\pi_n'} (t;Q^n) > x \}].
\end{aligned}
\end{equation*}
As $\liminf _{n \rightarrow \infty}  \TJ^n_{\pi_n'}(t; Q^n) \geq  \TJ^*(t; Q)$ 
$\mathbb{P}_\text{copy}$-a.s. by~\eqref{eq: definition of J*}, we have that
\begin{equation*}
  \mathbb{E}_{\mathbb{P}_\text{copy}}
  [ \lim\inf_n \mathbb{I}\{\TJ^n_{\pi_n'} (t;Q^n) > x \}]
  \geq \mathbb{E}_{\mathbb{P}_\text{copy}}
  [ \mathbb{I} \{\lim\inf_n \TJ^n_{\pi_n'} (t;Q^n) > x\} ]
  \geq \mathbb{P}_{\text{copy}} [\TJ^* (t;Q) > x].
\end{equation*} 
Combining equations above yields~\eqref{eq: HT lower bound in stochastic sense}
for any feasible p-FCFS policies. We can further extend~\eqref{eq: HT 
 lower bound in stochastic sense} to any feasible policies using 
Lemma~\ref{lemma: p fcfs}. This completes our proof.

\paragraph{Proof of Claim~\ref{claim:cost-lb}} Since
$n^{-1} \fA^n_{kl}(n \cdot)\rightarrow \tlambda \tp_k \q_{kl} e$ by Proposition~\ref{prop:
  joint conv. of predicted class A, U, S, and V},
\begin{equation}\label{eq: limit of cumulative cost function eq2}
\begin{aligned}
    n^{-1}[\fA_{kl}^n(nt_{i+1}) - \fA_{kl}^n(nt_i)]
= \lambda p_k \q_{kl} (t_{i+1} - t_i)  + o_n(1).
\end{aligned}
\end{equation}
Apply the convergence~\eqref{eq: limit of cumulative cost function eq2} to
rewrite $\mathbb{E}_{\xi_{kl, i}^n}[\ftau_l^n]$
\begin{equation}\label{eq: limit of cumulative cost function eq1}
\begin{aligned}
    \mathbb{E}_{\xi_{kl, i}^n}[\ftau_l^n]
  =&~ n^{-1} \big(n^{-1}[\fA_{kl}^n(nt_{i+1}) - \fA_{kl}^n(nt_i)]\big)^{-1}
     \int_{nt_i}^{nt_{i+1}} \ftau_l^n \mathrm{d}\fA^n_{kl}, \\
    =&~ n^{-1} \big[[\lambda p_k \q_{kl}(t_{i+1} - t_i)]^{-1} + o_n(1)\big]
       \int_{nt_i}^{nt_{i+1}} \ftau_l^n \mathrm{d}\fA^n_{kl},
  \end{aligned}
\end{equation}
where the last line holds since
$(x+\Delta x)^{-1}=x^{-1}-\Delta x+o(\Delta x)$.

We approximate $\int_{na}^{nb} \ftau_{l}^n\mathrm{d}\fA^n_{kl}$ using a
variant of Little's Law that we prove in
Section~\ref{section:proof-little-law}.
\begin{proposition}[Little's law]\label{prop: little's law}
  Given a classifier $\model$, suppose
  Assumptions~\ref{assumption: data generating process}, ~\ref{assumption:
    heavy traffic}, and~\ref{assumption: second order moments} hold. Then, for
  any $0\leq a< b \leq 1$
  \begin{subequations}
    \label{eq: little's law}
    \begin{align}
      \label{eq: little's law eq1}
    \frac{n^{-3/2}}{ \lambda^n \tp_k^n \q_{kl}^n (b - a)} 
    \int_{na}^{nb} \ftau_{l}^n\mathrm{d}\fA^n_{kl} 
    - \frac{1}{ \lambda^n \tp_k^n \q_{kl}^n (b - a)}
      \int_a^b \TfN_{kl}^n(t) \mathrm{d} t
      & = o(1),~\forall~k, l\in[K], \\
      \label{eq: little's law eq2}
    \frac{n^{-3/2}}{ \lambda^n \fp_l^n (b - a)} 
    \int_{na}^{nb} \ftau_{l}^n\mathrm{d}\fA^n_{l} 
    - \frac{1}{ \lambda^n \fp_l^n (b - a)}
      \int_a^b \TfN_{l}^n(t) \mathrm{d} t
      & = o(1),~\forall~ l\in[K].
  \end{align}
  \end{subequations}
Moreover, if either
  $\Tftau^n_l \rightarrow \Tftau_l \in \mathcal{C}$ or
  $\TfN_l^n\rightarrow \TfN_l \in \mathcal{C}$ holds, then
    $\tlambda\fp_l \Tftau_l = \TfN_l$.
\end{proposition}
\noindent Applying the proposition
$n^{-3/2} \int_{na}^{nb} \ftau_{l}^n\mathrm{d}\fA^n_{kl} - \int_a^b
\TfN_{kl}^n(t) \mathrm{d} t = o_n(1) O(|b-a|)$ to Eq.~\eqref{eq: limit of
  cumulative cost function eq1},
\begin{equation}\label{eq: cost in terms of Nkl}
  \begin{aligned}
    \mathbb{E}_{\xi_{kl, i}^n}[\ftau_l^n]
    =&~ n^{1/2} \big[[\lambda p_k \q_{kl} (t_{i+1} - t_i)]^{-1} + o_n(1)\big] \Big(\int_{t_i}^{t_{i+1}} \TfN_{kl} ^n(s) \mathrm{d} s + o_n(1)O(t_{i+1} - t_{i})\Big)\\
    =&~ n^{1/2} \Big[[\lambda_k p_k \q_{kl}(t_{i+1} - t_i)]^{-1} \int_{t_i}^{t_{i+1}} \TfN^n_{kl}(s) \mathrm{d}s + o_n(1) + o_n(1)O(\varepsilon)\Big],
  \end{aligned}
\end{equation}
since $\sup_i (t_{i+1} - t_i) = O(\varepsilon)$ and
$\limsup_n \|\TfN_{kl}\| \leq \limsup_n \|\TfN_l\| < \infty$ by
Proposition~\ref{prop: convergence and approximation of predicted class N,
  tau, T, and W}.

To rewrite $\int_a^b \TfN_{kl}^n(s) \mathrm{d} s$ in terms of the workload,
recall the key relation
$\TfN_{kl}^n = \frac{p_k^n \q_{kl}^n}{\sum_r p_{r}^n \q_{rl}^n} \TfN^n_l +
o_n(1)$ given in Proposition~\ref{proposition: prop of true class} (see
Section~\ref{subsection: proportion of true classes} for its proof).  We can
further approximate the queue length process $\TfN^n_l$ using the service rate
$\fmu_l$ and the remaining workload process $\TfW^n_l$.
\begin{lemma}[Relation between $\TfW^n_l$ and $\TfN^n_l$]
  \label{lemma: relation between N and W}
  Given a classifier $\model$, suppose Assumptions~\ref{assumption: data
    generating process}, ~\ref{assumption: heavy traffic},
  and~\ref{assumption: second order moments} hold. Then, for p-FCFS policies
  $\fmu_l \TfW^n_l - \TfN^n_{l} \rightarrow 0$ for all $l \in [K]$.
\end{lemma}
\noindent See Section~\ref{section:proof-relation between N and W} for the
proof. Applying Proposition~\ref{proposition: prop of true class} and
Lemma~\ref{lemma: relation between N and W},
\begin{equation*}
\begin{aligned}
    \int_a^b  \TfN_{kl}^n(s) \mathrm{d} s 
  =&~ \int_a^b  \fmu_l \frac{p_k \q_{kl}}{\sum_r p_{r} \q_{rl}}\TfW^n_l(s) \mathrm{d}s + o_n(1) O(|b-a|) .
\end{aligned}
\end{equation*}
Plugging this into the expression~\eqref{eq: cost in terms of Nkl} for
$\mathbb{E}_{\xi_{kl, i}^n}[\tau_l^n]$ 
\begin{equation}\label{eq: Nkl and Wl for lower bound}
\begin{aligned}
  \mathbb{E}_{\xi_{kl, i}^n}[\ftau_l^n]
  =&~ n^{1/2} \Big[[\lambda_k p_k \q_{kl}(t_{i+1} - t_i)]^{-1} 
     \Big(\int_{t_i}^{t_{i+1}}  
     \fmu_l\frac{p_k \q_{kl}}{\sum_r p_{r} \q_{rl}}\TfW^n_l(s) \mathrm{d}s 
     + o_n(1) O(t_{i+1} - t_{i})\Big)   +  o_n(1)  \Big]\\
  =&~ n^{1/2} \Big[[\frho_l (t_{i+1} - t_i)]^{-1} 
     \int_{t_i}^{t_{i+1}} \TfW^n_l(s) \mathrm{d}s 
     + o_n(1) \Big],
\end{aligned}
\end{equation}
where we use the shorthands $\fp_l = \sum_r p_{r} \q_{rl}$ and
$\frho_l = \frac{\lambda\fp_l}{\fmu_l}$ in the final line.

\subsection{Proof of Lemma~\ref{prop: properties of h and opt} (Gap A2)}
\label{subsection: continuity of the optimal workload allocation}  

For any $l\in [K]$, Assumption~\ref{assumption: on cost functions for showing
  the lower bound} implies $\fC'_l$ is continuous and strictly
increasing. Hence,
\begin{equation} \label{eq: workload as a function of x1}
  g(x) : = x + \sum_{l=2}^K \frho_l\cdot (\fC_l')^{-1}\Big( \frac{\fmu_1}{\fmu_l} \fC'_1(\frho_1^{-1}x)\Big)
\end{equation}
is continuous and strictly increasing with $g(0) = 0$, $g(r) \geq r$.  Let
$x_1(r)$ be a unique solution to $g(x) = r$ and
\begin{equation*}
  x_l(r) : =\frho_l\cdot (\fC_l')^{-1}\Big( \frac{\fmu_1}{\fmu_l} \fC'_1(\frho_1^{-1}x_1)\Big),~\forall~l \geq 2.
\end{equation*}
Evidently, $\underline{h}(r) = (x_1(r), \ldots, x_K(r))$ is a unique solution
to Eq.~\eqref{eq: linear systems for KKT}. To see (ii), note that $x_1(r)$ is
continuous with $x_1(0) = 0$ since $g^{-1}$ is continuous.
For (iii), for $r = 0$, $\underline{h}(0) = \mathbf{0}$ is clearly an optimal
solution to $\opt(0)$. When $r > 0$, we verify $\underline{h}(r)$ satisfies
the KKT conditions for $\opt(r)$. This is evident from the fact that
$\fC'_l(0) = (\fC'_l)^{-1}(0)=0$ and $\fC'_l$ and $(\fC'_l)^{-1}$ are strictly
increasing (Assumption~\ref{assumption: on cost functions for showing the
  lower bound}).

\subsection{Proof of Proposition~\ref{prop: little's law}}
\label{section:proof-little-law}

Our proof for Eqs.~\eqref{eq: little's law eq1} and~\eqref{eq: little's law
  eq2} is similar to the proof for \citet[Proposition 4]{VanMieghem95}. To
show Eqs.~\eqref{eq: little's law eq1}, consider the cumulative cost during
$t\in [a, b]$ where each job incurs a unit cost per unit time spent in the
system. We study three different cost charging schemes
\begin{equation*}
\begin{aligned}
  \text{Cost}_1^n(a, b) 
  & = \frac{1}{ \lambda^n \tp_k^n \q_{kl}^n (b - a) }
  \sum_{i = \fA^n_{kl}(na)}^{\fA^n_{kl}(nb)} \ftau^n_{li}, \\
   \text{Cost}_2^n(a, b)  
  & = \frac{1}{ \lambda^n \tp_k^n \q_{kl}^n (b - a)  }
  \int_{na}^{nb} \fN^n_{kl}(t) \mathrm{d}t,\\
  \text{Cost}_3^n(a, b)
  &= \frac{1}{ \lambda^n \tp_k^n \q_{kl}^n (b - a) }
  \sum_{i = \fA^n_{kl}(na)}^{\fA^n_{kl}(nb) - \fN^n_{kl}(nb)} \ftau^n_{li}.
\end{aligned}
\end{equation*}
$\text{Cost}_1^n(a, b)$ charges the entire cost at the job's arrival,
$\text{Cost}_3^n(a, b)$ at the job's departure, and $\text{Cost}_2^n(a, b)$
continuously. It is easy to verify 
\begin{equation*}
    \text{Cost}_3^n(a, b) \leq \text{Cost}_2^n(a, b) \leq \text{Cost}_1^n(a, b),
\end{equation*}
and
\begin{equation*}
\begin{aligned}
    n^{-3/2} (\text{Cost}_1^n(a, b) - \text{Cost}_3^n(a, b))
=&~     \frac{n^{-3/2}}{ \lambda^n \tp_k^n \q_{kl}^n (b - a) }  
        \sum_{i = \fA^n_{kl}(nb) - \fN^n_{kl}(nb)+1}^{\fA^n_{kl}(nb)} \ftau^n_{li}\\ 
\leq&~  \frac{n^{-1/2}}{ \lambda^n \tp_k^n \q_{kl}^n (b - a) } 
        \|\TfN^n_{kl}\|  \| \Tftau^n_l\| \rightarrow 0,
\end{aligned}
\end{equation*}
since $ \lambda^n \tp_k^n \q_{kl}^n \rightarrow  \lambda \tp_k \q_{kl} $ by Assumption~\ref{assumption: heavy
  traffic}, and $\|\TfN^n_{kl}\|$ and $\| \Tftau^n_l\|$ are bounded
(Proposition~\ref{prop: convergence and approximation of predicted class N,
  tau, T, and W}). Conclude
\begin{equation*}
\begin{aligned}
    o_n(1)
=&~ \frac{n^{-3/2}}{ \lambda^n \tp_k^n \q_{kl}^n (b - a)} 
    \int_{na}^{nb}  \ftau_{l}^n\mathrm{d}\fA^n_{kl} 
    -  \frac{n^{-3/2}}{ \lambda^n \tp_k^n \q_{kl}^n (b - a)}  
    \int_{na}^{nb} \fN_{kl}^n(t) \mathrm{d} t \\
=&~ \frac{n^{-3/2}}{ \lambda^n \tp_k^n \q_{kl}^n (b - a)}  
     \int_{na}^{nb}  \ftau_{l}^n\mathrm{d}\fA^n_{kl} 
     -\frac{1}{ \lambda^n \tp_k^n \q_{kl}^n (b - a)}  
    \int_{a}^{b} \TfN_{kl}^n(t) \mathrm{d} t.
\end{aligned}
\end{equation*}
The proof for Eq.~\eqref{eq: little's law eq2} can be established similarly
and we omit the details.

For the second part of the proposition, recall from Proposition~\ref{prop: convergence and approximation of predicted class N, tau, T, and W}(ii) that for
all $l\in [K]$, $\Tftau^n_l \rightarrow \Tftau_l \in \mathcal{C}$ and $\TfN^n_l \rightarrow \TfN_l \in \mathcal{C}$ are equivalent. Recall that by Eq.~\eqref{eq:
  little's law eq2},
\begin{equation*}
    \frac{n^{-3/2}}{ \lambda^n \fp_l^n (b - a)} 
    \int_{na}^{nb} \ftau_{l}^n\mathrm{d}\fA^n_{l} 
    - \frac{1}{\lambda^n \fp_l^n (b - a)}
    \int_a^b \TfN_{l}^n(t) \mathrm{d} t   = o_n(1),~\forall~ l\in[K].
\end{equation*}
For simplicity, for fixed $[a, b]$, let $\xi^n_l$ be the Lebesgue-Stieltjes
measure on $[0, 1]$ induced by $ n^{-1} \fA_l^n(n \cdot)$ and $\xi_l$ be the
Lebesgue-Stieltjes measure on $[0, 1]$ induced by $ \lambda \fp_l e$.  It is
easy to verify
\begin{equation*}
    n^{-3/2}\int_{na}^{nb} \ftau^n_l(t) \mathrm{d} \fA_l^n(t)
    = \int_a^b \Tftau^n_l~\mathrm{d}\xi^n_l.
\end{equation*}
Since $\xi^n_l \rightarrow \xi_l$ on sets and
$\| \Tftau^n_l \| \leq \limsup_n \|\Tftau^n_l\| < +\infty$ eventually
(Proposition~\ref{prop: convergence and approximation of predicted class N,
  tau, T, and W}), generalized Lebesgue convergence~\cite[Page 270]{Royden88}
implies
\begin{equation}\label{eq: little's law eq3}
n^{-3/2}\int_{na}^{nb} \ftau^n_l(t) \mathrm{d} \fA_l^n(t)
\rightarrow \int_a^b \Tftau_l(t) \mathrm{d} (\lambda \fp_l e) (t)
= \int_a^b \lambda \fp_l \Tftau_l(t) \mathrm{d}t.
\end{equation}
Next, we analyze the second term of Eq.~\eqref{eq: little's law
  eq2}. Dominated convergence gives
\begin{equation}\label{eq: little's law eq4}
\frac{1}{\lambda^n \fp_l^n (b - a)} \int_{a}^{b} \TfN^n_l(t) \mathrm{d} t  \rightarrow
\frac{1}{\lambda \fp_l (b - a)} \int_a^b \TfN_l(t) \mathrm{d} t.
\end{equation}

Combining Eqs.~\eqref{eq: little's law eq2},~\eqref{eq: little's law eq3}, and~\eqref{eq: little's law eq4} yields that for all $[a, b]\subset[0, 1]$,
\begin{equation}\label{eq: little's law eq5}
    \frac{1}{\lambda \fp_l  ( b - a)}
     \int_a^b \lambda  p_l \Tftau_l(t) \mathrm{d}t
    = \frac{1}{\lambda \fp_l  ( b - a)} \int_a^b \TfN_l(t) \mathrm{d} t.
\end{equation}
For a fixed $t\in [0, 1]$, inserting $a = t, b= t+\Delta t$ into  Eq.~\eqref{eq: little's law eq5} gives 
\begin{equation}\label{eq: little's law eq6}
\begin{aligned}
    \frac{1}{\lambda \fp_l (t+\Delta t) - \lambda \fp_l t} \int_t^{t + \Delta t} \TfN_l(s) \mathrm{d} s
&~ = \frac{1}{\lambda\fp_l} \cdot \frac{1}{\Delta t} \int_t^{t + \Delta t} \TfN_l(s) \mathrm{d} s  \rightarrow \frac{1}{\lambda \fp_l} \TfN_l(t),
\end{aligned}
\end{equation}
as $\Delta t \rightarrow 0$, where the convergence follows from continuity of $\TfN_l$ by Eqs.~\eqref{eq: approximation of tilde T} and~\eqref{eq: relation between W, tau, L}, and the mean value theorem for definite integrals.
Similarly, one can show as $\Delta t \rightarrow 0$,
\begin{equation}\label{eq: little's law eq7}
\frac{1}{\lambda \fp_l (t+\Delta t) - \lambda \fp_l t} \int_t^{t+ \Delta t} \lambda \fp_l \Tftau_l(s) \mathrm{d} s
\rightarrow \Tftau_l(t).
\end{equation}
Combining Eqs.~\eqref{eq: little's law eq5},~\eqref{eq: little's law eq6},
and~\eqref{eq: little's law eq7} yields the desired result
$\tlambda\fp_l \Tftau_l = \TfN_l,~\forall~l\in [K]$.

\subsection{Proof of Proposition~\ref{proposition: prop of true class}}
\label{subsection: proportion of true classes}

Recalling the definition~\eqref{eq:jobs}, for any $nt\in [0,n]$
\begin{equation*}
\begin{aligned}
\fN_{kl}^n (nt) 
=  \fA_{kl}^n (nt) - \sum_{i=1}^{(\fM_l^n \circ \fS_l^n \circ \fT_l^n)(nt)} \tY_{ik}^n \fY_{il}^n
= \fA_{kl}^n (nt) - \fZ^n_{kl}\Big((\fM_l^n \circ \fS_l^n \circ \fT_l^n)(nt)\Big).
\end{aligned}
\end{equation*}
By Lemma~\ref{lemma: uniform convergence of UZRV}, Proposition~\ref{prop:
  joint conv. of predicted class A, U, S, and V}, and Eq.~\eqref{eq:
  approximation of T},
\begin{equation*}
\begin{aligned}
\fZ_{kl}^n (nt) &= n \tp_k^n \q_{kl}^n t + n^{1/2} \TfZ_{kl}^n (t) + o(n^{1/2}),\\
\fS_l^n (nt) &= n\fmu_l^n t + n^{1/2} \TfS_l^n (t) + o(n^{1/2}),\\
\fT_l^n (nt) &= n \tlambda^n \fp_l^n (\fmu_l^n)^{-1} t +n^{1/2} \TfT_l^n (t) + o(n^{1/2}).
\end{aligned}
\end{equation*}

Recalling $\| \TfT^n_l\| <+\infty$ by Proposition~\ref{prop: convergence and 
approximation of predicted class N, tau, T, and W}, $(\fS_l^n \circ \fT_l^n)(nt)$ can be reformulated as 
\begin{equation*}
\begin{aligned}
    (\fS_l^n \circ \fT_l^n)(nt) 
=&~ \fmu^n_l \fT_l^n(nt)  + n^{1/2} \TfS^n_l(n^{-1} \fT_l^n(nt)) 
    + o(n^{1/2})\\
=&~  n \lambda^n \fp^n_l t + n^{1/2} \fmu^n_l \TfT^n_l(t) 
    + n^{1/2} \TfS^n_l(
     \lambda^n \fp^n_l (\fmu^n_l)^{-1} t + n^{-1/2} \TfT^n_l(t) + o(n^{-1/2}))
     + o(n^{1/2}) \\
    =& n \lambda^n \fp^n_l t + n^{1/2} \fmu^n_l \TfT^n_l(t) 
    + n^{1/2} \TfS^n_l(
        \lambda^n \fp^n_l (\fmu^n_l)^{-1} t + o_n (1))
    + o(n^{1/2}).
\end{aligned}
\end{equation*}
Therefore, we can rewrite $(\fM_l^n \circ \fS_l^n \circ \fT_l^n)(nt)$ as
\begin{equation} \label{eq: characterization of service completions in predicted classes}
\begin{aligned}
(\fM_l^n \circ \fS_l^n \circ \fT_l^n)(nt)
=&~n \tlambda^n t 
 + n^{1/2} (\fp_l^n)^{-1} \fmu_l^n \TfT_l^n (t) 
+ n^{1/2} (\fp_l^n)^{-1} \TfS_l^n ( \tlambda^n \fp_l^n (\fmu_l^n)^{-1} t 
 + o(1)) \\
&~ +n^{1/2} \TfM_l^n ( \tlambda^n \fp_l^n  t + o_n (1))
 + o(n^{1/2}).
\end{aligned}
\end{equation}
Since $\fA^n_{kl}(nt) = \lambda^n p_k^n \q^n_{kl} nt + n^{1/2} \TfA_{kl}^n(t) + o_n (1)$ and $\fZ_{kl}^n(nt) = p^n_k \q^n_{kl} nt + n^{1/2} \TfZ_{kl}^n(t) + o_n (1)$ by the proof of Proposition~\ref{prop: joint conv. of predicted
  class A, U, S, and V} and
Definition~\ref{def: concerned processes for arrival and service}, combining
equations above yields
\begin{equation*}
\begin{aligned}
    \TfN_{kl}^n (t) =&~ \TfA_{kl}^n (t) - \frac{\tp_k^n \q_{kl}^n}{\fp_l^n} \Big[ \fmu_l^n \TfT_l^n (t) 
    + \TfS_l^n (\tlambda^n \fp_l^n (\fmu_l^n)^{-1} t + o_n(1))\Big] \\
    &~- \tp_k^n \q_{kl}^n \TfM_l^n (\tlambda^n \fp_l^n t + o_n(1)) 
    - \TfZ_{kl}^n (\tlambda^n t + o_n (1)) + o_n(1).
\end{aligned}
\end{equation*}
Moreover, the proof of Proposition~\ref{prop: joint conv.  of predicted class
  A, U, S, and V} implies
\begin{equation*}
    \TfS_l^n \rightarrow \TfS_l, \quad 
    \TfA_{kl}^n  \rightarrow \TfA_{kl}: = \TfZ_{kl} \circ \tlambda e + \tp_k \q_{kl}\TtA_0, \quad 
    \TfM_l^n \rightarrow \TfM_l: = -\fp_l^{-1} \Big(\sum_{k=1}^K\TfZ_{kl} \Big)\circ \fp_l^{-1} e.
\end{equation*}
Since the limiting processes are continuous, by continuity of 
composition~\citep[Theorem 13.2.1]{Whitt02} and continuous mapping theorem, we 
have that 
\begin{equation*}
\begin{aligned}
&\TfS_l^n (\tlambda^n \fp_l^n (\fmu_l^n)^{-1} t +o_n(1)) = 
\TfS_l^n (\tlambda^n \fp_l^n (\fmu_l^n)^{-1} t) + o_n(1), 
\quad \TfA_{kl}^n  - \TfZ_{kl}^n ( \tlambda^n e + o_n(1)) 
= \tp_k^n \q_{kl}^n \TfA_0^n  + o_n(1),\\
&\TfM_{l}^n ( \tlambda^n \fp_l^n \cdot + o_n(1))= -(\fp_l^n)^{-1} \Big( \sum_{k=1}^K\TfZ_{kl}^n ( \tlambda^n \cdot ) \Big)  + o_n(1).
\end{aligned}
\end{equation*}
Thus, we can further rewrite $\TfN_{kl}^n (t) $ as 
\begin{equation*}
\begin{aligned}
&~\tp_k^n \q_{kl}^n \TfA_0^n(t)  - \frac{\tp_k^n \q_{kl}^n}{\fp_l^n}\Big[ \fmu_l^n \TfT_l^n (t) + \TfS_l^n (\tlambda^n \fp_l^n (\fmu_l^n)^{-1} t )\Big] + \frac{\tp_k^n \q_{kl}^n}{\fp_l^n} \sum_{k=1}^K\TfZ_{kl}^n ( \tlambda^n t) + o_n(1)\\
=&~\frac{\TfA_l^n (t) - \sum_{k=1}^K\TfZ_{kl}^n ( \tlambda^n t)}{\fp_l^n} \tp_k^n \q_{kl}^n- \frac{\tp_k^n \q_{kl}^n}{\fp_l^n}
\Big[ \fmu_l^n \TfT_l^n (t) + \TfS_l^n (\tlambda^n \fp_l^n (\fmu_l^n)^{-1} t )\Big] 
+ \frac{\tp_k^n \q_{kl}^n}{\fp_l^n} \sum_{k=1}^K\TfZ_{kl}^n ( \tlambda^n t) + o_n(1)\\
=&~\frac{\tp_k^n \q_{kl}^n}{\fp_l^n} \left[\TfA_l^n (t) - \fmu_l^n \TfT_l^n (t) - \TfS_l^n (\tlambda^n \fp_l^n (\fmu_l^n)^{-1} t )\right] + o_n (1),\\
\end{aligned}
\end{equation*}
where the second line follows from the identity
$\TfA_l^n (t) = \sum\limits_{k=1}^K \TfZ^n_{kl} (\tlambda^n t) + \fp^n_l
\TtA^n_0 (t) + o_n (1)$ we derived in the proof of Proposition~\ref{prop:
  joint conv. of predicted class A, U, S, and V}. Then, 
by~\eqref{eq: approximation of predicted N}, we have the desired result
$\TfN_{kl}^n (t) = \frac{\tp_k^n \q_{kl}^n}{\fp_l^n} \TfN_l^n (t) + o_n(1)$.

\subsection{Proof of Lemma~\ref{lemma: relation between N and W} (Gap A1)}
\label{section:proof-relation between N and W}

For any $nt\in [0, n]$, let $\fv^n_l(nt)$ be the amount of service, if any,
already given to the oldest predicted class $l$ job present in the system at
time $nt$. By definition,
$t\in[0, 1]$,
\begin{equation}\label{eq: relation between N and W eq1}
  \begin{aligned}
    \TfW^n_l(t) 
    =&~     n^{-1/2}\big[\fV_l^n( \fA_{l}^n(nt)) 
       - \fV_l^n\big( \fA_{l}^n(nt) - \fN^n_{l} (nt)\big) 
       - \fv^n_l(nt)\big]\\
    =&~      n^{1/2}      
       \big[ (\fmu^n_l)^{-1} (n^{-1}\fA^n_{l}(nt)) 
       -  (\fmu^n_l)^{-1} \big( n^{-1}\fA_{l}^n(nt) - n^{-1}\fN^n_{l} (nt)\big)\big] \\
     & + \big[
       \TfV_l^n( n^{-1}\fA_{l}^n(nt)) 
       - \TfV_l^n\big(n^{-1} \fA_{l}^n(nt) - n^{-1} \fN^n_{l} (nt)\big)\big]
       + o_n(1) - n^{-1/2}\fv^n_l(nt),
  \end{aligned}
\end{equation}
where the second equality follows from Proposition~\ref{prop: joint conv. of
  predicted class A, U, S, and V}, and $o_n(\cdot)$ is uniform over
$t\in [0, 1]$. Under the predicted-class GI/GI/1 structure,
$\BfV^n_l(s) = s/\fmu^n_l$ is linear, so the first term in~\eqref{eq:
  relation between N and W eq1} collapses to
\begin{align*}
  &n^{1/2}
    \big[ (\fmu^n_l)^{-1} (n^{-1}\fA^n_{l}(nt))
    -  (\fmu^n_l)^{-1} \big( n^{-1}\fA_{l}^n(nt) - n^{-1}\fN^n_{l} (nt)\big)\big]\\
  =&~     n^{-1/2} (\fmu^n_l)^{-1} \fN^n_l (nt)
    =    (\fmu^n_l)^{-1} \TfN^n_l (t).
\end{align*}
Since $\fmu^n_l \to \fmu_l$ (Assumption~\ref{assumption: heavy traffic}) and
$\limsup_n \|\TfN^n_l\| < +\infty$ (Proposition~\ref{prop: convergence and
  approximation of predicted class N, tau, T, and W}),
$(\fmu^n_l)^{-1} \TfN^n_l (t) = (\fmu_l)^{-1} \TfN^n_l (t) + o_n(1)$
uniformly in $t$. The collapse-first ordering means only convergence
$\fmu^n_l \to \fmu_l$ (no rate) is used in this step.

It remains to bound the second term in Eq.~\eqref{eq: relation between N and
  W eq1}. Notice that since $\TfV_l^n = \TfV_l + o_n(1)$ where $\TfV_l$ is
uniformly continuous on compact intervals by Proposition~\ref{prop: joint
  conv. of predicted class A, U, S, and V} and
$\limsup_n \|\TfN^n_l\| < +\infty $ by Proposition~\ref{prop: convergence
  and approximation of predicted class N, tau, T, and W},
\begin{align*}
  &  \TfV_l^n( n^{-1}\fA_{l}^n(nt))
    - \TfV_l^n(n^{-1} \fA_{l}^n(nt) - n^{-1} \fN^n_{l} (nt)) \\
   =&~ \TfV_l( n^{-1}\fA_{l}^n(nt))
    - \TfV_l(n^{-1} \fA_{l}^n(nt) - n^{-1} \fN^n_{l} (nt)) + o_n(1)
    = o_n(1).
\end{align*}

\subsection{Complementary proof for Proposition 6 in~\citet{VanMieghem95} (Gap A3)}
\label{subsubsection: complementary proof for the partition size in VM}

Compared to the proof of~\citet[Proposition 6]{VanMieghem95}, we adopt a
different partition of the time interval $[0, 1]$ to derive Eq.~\eqref{eq:
  using mean value theorem} using the mean-value theorem.  To show an
analogous result~\cite[Eq.~(94)]{VanMieghem95}, Van Mieghem picks a partition
using stopping times of $\sumTtW$ to ensure sufficiently small variation of
$\sumTtW$ over each subinterval. Without justification, \citet{VanMieghem95}
claims the partition size is small enough ($O(\varepsilon)$).

Despite best efforts, we found proving this claim difficult when the workload
$\sumTtW$ is a general reflected process.  When $\sumTtW$ is a reflected
Brownian motion, we give a proof that the partition size is still
$\sup_i (t_{i+1} - t_i) = O(\varepsilon)$ in Lemma~\ref{lemma: stopping time
  of tilde w} below; hence~\eqref{eq: using mean value theorem} would follow
even if $\{t_i\}$ is chosen as the stopping times as in~\citet{VanMieghem95}.
Our proof exploits the almost sure non-differentiablity of sample paths of
reflected Brownian motions.  (Alternatively, our previous proof provides a
simple justification for \cite[Eq.~(94)]{VanMieghem95} using our mean value
theorem result~\eqref{eq: using mean value theorem}.)

\begin{lemma}[Stopping times of $\sumTtW$]
  \label{lemma: stopping time of tilde w}
  Given $\varepsilon>0$, consider the sequence of stopping times
  $\{t_i(\varepsilon): i \in \mathbb{N}\}$ of $\sumTtW$
  \begin{equation*}
    \begin{aligned}
      &t_1(\varepsilon)=\min \{1, \inf \{0<t \leq 1:|\TtW_{+}(t)-\lfloor\TtW_{+}(0) / \varepsilon\rfloor \varepsilon| \geq \varepsilon\}\}, \\
      &t_{i+1}(\varepsilon)=\min
        \{1, \inf \{t_i(\varepsilon) <t \leq 1:|\TtW_{+}(t)-\TtW_{+}(t_i(\varepsilon))| \geq \varepsilon\}\}.
    \end{aligned}
  \end{equation*}
  Then, we have that
  \begin{equation*}
    \lim_{\varepsilon \rightarrow 0 }\sup_i (t_{i+1}(\varepsilon) - t_i(\varepsilon)) = 0
  \end{equation*}
\end{lemma}

\begin{proof}
  We prove by contradiction and will show that if Lemma~\ref{lemma: stopping
    time of tilde w} does not hold, then there exists $[a, b]\subset [0, 1]$
  such that $b-a > 0$ and $\sumTtW$ is a constant on $[a, b]$. We argue that
  the latter leads to a contraction using that $\sumTtW$ is a reflected
  Brownian motion as shown in Proposition~\ref{prop: convergence and
    approximation of predicted class N, tau, T, and W}. If $\sumTtW (t)=0$ for
  $t\in [a,b]$, then the underlying Brownian motion must stay at its running
  minimum throughout $[a,b]$ because of the definition of the reflection
  mapping~\citep{Whitt02}, but this is a zero probability
  event~\citep{MortersPe10}. If $\sumTtW (t)=c$ for some positive constant $c$
  and $t\in [a,b]$, it is contradictory to the nondifferentiability of
  Brownian motion~\citep{MortersPe10}.

  Suppose for the purpose of contradiction that there exists some
  $\delta > 0$, a sequence of $\varepsilon_k \rightarrow 0$, and a sequence of
  $\{i_k\}_{k=1}^{\infty}$ satisfying
  \begin{equation*}
    \begin{aligned}
      t_{i_k + 1}(\varepsilon_k) - t_{i_k}(\varepsilon_k) \geq \delta, \text{ and }
      |\TtW_+ (t) - \TtW_+(t_{i_k}(\varepsilon_k))| \leq \varepsilon_k
      ,~\forall~
      t \in [t_{i_k}(\varepsilon_k), t_{i_k}(\varepsilon_k) + \delta]\subset [0, 1].
    \end{aligned}
  \end{equation*}
  Let
  $I(k) =[t_{i_k}(\varepsilon_k), t_{i_k}(\varepsilon_k) + \delta]\subset[0,
  1]$ for all $k \geq 1$. We claim that there exists $b - a \geq \delta_0 > 0$
  and a subsequence $\{k_l\}_{l=1}^\infty$ such that $[a, b] \subset I(k_l)$
  for all $l \geq 1$.  Let $M = \lceil 2/\delta\rceil$ . Partition $[0, 1]$
  into $$0=a_0 < a_1 < \cdots < a_M = 1$$ with
  $a_{r + 1} - a_r = \delta / 2 > 0$, possibly except the last
  interval. Evidently, there exists some $r_0 \in \{0,1,..., M - 1\}$ such
  that $[a_{r_0}, a_{r_0 + 1}] \cap I(k) \neq \emptyset,$ for infinitely many
  $k$'s; otherwise
  $\sum_{m=0}^{M-1} \# \{k : [a_{m}, a_{m+1}] \cap I(k) \neq \emptyset, k\in
  \mathbb{N}\} < + \infty$, so that
  $\sum_{m=0}^{M-1} \# \{k : [a_{m}, a_{m+1}] \cap I(k) \neq \emptyset, k\in
  \mathbb{N}\} \geq \#\{k: k\in \mathbb{N}_+\} = \infty$ gives a
  contradiction.

  We next construct the aforementioned interval $[a, b]$ and subsequence $\{k_l\}_{l=1}^\infty$.
  Since $[a_{r_0}, a_{r_0 + 1}] \cap I(k) \neq \emptyset$ for infinitely many $k$, at least one of the following statement hold:
\begin{enumerate}[(i)]
\item there exists a subsequence $\{k_l\}_{l=1}^\infty$ such that
$t_{i_{k_l}} (\varepsilon_{k_l})>a_{r_0} $ for all $l$;
\item there exists a subsequence $\{k_l\}_{l=1}^\infty$ such that
$t_{i_{k_l}}(\varepsilon_{k_l}) + \delta < a_{r_0 + 1} $ for all $l$;
\item there exists a subsequence $\{k_l\}_{l=1}^\infty$ such that
$ t_{i_{k_l}}(\varepsilon_{k_l}) \leq a_{r_0} < a_{r_0 + 1} \leq t_{i_{k_l}}(\varepsilon_{k_l}) + \delta  $ for all $l$. 
\end{enumerate}
For the case of (i), by definition we have that for all $l$, $a_{r_0} < t_{i_{k_l}}(\varepsilon_{k_l}) \leq a_{r_0+1}$ since $I(k_l) \cap [a_{r_0}, a_{r_0+ 1} ] \neq \emptyset$. Therefore, for all $l$, we have that $a_{r_0} + \delta < t_{i_{k_l}} (\varepsilon_{k_l})+ \delta \leq a_{r_0+1} + \delta$. In other words, 
$[a_{r_0 + 1}, a_{r_0}+\delta] \subset I_{k_l},~\forall~ l\geq 1$.
Hence, we can set $a = a_{r_0+1}, b= a_{r_0} + \delta$, where $b-a \geq \delta/2$. For (ii) and (iii),
we can construct $a$ and $b$ similarly and we skip the details here.

Then, by $[a, b]\subset I(k_l)$, we have that
$\sup_{a \leq t, t' \leq b}|\TtW_+(t) - \TtW_+(t')| \leq 2\varepsilon_{k_l},~\forall~l\geq 1$,
which implies that $\TtW_+(t)$ is a constant on $[a, b]$. This completes our proof.
\end{proof} 


\section{Proof of the heavy traffic optimality of \ourmethod~(Theorem~\ref{theorem: optimality of our policy})}
\label{section: proof for optimality of our policy}

\subsection{Overview of the proof}
\label{subsection: proof overview of the optimality}

Our goal is to show convergence of $\Tftau^n$ to a limit satisfying~\eqref{eq: KKT condition for tau} under \ourmethod, from which
Lemma~\ref{lemma: convergence of tildeJ} will imply Theorem~\ref{theorem:
  optimality of our policy}.

\paragraph{Relationships between $(\Tftau_l^n, \TfN_l^n, \TfT_l^n, \TfW_l^n)$ and $\Tfa_l^n$} 
Since the \ourmethod~uses observable \textit{ages}, we need to connect $\Tfa_l^n$
and the endogenous processes
$(\Tftau_l^n, \TfN_l^n, \TfT_l^n, \TfW_l^n),~\forall~l\in [K]$. We prove the
equivalence between the original KKT conditions~\eqref{eq: KKT condition for
  tau} and the modified version for age~\eqref{eq: KKT condition for age},
provided that either $\Tftau_l^n \rightarrow \Tftau_l$ or
$\Tfa_l^n \rightarrow \Tfa_l$. For predicted class $l\in [K]$, $\flambda_l$ is the limiting arrival rate (see Definition~\ref{def: concerned processes for arrival and service}).
\begin{proposition}[Relationship between $\Tfa^n_l$ and $\Tftau^n_l$]
  \label{proposition: relation between a and tau}
  Given a classifier $f_\theta$ and a sequence of queueing systems, suppose that 
  Assumptions~\ref{assumption: data generating process},
  ~\ref{assumption: heavy traffic}, and~\ref{assumption:
    second order moments} hold. Under p-FCFS feasible policies, for any
  predicted class $l\in [K]$, (i)
  $\flambda_l \Tfa_l^n - \TfN_l^n \rightarrow{0}$, (ii)
  $\limsup_n \|\Tfa_l^n \| < \infty$; (iii) $\{\Tfa_l^n\}_n$ converges
  iff $\{\Tftau_l^n\}_n$ converges; (iv) their limits coincide: if there exist
  $\Tfa_l, \Tftau_l\in \mathcal{C}$ such that $\Tfa_l^n \rightarrow \Tfa_l$
  and $\Tftau_l^n \rightarrow \Tftau_l$, then $\Tfa_l = \Tftau_l$.
\end{proposition}
\noindent The proof of Proposition~\ref{proposition: relation between a and tau} requires the arrival rates of the predicted classes to converge to $\{\flambda_l \}_{l\in [K]}$, which is implied by Assumption~\ref{assumption: heavy traffic}. We also characterize the relationship between
$\Tfa_l^n$ and the policy process $\TfT_l^n$ in Corollary~\ref{corr: relation
  between age and T}, which allows for directly analyzing the dynamics of the scaled age
process $\Tfa_l^n$.

\paragraph{Convergence of the P$c\mu$ indices and the scaled age processes}
Since the \ourmethod~ serves the job that has the highest index value, the gap
between the class indices becomes small and the convergence~\eqref{eq: KKT
  condition for age} holds.
\begin{proposition}[Convergence of max difference of the P$c\mu$ indices]\label{prop: convergence of max gcmu difference}
  Given a classifier $f_\theta$, suppose that Assumptions~\ref{assumption: data
    generating process},~\ref{assumption: heavy traffic},~\ref{assumption: second order moments},
    ~\ref{assumption: on cost functions for showing the lower bound}, and~\ref{assumption: on cost
    functions for optimality} hold. Under the \ourmethod, 
\begin{equation}\label{eq: max gcmu index prelimit}
   \sup\limits_{t\in [0,1]}   \, \max\limits_{l, m \in [K]} \big| \fIpre^n_l (t) - \fIpre^n_m (t) \big|\rightarrow{0}. 
\end{equation}
\end{proposition}
\noindent By the continuity of the inverse cost function $(\fC_l')^{-1}$,
convergence of $\{\Tfa_l^n\}_{l\in [K]}$ follows (Lemma~\ref{lemma:
  convergence of age}) and we will have the desired final result.

We prove Proposition~\ref{prop: convergence of max gcmu difference} in
Section~\ref{section: proof of convergence of max gcmu difference}.
Specifically, we partition $[0, 1]$, the domain of the diffusion-scaled 
processes, into intervals of size $O(n^{-1/2})$ and show that
$\max_{l, m \in [K]}| \fIpre^n_l (t) - \fIpre^n_m (t) |$ do not exhibit
substantial growth within each interval given the carefully chosen size of the interval. The main technical challenge is to demonstrate that such 
growth does not accumulate over time. Since
$\max_{l, m \in [K]}| \fIpre^n_l (0) - \fIpre^n_m (0) | = 0$, we proceed via
induction over the intervals within the partition: for a fixed $\varepsilon > 0$, we show that 
\begin{enumerate}[(i)]
    \item at each endpoint $t$ of every interval, $\max_{l, m \in [K]}| \fIpre^n_l (t) - \fIpre^n_m (t) | \leq \varepsilon$ (Propositions~\ref{prop: max gcmu case i} and~\ref{prop: max gcmu case ii});
    \item within each interval $I$, $\sup_{t\in I} \max_{l, m \in [K]}| \fIpre^n_l (t) - \fIpre^n_m (t) | \leq 3\varepsilon/2$ (Proposition~\ref{prop: gcm difference within intervals}).
\end{enumerate}

We outline the proof for part (i) (part (ii) can be shown similarly).  
Given an interval $[t_1, t_2]$ and an induction hypothesis assumed for $t_1$, by symmetry it suffices to show
$\fIpre^n_{l}(t_2) - \fIpre^n_{m}(t_2) \leq \varepsilon$
for any $l,m\in [K]$.
First, for the case that predicted class $m$ is selected 
by the \ourmethod~at some time $ns\in [nt_1, nt_2]$, 
we use definition of the \ourmethod~to bound growth of $\fIpre^n_{l}(t_2) - \fIpre^n_{m}(t_2)$. In particular, 
\begin{equation}\label{eq: decomposition of max pcmu diff}
    \fIpre^n_{l}(t_2) - \fIpre^n_{m}(t_2) 
\leq    \underbrace{[\fIpre^n_{l}(t_2) - \fIpre^n_{l}(s)]}_{\text{bounded increase, }\leq \varepsilon/2}
    + \underbrace{[\fIpre^n_{l}(s) - \fIpre^n_{m}(s)]}_{\text{by the \ourmethod, $\leq 0$}}
    + \underbrace{[\fIpre^n_{m}(s) - \fIpre^n_{m}(t_2)]}_{\text{bounded increase, }\leq \varepsilon/2},
\end{equation}
where the first and the last term are bounded by $\varepsilon / 2$ due to
our choice of $t_2 - t_1 = O(n^{-1/2})$ and the smoothness of the cost
functions in Assumption~\ref{assumption: on cost functions for showing the
lower bound}, and the second term is non-positive since predicted class $m$ is chosen by the \ourmethod~at time $ns$.

For the other case that predicted class $m$ is never selected
by the \ourmethod~during the interval $[nt_1, nt_2]$, the analysis is
more involved and requires development of novel analysis techniques. 
If the server is idling at some  $ns \in [nt_1, nt_2]$, our analysis
is similar to~\eqref{eq: decomposition of max pcmu diff} and the second term becomes zero since the \ourmethod~is work-conserving. Otherwise, if there is no idling during $[nt_1, nt_2]$, then intuitively, the server is busy serving other $K-1$ predicted classes. 
By heavy traffic assumption $\sum_l \frho_l = 1 $
(Assumption~\ref{assumption: heavy traffic}), there exists at least one 
predicted class $k_0^n$ that receives sufficient service from the server
(See \eqref{eq: additional service in T}) and incurs sufficient descent 
in the age process (Lemma~\ref{lemma: age descent}) in $[nt_1, nt_2]$. 
Then, by strong convexity of the \ourpolicy~cost $\fC_{k_0^n}$ (Assumption~\ref{assumption: on cost functions for optimality}), the \ourpolicy~index of class $k^n_0$, say $\fIpre^n_{k^n_0}$, 
also incurs sufficient descent (Lemma~\ref{lemma: gcmu descent}). Such 
descent in the \ourpolicy~index enables us to bound the growth of $\fIpre^n_{l}$ 
and derive the desired result in Proposition~\ref{prop: max gcmu case ii}.

\subsection{Comparison to the optimality result in~\citet{VanMieghem95}}
\label{subsection: comparison to the optimality result in VM}

Plugging $Q^n= I$, our proof gives the optimality of the well-known \gcmu~where true class labels are known~\citep{VanMieghem95}. In this special case,
our analysis identifies missing arguments in~\citet{VanMieghem95}'s original
proof and provides conditions under which their original claims hold.  For
example, the derivatives of the first-order approximations for arrival time partial sum processes, $\{ (\flambda^n_l)^{-1} e\}_{l\in [K]}$, converge to $\{(\flambda_l)^{-1} e\}_{l\in [K]}$ in our framework. We find the same condition on the analogous processes, $\{\bar{U}_k^n\}_{k\in [K]}$ in~\citet{VanMieghem95}, to be one of the two alternative conditions that suffice to correctly connect the age and sojourn time processes in the G/G/1 systems. We refer the reader to Section~\ref{section: adding to the exposition in classical heavy traffic queueing analysis} for detailed explanations.

The first missing piece is that the \gcmu~uses the ages of waiting jobs for
scheduling, but~\citet{VanMieghem95} does not prove the \gcmu~achieves
optimality conditions defined in terms of the sojourn times~\cite[Eq
(54)]{VanMieghem95}. In Section~\ref{section: adding to the exposition in classical heavy traffic queueing analysis},  using our proof for Proposition~\ref{proposition: relation between a and tau}, we provide sufficient conditions to establish the equivalence between convergence of the scaled age and sojourn time processes in that paper, and complement the proof for the optimality of the \gcmu.

To prove the optimality of the \gcmu~rigorously, our analysis also shows that the age processes converge to a limit satisfying the optimality condition under the \gcmu~using
Proposition~\ref{prop: convergence of max gcmu difference} with $Q^n=I$, and thus
condition~\eqref{eq: KKT condition for tau} (extension of~\citet[Eq
(54)]{VanMieghem95}) is satisfied. This missing justification was nontrivial
(to us), and we hope our rigorous arguments provide analytical value to
subsequent works.

We found the proof of Proposition~\ref{prop: convergence of max gcmu
  difference} to be novel and involved. Our analysis of the index dynamics with the
particular choice of the partition size of the time horizon entails carefully handling errors of
diffusion approximations for predicted classes (Proposition~\ref{prop: convergence and approximation of predicted class N, tau, T, and W} and~\ref{prop: joint conv. of predicted class A, U, S, and V}). We control the
evolution of $\{\Tfa_l^n\}_{l\in [K]}$ under the \ourmethod, which requires
formally establishing relationships between $\Tfa_l^n$, $\Tftau_l^n$, and
$\TfT_l^n$.

In particular, our proof of Proposition~\ref{prop: convergence of max gcmu
  difference} identifies a previously unstated sufficient condition: strong convexity of the cost functions in Assumption~\ref{assumption:
  on cost functions for optimality}. The curvature ensures that if some
predicted class is not served and its index increases in a subinterval of the
partition, then there is another predicted class $k^n_0$ that receives ample
service so that the index $\fIpre^n_{k^n_0}$ decreases enough
(Lemma~\ref{lemma: gcmu descent}), implying that
the gap between the indices remains small. On the other hand, under the
\textit{strict} convexity~\citet{VanMieghem95} assumes, we were unable to show
the desired convergence he claims (either \cite[Eq (54)]{VanMieghem95} or a
more general version in Proposition~\ref{prop: convergence of max gcmu
  difference}).

\subsection{Comparison to the optimality result in~\citet{MandelbaumSt04}}
\label{subsection: comparison to the optimality result in MS}

Our analysis with perfect classification ($Q^n=I$) also identifies missing pieces in the optimality proof by~\citet{MandelbaumSt04} for the \gcmu~with sojourn time cost (called $D\text{-}Gc\mu$ in~\citep{MandelbaumSt04}) in single-server systems and provides conditions for the claims to hold.

First, similarly to~\citep{VanMieghem95},~\citet{MandelbaumSt04} suggests 
using age processes for the $D\text{-}Gc\mu$ rule, but they did not prove that their 
$D\text{-}Gc\mu$ rule satisfies optimality conditions they adopted, which are based 
on sojourn time processes and identical to~\eqref{eq: KKT condition for tau} in 
the single-server case. 
Specifically, we find that~\cite[Eq (66)]{MandelbaumSt04} that connects $D\text{-}Gc\mu$ to the preceding analysis in~\citep{MandelbaumSt04} should have been shown in terms of the queue length and age processes similarly to our Proposition~\ref{proposition: relation between a and tau} (i). Using an equivalence between the age and sojourn time processes analogous to Proposition~\ref{proposition: relation between a and tau} (iii) and (iv), the optimality of $D\text{-}Gc\mu$ could be obtained. 


Our analysis shows the optimality of the $D\text{-}Gc\mu$ in the single-server case 
requires weaker assumptions on cost functions than those adopted 
in~\citet{MandelbaumSt04}. The optimality in~\cite[Theorem 2]{MandelbaumSt04} is built on the attraction propery of the fluid-scaled queue length limit~\cite[Theorem 3]{MandelbaumSt04}. In the single-server case, the key implications of the attraction property are the small gaps between the class indices over subintervals~\cite[Eqs. (55), (56)]{MandelbaumSt04}, which are analogous to our Propositions~\ref{prop: max gcmu case i},~\ref{prop: max gcmu case ii}, and~\ref{prop: gcm difference within intervals}.~\citet[Theorem 3]{MandelbaumSt04} requires the cost functions to be twice continuously differentiable (and strongly convex) in order for the workload and queue length limits to be amenable to analysis in the multi-server setting. In contrast, our analysis directly identifies the dynamics of the age processes under the $D\text{-}Gc\mu$ in the single-server case, and prove the counterpart propositions under the weaker conditions, namely Assumptions~\ref{assumption: on cost functions for showing the lower bound} and~\ref{assumption: on cost functions for optimality}.

\subsection{Detailed proof of Theorem~\ref{theorem: optimality of our policy}}

We begin by showing the convergence of the age process, whose proof we give in
Section~\ref{section:proof-convergence of age}. We write
$\VTfa^n:=\{\Tfa_l^n\}_{l\in [K]}$ for the vector-valued diffusion-scaled age process.
\begin{lemma}[Convergence of $\Tfa^n$]
  \label{lemma: convergence of age}
  Given a classifier $\model$, suppose that Assumptions~\ref{assumption: data
    generating process},~\ref{assumption: heavy traffic},~\ref{assumption: second order moments},~\ref{assumption: on cost functions for showing the lower bound}, and~\ref{assumption: on cost
    functions for optimality} hold. Under the \ourmethod, there exists $\VTfa\in \mathcal{C}^K$ such that
  $\VTfa^n\rightarrow \VTfa$ in $(\mathcal{D}^K, \| \cdot \|)$
  $\mathbb{P}_\text{copy}$-a.s..
\end{lemma}
\noindent From the above lemma and Proposition~\ref{prop: convergence and
  approximation of predicted class N, tau, T, and W}, the relation between
$\Tfa^n_l$ and $\Tftau_l^n$ in Proposition~\ref{proposition: relation between
  a and tau} implies convergence of $\Tftau^n_l \rightarrow \Tftau_l$,
$\Tfa^n_l \rightarrow \Tfa_l$, $\TfN^n_l \rightarrow \TfN_l$,
$\TfT^n_l \rightarrow \TfT_l$, and $\TfW^n_l \rightarrow \TfW_l$.

If $\Tftau_l, \Tfa_l \in \mathcal{C}$, Proposition~\ref{proposition: relation
  between a and tau} implies $\Tftau_l = \Tfa_l$. Since
$\Tfa_l\in \mathcal{C}$ by Lemma~\ref{lemma: convergence of age}, it is easy
to verify $\TfN_l, \TfT_l, \TfW_l \in \mathcal{C}$ from the relation between
$\Tfa_l$ and $\TfN_l$ in Proposition~\ref{proposition: relation between a and
  tau}, relation~\eqref{eq: approximation of predicted N} between $\TfN_l$ and
$\TfT_l$, and relation~\eqref{eq: W in terms of L and T} between $\TfT_l$ and
$\TfW_l$. Using the relation~\eqref{eq: relation between W, tau, T} between
$\Tftau_l$, $\TfW_l$, and $\TfT_l$, one can check that
$\Tftau_l\in \mathcal{C}$.

By Lemma~\ref{lemma: relation between N and W} and Proposition~\ref{prop:
  little's law}, we have $\Tftau_l = \TfW_l / \frho_l,~\forall~l\in[K]$.
Proposition~\ref{prop: convergence of max gcmu difference} then implies
\begin{equation} \label{eq: chracterization of the sojourn time under pcmu}
  \Tftau_l = \TfW_l / \frho_l,~\forall~l\in[K],
  \quad \sum_{l\in [K]} \TfW_l = \sumTtW,
  \quad \fmu_l \fC_l'(\Tftau_l) = \fmu_m \fC_m'(\Tftau_m),~\forall~l,m\in[K].
\end{equation}
By Lemma~\ref{prop: properties of h and opt}, it follows
$\frho_l\Tftau_l = [h(\sumTtW)]_l,~\forall~l\in [K]$. This yields 
$\TJ^n_\ourpolicy (\cdot;Q^n) \rightarrow \TJ^* (\cdot;Q)$ 
$\mathbb{P}_{\text{copy}}\text{-}a.s.$ 
according to Theorem~\ref{theorem: HT lower bound} and 
Lemma~\ref{lemma: convergence of tildeJ}.

The weak convergence on the original systems $\TJ^n_\ourpolicy (\cdot;Q^n) \Rightarrow \TJ^* (\cdot;Q)$ in 
$(\mathcal{D}, \| \cdot \|)$ follows from~\cite[Lemma 3.2, Lemma 3.7]{Kallenberg97}.
Moreover, for any $x\in\mathbb{R},~t\in [0,1],$ by reverse Fatou's lemma and the $\mathbb{P}_{\text{copy}}\text{-}a.s.$ convergence of $\TJ^n_\ourpolicy (\cdot;Q^n)$, we have
\begin{equation*}
\begin{aligned}
\lim\sup_n \mathbb{P}^n [\TJ^n_\ourpolicy(t; Q^n) > x]  
\leq& \mathbb{E}_{\mathbb{P}_\text{copy}}
[ \lim\sup_n \mathbb{I}\{\TJ^n_{\ourpolicy} (t;Q^n) > x \}]\\
=&  \mathbb{E}_{\mathbb{P}_\text{copy}} 
[ \mathbb{I}\{\TJ^* (t;Q^n) > x \}]\\
=& \mathbb{P}_{\text{copy}} [\TJ^* (t;Q) > x]. 
\end{aligned}
\end{equation*}
Combining this with $\lim\inf_n \mathbb{P}^n [\TJ^n_\ourpolicy(t; Q^n) > x]  
\geq \mathbb{P}_{\text{copy}} [\TJ^* (t;Q) > x]$ from Theorem~\ref{theorem: HT lower bound} gives the desired result:
$\mathbb{P}^n [\TJ^n_\ourpolicy (t;Q^n) > x] \rightarrow 
\mathbb{P}_{\text{copy}} [\TJ^* (t;Q) > x]$.

\subsection{Proof of Lemma~\ref{lemma: convergence of age}}
\label{section:proof-convergence of age}

By Proposition~\ref{prop: convergence of max gcmu difference} and
Lemma~\ref{lemma: gcmu pre limit and asymptotic limit}, the \ourmethod~gives
$\max_{l, s \in [K]}\| \fmu_l \fC_l'(\Tfa^n_l) - \fmu_s \fC_l'(\Tfa^n_s)\|
\rightarrow 0$.  Given $s\in [K]$, for any $l\in [K]$, since $\fmu_l > 0$ by
Assumption~\ref{assumption: data generating process} and Definition~\ref{def:
  concerned processes for arrival and service}, we have
$\fC_l'(\tilde{a}_l^n) - \frac{\fmu_s}{\fmu_l } \fC_s' (\Tfa_s^n) \rightarrow
0.$ Letting
$f_s(\cdot):=\sum_{l = 1}^K \frho_l \cdot (\fC_l')^{-1}
\big(\frac{\fmu_s}{\fmu_l} \fC_s' (\cdot) \big)$, note that
$\sum_{l=1}^K\frho_l \Tfa_l^n - \big(f_s \circ \Tfa_s^n\big) \rightarrow 0$
from continuity of $(\fC_l')^{-1}$ (Assumption~\ref{assumption: on cost
  functions for showing the lower bound}). Under p-FCFS feasible policies,
Lemma~\ref{lemma: relation between N and W} and Proposition~\ref{proposition:
  relation between a and tau} implies
$\TfW^n_l - \frho_l \Tfa^n_l \rightarrow 0$.  Applying Proposition~\ref{prop:
  convergence and approximation of predicted class N, tau, T, and W}, there
exists $\TtW_+ \in \mathcal{C}([0, 1], \mathbb{R})$ such that
$\sum_{l = 1}^K \frho_l \Tfa^n_l \rightarrow \sumTtW$.  Hence,
$f_s \circ \Tfa_s^n \rightarrow \sumTtW$.

Since $f_s$ is continuous and strictly increasing, $f_s^{-1}$ is
well-defined and also continuous. Conclude
$(f_s^{-1}, f_s \circ \Tfa_s^n) \rightarrow (f_s^{-1}, \sumTtW)$ in
$\mathcal{C}^2$ under the product topology induced by $\Vert\cdot\Vert$.  By
the continuity of composition (e.g.,~\cite[Theorem 13.2.1]{Whitt02}),
$\Tfa_s^n = f_s^{-1} \circ \big(f_s \circ \Tfa_s^n \big) \rightarrow \Tfa_s
: = f_s^{-1} \circ \sumTtW$ where $\Tfa_s\in \mathcal{C}$ by continuity of
$f_s^{-1}$ and $\sumTtW$. This completes our proof.

\subsection{Proof of Proposition~\ref{proposition: relation between a and
    tau} (Gap B2)}
\label{subsection: proof of the relationship betwen age and sojourn time}

We show the asymptotically linear relationship (i) between $\Tfa_l^n$ and
$\TfN_l^n$. Other results immediately follow from (i) and
Propositions~\ref{prop: convergence and approximation of predicted class N,
  tau, T, and W} and Proposition~\ref{prop: little's law}. We use a
reformulation of the age process.
\begin{claim}
\begin{equation}\label{eq: age_proof_eq1}
    \fa_l^n(nt) = nt - \fU_l^n( \fA_l^n(nt) - \fN_l^n(nt) + 1) + o_n(n^{1/2}).
\end{equation}
\end{claim}
\noindent Since
$ \fU^n_l(nt) = n (\flambda_l^n)^{-1} t + n^{1/2} \TfU^n_l(t) + o_n(n^{1/2})$ by
Proposition~\ref{prop: joint conv. of predicted class A, U, S, and V}, we can
further rewrite~\eqref{eq: age_proof_eq1} as
\begin{equation*}
\begin{aligned}
    \Tfa_l^n(t) 
    = n^{1/2} [t - (\flambda^n_l)^{-1} (n^{-1} (\fA_l^n(nt) - \fN_l^n(nt) + 1))] - \TfU^n_l(n^{-1} (\fA_l^n(nt) - \fN_l^n(nt) + 1)) + o_n(1).
\end{aligned}
\end{equation*}
Recall $\fA^n_l(nt) = n \flambda^n_l t + n^{1/2} \TfA^n_l(t) + o_n(n^{1/2})$ by
Proposition~\ref{prop: joint conv. of predicted class A, U, S, and V},
$\TfN_{kl}^n:=n^{-\frac{1}{2}}\fN_{kl}^n$ by Definition~\ref{definition:
  concerned process}, and
$\limsup_n\|\TfN_{kl}^n\| \leq \limsup_n\|\TfN_{l}^n\| < +\infty$ by
Proposition~\ref{prop: convergence and approximation of predicted class N,
  tau, T, and W}.  Evidently,
\begin{equation*}
\begin{aligned}
    n^{-1} (\fA_l^n(nt) - \fN_l^n(nt) + 1)
    =&~ \flambda^n_l t + n^{-1/2}\TfA^n_l(t) 
    - n^{-1/2}\TfN^n_l(t) + o_n(n^{-1/2})\\
    \stackrel{(a)}{=}&~ \flambda_l t + n^{-1/2}\TfA_l(t) 
    - n^{-1/2}\TfN^n_l(t) + o_n(n^{-1/2})
    =\flambda_l t + o_n(1), \\
      \TfU^n_l(n^{-1} (\fA_l^n(nt) - \fN_l^n(nt) + 1))  
\stackrel{(b)}{=} &~        \TfU_l(\flambda_l t + o_n(1))  + o_n(1)\\
\stackrel{(c)}{=} &~  - (\flambda_l)^{-1} \TfA_l(t + o_n (1)) + o_n(1)
\stackrel{(d)}{=} - (\flambda_l)^{-1} \TfA_l(t) + o_n(1)
\end{aligned}
\end{equation*}
where we used $n^{1/2} |\flambda_l - \flambda_l^n |= o_n(1)$ by Assumption~\ref{assumption: heavy traffic} and Proposition~\ref{prop: joint
  conv. of predicted class A, U, S, and V} in step (a), $\TfU^n_l \rightarrow \TfU_l$ in step (b) and $\TfU_l(t) = - (\flambda_l)^{-1} \TfA_l(\flambda_l^{-1} t)$ in step (c) by the proof of Proposition~\ref{prop: joint conv. of predicted class A, U, S, and V}, and the uniform continuity of $\TfA_l$ on compact intervals in step (d).
Using Proposition~\ref{prop: joint conv. of predicted class A, U, S, and V}, as well as the right-hand side of step (a), we have
\begin{equation*}
\begin{aligned}
     & n^{1/2} [t - (\flambda^n_l)^{-1} (n^{-1} (\fA_l^n(nt) - \fN_l^n(nt) + 1))]
\\=&~     n^{1/2}t 
        - n^{1/2} (\flambda_l^n)^{-1} \big(\flambda_l t + n^{-1/2} \TfA_l^n (t) - n^{-1/2} \TfN_l^n (t) + o(n^{-1/2})\big)\\
=&~     - (\flambda_l)^{-1}[\TfA_l(t) - \TfN^n_l(t)]  + o_n(1).
\end{aligned}
\end{equation*}  

\paragraph{Proof of claim} For fixed $t\in[0, 1]$, we first consider the case that $\fA^n_l(nt) =
0$. Since there is no arrival to the predicted class $l$ at time $nt$, it is
easy to verify that $\fa^n_l(nt)=0$, $\fN_l^n(nt) = 0$, and
$nt \leq \fu_{l1}$. Therefore, we obtain that
\begin{equation*}
    \begin{aligned}
    \Big|\fa_l^n(nt) -[nt - \fU_l^n( \fA_l^n(nt) - \fN_l^n(nt) + 1)]\Big|
    =|0 -[nt - \fU_l^n( 1)]|
    \leq&~ |\fu_{l1}| = o_n(n^{1/2}),
\end{aligned}
\end{equation*}
where the last equality follows from Propositon~\ref{prop: joint conv. of
  predicted class A, U, S, and V}.  When $\fA_l^n(nt) \geq 1$,
$\fA_l^n(nt) - \fN_l^n(nt)$ jobs from the predicted class $l$ have completed
service and exited the queue. Under a p-FCFS policy, the oldest customer from
the predicted class $l$ at time $nt$ corresponds to the
$[\fA_l^n(nt) - \fN_l^n(nt) + 1]$th arrival of predicted class $l$. From the
definition of $\fa^n_l(nt)$ as the time difference between $nt$ and the
arrival time of the oldest job in predicted class $l$, the exact formulation
$\fa^n_l(nt) = nt - \fU_l^n( \fA_l^n(nt) - \fN_l^n(nt) + 1)$ follows. This
completes our proof of~\eqref{eq: age_proof_eq1}.


\section{Proof for convergence of the P$c\mu$ indices (Proposition~\ref{prop: convergence of max gcmu difference}, Gap B1)}
\label{section: proof of convergence of max gcmu difference}

Recalling the strong convexity of $\fC_l$, 
\begin{align*}
  \fC_l(y) \geq \fC_l (x) +  \fC_l' (x) (y-x) + \frac{m}{2}(y-x)^2, ~\forall~x, y,~\forall~l\in [K].
\end{align*}
for some $m>0$, we use the following constants
\begin{equation}
  \label{eq: special constants}
\begin{aligned}
&   \fmu_{\min} = \min_{l\in [K]} \fmu_l,\quad 
    \fmu_{\max} = \max_{l\in[K]} \fmu_l, \quad 
    \frho_{\min}= \min_{l\in[K]} \frho_l,\quad
    \frho_{\max}=\max_{l\in[K]} \frho_l,\\
&   \alpha_0 : = \frac{\frho_{\min}}{3(K-1)\frho_{\max}}, \quad 
    \beta_0 : = \frac{\fmu_{\min} \alpha_0}{2},
    \quad \gamma_0 : =  \frac{\frho_{\min}\alpha_0}{1-\frho_{\min}}.
\end{aligned}
\end{equation}
Since $\fC_l'$ is uniformly continuous on the compact set
$[0, \limsup_n \|\Tfa_l^n\|]$ where $ \limsup_n \|\Tfa_l^n\| < +\infty$
according to Proposition~\ref{proposition: relation between a and tau}, we
have the following result.
\begin{lemma}[Continuity of $\fC_l'$]\label{lemma: continuity of C}
  Given a classifier $f_\theta$, suppose that Assumption~\ref{assumption: on cost
    functions for showing the lower bound} holds. For any $\varepsilon>0$,
  there exists $\delta_1(\varepsilon), \delta_2(\varepsilon)>0$ such that for
  any $a_1, a_2 \in [0, \limsup_n \|\Tfa_l^n\|]$,
  \begin{enumerate}[(i)]
  \item if $|a_2 - a_1| \leq \delta_1(\varepsilon)$, then
    $|\fC_l'(a_2) - \fC_l' (a_1)| < \frac{\varepsilon}{8
      \fmu_{\max}},~\forall~l\in[K]$;
  \item if $|a_2 - a_1| \leq \delta_2(\varepsilon)$, then $|\fC_l'(a_2) - \fC_l' (a_1)| < \frac{m\beta_0}{2\fmu_{\max}}\delta_1(\varepsilon),~\forall~l\in[K]$.
  \end{enumerate}
\end{lemma}

Our proof is separated into three propositions. Below, we suppose
Assumptions~\ref{assumption: data generating process},~\ref{assumption: heavy
  traffic},~\ref{assumption: on cost functions for showing the lower
  bound},~\ref{assumption: on cost functions for optimality},
and~\ref{assumption: second order moments} hold.  For any $\varepsilon>0$, let
$\delta_1(\varepsilon), \delta_2(\varepsilon)$ be constants defined in
Lemma~\ref{lemma: continuity of C} and define
$\delta_1^n (\varepsilon) \defeq n^{-1/2}\delta_1 (\varepsilon), \delta_2^n
(\varepsilon) \defeq n^{-1/2}\delta_2 (\varepsilon)$.  Partition the time
interval $[0,n]$ into subintervals of length no more than
$n\delta_1^n (\varepsilon)$. Let $N(\varepsilon)$ be large enough that the three
propositions below hold for $n \ge N(\varepsilon)$. We prove the proposition by
an inductive argument over these subintervals, given after the three
propositions are stated.

See Section~\ref{section:proof-max gcmu case i} for the proof of the first proposition.
\begin{proposition}[Max difference of the P$c\mu$ indices at endpoints: Case I]
  \label{prop: max gcmu case i}
  Let $t_1 \in [0,1- \delta_1^n(\varepsilon)]$ be such that
  $ \max_{l, m \in [K]}| \fIpre_l^n (t_1) - \fIpre_m^n (t_1)| < \varepsilon$
  and all predicted classes are selected by the \ourmethod~in
  $[nt_1, n(t_1 + \delta_1^n(\varepsilon))]$. Then, there exists
  $N(\varepsilon) > 0$ such that for any $n>N(\varepsilon)$
  \begin{align*}
    \max_{l_1, l_2\in [K]}
    |\fIpre_{l_1}^n (t_1+ \delta_1^n(\varepsilon))
    - \fIpre_{l_2}^n (t_1 + \delta_1^n(\varepsilon))|  < \varepsilon.
  \end{align*}
\end{proposition}

We prove the second proposition in Section~\ref{section:proof-max gcmu case
  ii}.
\begin{proposition}[Max difference of the \ourmethod~indices at endpoints: Case
  II]
  \label{prop: max gcmu case ii}
  Let $t_1 \in [0,1- \delta_1^n(\varepsilon)]$ be such that
  $ \max_{l, m \in [K]}| \fIpre_l^n (t_1) - \fIpre_m^n (t_1)| < \varepsilon$
  and some predicted class is NOT selected for service under the \ourmethod~in
  $[nt_1, n(t_1 + \delta^n_1(\varepsilon))]$. Then, there exists
  $N(\varepsilon) > 0$ such that for any $n>N(\varepsilon)$
  \begin{enumerate}[(i)]
  \item (No Idling) if there is no server idle time in $[nt_1, n(t_1+\delta
    _1^n(\varepsilon))]$, then there exists
    $s^n_1 \in [t_1 + \gamma_0 \delta_1^n(\varepsilon),
    t_1 + \delta_1^n(\varepsilon)]$ such that 
    $\max_{l_1, l_2\in [K]}|\fIpre_{l_1}^n (s^n_1) - \fIpre_{l_2}^n (s^n_1)|
    < \varepsilon$;
  \item (Idling) if server idling occurs in $[nt_1, n(t_1+\delta_1^n(\varepsilon))]$,
    then $ \max_{l_1, l_2\in [K]} |\fIpre_{l_1}^n (t_1+ \delta_1^n(\varepsilon))
    - \fIpre_{l_2}^n (t_1 + \delta_1^n(\varepsilon))|  < \varepsilon.$
  \end{enumerate}
\end{proposition}

Finally, see Section~\ref{section:proof-gcm difference within intervals} for
the proof of the third proposition.
\begin{proposition} [Max difference of the \ourmethod~indices within intervals]
  \label{prop: gcm difference within intervals}
  Let $t_1 \in [0,1]$ be such that
  $\max_{l, m \in[K]}\big| \fIpre_l^n (t_1) -\fIpre_m^n (t_1) \big| <
  \varepsilon$.  Then, there exists $N(\varepsilon) > 0$ such that for any
  $n>N(\varepsilon)$
  \begin{align*}
    \max_{l_1, l_2\in [K]} \,
    \sup_{t\in [t_1,(t_1 + \delta_1^n(\varepsilon))\wedge 1 ]}
    \big| \fIpre_{l_1}^n (t) -\fIpre_{l_2}^n (t) \big| < 3\varepsilon/2.
  \end{align*}
\end{proposition}

\paragraph{Proof of Proposition~\ref{prop: convergence of max gcmu difference}.}
Fix $\varepsilon>0$ and $n\geq N(\varepsilon)$. Call a point $s\in[0,1]$
\emph{$\varepsilon$-good} if
$\max_{l,m\in[K]}|\fIpre_l^n(s)-\fIpre_m^n(s)|<\varepsilon$. We construct an
increasing sequence of $\varepsilon$-good points $0=s_0<s_1<\cdots<s_J=1$ whose
consecutive gaps satisfy
\begin{equation}\label{eq: good point spacing}
  \gamma_0\delta_1^n(\varepsilon) \leq s_{j+1}-s_j \leq \delta_1^n(\varepsilon),
  \qquad j=0,\dots,J-1.
\end{equation}

\noindent\emph{Base case.} At $t=0$ all ages are zero, so all indices coincide
and $s_0=0$ is $\varepsilon$-good.

\noindent\emph{Recursive step.} Suppose $s_j<1$ is $\varepsilon$-good. We obtain
the next $\varepsilon$-good point by applying the relevant proposition on
$[s_j,(s_j+\delta_1^n(\varepsilon))\wedge1]$.
\begin{enumerate}[(a)]
\item If all predicted classes are served on this interval, then
  Proposition~\ref{prop: max gcmu case i} makes
  $s_{j+1}:=(s_j+\delta_1^n(\varepsilon))\wedge1$ $\varepsilon$-good.
\item If some predicted class is not served and either idling occurs or
  $s_j+\delta_1^n(\varepsilon)\geq1$, then
  Proposition~\ref{prop: max gcmu case ii}(ii) makes the same
  $s_{j+1}:=(s_j+\delta_1^n(\varepsilon))\wedge1$ $\varepsilon$-good.
\item Otherwise some predicted class is not served and no idling occurs, and
  Proposition~\ref{prop: max gcmu case ii}(i) provides an $\varepsilon$-good
  point
  $s_{j+1}\in[s_j+\gamma_0\delta_1^n(\varepsilon),\,s_j+\delta_1^n(\varepsilon)]$.
\end{enumerate}
In every case the maximal index gap is restored to below $\varepsilon$ at
$s_{j+1}$, so the error does not accumulate across steps, and~\eqref{eq: good
  point spacing} holds because $\gamma_0\in(0,1)$. Each step advances time by at
least $\gamma_0\delta_1^n(\varepsilon)=\gamma_0 n^{-1/2}\delta_1(\varepsilon)>0$,
so the recursion reaches $s_J=1$ after at most
$\lceil n^{1/2}/(\gamma_0\delta_1(\varepsilon))\rceil$ steps.

\noindent\emph{From the good points to the uniform bound.} Each
$\varepsilon$-good $s_j$ meets the hypothesis of
Proposition~\ref{prop: gcm difference within intervals}, which gives
\begin{equation*}
  \sup_{t\in[s_j,(s_j+\delta_1^n(\varepsilon))\wedge1]}\,
  \max_{l,m\in[K]}|\fIpre_l^n(t)-\fIpre_m^n(t)|<3\varepsilon/2.
\end{equation*}
Since $s_{j+1}\leq s_j+\delta_1^n(\varepsilon)$, the intervals
$[s_j,(s_j+\delta_1^n(\varepsilon))\wedge1]$ cover $[0,1]$, hence
$\sup_{t\in[0,1]}\max_{l,m\in[K]}|\fIpre_l^n(t)-\fIpre_m^n(t)|<3\varepsilon/2$ for
all $n\geq N(\varepsilon)$. As $\varepsilon>0$ is arbitrary,~\eqref{eq: max gcmu
  index prelimit} follows.

\subsection{Preliminaries}
\label{subsec: intermediate results}

\paragraph{Facts about limiting diffusion processes} We use the following
basic facts to analyze the dynamics of $\VTfa^n$.
\begin{lemma}[Continuity of $\TfA_l$ and $\TfS_l$]\label{lemma: continuity of
    tilde A and S}
  There exists $N(\varepsilon)$ such that for $n> N(\varepsilon)$ and
  $t_1, t_2\in [0, 1]$,
  \begin{enumerate}[(i)]
  \item  if $|t_2 - t_1| < \delta^n_1(\varepsilon)$, then
    $|\TfA_l(t_2) - \TfA_l(t_1) |
    < \alpha_0 \delta_1(\varepsilon)/3,~\forall~l\in[K]$;
  \item if $|t_2 - t_1| < \delta^n_1(\varepsilon)$, then
    $|\TfS_l(n^{-1}\fT^n_l(nt_2)) - \TfS_l (n^{-1}\fT^n_l(nt_1)) |
    <\alpha_0\delta_1(\varepsilon)/3,~\forall~l\in[K]$.
  \end{enumerate}
\end{lemma}
\begin{proof}
  By Proposition~\ref{prop: joint conv. of predicted class A, U, S, and V}, we
  have $\sup_{t\in[0, 1]}|\TfA_l(t + o_n(1)) - \TfA_l(t)| = o_n(1)$ by uniform
  continuity of $\TfA_l$ over a closed interval of which $[0,1]$ is a proper
  subset for all $l\in [K]$. (i) is a direct consequence of
  $|t_2 - t_1| = o_n(1)$. To see (ii), we have
  $\sup_{t\in[0, 1]}|\TfS_l(t + o_n(1)) - \TfS_l(t)| = o_n(1)$ similarly, and
  $\sup_{t\in[0, 1]}|n^{-1}\fT_l(n(t + o_n(1))) - n^{-1}\fT_l(nt)| = o_n(1)$
  by~\eqref{eq: approximation of T}.
\end{proof}

\begin{lemma}[Relation between $\Tfa^n_l$ and $\TfT^n_l$] \label{corr:
    relation between age and T} Given a classifier $\model$, suppose
  Assumptions~\ref{assumption: data generating process}, ~\ref{assumption:
    heavy traffic}, and~\ref{assumption: second order moments} hold. Under
  p-FCFS feasible policies,
  \begin{equation}
    \label{eq: relation between age and T}
    \Tfa_l^n (t) = 
    n^{1/2}  t -  n^{-1/2} \frho_l^{-1} \fT^n_l(nt) 
    + \flambda_l^{-1} \TfA_l(t) - \flambda_l^{-1}  \TfS_l (n^{-1}\fT^n_l(nt))
    + o_n(1).
  \end{equation}
\end{lemma}
\begin{proof}
  Recalling $\fA^n_l(nt) = n \BfA^n_l(t) + n^{1/2} \TfA^n_l(t) + o_n(n^{1/2})$,
  $\fS^n_l(nt) = n \BfS^n_l(t) + n^{1/2} \TfS^n_l(t) + o_n(n^{1/2})$
  (Proposition~\ref{prop: joint conv. of predicted class A, U, S, and V}),
  \begin{equation}\label{eq: relation between N and T eq2}
    \begin{aligned}
      \TfN^n_l(t) 
      =&~ n^{1/2} \BfA^n_l(t) + \TfA^n_l(t)
      - n^{1/2} \BfS^n_l(n^{-1}\fT^n_l(nt))
      - \TfS^n_l (n^{-1}\fT^n_l(nt)) +  o_n(1)\\
      =&~ n^{1/2} \BfA_l(t) + \TfA_l(t)
      - n^{1/2} \BfS_l(n^{-1}\fT^n_l(nt))
      - \TfS_l (n^{-1}\fT^n_l(nt))
      + o_n(1)\\
      =&~ n^{1/2} \flambda_l t -  n^{-1/2} \fmu_l \fT^n_l(nt) 
      + \TfA_l(t) - \TfS_l (n^{-1}\fT^n_l(nt))
      + o_n(1),
    \end{aligned}
  \end{equation}
  where we used $\fN^n_l(nt) = \fA^n_l(nt) - \fS^n_l(\fT^n_l(nt))$,
  $n^{1/2} (\BfA^n_l - \BfA_l) = o_n(1)$, $n^{1/2} (\BfS^n_l - \BfS_l) = o_n(1)$
  from Assumption~\ref{assumption: heavy traffic}, and boundedness of
  $n^{-1}\fT^n_l(n\cdot)$~\eqref{eq: approximation of T}.  Noting
  $\Tfa_l^n (t) = \flambda_l^{-1} \TfN_l^n (t)+ o_n(1)$ by
  Proposition~\ref{proposition: relation between a and tau}, we have the desired result.
\end{proof}

\paragraph{Asymptotic P$c\mu$ index}
For any predicted class $l\in[K]$ and $t\in [0, 1]$, the P$c\mu$ index and its
asymptotic counterpart is
\begin{equation}
  \label{eq: def gcmu index asymptotic limit}
  \fIpre^n_l(t):=
  \fmu^n_l\cdot n^{1/2} (\fC^n_l)'(\fa_l^n (nt)),
  ~~~~~~\fIlim_l^n (t):= \fmu_l \cdot \fC_l' ( \Tfa_l^n (t))
\end{equation}
Their difference can be bounded by
\begin{equation*}
  |\fIlim_l^n (t) - \fIpre_l^n (t) | 
  \leq  \fC_l'(\Tfa_l^n (t)) \cdot |\fmu^n_l - \fmu_l| 
  + \fmu^n_l \cdot |n^{1/2} (\fC^n_l)'(n^{1/2}\Tfa_l^n(t)) 
  - \fC_l' (\Tfa_l^n (t))|.
\end{equation*}
Note that $\limsup_n \fmu_l^n < +\infty$ from
$n^{1/2}(\fmu_l^n - \fmu_l) \rightarrow 0 $ (Assumption~\ref{assumption: heavy
  traffic}), $\limsup_{n} \|\fC'_l (\Tfa_l^n(\cdot)) \| < +\infty$ since
$\fC'_l$ is continuous, and $\limsup_n \| \Tfa^n_l \| < +\infty$ by
Proposition~\ref{proposition: relation between a and tau}.  Since
$n^{1/2}(\fC^n_l)'(n^{1/2}\cdot) \rightarrow \fC'_l$ by
Assumption~\ref{assumption: on cost functions for showing the lower bound}, we
can conclude $\sup_{t\in[0, 1]}|\fIlim_l^n (t) - \fIpre_l^n (t) | = o_n(1)$.
\begin{lemma}\label{lemma: gcmu pre limit and asymptotic limit}
  There exists $N(\varepsilon) > 0$ such that for any
  $n \geq N(\varepsilon)$,
  \begin{equation*}
    \max_{l\in[K]}\sup_{t\in [0,1]}\big|\fIlim_l^n (t) - \fIpre_l^n (t) \big| 
    \leq \min\Big\{\frac{\varepsilon}{16},  \frac{m\beta_0}{4} \delta_1(\varepsilon)\Big\}.
  \end{equation*}
\end{lemma}

\paragraph{Bounding the difference between P$c\mu$ indices}
When the difference of the scaled ages
$\{\Tfa_l^n\}_{l\in [K]}$ is bounded, we demonstrate bounded differences of
the indices over sufficiently small intervals.
\begin{lemma}[P$c\mu$ index: Continuity I]\label{lemma: continuity of gcmu index} 
  There exists $N(\varepsilon) > 0$ such that for any $n \geq N(\varepsilon)$,
  $l\in[K]$, and $0 \leq t_1 < t_2 \leq 1$,
  \begin{enumerate}[(i)]
  \item  if $\Tfa^n_l(t_2) - \Tfa^n_l(t_1) \leq \delta_1(\varepsilon)$, then
    $\fIpre_l^n(t_2) - \fIpre_l^n(t_1)  \leq \frac{\varepsilon}{4}$;
  \item if $\Tfa^n_l(t_2) - \Tfa^n_l(t_1) \geq 0$,
    then $\fIpre_l^n(t_2) - \fIpre_l^n(t_1) \geq \max\{-\frac{\varepsilon}{4}, -m \beta_0 \delta_1(\varepsilon)\}.$
  \end{enumerate}
\end{lemma}
\begin{proof}
  By Lemma~\ref{lemma: gcmu pre limit and asymptotic limit}, it suffices to
  show $\fIlim_l^n (t_2) - \fIlim_l^n (t_1) \leq \varepsilon/8$ for (i) and
  $ \fIlim_l^n(t_2) - \fIlim_l^n(t_1) \geq 0$ for (ii).  Noting $\fC_l'$
  is non-decreasing, we have
  \begin{equation*}
    \begin{aligned}
      \fIlim_l^n (t_2) - \fIlim_l^n (t_1) 
      =       \fmu_l  [\fC_l'(\Tfa_l^n (t_2)) 
      - \fC_l'(\Tfa_l^n (t_1))]
      \leq   \fmu_l  [\fC_l'(\Tfa_l^n (t_1) + \delta_1(\varepsilon)) 
      -  \fC_l'(\Tfa_l^n (t_1))]
    \end{aligned}
  \end{equation*}
  which yields (i) by continuity of $\fC'_l$ from Lemma~\ref{lemma: continuity
    of C}.  For (ii), non-decreasing $\fC'_l$ again implies
  $\fIlim_l^n(t_2) - \fIlim_l^n(t_1) = \fmu_l [\fC_l'(\Tfa_l^n (t_2)) -
  \fC_l'(\Tfa_l^n (t_1))] \geq 0.$ 
\end{proof}

Next, we bound the differences when the age process has a negative jump due to
a job's departure following service completion.
\begin{lemma}[Size of a negative jump of P$c\mu$ index]\label{lemma: size of gcmu's negative jump}
  There exists $N(\varepsilon) > 0$ such that for $n \geq N(\varepsilon)$
\begin{equation*}
        \sup_{t\in[0, 1]}|\fIpre_l^n (t^-) - \fIpre_l^n (t)|
\leq    \min\Big\{\frac{\varepsilon}{4}, m \beta_0 \delta_1(\varepsilon)\Big\}.
\end{equation*}
\end{lemma}

\begin{proof}
  By Lemma~\ref{lemma: gcmu pre limit and asymptotic limit}, it suffices to
  show
  $\fIlim_l^n (t^-) - \fIlim_l^n (t) \leq \min\{\frac{\varepsilon}{8}, \frac{m
    \beta_0}{2} \delta_1(\varepsilon)\}.$ $\Tfa^n_l (t) \neq \Tfa^n_l (t^-) $
  only arises when a job from the predicted class $l$ completes service and
  leaves the system at time $t$. By definition, the age process will incur a
  negative jump that corresponds to the interarrival time of two consecutive
  jobs. It follows from Proposition~\ref{prop: joint conv. of predicted class
    A, U, S, and V} that
  \begin{equation*}
    |\Tfa_l^n (t) - \Tfa_l^n(t^-)|  
    \leq  n^{-1/2}\sup_{1\leq i \leq \fA_l^n (n)} \fu_{li}^n 
    \leq \min\{\delta_1(\varepsilon), \delta_2(\varepsilon)\}.
  \end{equation*}
  for all sufficently large $n$. Combining the above and continuity of
  $\fC'_l$ from Lemma~\ref{lemma: continuity of C},
  \begin{equation*}
    |\fIlim_l^n (t) - \fIlim_l^n (t^-) |
    =\fmu_l  |\fC_l' ( \Tfa_l^n (t)) - \fC_l' (\Tfa_l^n (t^-))|
    \leq \min\Big\{ \frac{\varepsilon}{8}, \frac{m\beta_0}{2} \delta_1(\varepsilon)\Big\}.
  \end{equation*}
\end{proof}

\subsection{Proof of Proposition~\ref{prop: max gcmu case i}}
\label{section:proof-max gcmu case i}
Without loss of generality, we fix $\varepsilon > 0$, $n > N(\varepsilon)$,
and $t_1 \in [0, 1 - \delta^n_1(\varepsilon)]$. For simplicity, let
$t_2 = t_1 + \delta_1^n(\varepsilon)$. Choose any $l_1, l_2 \in[K]$. By
symmetry, it suffices to show that
$\fIpre_{l_1}^n (t_2) - \fIpre_{l_2}^n (t_2) < \varepsilon$.  Let $s^n_0$
denote the largest (scaled) time point in $[t_1, t_2]$ at which the predicted
class $l_2$ is selected by \ourmethod
\begin{equation*}
    s^n_0 : = \sup\left\{t \mid t \in [t_1, t_2],~\fIpre^n_{l_2}(t) = \max_{l\in [K]} \fIpre^n_l(t)\right\}.
\end{equation*}
We can obtain from the definition of \ourmethod~that
\begin{equation*}
  \fIpre_{l_1}^n (t_2) - \fIpre_{l_2}^n (t_2) 
  =   \underbrace{[\fIpre_{l_1}^n (t_2) 
    - \fIpre_{l_1}^n (s^n_0)]}_{\text{by Lemmas~\ref{lemma: continuity of gcmu index} and~\ref{lemma: size of gcmu's negative jump}, }\leq \varepsilon/2} 
  + \underbrace{[\fIpre_{l_1}^n (s^n_0) - \fIpre_{l_2}^n (s^n_0)]}_{\text{by \ourmethod, } \leq 0}  
  + \underbrace{[\fIpre_{l_2}^n (s^n_0) - \fIpre_{l_2}^n (t_2)]}_{\text{by Lemmas~\ref{lemma: continuity of gcmu index}, }\leq \varepsilon/2} .
   \end{equation*}

The second term satisfies
$\fIpre_{l_1}^n (s^n_0) - \fIpre_{l_2}^n (s^n_0)\leq 0$ since predicted class
$l_2$ is selected for service by \ourmethod~at time $s^n_0$. The other two
terms can be bounded by $\varepsilon / 2$ due to our selection of $t_2$ and
continuity of P$c\mu$ index, as show in Lemmas~\ref{lemma: continuity of gcmu
  index} and~\ref{lemma: size of gcmu's negative jump}. In particular, the
first term can be bounded by
\begin{equation*}
\begin{aligned}
        \fIpre^n_{l_1}(t_2) - \fIpre^n_{l_1}(s^n_0)
=       \underbrace{[ \fIpre^n_{l_1}(t_2) 
        - \fIpre^n_{l_1}((s^n_0)^-)]}_{\text{by Lemma~\ref{lemma: continuity of gcmu index}, } \leq \varepsilon/4}
        +  \underbrace{[ \fIpre^n_{l_1}((s^n_0)^-) 
        - \fIpre^n_{l_1}(s^n_0)]}_{\text{by Lemma~\ref{lemma: size of gcmu's negative jump}, }\leq \varepsilon/4}
\leq   \varepsilon/2,
\end{aligned}
\end{equation*}
since $\Tfa^n_{l_1}(t_2) - \Tfa^n_{l_1}(s^n_0) \leq n^{1/2}(t_2 - s^n_0) \leq \delta_1(\varepsilon)$. Similarly, the third term satisfies
\begin{equation*}
\begin{aligned}
    \fIpre^n_{l_2} (s^n_0) - \fIpre^n_{l_2} (t_2)
    \leq \varepsilon/2,
\end{aligned}
\end{equation*}
by Lemma~\ref{lemma: continuity of gcmu index}, since $l_2$ is not served
on the scaled interval $[s^n_0, t_2]$, and thus
$\Tfa^n_{l_2}(t_2) - \Tfa^n_{l_2}(s^n_0) \geq 0$.

\subsection{Proof of Proposition~\ref{prop: max gcmu case ii}}
\label{section:proof-max gcmu case ii}

Let $t_2 = t_1 + \delta_1^n(\varepsilon)$. By symmetry, it suffices to show
$\fIpre_{l_1}^n (t_2) - \fIpre_{l_2}^n (t_2) < \varepsilon$ for
$l_1, l_2 \in [K]$. First, consider the scenario (ii) where idling occurs in
$[nt_1, nt_2]$, i.e.,
$\sum_l \fT^n_l(nt_2) - \sum_l \fT^n_l(n t_1) < n \delta_1^n(\varepsilon)$.
Since we only consider work conserving policies, idling implies that there is
no job in queue at some time $ns^n_2 \in [nt_1, nt_2]$. Consequently, the age
of all predicted classes is zero $\Tfa^n_l( s^n_2) = 0, ~\forall~l \in
[K]$. Then,
\begin{equation*}
  \fIpre_{l_1}^n(t_2) -\fIpre_{l_2}^n(t_2)
  = \underbrace{[\fIpre_{l_1}^n(t_2) - \fIpre_{l_1}^n(s^n_2) ]}_{\text{by Lemma~\ref{lemma: continuity of gcmu index}, }\leq \varepsilon/2 }
  + \underbrace{[\fIpre_{l_1}^n(s^n_2) -\fIpre_{l_2}^n(s^n_2)]}_{\text{by definition, }= 0}
  + \underbrace{[\fIpre_{l_2}^n(s^n_2) - \fIpre_{l_2}^n(t_2)]}_{\text{by $\Tfa^n_{l_2}( s^n_2) = 0$, }\leq 0}
  \leq \varepsilon,
\end{equation*}
since $\fIpre_{l_2}^n(t_2) \geq 0$. 

The case (i) where no idling occurs is more complicated. We begin by showing
that the age and the P$c\mu$ index decrease sufficiently. See
Section~\ref{section:proof-age-descent} for the proof of the following result.
\begin{lemma} [Sufficient descent in age process]
  \label{lemma: age descent}
  For all $t_1 \in [0,1 - \delta_1^n(\varepsilon)]$, assume 
  \begin{enumerate}[(i)]
  \item (Non-Selected Class) at least one predicted class, say $l^n_0$, is not
    selected by \ourmethod~in time interval
    $[nt_1, n(t_1 +\delta^n_1(\varepsilon))]$;
  \item (No Idling)
    $\sum_l \fT^n_l(n (t_1 + \delta_1^n(\varepsilon))) - \sum_l \fT^n_l(n t_1)
    = n \delta_1^n(\varepsilon)$.
  \end{enumerate}
  There exists $N(\varepsilon)$ such that for all $n > N(\varepsilon)$, there
  is a predicted class $k^n_0$ whose age process decreases sufficiently:
    $\Tfa_{k^n_0}^n(t_1 + \delta^n_1(\varepsilon)) - \Tfa_{k^n_0}^n (t_1) 
    \leq - 2\alpha_0 \delta_1(\varepsilon)$.
\end{lemma}     
\noindent Let $k^n_0$ be the predicted class with
$\Tfa_{k^n_0}^n(t_2) - \Tfa_{k^n_0}^n (t_1) \leq - 2\alpha_0
\delta_1(\varepsilon).$ Let $s^n_1$ denote the smallest scaled time in
$[t_1, t_2]$ at which $\Tfa_{k^n_0}^n$ experience such decrease
\begin{equation*}
    s^n_1 : = \inf\{t \mid t \in [t_1, t_2],~\Tfa_{k^n_0}^n(t) - \Tfa_{k^n_0}^n (t_1) \leq - 2\alpha_0 \delta_1(\varepsilon)\}.
\end{equation*}

If predicted class $l_2$ is selected for service by \ourmethod~in
$[nt_1, ns^n_1]$, we can show
$\fIpre_{l_1}^n(s_1^n) - \fIpre_{l_2}^n(s^n_1) \leq \varepsilon$ by a similar
analysis as the proof of Proposition \ref{prop: max gcmu case i}. The crux of
our proof lies in the scenario where $l_2$ is not selected in
$[nt_1, ns^n_1]$.  At $s^n_1$, $\Tfa_{k^n_0}^n$ has a negative jump by a
service completion in predicted class $k^n_0$ and the P$c\mu$ index decreases
sufficiently.
\begin{lemma}[Sufficient descent in P$c\mu$ index]\label{lemma: gcmu descent}
  For any $0 \leq t_1 < t_2 \leq 1$, assume there exists some predicted class
  $k^n_0$ satisfying
    $\Tfa_{k^n_0}^n(t_2) - \Tfa_{k^n_0}^n (t_1) 
    \leq - 2\alpha_0 \delta_1(\varepsilon)$.
  There exists $N(\varepsilon) > 0$ such that for any $n \geq N(\varepsilon)$,
  the P$c\mu$ index for this predicted class decreases sufficiently 
\begin{equation*}
  \fIpre_{k^n_0}^n (t_2) 
  - \fIpre_{k^n_0}^n (t_1) \leq  - 3m\beta_0\delta_1(\varepsilon).
\end{equation*}
\end{lemma}
\begin{proof}
  By Lemma~\ref{lemma: gcmu pre limit and asymptotic limit}, it suffices to
  show
\begin{align*}
  \fIlim_{k^n_0}^n (t_2) - \fIlim_{k^n_0}^n (t_1)
  = \fmu_{k^n_0}
  \left(\fC'_{k^n_0}( \Tfa^n_{k^n_0} (t_2))
  - \fC'_{k^n_0}( \Tfa^n_{k^n_0} (t_1)) \right)
  \leq -4m \beta_0 \delta_1(\varepsilon).
\end{align*}
Since $\fC_{k^n_0}$ is strongly convex,
\begin{align*}
[\fC_{k^n_0}' (\Tfa_{k^n_0}^n (t_2)) - \fC_{k^n_0}' (\Tfa_{k^n_0}^n (t_1))] 
[ \Tfa_{k^n_0}^n (t_2) - \Tfa_{k^n_0}^n (t_1)] 
\geq m [\Tfa_{k^n_0}^n (t_2) - \Tfa_{k^n_0}^n (t_1)]^2.
\end{align*}
Then, $\Tfa_{k^n_0}^n(t_2) - \Tfa_{k^n_0}^n(t_1)\leq -2\alpha_0 \delta_1(\varepsilon)$ yields 
\begin{equation*}
    \fC_{k^n_0}' (\Tfa_{k^n_0}^n (t_2)) - \fC_{k^n_0}' (\Tfa_{k^n_0}^n (t_1))
    \leq m [\Tfa_{k^n_0}^n (t_2) - \Tfa_{k^n_0}^n (t_1)]
    \leq -2m \alpha_0 \delta_1(\varepsilon),
\end{equation*}
and
\begin{equation*}
\begin{aligned}
        \fIlim_{k^n_0}^n (t_2) - \fIlim_{k^n_0}^n (t_1) 
=     \fmu_{k^n_0} [\fC'_{k^n_0}( \Tfa^n_{k^n_0} (t_2)) - \fC'_{k^n_0}( \Tfa^n_{k^n_0} (t_1))] 
\leq -2\fmu_{\min} m \alpha_0 \delta_1(\varepsilon)
=  - 4m \beta_0 \delta_1(\varepsilon).
\end{aligned}
\end{equation*}
\end{proof}

By Lemma~\ref{lemma: gcmu descent}, we have
$\fIpre_{k^n_0}^n (s^n_1) - \fIpre_{k^n_0}^n (t_1) \leq -
3m\beta_0\delta_1(\varepsilon).$ Assume first that $l_1 \neq k^n_0$. Then
$\fIpre_{l_1}^n ((s^n_1)^-) = \fIpre_{l_1}^n (s^n_1) $ because predicted class
$l_1$ is not served at $s^n_1$. Consequently, there is no negative jump of the
index and
\begin{equation}\label{eq: upper bound of l1}
\begin{aligned}
        \fIpre_{l_1}^n(s_1^n) 
=&~     \underbrace{[\fIpre_{l_1}^n(s^n_1) - \fIpre_{l_1}^n((s^n_1)^-)]}
        _{\text{no negative jump at } s^n_1,~=0}
        + \underbrace{[\fIpre_{l_1}^n((s^n_1)^-) - \fIpre_{k^n_0}^n((s^n_1)^-)]}
        _{\text{by \ourmethod, } \leq 0} \\
        & + \underbrace{[ \fIpre_{k^n_0}^n((s^n_1)^-) 
        -  \fIpre_{k^n_0}^n(s^n_1)]}_{\text{by Lemma~\ref{lemma: size of gcmu's negative jump}, }\leq m \beta_0 \delta_1(\varepsilon)}
      + \underbrace{[\fIpre_{k^n_0}^n(s^n_1) - \fIpre_{k^n_0}^n(t_1)]}_{\leq -3m \beta_0\delta_1(\varepsilon)} 
        + \fIpre_{k^n_0}^n(t_1)\\
\leq&~  \fIpre_{k^n_0}^n(t_1) - m \beta_0 \delta_1(\varepsilon).
\end{aligned}
\end{equation}
For $\fIpre_{l_2}^n (s^n_1)$, since $l_2$ is NOT selected by the \ourmethod~in $[nt_1, ns^n_1]$, $\Tfa^n_{l_2}(s^n_1) - \Tfa^n_{l_2}(t_1) \geq 0$, which yields
\begin{equation}\label{eq: lower bound of l2}
    \fIpre_{l_2}^n (s^n_1) \geq \fIpre_{l_2}^n (t_1) - m\beta_0 \delta_1(\varepsilon)
\end{equation}
by Lemma~\ref{lemma: continuity of gcmu index}. By the condition in the proposition, subtracting~\eqref{eq: lower bound of l2} from~\eqref{eq: upper bound of l1} yields
\begin{equation*}
    \fIpre_{l_1}^n(s^n_1)  - \fIpre_{l_2}^n (s^n_1)
    \leq \fIpre_{k^n_0}^n(t_1) - \fIpre_{l_2}^n (t_1) \leq \varepsilon.
\end{equation*}
If instead $l_1 = k^n_0$, then combining the descent
$\fIpre_{k^n_0}^n(s^n_1) - \fIpre_{k^n_0}^n(t_1) \leq -3m\beta_0\delta_1(\varepsilon)$
of Lemma~\ref{lemma: gcmu descent} with~\eqref{eq: lower bound of l2} yields the
same bound with additional slack,
\begin{equation*}
    \fIpre_{l_1}^n(s^n_1)  - \fIpre_{l_2}^n (s^n_1)
    \leq \fIpre_{l_1}^n(t_1) - \fIpre_{l_2}^n (t_1) - 2m\beta_0\delta_1(\varepsilon)
    \leq \varepsilon.
\end{equation*}

To show $s^n_1 \geq t_1 + \gamma_0\delta^n_1(\varepsilon)$, recall from the
choice of $s^n_1$ that
\begin{equation*}
    \Tfa_{k^n_0}^n(s^n_1) - \Tfa_{k^n_0}^n (t_1) \leq - 2\alpha_0 \delta_1(\varepsilon).
\end{equation*}
Then, by Lemmas~\ref{lemma: continuity of tilde A and S},~\ref{corr: relation
  between age and T}, it is easy to verify that for sufficiently large $n$
\begin{equation*}
\begin{aligned}
    n(s^n_1 - t_1)
\geq  \fT^n_{k^n_0} (ns^n_1) - \fT^n_{k^n_0} (nt_1)
\geq&~  n^{1/2} \cdot \frho_{k^n_0}[ n^{1/2} (s^n_1 - t_1) +   2\alpha_0 \delta_1(\varepsilon) - \alpha_0 \delta_1(\varepsilon)]\\
\geq&~ n^{1/2} \cdot \frho_{\min} [n^{1/2} (s^n_1 - t_1)  + \alpha_0 \delta_1(\varepsilon)],
\end{aligned}
\end{equation*}
where $\fT^n_{k^n_0} (ns^n_1) - \fT^n_{k^n_0} (nt_1) \leq n(s^n_1 - t_1)$
follows from the definition of the policy process $\fT^n_{k^n_0}$. This yields
the desired result that
$s^n_1 - t_1 \geq \frac{\frho_{\min}\alpha_0}{1- \frho_{\min}}
\delta^n_1(\varepsilon) = \gamma_0\delta^n_1(\varepsilon)$. Note that
$\gamma_0 \in (0, 1)$ because the critical load condition
$\sum_k \rho_k=1$ in Assumption~\ref{assumption: heavy traffic} implies
$\rho_{\text{min}} \leq \frac{1}{K}$, $\rho_{\text{max}} \geq \frac{1}{K}$,
and $1 - \rho_{\min} \geq (K-1)/K$, giving
$\gamma_0 \leq 1/[3(K-1)^2] \leq 1/3 < 1$ for $K \geq 2$.

\subsubsection{Proof of Lemma~\ref{lemma: age descent}}
\label{section:proof-age-descent}

  By condition (i), it is clear that 
  $\fT_{l^n_0}^n (nt_1 + n \delta_1^n (\varepsilon)) - \fT_{l^n_0}^n (nt_1) =0$,
  since the predicted class $l^n_0$ is not selected by \ourmethod~in $[nt_1, n(t_1 + \delta_1^n(\varepsilon))]$. Intuitively, the server is busy for serving other predicted classes, implying positive stochastic fluctuations of the policy processes dedicated to the other predicted classes, and there must be at least one predicted classes that absorbs the additional service. In particular, we claim that there exists some predicted class $k^n_0$ such that 
  \begin{equation}\label{eq: additional service in T}
    \fT_{k^n_0}^n (n(t_1 + \delta_1^n (\varepsilon))) - \fT_{k^n_0}^n (nt_1) \geq \Big(\frac{\frho_{\min}}{ K - 1} + \frho_{k^n_0} \Big) \cdot n \delta^n_1(\varepsilon).
  \end{equation}
  For simplicity, for all $l\in [K]$, let
  \begin{equation*}
    \Delta \fT^n_{l}(nt_1): = \fT_{l}^n (n(t_1 + \delta_1^n (\varepsilon))) - \fT_{l}^n (nt_1),
    \quad w^n_l: =   \Delta \fT^n_{l}(nt_1) / (n \delta^n_1(\varepsilon)),
  \end{equation*}
  where $\Delta \fT^n_{l}(nt_1)$ represents the service time allocated to predicted class $l$ during $[nt_1, n(t_1 + \delta_1^n(\varepsilon))]$, and $w^n_l$ denotes its proportion in the time interval. According to conditions (i) and (ii) and Assumption~\ref{assumption: heavy traffic}, it is easy to verify that
  \begin{equation*}
    \sum_{l \neq l^n_0} w^n_l = 1, \quad \sum_{l \neq l^n_0} \frho_l = 1 - \frho_{l^n_0} \leq 1 - \frho_{\min}.
  \end{equation*}
  Rearranging the terms, we can claim that there exists some predicted class
  $k^n_0$ satisfying
  $w^n_{k^n_0} - \frho_{k^n_0} \geq \frac{\frho_{\min}}{ K - 1}$, which is
  equivalent to~\eqref{eq: additional service in T}.

  Combining~\eqref{eq: additional service in T} and Lemmas~\ref{lemma:
    continuity of tilde A and S},~\ref{corr: relation between age and T}, for
  sufficient large $n$
  \begin{equation*}
    \begin{aligned}
      \Tfa_{k^n_0}^n (t_1 + \delta^n_1(\varepsilon)) - \Tfa_{k^n_0}^n (t_1)
      =&~     n^{1/2} \delta^n_1(\varepsilon) - n^{-1/2} \frho_{k^n_0}^{-1} \Delta \fT^n_{k^n_0}(nt) 
      + \alpha_0\delta_1(\varepsilon)\\
      \leq&~  \delta_1(\varepsilon) 
      - \frho_{k^n_0}^{-1} \Big(\frac{\frho_{\min}}{ K - 1} + \frho_{k^n_0} \Big)
      \cdot \delta_1(\varepsilon) 
      + \alpha_0\delta_1(\varepsilon)\\
      \leq&~  - 3 \alpha_0 \delta_1(\varepsilon) + \alpha_0 \delta_1(\varepsilon)
      = - 2\alpha_0 \delta_1(\varepsilon).
    \end{aligned}
  \end{equation*}

\subsection{Proof of Proposition~\ref{prop: gcm difference within intervals}}
\label{section:proof-gcm difference within intervals}

Fix $t_2 \in [t_1,(t_1 + \delta_1^n(\varepsilon))\wedge 1 ]$. By symmetry, it
suffices to show
$\fIpre_{l_1}^n (t_2) - \fIpre_{l_2}^n (t_2) < 3\varepsilon/2$ for any
$l_1, l_2 \in[K]$.  When $l_2$ is selected for service by \ourmethod~in
$[nt_1, nt_2]$, we can employ a similar analysis as in the proof of
Proposition~\ref{prop: max gcmu case i} to show
$\fIpre_{l_1}^n (t_2) - \fIpre_{l_2}^n (t_2) < \varepsilon$. For the other
case, we have from Lemma~\ref{lemma: continuity of gcmu index} that
$\fIpre_{l_2}^n (t_2) - \fIpre_{l_2}^n (t_1) \geq - \varepsilon/4$, since
$\Tfa^n_{l_2}(t_2) - \Tfa^n_{l_2}(t_1) \geq 0$. Also, once again by
Lemma~\ref{lemma: continuity of gcmu index}, one can check
$\fIpre_{l_1}^n (t_2) - \fIpre_{l_1}^n (t_1) \leq \varepsilon/4$ since
$t_2 - t_1 \leq \delta^n_1(\varepsilon)$. Combining equations above yields the
desired result.


\section{Adding to the exposition in classical heavy traffic queueing analysis}
\label{section: adding to the exposition in classical heavy traffic queueing analysis}

In this section, we use our analysis to identify sufficient conditions that address gaps in the original arguments for the lower bound in the classical heavy traffic analysis~\citet{VanMieghem95}. In this sense, we observe that the original proofs in~\citet{VanMieghem95} leave certain gaps, which our analysis addresses by providing sufficient conditions for the claims to hold rigorously. 

When interarrival and service times are i.i.d, the arguments in~\citet{VanMieghem95} hold under their original assumptions. However, when interarrival or service times are not i.i.d in the G/G/1 systems, we identify specific sufficient conditions either on $\BtU_{k}^n$ or $\BtV_{k}^n$ under which the original arguments can be rigorously complemented. Based on our analysis, we found it challenging to rigorously address these gaps without using such conditions.

\subsection{Conditions for complementing the lower bound proof in~\citet{VanMieghem95}}
\label{section: sufficient conditions for the lower bound proof in Van Mieghem}
We provide conditions under which our analysis complements the original proof for the lower bound in~\citet[Proposition 6]{VanMieghem95} for the G/G/1 systems studied in that work. Specifically, these conditions suffice to prove a key intermediate step,~\citet[Proposition 3]{VanMieghem95}, which is analogous to our Lemma~\ref{lemma: relation between N and W}. The convergence is in the almost-sure sense in the uniform metric as in the sample path analysis in~\citet{VanMieghem95}.

\begin{claim} 
If for each class $k \in [K]$, either (i) $(\BtV_k^n)' = (\BtV_k^*)' + o_n(1)$ or (ii) $\BtV_k^n = \BtV_k^*+ o_n(n^{-1/2})$ is satisfied, then the proof of~\citet[Proposition 3]{VanMieghem95} can be rigorously complemented, allowing the queue lengths to be translated to workloads for deriving the lower bound result in~\citet[Proposition 6]{VanMieghem95}. 
\end{claim}
\noindent When the service times are i.i.d in all systems $n$, $\BtV_k^n$ and $\BtV_k^*$ are deterministic and linear. In this case, the first condition is satisfied by the convergence of $\BtV_k^n$ to $\BtV_k^*$ as argued in~\citet[Proposition 1]{VanMieghem95}, and the proofs for~\citet[Proposition 3, 6]{VanMieghem95} apply. Our focus, however, is on G/G/1 systems where the service times are not i.i.d., and thus $\BtV_k^n$ may not be linear processes.

\paragraph{Proof of the claim} 

We argue that having at least one of the two stated conditions is sufficient to prove~\citet[Proposition 3]{VanMieghem95}, a key intermediate result for the lower bound in~\citet[Proposition 6]{VanMieghem95}. Specifically, either condition can be used to derive a variant of Eq (82) from Eq (81) in the proof of~\citet[Proposition 3]{VanMieghem95}. Using the Mean-Value Theorem and the uniform continuity of $(\BtV_k^*)'$ on a compact interval, if condition (i) holds, 

\begin{equation*}
\begin{aligned}
n \BtV_k^n \big(n^{-1} A_k^n (t)\big) - n \BtV_k^n \big( n^{-1} A_k^n (t) - n^{-1} N_k^n (t) \big)&= (\BtV_k^n)' \big(\xi_k^n (t)) \cdot  N_k^n (t)\\
&=(\BtV_k^n)' \big(n^{-1} A_k^n (t)) \cdot  N_k^n (t) + o(N_k^n (t))\\
&=(\BtV_k^*)' \big(n^{-1} A_k^n (t)) \cdot  N_k^n (t) + o(N_k^n (t)) + o_n (1),
\end{aligned}
\end{equation*}
where $\xi_k^n (t) \in [n^{-1} A_k^n (t) - n^{-1} N_k^n (t), n^{-1} A_k^n (t)]$. Similarly, if condition (ii) is guaranteed, we have
\begin{equation*}
\begin{aligned}
&~ n \BtV_k^n \big(n^{-1} A_k^n (t)\big) - n \BtV_k^n \big( n^{-1} A_k^n (t) - n^{-1} N_k^n (t) \big)\\
=&~ n \BtV_k^* \big(n^{-1} A_k^n (t)\big) - n \BtV_k^* \big( n^{-1} A_k^n (t) - n^{-1} N_k^n (t) \big) + o(n^{1/2})\\
=&~ (\BtV_k^*)' \big(n^{-1} A_k^n (t)) \cdot  N_k^n (t) + o(N_k^n (t)) + o (n^{1/2}).
\end{aligned}
\end{equation*}
In either case, the right-hand side can be plugged into~\citet[Eq (81)]{VanMieghem95}, ensuring the proof of~\citet[Proposition 3]{VanMieghem95} is complete. Without either condition, it remains challenging to rigorously derive~\citet[Eq (82)]{VanMieghem95} or its equivalent form based on our analysis.

\subsection{Conditions for complementing the optimality proof in~\citet{VanMieghem95,MandelbaumSt04}}
\label{section: sufficient conditions for the optimality proof in Van Mieghem}
As mentioned in Sections~\ref{section:introduction},~\ref{section: heavy-traffic optimality of pcmu}, and~\ref{subsection: comparison to the optimality result in VM}, the optimality proofs in the classical heavy traffic queueing literature~\citep{VanMieghem95, MandelbaumSt04} are incomplete because they did not establish the equivalence between convergence of the scaled age and sojourn time processes. With $Q^n=I$, our analysis provides a proof for the equivalence even in the traditional setting with known true classes and fills the gap in the literature. We find two alternative conditions that can guarantee the asymptotic equivalence between age and sojourn time in~\citet{VanMieghem95, MandelbaumSt04}.

\begin{claim}
If for each class $k\in [K]$, either (i) $(\BtU_k^n)' = (\BtU_k^*)' + o_n (1)$ or (ii) $\BtU_k^n = \BtU_k^* + o_n (n^{-1/2})$, then the age and sojourn time processes in~\citet{VanMieghem95, MandelbaumSt04} can be shown asymptotically equivalent in a manner analogous to our Proposition~\ref{proposition: relation between a and tau}. These conditions, combined with our proof for Proposition~\ref{prop: convergence of max gcmu difference} with $Q^n=I$, complement the original optimality proofs for the \gcmu~and $D\text{-}Gc\mu$.
\end{claim}
\noindent Note that the conditions resemble those for the lower bound result in the previous section. In~\citet{MandelbaumSt04}, the interarrival times are i.i.d, so $\BtU_k^n$ is linear, and the $o(n^{-1/2})$ rate of arrival-rate convergence assumed in~\citet[Eq (2)]{MandelbaumSt04} then yields $\BtU_k^n = \BtU_k^* + o_n(n^{-1/2})$ uniformly on compacts, satisfying condition (ii), and hence condition (i) as well.

\paragraph{Proof of the claim} Plugging $Q^n=I$ in our framework, Assumptions~\ref{assumption: heavy traffic},~\ref{assumption: on cost functions for showing the lower bound},~\ref{assumption: second order moments}, and~\ref{assumption: on cost functions for optimality} and the result of 
Proposition~\ref{prop: joint conv. of predicted class A, U, S, and V} 
satisfy all assumptions in~\citet{VanMieghem95}. Thus, the relationship between the age and sojourn time processes in Proposition~\ref{proposition: relation between a and tau} and the reformulation of the age processes in~\eqref{eq: age_proof_eq1}, which stemmed from Proposition~\ref{prop: joint conv. of predicted class A, U, S, and V}, apply to the classical setting in~\citet{VanMieghem95} with known true classes, where the age processes can be written as $a_k^n(nt) = nt - U_k^n( A_k^n(nt) - N_k^n(nt) + 1) + o_n(n^{1/2})$.
    
 Following the proof of Proposition~\ref{proposition: relation between a and tau} with $Q^n = I$, we would have under condition (i) that
\begin{equation*}
\begin{aligned}
      \BtU^n_k(n^{-1} (A_k^n(nt) - N_k^n(nt) + 1))
=&~   \BtU^n_k \big( \BtA_k^n (t) + n^{-1/2} \tilde{A}_k^n (t)- n^{-1/2} \tilde{N}_k^n (t) + o(n^{-1/2})  \big)\\
=&~   t + (\BtU_k^*)' (\BtA^n_k (t))\cdot n^{-1/2} \big(\TtA_k^n (t)-\TtN_k^n (t)  \big) + o(n^{-1/2}),
\end{aligned}
\end{equation*} 
where we use $\BtU^n_k =(\BtA_k^n)^{-1}$ in the proof of~\citet[Proposition 2]{VanMieghem95}, the Mean Value Theorem, and the uniform continuity of $(\BtU_k^*)'$. With condition (ii), we can write similarly that
\begin{equation*}
\begin{aligned}
  \BtU^n_k(n^{-1} (A_k^n(nt) - N_k^n(nt) + 1))
=&~\BtU^*_k(n^{-1} (A_k^n(nt) - N_k^n(nt) + 1)) +  o(n^{-1/2}) \\
=&~\BtU^*_k \big( \BtA_k^n (t) + n^{-1/2} \tilde{A}_k^n (t)- n^{-1/2} \TtN_k^n (t) + o(n^{-1/2})  \big) + o(n^{-1/2}) \\
=&~   t + (\BtU_k^*)' (\BtA^n_k (t))\cdot n^{-1/2} \big(\tilde{A}_k^n (t)-\TtN_k^n (t)  \big) + o(n^{-1/2}).
\end{aligned}
\end{equation*} 
Plugging the right-hand side of either case into the remaining proof for Proposition~\ref{proposition: relation between a and tau}, the desired asymptotic equivalence between age and sojourn time in~\citet{VanMieghem95} is established. Under the assumptions in~\citet{VanMieghem95} without condition (i) or (ii), it is not evident how to approximate $\BtU^n_k(n^{-1} (A_k^n(nt) - N_k^n(nt) + 1))$ and establish the asymptotic equivalence between the age and sojourn time processes.


\section{Proofs for Section~\ref{section: characterization of the optimal cost with q}}
\label{section: proof for characterization of the optimal cost with q}

\subsection{Proof for Proposition~\ref{proposition: information gap; quadratic cost}}
\label{section: proof of cumulative cost rate}

From Theorem~\ref{theorem: HT lower bound},
$\TJ^*(t; Q) = \int_0^{t} \sum_{l=1}^K \sum_{k=1}^K \lambda p_k \q_{kl}C_k
(\Tftau_l (s))\mathrm{d} s$ where $\{\Tftau_l\}_{l\in [K]}$ is characterized
as
\begin{equation}\label{eq: quadratic characterization}
    \Tftau_l = \TfW_l / \frho_l,~\forall~l\in[K],\quad \sum_l \TfW_l = \sumTtW, \quad \fmu_l \fC_l'(\Tftau_l) = \fmu_m \fC_m'(\Tftau_m),~\forall~l,m\in[K],
\end{equation}
by our argument in Section~\ref{section: heavy-traffic optimality of pcmu}. 
According to Assumption~\ref{assumption: quadratic cost functions}, 
we can equivalently reformulate~\eqref{eq: quadratic characterization} as
\begin{equation*}
    \frho_l \Tftau_l(t; Q) = \frac{(\beta_l (Q))^{-1}}
    {\sum_{m=1}^K (\beta_m (Q))^{-1}} \sumTtW(t),\quad
    \beta_l (Q) = \frac{\fmu_l\fc_l}{\frho_l},~\forall~t\in [0,1],~\forall~l\in[K].
\end{equation*}
For any $s\in [0, t]$, the integrand 
$\sum_{l=1}^K \sum_{k=1}^K \lambda p_k \q_{kl}C_k (\Tftau_l (s))$ 
can be written as
\begin{equation*}
\begin{aligned}
 \sum\limits_{l=1}^K \sum\limits_{k=1}^K \lambda p_k \q_{kl} 
    \frac{c_k}{2} \frac{1}{\frho_l^2} \Big( \frac{\sumTtW (s)}
    {\sum_{m=1}^K \frac{\beta_l (Q)}{\beta_m (Q)}}\Big)^2 
=&~   \frac{1}{2}\sumTtW^2(s)  
    \sum\limits_{l=1}^K \frac{1}{\frho_l^2} 
    \Big( \frac{1}{\sum_{m=1}^K \frac{\beta_l (Q)}{\beta_m (Q)}}\Big)^2 \sum\limits_{k=1}^K \lambda p_k \q_{kl} c_k   \\
=&~ \frac{1}{2}\sumTtW^2 (s)\sum\limits_{l=1}^K \frac{\beta_l (Q)}{\big(\sum_{m=1}^K \frac{\beta_l (Q)}{\beta_m (Q)}\big)^2},
\end{aligned}
\end{equation*}
where  the last equality holds since $\beta_l (Q)= \fmu_l \fc_l / \frho_l$, $\fc_l = \frac{\sum \lambda p_k \q_{kl} c_k}{\sum \lambda p_k \q_{kl}}$ by definition. The summation term can be further written as
\begin{equation*}
    \begin{aligned}
    \sum\limits_{l=1}^K \frac{\beta_l (Q)}{\big(\sum_{m=1}^K \frac{\beta_l (Q)}{\beta_m (Q)}\big)^2}
=&~  \sum_{l=1}^K \frac{\beta_l (Q) \big(\prod_{r\neq l} \beta_r (Q)\big)^2}
    {\big(\sum_{m=1}^K \frac{\beta_l (Q)}{\beta_m (Q)}\big)^2\big(\prod_{r\neq l} \beta_r (Q)\big)^2}\\
=&~\sum_{l=1}^K \frac{\beta_l (Q) \big(\prod_{r\neq l} \beta_r (Q)\big)^2}
    {\big(\sum_{m=1}^K \prod_{r\neq m} \beta_m (Q)\big)^2}
= \frac{\prod_{r=1}^K \beta_r (Q)}
    {\sum_{m=1}^K \prod_{r\neq m} \beta_r (Q)} =  \frac{1}{\sum_{m=1}^K  (\beta_m (Q))^{-1}}.
\end{aligned}
\end{equation*}

Applying a similar approach to the proof of Proposition~\ref{prop: convergence of max gcmu difference} with the corresponding class indices, the age process converges under the
\naivemethod, and by Lemma~\ref{lemma: convergence of tildeJ}, the cumulative
cost converges to
$$\TJ_\text{Naive}(t; Q) = \sum_{l=1}^K \sum_{k=1}^K \int_0^{t}
\lambda p_k \q_{kl}C_k (\Tftau_{l, \text{Naive}} (s))\mathrm{d} s,$$ 
where
$\{\Tftau_{l, \text{Naive}}\}_{l\in [K]}$ is the limit of the sojourn time
process under the \naivemethod. By similar analysis as in the proof of
Theorem~\ref{theorem: optimality of our policy}, the limit
$\{\Tftau_{l, \text{Naive}}\}_{l\in [K]}$ is characterized by
\begin{equation*}
    \Tftau_{l, \text{Naive}} = \TfW_l / \frho_l,~\forall~l\in[K],\quad \sum_l \TfW_l = \sumTtW, \quad \fmu_l C_l'(\Tftau_{l, \text{Naive}}) = \fmu_m C_m'(\Tftau_{m, \text{Naive}}),~\forall~l,m\in[K].
\end{equation*}
In contrast to Eq.~\eqref{eq: quadratic characterization}, each predicted class $l\in [K]$ is associated with the
\textit{original cost function} $C_l$ in the above characterization, which
does \textit{not} take into account misclassification errors in the marginal
cost rate of the class. It follows that
\begin{equation*}
  \frho_l \Tftau_{l, \text{Naive}}(t; Q) = \frac{(\beta_{l, \text{Naive}} (Q))^{-1}}
  {\sum_{m=1}^K (\beta_{m, \text{Naive}} (Q))^{-1}} \sumTtW(t),\quad
  \beta_{l, \text{Naive}} (Q) = \frac{\fmu_l c_l}{\frho_l},~\forall~l\in[K].
\end{equation*}
Combining the equations above and noting $\beta_l(Q) = \fmu_l \fc_l / \frho_l$, we have
\begin{equation}\label{eq: naive gcmu quadratic cost}
\begin{aligned}
    \TJ_\text{Naive}(t; Q)
=&~ \sum\limits_{l=1}^K \sum\limits_{k=1}^K  \int_0^t \lambda p_k \q_{kl} \frac{c_k}{2 \frho_l^2} \sumTtW^2(s) 
    \Big(\frac{(\beta_{l, \text{Naive}} (Q))^{-1}}
    {\sum_{m=1}^K (\beta_{m, \text{Naive}} (Q))^{-1}}\Big)^2 \mathrm{d}s\\
=&~ \sum\limits_{l=1}^K  
\frac{\beta_l(Q)}{\big(\sum_m \frac{\beta_{l, \text{Naive}}(Q)}{\beta_{m, \text{Naive}}(Q)}\big)^2 } \cdot 
\frac{1}{2}\int_0^t \sumTtW^2(s) \mathrm{d}s.
\end{aligned}
\end{equation}


\section{Proof of results in Section~\ref{section: design of an AI triage system}}
\label{section: appendix triage}

\subsection{Joint convergence of the AI-based triage system}
\label{section: joint convergence of triage system}

We define the concerned processes below to analyze the AI triage system.
\begin{definition}[Arrival processes of the AI-based triage system] 
\label{definition: arrival processes of the AI Triage system}
Given a classifier $\model$, filtering level $\fl$, 
the number of hired reviewers $\Gamma(\fl)$, 
and a sequence of queueing systems, suppose that 
Assumptions~\ref{assumption: data generating process for triage system}
and~\ref{assumption: AI triage system} hold. We define the following for any 
system $n$, reviewer $r\in \Gamma(\fl)$, and time $t\in [0, n]$:
\begin{enumerate}[(i)]
\item (Arrival process of the triage system) Let 
$\tU^n_0(t): = \sum_{i=1}^{\lfloor t \rfloor} u_i^n$ be the partial sum of 
interarrival times among the first 
$\lfloor t \rfloor$ jobs arriving at the triage system, and 
$\tA^n_0(t)$ be the number 
of jobs that arrive at the triage system $n$ up to time $t$. Moreover, let 
$\TtU^n_0(t)$, $\TtA_0^n(t)$ be the corresponding diffusion-scaled process, 
defined as 
\begin{equation*}
\TtU^n_0(t) = n^{-1/2}[ \tU^n_0(nt) - \Lambda_n^{-1} \cdot n t ],
\quad 
\TtA^n_0(t) = n^{-1/2}[ \tA^n_0(nt) - \Lambda_n \cdot n t ],
~\forall~t\in[0, 1];  
\end{equation*}
\item (Arrival process of jobs filtered out) 
For each class $k\in\{1, 2\}$, 
let $\tUflk^n(t)$ be the partial sum of interarrival times among the first 
$\lfloor t \rfloor$ class $k$ jobs that are filtered out, and 
$\tAflk^n(t)$ be the number of class $k$ jobs that 
are filtered out by the filtering system up to time $t$, i.e.,
$\tAflk^n(t) = \sum_{i=1}^{\tA^n_0(t)} \mathbb{I} (\model(X^n_i) < \fl) \cdot Y_{ik}^n,
~\forall~k\in\{1, 2\}$. Moreover, let 
$\TtUfl^n(t)$ and $\TtAflk^n(t)$ be the 
corresponding diffusion-scaled processes, defined as
\begin{equation*}
\begin{aligned}        
\TtUflk^n(t) 
=&~ n^{-1/2}\Big[ \tUflk^n(nt) 
- (\Lambda_n p^n_k (1-g^n_k (\fl)))^{-1} \cdot  n t\Big], 
~\forall~t\in[0, 1],~\forall~k\in\{1, 2\}\\
\TtAflk^n(t) 
=&~ n^{-1/2}\Big[ \tAflk^n(nt) 
- \Lambda_n p^n_k (1-g^n_k (\fl)) \cdot n t \Big],
~\forall~t\in[0, 1],~\forall~k\in\{1, 2\};
\end{aligned}
\end{equation*}
\item (Arrival process of the queueing system) 
Let $\tUps^n(t)$ be the partial sum of interarrival times among the first 
$\lfloor t \rfloor$ jobs that pass through the filtering system and arrive at 
the queueing system, and 
$\tAps^n(t)$ be the number 
of jobs that pass through the filtering system and arrive at the 
queueing system up to time $t$, i.e., 
$\tAps^n(t) = \sum_{i=1}^{\tA^n_0(t)} \mathbb{I} (\model(X^n_i) \geq \fl)$. Also, let 
$\TtUps^n$ and $\TtAps^n$ be the corresponding diffusion-scaled arrival process, 
defined as 
\begin{equation*}
\begin{aligned}
\TtUps^n(t) =&~ 
n^{-1/2}[ \tUps^n(nt) 
- n t \cdot (\Lambda_n \sum_{k=1}^2 p^n_k g^n_k (\fl))^{-1}],
~\forall~t\in[0, 1],\\
\TtAps^n(t) =&~ 
n^{-1/2}[ \tAps^n(nt) 
- n t \cdot \Lambda_n  \sum_{k=1}^2 p^n_k g^n_k (\fl)],
~\forall~t\in[0, 1];
\end{aligned}
\end{equation*}
\item (Arrival process of each reviewer) 
Let $\tUpsr^n(t): = \sum_{s=1}^{\lfloor t \rfloor} u^n_{s,r}$ 
be the partial sum of interarrival times among the first 
$\lfloor t \rfloor$ jobs that are assigned to reviewer $r$, and 
$\tApsr^n(t)$ be the number 
of jobs that are assigned to reviewer $r$ up to time $t$, i.e., 
$\tApsr^n(t) = \sum_{j=1}^{\tAps^n(t)} B_{jr}$. Moreover, 
let $\VTtUpsr^n(t) = \{\TtUpsr^n(t)\}_{r\in \Gamma(\fl)}$,
$\VTtApsr^n(t) = \{\TtApsr^n(t)\}_{r\in \Gamma(\fl)}$ 
be the corresponding diffusion-scaled arrival process, defined as 
\begin{equation*}
\begin{aligned}
  \TtUpsr^n(t) =&~ 
  n^{-1/2}\Big[ \tUpsr^n(nt)
  - nt \cdot \frac{\Gamma(\fl)}{\Lambda_n \sum_{k=1}^2 p^n_k g^n_k (\fl)} \Big],
  ~\forall~t\in[0, 1],\\
\TtApsr^n(t) =&~ 
n^{-1/2}\Big[ \tApsr^n(nt)
- nt \cdot \frac{\Lambda_n}{\Gamma(\fl)} \sum_{k=1}^2 p^n_k g^n_k (\fl) \Big],
~\forall~t\in[0, 1];
\end{aligned}
\end{equation*}
\item (Split probability) Let $\pflk^n$ be the probability that 
a job arriving at the triage system is of class $k$ and is filtered out by the
filtering system, i.e., 
$\pflk^n = p^n_k (1-g^n_k (\fl))$, 
and $\pps^n$ be the probability that a job arriving at the
triage system passes through the filtering system, i.e., 
$\pps^n = \sum_{k=1}^2 p^n_k g^n_k (\fl)$. Moreover, let $\pflk$ and $\pps$ be 
the corresponding limiting probability defined as 
$\pflk =p_k (1-g_k (\fl))$ and
$\pps = \sum_{k=1}^2 p_k g_k (\fl)$.
\item (Spliting process) Let $\tSpflk(t)$ be the number of jobs that
are of class $k$ and filtered out by the filtering system among the first $\lfloor t \rfloor$ jobs 
arriving at the triage system,
and $\tSppsr(t)$ be the number of jobs that are 
assigned to reviewer $r$ among the first $\lfloor t \rfloor$ jobs arriving at the 
triage system. 
Moreover, let $\VTtSpfl(t) = \{\TtSpflk(t)\}_{k\in \{1, 2\}}$, 
$\VTtSpps(t) = \{\TtSppsr(t)\}_{r\in \Gamma(\fl)}$ be the 
corresponding diffusion-scaled splitting process, defined as 
\begin{equation*}
\begin{aligned}
\TtSpflk(t)=&~n^{-1/2} [\tSpflk(nt) - \pflk^n\cdot nt],
~\forall~t\in [0, 1],~k\in \{1, 2\}\\
\TtSppsr(t)=&~n^{-1/2} [\tSppsr(nt) - \frac{\pps^n \cdot nt}{\Gamma(\fl)}],
~\forall~t\in [0, 1].
\end{aligned}
\end{equation*}
\end{enumerate}
\end{definition}
Similar to Definition~\ref{definition: U, Z, R, V}, we define processes above 
on $[0, n]$ or $[0, 1]$ for analysis simplicity. These processes can be 
naturally extended to $[0, +\infty)$ to apply the martingale FCLT 
(Lemma~\ref{lemma: Martingale FCLT})
and FCLT for split processes from~\citep[Theorem 9.5.1]{Whitt02}, 
which yields the joint convergence result below. With a slight abuse of notation,
we adopt Assumption~\ref{assumption: second order moments} to guarantee 
unform integrability of quantities associated with the triage system.  
\begin{lemma}[Joint convergence of the AI-based triage system]
\label{lemma: joint convergence of triage system}
Given a classifier $\model$, filtering level $\fl$, the number of hired reviewers  
$\Gamma(\fl)$, and a sequence of queueing systems, suppose that 
Assumptions~\ref{assumption: 
second order moments},~\ref{assumption: data generating process for triage 
system}, and~\ref{assumption: AI triage system} hold. 
Then, we have that: 
(i) there exists Brownian motion $(\TtA_0, \VTtSpfl, \VTtSpps)$ such that 
$(\TtA^n_0, \VTtSpfl^n, \VTtSpps^n) \Rightarrow (\TtA_0, \VTtSpfl, \VTtSpps)$ 
in $(D^{\Gamma(\fl) + 3}, WJ_1)$; 
(ii) there exist continuous stochastic processes $(\TtAflp, \TtAfln, \VTtApsr)$ 
such that
\begin{equation*}
(\TtAflp^n, \TtAfln^n, \VTtApsr^n) \Rightarrow 
(\TtAflp, \TtAfln, \VTtApsr), \quad\text{in } (D^{\Gamma(\fl)+2}, WJ_1),
\end{equation*}
where $\TtAflk(t) = \pflk\TtA_0(t) + \TtSpflk(\Lambda t)$
and $\TtApsr(t) = \frac{\pps \TtA_0(t)}{ \Gamma(\fl)} 
+ \TtSppsr(\Lambda t)$;
(iii) there exists continuous stochastic processes 
$(\TtUflp, \TtUfln, \VTtUpsr)$ such that 
\begin{equation*}
(\TtUflp^n, \TtUfln^n, \VTtUpsr^n) \Rightarrow 
(\TtUflp, \TtUfln, \VTtUpsr), \quad\text{in } (D^{\Gamma(\fl)+2}, WJ_1).
\end{equation*}

\end{lemma}
\begin{proof}
As for (i), according to Assumption~\ref{assumption: second order moments}, 
we have that  $\text{Var}[u^n_1]< +\infty$ for each $n$,
and $\text{Var}[u^n_1]$ converges to some constant $\sigma^2_u$.
Then, by martingale FCLT (Lemma~\ref{lemma: Martingale FCLT}), 
it is easy to show that 
$(\TtU^n_0, \VTtSpfl^n, \VTtSpps^n)$ jointly convrges to 
$(\TtU_0, \VTtSpfl, \VTtSpps)$. Here, 
$\TtU_0$ is a zero-drift Brownian motion with variance being some 
$\sigma^2_u$,
and $\VTtSpps$  is a zero-drift Bronian motion with covariance matrix being
$\Sigma = (\sigma^2_{r_1, r_2})$, where 
$\sigma^2_{r_1, r_1} = \frac{\Gamma(\fl)-1}{\Gamma^2(\fl)}$ and 
$\sigma^2_{r_1, r_2} = - \frac{1}{\Gamma^2(\fl)},~\forall~r_1\neq r_2$.
According to~\citep[Corollary 13.8.1]{Whitt02}, the joint convergence 
of $(\TtA_0, \VTtSpfl, \VTtSpps)$ follows immediately. (ii) is a direct consequence 
of (i) and~\citep[Theorem 9.5.1]{Whitt02}. Then, by~\citep[Corollary 13.8.1]{Whitt02}, 
(iii) is a corollary of (ii).
\end{proof}

With a slight abuse of notation, we extend from Definition~\ref{definition: U, Z, R, V} 
and Section~\ref{section: model} in order to define $\fZ^n_{kl, r}$, $\fR^n_{l, r}$, 
$\tVpsr^n$ on the jobs that are assigned to each reviewer $r$. 
Let $\VTfZ^n: = \{\TfZ^n_{kl, r}\}_{k,l\in \{1, 2\}, r\in [\Gamma(\fl)]}$, 
$\VTfR^n: = \{\TfR^n_{l, r}\}_{l\in \{1, 2\}, r\in [\Gamma(\fl)]}$, $\VTtVps^n: = \{\TtVpsr^n\}_{r\in [\Gamma(\fl)]}$ 
be the corresponding diffusion-scaled processes. 
As the job assignment process is independent of any other random objects by Assumption~\ref{assumption: AI triage system}, it is easy to show
that $\{(\TfZ^n_{kl, r}, \TfR^n_{l, r}, \TtVpsr^n)\}$ are i.i.d. processes 
across all reviewers. Therefore, 
by independence and Lemma~\ref{lemma: joint convergence under indep}, we can extend 
Lemma~\ref{lemma: individual weak convergence} to achieve joint convergence of 
$(\VTfZ^n, \VTfR^n, \VTtVps^n)$ over all reviewers. 
\begin{lemma}[Joint weak convergence of the AI-based triage system I]
\label{lemma: joint weak convergence of triage system I}
Suppose that Assumptions~\ref{assumption: 
second order moments},~\ref{assumption: data generating process for triage 
system}, and~\ref{assumption: AI triage system} hold. 
Then, there exist Brownian motions $(\VTfZ, \VTfR, \VTtVps)$ such that
\begin{equation*}
(\VTfZ^n, \VTfR^n, \VTtVps^n) \Rightarrow 
(\VTfZ, \VTfR, \VTtVps),\quad \text{in } (D^{7\Gamma(\fl)}, WJ_1).
\end{equation*}
\end{lemma}

Next, we claim that $(\TtUflp^n, \TtUfln^n, \VTtUpsr^n)$ and 
$(\VTfZ^n, \VTfR^n, \VTtVps^n)$ are independent processes under 
Assumption~\ref{assumption: data generating process for triage system}.
Recall that by Definition~\ref{definition: concerned process of triage system},
$\tUps^n(t)$ denotes the partial sum of interarrival times among the first
$\lfloor t \rfloor$ jobs that pass throught the filtering system. 
Let $\{u^n_j: j\in \mathbb{N}\}$ and 
$\{(X^n_{j}, v^n_{j}, Y^n_{j}): j\in \mathbb{N}\}$ be the interarrival time 
and tuples for jobs that \emph{pass through the filtering system}. Then, we have that 
$\tUps^n(t): = \sum_{j= 1}^{\lfloor t\rfloor} u^n_j$.

We first argue that $\{u^n_j: j\in \mathbb{N}\}$ and 
$\{(X^n_{j}, v^n_{j}, Y^n_{j}): j\in \mathbb{N}\}$ are independent.
Let $\{u_i^n: i \in \mathbb{N}\}$ and  $\{(X^n_i, v^n_i, Y^n_i): i \in \mathbb{N}\}$
be the interarrival times and tuples for all jobs arriving at the triage system. 
Note that the primitive sequences $\{u_i^n: i \in \mathbb{N}\}$ and  
$\{(X^n_i, v^n_i, Y^n_i): i \in \mathbb{N}\}$ are independent by 
Assumption~\ref{assumption: data generating process for triage system} (ii).
Therefore, by construction, $\{u^n_j: j\in \mathbb{N}\}$ are the thinned
interarrival times from $\{u_i^n: i \in \mathbb{N}\}$, where each arriving job is 
retained independetly with equal probability $\pps^n$. 
Moreover, since $\{(X^n_i, v^n_i, Y^n_i): i \in \mathbb{N}\}$ are i.i.d. by 
Assumption~\ref{assumption: data generating process for triage system} (i), 
$\{(X^n_{j}, v^n_{j}, Y^n_{j}): j\in \mathbb{N}\}$ are also i.i.d., following the 
conditional distribution $(X^n_{1}, v^n_{1}, Y^n_{1}) \mid \model(X^n_1) \geq \fl$. 
It is important to note that 
the realization of $u^n_j$ can not provide additional information on 
a job that is known to have been retained, i.e., $\model(X^n_i) \geq \fl$), and the tuple $(X^n_{j}, v^n_{j}, Y^n_{j})$ by 
Assumption~\ref{assumption: data generating process for triage system} (ii).

According to analysis above, $\tUps^n(t)$ and 
$\{(X^n_{j}, v^n_{j}, Y^n_{j}): j\in \mathbb{N}\}$ are independent, as the 
former is a function of $\{u^n_j : j\in \mathbb{N}\}$. 
Let $\{(X^n_{s, r}, v^n_{s, r}, Y^n_{s, r}): s \in \mathbb{N}\}$
be the tuples for jobs assigned to some reviewer $r$, which is split from $\{(X^n_{j}, v^n_{j}, Y^n_{j}): j\in \mathbb{N}\}$ 
according to the reviewer assignment $\{\mathbf{B}^n_j: j \in \mathbb{N}\}$.
Since $\{\mathbf{B}^n_j: j \in \mathbb{N}\}$
is independent of any other random objects by 
Assumption~\ref{assumption: data generating process for triage system}, 
we can adopt a similar approach as above to establish 
independence between $(\TtUflp^n, \TtUfln^n, \VTtUpsr^n)$ and
$\{(X^n_{s, r}, v^n_{s, r}, Y^n_{s, r}): s \in \mathbb{N}, r\in [\Gamma(\fl)]\}$, which further yields independence between
$(\TtUflp^n, \TtUfln^n, \VTtUpsr^n)$ and $(\VTfZ^n, \VTfR^n, \VTtVps^n)$. Finally, according to 
Lemmas~\ref{lemma: joint convergence under indep},~\ref{lemma: joint convergence of triage system}, and~\ref{lemma: joint weak convergence of triage system I}, 
such independence leads to 
the joint weak convergence of the AI triage system below
(Lemma~\ref{lemma: joint weak convergence of triage system}),
which extends Lemma~\ref{lemma: joint weak convergence}.

\begin{lemma}[Joint weak convergence of the AI-based triage system II]
\label{lemma: joint weak convergence of triage system}
Suppose that Assumptions~\ref{assumption: 
second order moments},~\ref{assumption: data generating process for triage 
system}, and~\ref{assumption: AI triage system} hold.  Then, we have that 
\begin{equation*}
(\TtUflp^n, \TtUfln^n, \VTtUpsr^n, 
\VTfZ^n, \VTfR^n, \VTtVps^n) \Rightarrow 
(\TtUflp, \TtUfln, \VTtUpsr, \VTfZ, \VTfR, \VTtVps),
\quad \text{in } (D^{8\Gamma(\fl)+2}, WJ_1).
\end{equation*}
\end{lemma}
Similarly to Lemma~\ref{lemma: uniform convergence of UZRV}, we can then 
strengthen the convergence to uniform topology and conduct 
sample path analysis on copies of the original processes. With a slight abuse
of notation, we still use 
$(\Omega_\text{copy}, \mathcal{F}_\text{copy}, \mathbb{P}_\text{copy})$ to denote the 
common probability space. 

\begin{lemma}[Uniform Convergence of the AI Triage System]
\label{lemma: uniform convergence of triage system}
Suppose that Assumptions~\ref{assumption: data generating process for triage 
system},~\ref{assumption: AI triage system}, and~\ref{assumption: 
second order moments} hold.  Then, there 
exist stochastic processes  $(\TtUflp^n, \TtUfln^n, \VTtUpsr^n, 
\VTfZ^n, \VTfR^n, \VTtVps),~\forall~n \geq 1$ and 
$(\TtUflp, \TtUfln, \VTtUpsr, \VTfZ, \VTfR, \VTtVps)$ defined on a common probability 
space $(\Omega_\text{copy}, \mathcal{F}_\text{copy}, \mathbb{P}_\text{copy})$ 
such that 
$(\TtUflp^n, \TtUfln^n, \VTtUpsr^n, \VTfZ^n, \VTfR^n, \VTtVps), ~\forall~n \geq 1$ and 
$(\TtUflp, \TtUfln, \VTtUpsr, \VTfZ, \VTfR, \VTtVps)$ are identical in distribution to 
their original counterparts and 
\begin{equation*}
    (\TtUflp^n, \TtUfln^n, \VTtUpsr^n, 
    \VTfZ^n, \VTfR^n, \VTtVps) \rightarrow 
    (\TtUflp, \TtUfln, \VTtUpsr, \VTfZ, \VTfR, \VTtVps),
    \quad \text{in } (D^{8\Gamma(\fl)+2}, \| \cdot \|), 
    \quad \mathbb{P}_\text{copy}-a.s..
    \end{equation*}
\end{lemma}

\subsection{Sample path analysis of each reviewer}
\label{section: sample path analysis of each reviewer}

In this section, we conduct sample path analysis for each reviewer. 
We adopt a similar analysis approach as in
Section~\ref{section: convergence and lower bound} 
and~\ref{section: heavy-traffic optimality of pcmu}. 
In particular, we consider copies of the original processes 
defined on the common probability space $\mathbb{P}_\text{copy}$, as shown 
in Lemma~\ref{lemma: uniform convergence of triage system}. 
We establish all subsequent results regarding almost sure convergence for 
the copied processes, which can then be converted back into corresponding 
weak convergence results for the original processes.

\paragraph{Heavy Traffic Condition for Each Reviewer}
We first show that our 
Assumptions~\ref{assumption: data generating process for triage system}
and~\ref{assumption: AI triage system} are compatible with 
Assumptions~\ref{assumption: data generating process} 
and~\ref{assumption: heavy traffic}
we adopt for each single-server queueing system.
\begin{definition}
Given a classifier $\model$, filtering level $\fl$, toxicity level $\tx$, the number of hired reviewers $\Gamma(\fl)$, and a sequence of queueing systems, suppose that 
Assumptions~\ref{assumption: data generating process for triage system}
and~\ref{assumption: AI triage system} hold. We define the following for 
any system $n$ and reviewer $r$:
\begin{enumerate}[(i)]
\item (Class prevalence) Let $\tp^n_{k, r}(\fl)$ be the conditional probability 
that a job that passes through the filtering system and is assigned to reviewer $r$ 
is of class $k$, i.e.,
$\tp^n_{k, r}(\fl): = \P^n[\tY_{1k, r}^n=1 \mid \model(X^n_{1, r}) \geq \fl]$. 
Moreover, let $\tp_k(\fl)$ be the limiting probability, defined as
$\tp_k(\fl): = \frac{\tp_k g_k (\fl)}{\tp_1 g_1 (\fl)+\tp_2 g_2 (\fl)}$;
\item (Confusion matrix) Let $\q^n_{kl, r}(\mathbf{z})$ be the conditional 
probability that a class $k$ job arriving at reviewer $r$ is predicted as class 
$l$, i.e., 
$\q^n_{kl, r}(\mathbf{z}): = 
\P^n [ \fY^n_{1l, r} = 1 \mid \model(X^n_{1, r}) \geq \fl, \tY_{1k, r}^n=1]$. 
Moreover, let $\q_{kl}(\mathbf{z})$ be the limiting probability, defined as
$\q_{k1}(\fl, \tx) = \frac{g_k (\tx)}{g_k (\fl)}$, 
$\q_{k2}(\fl, \tx)
= \frac{g_k (\fl) - g_k (z_{\text{tx}})}{g_k (\fl)},~\forall~k\in\{1, 2\}$;
\item (Arrival rate) Let 
$\lambda^n_r = \frac{\Lambda_n}{\Gamma(\fl)}[p^n_1 g^n_1(\fl) + p^n_2 g^n_2(\fl)]$ be the arrival rate of jobs assigned to reviewer $r$. Moreover, let
$\lambda = \frac{\Lambda}{\Gamma(\fl)}[p_1 g_1(\fl) + p_2 g_2(\fl)]$ be the limiting arrival rate.
\end{enumerate}
\end{definition}
\noindent 
We define the arrival rate $\lambda^n_r$ based on 
Lemma~\ref{lemma: joint convergence of triage system}, which shows that 
$n^{-1}\tApsr^n(nt) = \frac{\Lambda_n t}{\Gamma(\fl)} 
\cdot \sum_{k=1}^2 p_k^n g^n_k(\fl) + o_n(1)$. 
According to Assumptions~\ref{assumption: data generating process for triage system}
and~\ref{assumption: AI triage system}, it is easy to verify that 
class prevalence, confusion matrix, and arrival rate all converges to 
their limiting values at the rate of $n^{1/2}$. 
\begin{lemma}\label{lemma: conv of lambda p q for each reviewer}
Given a classifier $\model$, filtering level $\fl$, toxicity level $\tx$, the number of hired reviewers $\Gamma(\fl)$, and a sequence of queueing systems, suppose that 
Assumptions~\ref{assumption: data generating process for triage system}
and~\ref{assumption: AI triage system} hold. Then, 
for any $k, l\in\{1, 2\}$, an reveiwer $r\in [\Gamma(\fl)]$, we have that
\begin{equation*}
    n^{1/2}(\lambda^n_r - \lambda) \rightarrow 0, \quad
    n^{1/2}(\tp^n_{k, r}(\fl) - \tp_{k}(\fl)) \rightarrow 0, \quad
    n^{1/2}(\q^n_{kl, r}(\mathbf{z}) - \q_{kl}(\mathbf{z})) \rightarrow 0.
\end{equation*}
\end{lemma}
\noindent 
As a direct corollary of Lemma~\ref{lemma: conv of lambda p q for each reviewer}
and Assumption~\ref{assumption: AI triage system} (ii), for each reviewer $r$, their limiting traffic intensity 
satisfies 
\begin{equation*}
    \lambda \sum_{k=1}^2 \frac{p_k(\fl)}{\mu_k} 
    = \frac{\Lambda}{\Gamma(\fl)}[p_1 g_1(\fl) + p_2 g_2(\fl)] \cdot 
    \sum_{k=1}^2 \frac{p_k g_k(\fl)}{\mu_k(p_1 g_1(\fl) + p_2 g_2(\fl))} = 1.
\end{equation*}
Therefore, all reviewers operate under heavy traffic conditions, each of which satisfies
Assumption~\ref{assumption: data generating process} 
and~\ref{assumption: heavy traffic}. This enables us to directly apply results 
for the single-server queueing system to each reviewer. Since the analysis is 
similar, we only present the main results below and skip proof details.

\paragraph{Endogenous Processes of the AI-based Triage System}
We define the concerned endogenous processes below to analyze the 
AI-based triage system following Definition~\ref{definition: concerned process}.
\begin{definition}[Endogenous processes of the AI-based triage system]
\label{definition: concerned process of triage system}
Given the filtering level $\fl$, toxicity level $\tx$,  
and the number of hired reviewers $\Gamma(\fl)$, we define the following processes
associated with reviewer $r$ in system $n$:
\begin{enumerate}[(i)]
\item (Input process for predicted classes) Let $\fL_{l, r}^n (t)$ be the total service time
requested by all jobs predicted as class $l$ by time
$t\in[0, n]$, i.e., $\fL_{l, r}^n (t) =\sum_{s=1}^{\tApsr^n (t)} \fY_{sl,r}^n \tv_{s,r}^n$, $t\in [0,n]$. Moreover, let $\TfL_{l,r}^n (t)$ be the corresponding diffusion-scaled process, defined as 
\begin{equation*}
\TfL_{l, r}^n (t) = n^{-1/2}
\Big[\fL_{l,r}^n (nt) - \frac{\Lambda^n}{\Gamma(\fl)} 
\sum_{k=1}^K \frac{p^n_k g^n_k(\fl)}
{\mu^n_k}\q_{kl}^n(\mathbf{z}) \cdot nt \Big],~t\in [0,1].
\end{equation*}
\item (Cumulative total input process) Let
$\sumtL^n(t; \mathbf{z}, r) = \sum_l \fL_{l,r}^n (t), t\in[0, n]$ be the
cumulative total input process and 
$\sumTtL^n(t; \mathbf{z}, r) := \sum_{l=1}^K \TfL_{l,r}^n (t),~t\in[0, 1]$ be the
corresponding diffusion-scaled process, i.e.,
\begin{equation*}
\sumTtL^n (t; \mathbf{z}, r) = n^{-1/2} \Big[\sumtL^n (nt; \mathbf{z}, r) - \frac{\Lambda^n}{\Gamma(\fl)} 
\sum_{k=1}^K \frac{p^n_k g^n_k(\fl)}{\mu^n_k}\cdot nt\Big],~\forall~t\in [0,1].
\end{equation*}
\item (Policy process) Let $\fT_{l, r}^n (t)$ be total amount of service time dedicated to predicted class $l$ in $[0, t]$.
\item (Remaining workload process) Let $\fW_{l,r}^n (t)$ be the remaining service
time requested by jobs predicted as class $l$ and present---waiting for
service or being served---at time $t \in [0, n]$,
\begin{equation*}
\fW_{l,r}^n (t) = \fL_{l,r}^n (t) - \fT_{l,r}^n (t),\quad t\in [0,n].
\end{equation*}
and $\TfW^n_{l,r}(t): = n^{-1/2} \fW^n_{l,r}(nt),~\forall~t\in [0, 1]$ be the
corresponding diffusion scaled process.
\item (Total remaining workload process) Let
$\sumtW^n(t; \mathbf{z}, r) = \sum_l \fW_{l,r}^n (t)$ be the total remaining workload p
rocess and $\sumTtW^n(t; \mathbf{z}, r): = n^{-1/2} \sum_{l=1}^K \fW^n_l(nt; \mathbf{z}, r),
~\forall~t\in [0, 1]$
be the corresponding diffusion scaled process.
\end{enumerate}
\end{definition}
\noindent
Extending Lemma~\ref{lemma: conv. of predicted L} and Proposition~\ref{prop: 
convergence and approximation of predicted class N, tau, T, and W},
we have the following results for the endogenous processes of each reviewer $r$.
\begin{lemma}[Convergence of $\sumTtL^n(t; \mathbf{z}, r)$ and $\sumTtW^n(t; \mathbf{z}, r)$] 
\label{lemma: conv. of L and W triage system} 
Suppose that Assumptions~\ref{assumption: 
second order moments},~\ref{assumption: data generating process for triage 
system}, and~\ref{assumption: AI triage system} hold. 
\begin{enumerate}[(i)]
\item for a sequence of feasible policies $\{\policyn\}$, 
we have that for each reviewer $r$,
\begin{equation*}
\begin{aligned}
&\sumTtL^n(\cdot; \mathbf{z}, r) \rightarrow \sumTtL(\cdot; \mathbf{z}, r)~\text{in } (\mathcal{D}, \|\cdot \|)~\mathbb{P}_\text{copy}\text{-a.s.},~\text{where}\\
&\sumTtL(t; \mathbf{z}, r):=
\TtVpsr\Big(\frac{\Lambda t}{\Gamma(\fl)}[p_1 g_1(\fl) + p_2 g_2(\fl)]\Big) 
+ \sum\limits_{k=1}^K \frac{\tp_k(\fl)}{\tmu_k} \TtApsr(t),~t\in [0,1].
\end{aligned}
\end{equation*}
\item for a sequence of work-conserving p-FCFS feasible policy, 
we have that for each reviewer $r$,
$\sumTtW^n(\cdot; \mathbf{z}, r) \rightarrow \sumTtW(\cdot; \mathbf{z}, r) 
:= \phi(\sumTtL (\cdot; \mathbf{z}, r))$ in $(\mathcal{D}, \|\cdot \|)~\mathbb{P}_\text{copy}-a.s.$, where $\phi$ is the reflection mapping.
\end{enumerate}
\end{lemma}
\noindent
Starting from Lemma~\ref{lemma: conv. of L and W triage system}, we 
can then follow the sample path analysis and establish 
Theorem~\ref{theorem: total cost of triage queue system} in reviewer-level; we skip the detailed 
proof here.

\subsection{Simulation of the total cost of the AI-based Triage System}
\label{section: simulation of the total cost of the AI-based Triage System}
As shown in Theorem~\ref{theorem: total cost of triage queue system}, 
the limiting total cost is solely determined by (i) the limiting exogenous quantities, 
such as arrival rate $\Lambda$, 
class prevalence $p_k(\fl)$, confusion matrix $\q_{kl}(\fl)$, etc; and 
(ii) the limiting total workload process $\sumTtW(\cdot; \mathbf{z}, r)$. 
Though (i) can be easily estimated, (ii) requires a more detailed analysis to 
assist a practical estimation. 

According to Lemma~\ref{lemma: conv. of L and W triage system}, 
$\sumTtW(\cdot; \mathbf{z}, r)$ is a continuous stochastic process. Therefore, it 
suffices to approximate the integral by a Riemann sum, using the explicit characterization of $ \sumTtW(t; \mathbf{z}, r)$ in the lemma.

To that end, we analyze $\TtApsr(t)$ and $\TtVpsr$ separately.
By Lemma~\ref{lemma: joint convergence of triage system},
we have that 
$\TtApsr(t) = \frac{\pps \TtA_0(t)}{ \Gamma(\fl)} 
+ \TtSppsr(\Lambda t)$.
Note that $\TtU_0$ is a zero-drift Brownian motion with variance being
$\sigma^2_u < +\infty$ by Assumption~\ref{assumption: second order moments}, which can be estimated similarly as in Section~\ref{subsection: proof of individual weak convergence}. Then, by~\cite[Corollary 13.8.1]{Whitt02}, we have that 
$\TtA_0(t) = - \TtU_0(\Lambda t)$ and 
\begin{equation*}
  \TtApsr(t) = -\frac{\Lambda \pps}{\Gamma(\fl)} \TtU_0(\Lambda t) + 
  \TtSppsr(\Lambda t).
\end{equation*}
Here, $\VTtSpps$  is a zero-drift Bronian motion with covariance matrix being
$\Sigma = (\sigma^2_{r_1, r_2})$, where 
$\sigma^2_{r_1, r_1} = \frac{\Gamma(\fl)-1}{\Gamma^2(\fl)}$ and 
$\sigma^2_{r_1, r_2} = - \frac{1}{\Gamma^2(\fl)},~\forall~r_1\neq r_2$; see discussion following~\cite[Theorem 9.5.1]{Whitt02}. 
For $\TtVpsr$, according to Assumption~\ref{assumption: second order moments}, it 
is easy to verify that $\text{Var}[v^n_{s,r}]< +\infty$ for each $n$ and 
converges to some constant $\sigma^2_{v}(\fl)= 
 \alpha_{v,1}  p_1(\fl)+ \alpha_{v,2} p_2(\fl) 
- \big(\frac{1}{\mu_1} p_1(\fl) + \frac{1}{\mu_2} p_2(\fl)\big)^2$.
Then, by martingale FCLT (Lemma~\ref{lemma: Martingale FCLT}), we have that 
$\TtVpsr$ is a zero-drift Brownian motion with variance being $\sigma^2_v(\fl)$.

In this way, we rewrite $\sumTtW(t; \mathbf{z}, r)$ as a function of 
(multi-dimensional) Brownian motion, whose Riemann sum can be easily simulated.

\section{Numerical Experiments for the AI-based triage system}
\label{section: experiment triage system}

Our formulation trades off multiple desiderata, in contrast to the standard
industry practice that choose $z$ solely based on prediction metrics,
such as maximizing recall subject to a fixed high precision
level~\citep{ChandakAb23}. To compare our proposed approach with such standard
triage design approaches, we consider the 2-class content moderation problem
described in Section~\ref{section: characterization of the optimal cost with q}. 
We assume the covariates
for positive and negative classes are generated in the same fasion as in the
2d logistic regression problem in 
Section~\ref{section: characterization of the optimal cost with q}, and
consider the logistic regression classifier $\model$ developed by minimizing 
the equally-weighted cross-entropy loss ($w_1 = 1$, $w_2 = 1$). 
For simplicity, we fix the toxicity level
$\tx = 0.5$ and only study how filtering level $\fl$ affects the total cost.

We examine the setting where the positive class has a relatively high arrival rate
to mimic the setting where only flagged content is sent to the triage system, which
results in a relatively high proportion of positive class; 
recall Figure~\ref{fig:diagram}. 
In particular, we set $[\Lambda_1, \Lambda_2] = [10000, 40000]$, 
$[\mu_1, \mu_2] = [50, 200]$, where 
$[\Lambda_1, \Lambda_2]$ is the arrival rate of positive and negative classes to 
the triage system, and $[\mu_1, \mu_2]$ is the common service rate for the positive 
and negative classes across all reviewers. We consider three cases: 
(i) filtering costs dominate, (ii) hiring costs dominate,
and (iii) a trade-off between filtering cost and hiring costs.
The filtering costs and hiring costs are set as follows:
(i) $[\ctrp, \ctrn] = [200, -3]$, $c_r = 800$, 
(ii) $[\ctrp, \ctrn] = [20, -3]$, $c_r = 8000$, 
and (iii) $[\ctrp, \ctrn] = [20, -3]$, $c_r = 800$.
In all cases, the misclassification costs are set as 
$[\cfp, \cfn, \ctp, \ctn] = [3, 3, -3, -3]$, 
and the delay costs are set as
$C_{\cdot}(t) = c_{\cdot} t^2/2$ with $c_1 = 15, c_2 = 1$. 

Our goal is to find the best filtering level $\fl$ that minimizes the total
cost.  We compare our method that minimizes~\eqref{eq: total cost of triage
  queue system} to the following classical method from~\citet{ChandakAb23}, 
which finds the filtering level $\fl$ by maximizing recall subject to a high 
precision level lower bound $z_\text{prec}\in [0,1]$. \footnote{We follow notations 
used in~\citet{ChandakAb23}. For the filtering system, precision and recall 
are calculated by treating safe content as the positive class.}  
Both methods can be effectively implemented using small
set of validation data and a linear search. 
We set the search range for $\fl$ as $[0.05, 0.48]$.

\begin{figure}[t]
\vspace{-1em}
  \centering  
\begin{minipage}[t]{0.32\textwidth}
\centering
\includegraphics[width=0.985\columnwidth]{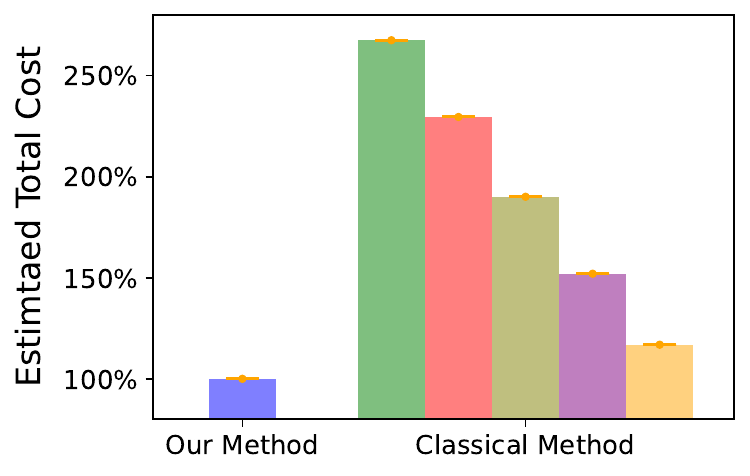} 
(i) Filtering costs dominate
\end{minipage}
\begin{minipage}[t]{0.32\textwidth}
  \centering
\includegraphics[width=0.985\columnwidth]{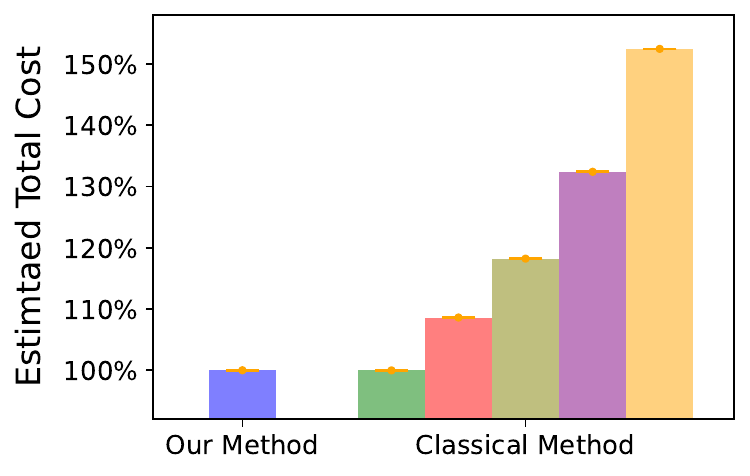} 
(ii) Hiring costs dominate
\end{minipage}
\begin{minipage}[t]{0.32\textwidth}
\centering
\includegraphics[width=0.985\columnwidth]{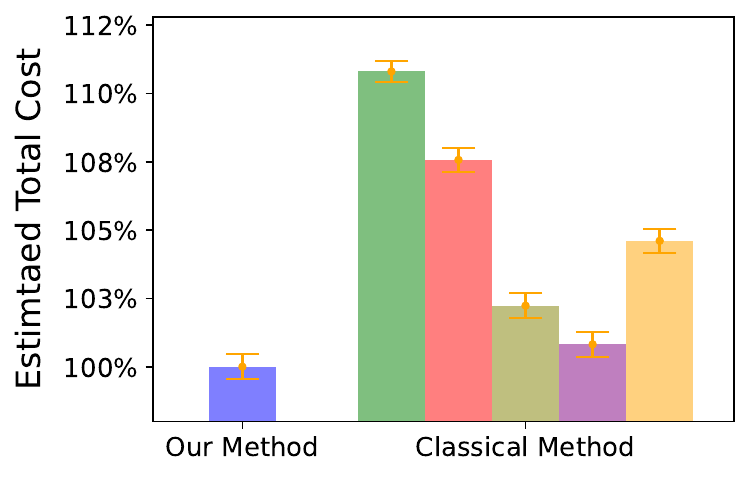} 
(iii) Trade-off between \\ filtering and hiring costs
\end{minipage}
\caption{For different methods, we consider the selected filtering level $\fl$
and present the associated estimated total cost of the AI-based triage system. 
The classical method maximizes recall subject to the
precision level $[0.93, 0.94, 0.95, 0.96, 0.97]$, positioned from left to right.
This method exhibits highly varying total cost even at high precision levels, making 
it hard to determine the best filtering level. In contrast, our method effectively
minimizes the total cost by cheap simulations of (reflected) Brownian motion.}
\label{fig: triage total cost}
\vspace{-.75em}
\end{figure}

In Figure~\ref{fig: triage total cost}, we present the average total cost 
over 10K sample paths of the simulated (reflected) Brownian motion, with 
$2\times$ the standard error encapsulated in the orange brackets. For the 
classical method, we set the precision level as
$[0.93, 0.94, 0.95, 0.96, 0.97]$, positioned from left to right. 
To facilitate comparison with our method, we normalize the estimated 
total cost by that of our method. We observe that the classical method
exhibits highly varying total cost (by $\sim 10\%-150\%$) even at high precision 
levels in Figure~\ref{fig: triage total cost}. This demonstrates the importance of 
selecting the right filtering level to minimize the total cost. In addition, 
for the classical method, it also shows the total cost is highly sensitive to 
the precision level. Therefore, the precision level serves as an important 
hyperparameter, and it is challenging to determine the best precision 
level that corresponds to optimal filtering level using the classical method. 

Such challenge arises since our method takes a holistic view of the entire 
triage system, yet the classical method only considers the prediction metrics. 
In our toy example, a higher precision level leads to a lower selected filtering 
level $\fl$, which results in lower filtering costs and higher hiring costs.
Figure~\ref{fig: triage total cost} (i)(ii) corresponds to 
simpler settings where the total cost aligns with prediction metrics. 
That is, when filtering costs or hiring costs dominate, the total cost is 
monotone with respect to the filtering level and thus the precision level, as shown 
in Figure~\ref{fig: triage total cost} (i)(ii). Therefore, when adopting the 
classical method, we can simply choose the precision level at the search boundary, 
which yields a filtering level near the search boundary that minimizes
the total cost. In contrast, in Figure~\ref{fig: triage total cost} (iii), where 
there is a trade-off between filtering and hiring costs, the total cost is 
non-monotone and ``U"-shaped with respect to the precision/filtering level. 
In this case, prediction metrics fails to capture the total cost. While 
hyperparameter (precision) tuning based on total costs is possible, 
the classical method merely shifts our search space to hyperparameters 
(precision levels). In other words, hyperparameter tuning is equivalent 
to a naive line search for the decision variable (filtering level $\fl$) 
based on simulated total cost and the classical method does not serve as
effective objectives/metrics. More importantly, without our 
Theorem~\ref{theorem: total cost of triage queue system}, 
the total cost can only be estimated through multiple costly simulations 
of the entire triage system. 
Our method, in contrast, effectively identifies the correct objective and 
finds the best filtering level through cheap simulations of 
(reflected) Brownian motion.

Our numerical experiments demonstrate the 
the effectiveness of our method and the importance of taking a holistic 
view of the entire content moderation system. We hope our method paves the 
way for more advanced system designs for complex AI-based triage systems in
practical use.

\end{document}